\documentclass[12pt]{amsart}
\usepackage{amsmath,bbm,amsfonts,amssymb,amsxtra,amscd,enumerate,amsthm,yfonts,verbatim}
\usepackage{color,xcolor}
\usepackage{latexsym}
\usepackage{epsfig}
\usepackage{fullpage}
\usepackage{hyperref}
\usepackage{comment}
\usepackage{enumitem}
\usepackage{esint}

\newcommand\e{\epsilon}
\renewcommand\l{\lambda}

\newcommand\les{\lesssim}
\newcommand\ges{\gtrsim}

\newcommand{\lmax}{\lambda^{max}}

\newcommand{\Mf}{\mathcal B}
\newcommand{\MF}{\mathcal M}
\newcommand{\gf}{\mathfrak{g}}

\newcommand\R{\mathbb{R}}
\newcommand\C{\mathbb{C}}
\newcommand\Z{\mathbb{Z}}
\newcommand\N{\mathbb{N}}
\renewcommand\S{\mathbb{S}}
\newcommand\K{\mathcal K}
\newcommand\FH{{\mathcal F}}
\newcommand\FtH{{\tilde{\mathcal F}}}
\newcommand\F{{\mathcal F}}
\newcommand{\calO}{\mathcal O}
\newcommand\W{\mathcal{W}}

\newcommand\cT{\mathcal{T}}
\newcommand\cK{\mathcal{K}}
\newcommand\tH{\tilde{H}}
\newcommand\tF{\tilde{F}}
\newcommand\tK{\tilde{K}}
\newcommand\tcK{\tilde{\mathcal{K}}}
\newcommand\tP{\tilde{P}}
\newcommand{\tm}{{\tilde m}}

\newcommand{\tX}{\tilde X}
\newcommand{\A}{\mathcal{C}^2}

\newcommand{\tpsi}{\tilde \psi}

\newcommand{\PP}{\mathbf{P}}

\newcommand{\lc}{\Lambda}

\newcommand\la{\langle}
\newcommand\ra{\rangle}
\newtheorem{t1}{Theorem}
\numberwithin{t1}{section} 
\newtheorem{l1}[t1]{Lemma}
\newtheorem{p1}[t1]{Proposition}
\newtheorem{c1}[t1]{Corollary}
\newtheorem{d1}[t1]{Definition}

\newtheorem{r1}{Remark}

\numberwithin{equation}{section}

\newcommand{\He}{H^1_e}
\newcommand{\dHe}{\dot{H}^1_e}
\newcommand{\Hde}{H^2_e}
\newcommand{\dHde}{\dot{H}^2_e}

\newcommand{\ixi}{\hat\xi}
\newcommand{\ieta}{\hat\eta}
\newcommand{\bX}{\bar{X}}
\newcommand{\LX}{{L\bar{X}}}

\usepackage{mathtools}
\mathtoolsset{showonlyrefs}

\begin{document}

\title[$2$-equivariant Schr\"odinger maps]{Near soliton evolution for $2$-equivariant Schr\"odinger Maps in two space dimensions}

\author{Ioan Bejenaru }
\address{Department of Mathematics, University of California, San Diego}
\email{bejenaru@ucsd.edu}

\author{ Mohandas Pillai}
\address{Department of Mathematics, University of California, San Diego}
\email{mkpillai@ucsd.edu}

\author{ Daniel Tataru}
\address{Department of Mathematics, University of California, Berkeley}
\email{tataru@math.berkeley.edu}

\begin{abstract}  
  We consider equivariant solutions for the Schr\"odinger Map equation in $2+1$ dimensions, with values into $\S^2$. Within each equivariance class $m \in \Z$ this admits a lowest energy nontrivial steady state $Q^m$, which extends to a two dimensional family of steady states by scaling and rotation. If $|m| \geq 3$ then these ground states are known to be stable in the energy space $\dot H^1$, whereas instability and even finite time blow-up along the  ground state family may occur if $|m| = 1$. 
  In this article we consider the most delicate case 
  $|m| = 2$. Our main result asserts that small $\dot H^1$ perturbations of the ground state $Q^2$ yield global in time solutions, which satisfy global dispersive bounds. Unlike the higher equivariance classes, here 
  we expect solutions to move arbitrarily far along the 
  soliton family; however, we are able to provide 
a time dependent bound on the growth of the scale modulation parameter. We also show that within the
equivariant class the ground state is stable in a slightly stronger topology $X \subset
  \dot H^1$.
\end{abstract}


\maketitle

\setcounter{tocdepth}{1}
\tableofcontents

\section{Introduction}

In this article we consider the Schr\"odinger map equation in $\R^{2+1}$
with values into $\S^2$,
\begin{equation}
u_t = u \times \Delta u, \qquad u(0) = u_0.
\label{SM}\end{equation}
This equation admits a conserved energy,
\[
E(u) = \frac12 \int_{\R^2} |\nabla u|^2 dx,
\]
and is invariant with respect to the dimensionless scaling
\[
 u(t,x) \to u(\lambda^2 t, \lambda x).
\]
The energy is invariant with respect to the above scaling,
therefore the Schr\"odinger map equation in $\R^{2+1}$ is 
{\em energy critical}.

In reviewing the literature about this problem, we note that \eqref{SM} can be generalized in several ways. The simplest one is by considering maps $u: \R^{n+1} \rightarrow \S^2$ (keeping in mind that $n=2$ renders the energy critical dimension). One can also generalize the problem by replacing the target manifold $\S^2$ with a K\"ahler manifold with a complex structure. 

Local solutions for regular large initial data have been constructed
in \cite{SuSuBa} and \cite{MG}. When the target manifold $\S^2$ is replaced  with a K\"ahler manifold with a compatible complex structure, local well-posedness for regular data has been established in \cite{DiWa-2}, see also \cite{Di} and \cite{KeLaPoStTo}. 

The global in time problem is a very difficult one, except in dimension $n=1$, that is for maps $u: \R^{1+1} \rightarrow \S^2$ (or some generalization of the base/target manifold), where it becomes energy subcritical;  large data global well-posedness in this case has been established in \cite{DiWa-1}.

The definitive result for the small data problem, 
for maps $u:\R^{n+1} \rightarrow \S^2$, with $n \geq 2$,
was obtained by two of the authors and collaborators in \cite{BIKT-1}, following earlier 
results in \cite{Be1},
\cite{Be2}, \cite{bik}, \cite{IoKe-1}, \cite{IoKe-2}, \cite{Ka},
\cite{KaKo}, \cite{KeNa}, \cite{NaStUh}, \cite{NaStUh2},
\cite{NaShVeZe}.  There global well-posedness and scattering were proved for initial data which
is small in the critical Sobolev space $\dot H^{\frac{n}2}$, which agrees with the  
energy space $\dot H^1$ when $n =2$.  The counterpart of these results when $\S^2$ is replaced with a more general Kh\"aler manifold was established more recently in \cite{Ze1} and \cite{Ze2}. 
The same problem with the base space $\R^2$ replaced with the hyperbolic plane was considered in \cite{LaLuOhSh}, where the authors establish the asymptotic stability of (a large class of) harmonic maps under the Schr\"odinger maps evolution.

We now return to the setup of maps $u:\R^{2+1} \rightarrow \S^2$ which is the energy critical one.  
The space of finite energy states for this problem is the space of $\dot H^1$ 
maps from $\R^2$ into $\S^2$, which separates into connected components according to the 
homotopy class $m \in \Z$. Within each homotopy class there exists an energy 
minimizer, called the ground state, which is unique up to symmetries, namely scaling and isometries of the base space $\R^2$ and of the target space $\S^2$. For each integer $m$ we  denote the corresponding ground state family by $\mathcal Q^m$. Then a natural question is whether the ground states are stable with respect 
to the Schr\"odinger map flow.

The ground state family for the trivial  $m=0$ homotopy class consists of constant maps, which have energy $0$. All the global results mentioned above (in the case $n=2$) required small energy; as a consequence those maps have trivial topology and, with respect to the energy norm, they can be seen a small perturbations of these constant maps.

For $m \neq 0$, the generator $Q^m$ for the ground state family can be taken to belong to the class of $m$-equivariant maps, which satisfy 
\[
u (R x) = R^{m} u(x),
\]
where $R$ stands for rotations around the origin in the plane, and around
the vertical axis on the sphere, with the same angle. The class of $m$-equivariant maps is closed with respect to the Schr\"odinger map flow, therefore it is natural to restrict the above stability question of $\mathcal Q^m$ to the class of $m$-equivariant maps. 

As it turns out, the answer to this question depends on the equivariance class $m$:
\begin{itemize}
\item If $|m| \geq 3$, then the ground state is stable in $\dot H^1$ within the equivariant class; this was proved in  \cite{GuKaTs-2} for $|m| \geq 4$, and in \cite{GuNaTs} for $|m|=3$.

\item If $|m| = 1$, then the  ground state is unstable in $\dot H^1$ within the equivariant class; this was proved in \cite{BeTa-1}. Furthermore, finite time blow-up may also occur along the ground state family, see \cite{MeRaRo} and \cite{Pe}.
\end{itemize}

This left open the case $|m|=2$, which is in some sense borderline and does not fit either of the two patterns above. 

The objective of the present article is to investigate the ground state stability exactly in this remaining case $|m|=2$ within the $\pm 2$ equivariance class.
In brief, our main results are as follows:

\begin{itemize}
    \item Small $\dot H^1$ perturbations of the ground state yield \emph{global in time} 
    solutions.
    \item The solutions satisfy universal dispersive and local energy bounds.
\item The solutions can move along the ground state manifold, but we prove 
quantitative time dependent bounds on the modulation parameters.
\item The ground state is stable in a slightly stronger topology.     
\end{itemize}

\subsection{ Homotopy classes and the ground state} \label{HCGS}
Given a $\dot H^1$ map $u : \R^2 \to \S^2$, its homotopy class 
$m \in \Z$ is defined as 
\[
m = \frac{1}{4\pi} \int_{\R^2} u \cdot (\partial_x u \times \partial_y u)\, dx
\]
and roughly speaking counts the number of times the plane wraps around the sphere, taking also the orientation into account.

The homotopy classes generate a partition of the family of all  $\dot H^1$ maps into connected components. Within each homotopy class, one 
may look for energy minimizers, which are called \emph{ground states}.
The Euler-Lagrange equation for ground states is the harmonic map equation, 
\[
-\Delta u  = u |\nabla u|^2,
\]
so the ground states are in particular harmonic maps and also steady states for the Schr\"odinger map equation.

For each nonzero integer $m$ there exists a  family $\mathcal Q^m$  of  ground states. To describe these families we begin with the maps $Q^m$ 
defined in polar coordinates by
\begin{equation} \label{FQm}
Q^m(r,\theta) = e^{m \theta R} \bar Q^m(r), \qquad
\bar Q^m(r)=\left( \begin{array}{ccc} h_1^m(r) \\ 0 \\ h_3^m(r) \end{array} \right),
\qquad m \in \Z \setminus\{0\}
\end{equation}
with 
\[
h_1^m(r)=\frac{2r^m}{r^{2m}+1}, \qquad h_3^m(r)=\frac{r^{2m}-1}{r^{2m} + 1}. 
\]
Here $R$ is the generator of horizontal rotations, which 
can be interpreted as a matrix or, equivalently, as the operator below
\[
 R = \left( \begin{array}{ccc} 0 & -1 & 0 \\ 1 & 0 & 0 \\ 0 & 0 & 0 \end{array}
\right)  , \qquad R u  = \overrightarrow{k} \times u.
\]
The families $\mathcal Q^m$ are constructed from $Q^{ m}$
via the symmetries of the problem, namely scaling and orientation preserving isometries
of the base space $\R^2$ and of the target space $\S^2$. $\mathcal Q^{-m}$ is equivalent to 
 $\mathcal Q^m$ modulo one reflection. 
The elements of $\mathcal Q^m$ are harmonic maps from $\R^2$ into
$\S^2$, and admit a variational characterization as the unique energy
minimizers, up to symmetries, among all maps $u: \R^2 \to \S^2$
within their homotopy class.  \medskip

In the above context, a natural question is to study Schr\"odinger
maps for which the initial data is close in $\dot H^1$ to one of the
$\mathcal Q^m$ families.  One may try to think of this as a small data
problem, but in some aspects it turns out to be closer to a large data
problem.  Studying this in full generality is very difficult.  In this
article we confine ourselves to a class of maps which have some extra
symmetry properties, namely the {\em equivariant} Schr\"odinger maps.
These are indexed by an integer $m$ called the equivariance class, and
consist of maps of the form
\begin{equation} \label{equiv}
u(r,\theta) = e^{m \theta R} \bar{u}(r).
\end{equation} 
In particular the maps $Q^m$ above are $m$-equivariant.
The case $m=0$ would correspond to spherical symmetry.
Restricted to equivariant functions the energy has the form
\begin{equation} \label{energy}
 E(u) = \pi \int_{0}^\infty \left( |\partial_r \bar{u}(r)|^2 +
\frac{m^2}{r^2}(\bar{u}_1^2(r)+\bar{u}_2^2(r)) \right) r dr.
\end{equation}
 In the case $m \ne 0$ this implies that the $m$-equivariant maps have the property that $\bar{u}(0):=\lim_{r \rightarrow 0} \bar{u}(r)$ and
  $\bar{u}(\infty):= \lim_{r \rightarrow \infty} \bar{u}(r)$ are well-defined distinct elements of the set $\{ \pm \vec{k} \}$.
In particular, if $m \neq 0$  this implies that m-equivariant maps
have a limit at $0$ and at infinity, both in the set $\{\pm \vec{k} \}$.
The homotopy class of $m$-equivariant maps depends on whether the two limits are equal:

\begin{enumerate}
    \item If $\bar{u}(0) = \bar{u}(\infty)$ then the map $u$ is topologically 
    trivial, i.e. has homotopy $0$.
    \item If $\bar{u}(0)=-\vec{k}$ and $\bar{u}(\infty) = \vec{k}$  then the map $u$ has homotopy $m$.
     \item If $\bar{u}(0)=\vec{k}$ and $\bar{u}(\infty) = -\vec{k}$  then the map $u$ has homotopy $-m$.
\end{enumerate}
To fix the notations in the sequel we will work with the connected component of $m$-equivariant maps with homotopy $m$, i.e. with  $m$-equivariant maps 
  satisfying
  \begin{equation} \label{eqvl} 
   \bar{u}(0)= -\vec{k}, \qquad \bar{u}(\infty)=\vec{k}.
  \end{equation}

Intersecting the full set $\mathcal Q^{m}$ with the $m$-equivariant
class and with the homotopy class of $Q^m$ we obtain the two parameter family $\mathcal Q^m_e$ generated from $Q^m$ by rotations and scaling,
\[
 \mathcal Q^m_e = \{ Q^m_{\alpha,\lambda}; \alpha \in \R/ 2\pi\Z, \lambda \in \R^+\}, 
\qquad  Q^m_{\alpha,\lambda}(r,\theta)= e^{m \alpha R}  Q^m (\lambda r,\theta)
\]
Here $Q^m_{0,1} = Q^m$. Their energy depends on $m$ as follows:
\[
E(Q^m_{\alpha,\l}) = 4\pi |m|:= E(\mathcal Q^m).
\]

 \subsection{ A brief history of the problem} The study of equivariant Schr\"odinger maps for $m$-equivariant
 initial data close to $\mathcal Q^m_e$ was initiated by Gustafson, Kang, Tsai
 in \cite{GuKaTs-1}, \cite{GuKaTs-2}, and continued by Gustafson, Nakanishi,
 Tsai in \cite{GuNaTs}. The energy conservation suffices to confine
 solutions to a small neighborhood of $\mathcal Q^m_e$ due to the
 inequality (see \cite{GuKaTs-1})
\begin{equation}
\text{dist}_{\dot H^1} (u,\mathcal Q^m_e)^2= \inf_{\alpha,\lambda}
\| Q^m_{\alpha,\lambda} - u\|_{\dot H^1}^2 \lesssim E(u) - E(\mathcal Q^m),
\label{stab}\end{equation}
which holds for all $m$-equivariant maps $u:\R^2 \to \S^2$ in the homotopy class of
$\mathcal Q^m_e$ with $0 \leq E(u) - E(\mathcal Q^m) \ll 1$.  One can
interpret this as an orbital stability result for $\mathcal Q^m_e$.
However, this does not say much about the global behavior of solutions
since these soliton families are noncompact; thus one might have even
finite time blow-up while staying close to a soliton family.

To track the evolution of an $m$-equivariant Schr\"odinger 
map $u(t)$ along $\mathcal Q^m_e$ we use  
functions $(\alpha(t),\lambda(t))$ describing trajectories 
in $\mathcal Q^m_e$. A natural choice is to 
choose them as minimizers for the infimum in \eqref{stab}; 
while this is feasible and useful for the local theory, see \cite{GuKaTs-1},
it is much less helpful for the purpose of the global theory.
Instead, we will allow ourselves more freedom, and be content
with any choice $(\alpha(t),\lambda(t))$ satisfying
\begin{equation}
\| u - Q_{\alpha(t),\lambda(t)}^m\|_{\dot H^1}^2 \lesssim
E(u) - E(\mathcal Q^m). 
\label{goodal}\end{equation}

An important preliminary step in this analysis is the next result
concerning both the local wellposedness in $\dot H^1$ and the
persistence of higher regularity:
\begin{t1} \label{th2}
Let $|m| \geq 1$. There exists $\delta >0$ and $\sigma >0$ such that the following holds true. Given  an $m$-equivariant initial data $u_0$ in the homotopy class of  $\mathcal Q^m_e$ and with $E(u_0) \leq 4\pi |m| + \delta^2$, let $\lambda_0=\lambda(u_0)$
be the  choice of $\lambda$ which (in combination with $\alpha$) achieves the infimum in \eqref{stab}. Then the equation \eqref{SM} with initial data $u_0$ is locally well-posed in $\dot H^1$ on a time interval $I=[-\frac{\sigma}{\lambda_0^2}, \frac{\sigma}{\lambda_0^2}]$, and 
$\lambda \approx \lambda_0$ in $I$. 

If, in addition, $u_0 \in \dot{H}^2$ then $u \in C_t (I:\dot{H}^2)$. 
Furthermore, the $\dot H^1 \cap \dot H^2$ regularity persists for as long as 
the function $\lambda(t)$ in \eqref{goodal} remains in a compact set. If $T_{max}$ is the maximal time 
for which there is a unique solution $u(t) \in C([0,T_{max}):\dot H^2)$, then we have that $T_{\max} \geq \frac{\sigma}{\lambda_0^2} $. In addition, if $T_{max}$ is finite then  $\lim_{t \nearrow T_{max}} \lambda(t)=\infty$. 
\end{t1}

This follows from Theorem 1.1 in \cite{GuKaTs-1} and Theorem $1.4$ in \cite{GuKaTs-2}.  Given the above
result, the main remaining  problem is to understand whether the steady
states $Q_{\alpha,\lambda}^m$ are stable or not; in the latter case,
one would like to understand the dynamics of the motion of the
solutions along the soliton family.  The steady states turn out to be stable
in the case of large $m$, which  was considered in prior work:

\begin{t1}[\cite{GuKaTs-2} for $|m| \geq 4$, \cite{GuNaTs} for $|m|=3$]
The solitons $Q_{\alpha,\lambda}^m$ are asymptotically stable in the $\dot H^1$
topology within the $m$-equivariant class.
\end{t1}

However, in the cases $|m|=1,2$ the behavior is expected to be different in nature.
The $|m|=1$ case was considered by two of the authors in \cite{BeTa-1}.  There 
it is shown that the ground states are in general unstable in the energy topology,
but on the other hand that they are stable with respect to perturbations which are 
small in a stronger topology. Further, Merle, Raphael and Rodnianski constructed solutions which blow up in finite time along the ground state family, see \cite{MeRaRo}; soon after, Perelman provided an alternative construction of singularity formation in \cite{Pe}. 

In this article we consider the case $m=2$, where there are no analogue results to the ones just mentioned for in the case $m=1$ or $m \geq 3$; in other words no blow-up construction is known, and no stability regimes have been established. 

Our first result  will be to prove  the global in time well-posedness of the equivariant Schr\"odinger Maps with data near $Q^2$;  as part of the proof, we establish global bounds on the potential growth of the noncompact modulation parameter $\lambda(t)$, and thus on the growth of higher regularity norms. Our result leaves open the possibility of blow-up/relaxation in infinite time, which will be considered in a forthcoming paper.

Our second result identifies stability regimes and it is similar in spirit to the one obtained in \cite{BeTa-1} in the case $m=1$. It essentially says that if the original data $u_0$ is close to $Q^2$ in a slightly stronger topology than the  one induced by the energy, then we obtain uniform global bounds and rule out the possibility of blow-up/relaxation in infinite time. There are two main differences between our result here and the one in \cite{BeTa-1} for $m=1$. The first one is that here the slightly stronger topology is rather close to the one induced by the energy. The second one is that while the smallness of  $u_0 - Q^2$ in the stronger norm guarantees stability, in effect our result still applies for any finite size of $u_0 - Q^2$ in the stronger norm, in which case it guarantees global bounds (no infinite-time blow-up) in this scenario.

\subsection{ The main results}
In order to state our first main result, we need to recall two important concepts: the modulation parameters $\alpha$ and $\lambda$, and the reduced field $\psi$.  

The modulation parameters describe the motion of our solution along the ground state family, and have been already discussed in the paragraph leading to \eqref{goodal}. As noted there, one has quite a bit of flexibility on what are good choices for these parameters. In practice one imposes some orthogonality conditions which uniquely identify these parameters in the regime described in \eqref{goodal}. Our \emph{orthogonality condition} is defined by the relation \eqref{ldef} in Section~\ref{seccoulomb}. 

Next, a common strategy in analyzing the SM equation \eqref{SM} is to write $\nabla u \in T_u \S^2$ in an appropriate gauge, that is a choice of an orthonormal frame in $T_u \S^2$. In complex notation, this representation gives rise to the complex valued differentiated field $\psi$, which is small in $L^2$ and solves a nonlinear Schr\"odinger type evolution.  This construction is developed in Section \ref{seccoulomb}; our gauge choice is the classical Coulomb gauge. The $L^2$ size of $\psi$
is provided by a very simple formula,
\begin{equation}\label{u=psi}
\|\psi(0)\|_{L^2}^2 = E(u_0)-8\pi,  
 \end{equation}
which in particular shows it is a conserved quantity.  One advantage of working with the
 differentiated field $\psi$ is that it is more amenable to a dispersive type analysis, that is $\psi$ can be measured in function spaces which contain information such as Strichartz estimates and local energy decay estimates. 

With these objects at hand, the analysis of the dynamics of the Schr\"odinger map $u$ near $\mathcal Q^2_e$ (and more generally near $\mathcal Q^m_e$) can be essentially reduced to two major tasks: 
\begin{enumerate}[label=\roman*)]
\item  a dispersive analysis of the nonlinear Schr\"odinger type  PDE which governs the evolution of the differentiated field $\psi$;

\item  an analysis of the \emph{modulation equation}, i.e. the nonlinear ODE that governs the evolution of the modulation parameters $\alpha$ 
and $\lambda$.
\end{enumerate}

One has to keep in mind that the PDE and the ODE mentioned above are coupled and thus the two tasks are not performed independently. In most of the prior works, the variation of $\lambda(t)$ was small and one could essentially freeze its value for the PDE analysis; at the same time, once the PDE yields enough information of dispersive type, one could show 
that $\lambda'$ is small in $L^1_t$. While this description is a bit oversimplified, it highlights the light coupling in the prior works between the PDE and ODE analysis discussed above. In our context, the coupling between the PDE and the ODE analysis is significantly more involved at all levels.

We also  emphasize another novel feature of our analysis: in all prior works on the Schr\"odinger Map equation \eqref{SM} in this context, the results seem to indicate that if blow-up occurs, that is $\lim_{t \nearrow T_{max}} \lambda(t)=+\infty$, then the dispersive properties of the field $\psi$ are lost on the (forward) maximal interval of existence $I_{max}=[0,T_{max})$; see for instance Lemma 3.1 in \cite{GuKaTs-1} and Lemma 2.6 in \cite{GuKaTs-2}, as well as Proposition 5.2 in \cite{GuNaTs} where one assumes that $\lambda$ stays close to $1$ in order to recover dispersive estimates. 

Surprisingly, our first result in this paper establishes uniform dispersive estimates (measured by Strichartz and local energy decay norms)
on the full interval time of existence even in the case of potential blow-up (whether in finite or infinite time).

\begin{t1} \label{tmain-L}
Assume that we have a $\pm 2$-equivariant initial data $u_0 \in \dot H^1$ in the homotopy class of $Q^{\pm 2}$, and with energy
\begin{equation}\label{small-energy}
E(u_0) \leq     8\pi+\delta^2, \qquad \delta \ll 1.
\end{equation}
Let $I_{max}=(T_{min},T_{max})$ be the maximal time of existence of the solution $u$ to the Schr\"odinger map flow \eqref{SM} with initial data $u_0$. 
Then the associated Coulomb gauge field $\psi$ and the associated modulation parameters $\lambda$ and $\alpha$ satisfy
\begin{equation} \label{l1est-m}
    \|\frac{\lambda'}{\lambda^2} \|_{L^2_t(I_{max})}
    +  \|\frac{\alpha'}{\lambda} \|_{L^2_t(I_{max})}
    \les \delta,
\end{equation}
respectively
\begin{equation} \label{psiest-m}
    \| \psi \|_{l^2 S(I_{max})} \les \delta. 
\end{equation}

\end{t1}

Importantly, the implicit constants in the above bounds do not depend on the length of the time interval. The space $l^2 S$, which is defined later in Section \ref{s:linear}, contains standard dispersive estimates, namely Strichartz bounds and local energy decay. In particular we have the embedding $l^2 S \subset L^4_{x,t}$.

The above theorem brings a new insight into the 
stability theory: even if blow-up occurs, the dispersive properties of the gauge field (as measured by Strichartz and local energy decay norms) hold true uniformly all the way up to the blow up time. 
However, we point out that the estimate in \eqref{l1est-m} is not strong enough to preclude the scenario that 
$T_{\max} < +\infty$ and $\lim_{t \nearrow T_{max}} \lambda(t) = \infty$ which corresponds to blow-up
in finite time; nor does it preclude a similar scenario happening at $T_{min} > -\infty$.

\medskip

The $l^2$ summation in \eqref{psiest-m} corresponds to a Littlewood-Paley dyadic frequency decomposition for $\psi_0 \in L^2$ or equivalently for $u_0 - Q^{\pm 2}_{\alpha,\lambda} \in \dot H^1$ 
with $(\alpha,\lambda)$ as in \eqref{goodal}.
It will also be of interest to consider the situation where the $l^2$ Besov bound is supplemented with an $l^1$ Besov bound.
For this purpose we introduce a stronger topology
\begin{equation}
 X = \dot B^{1}_{2,1} \subset \dot H^1.    
\end{equation}
This is a space of functions in $\R^2$, which, when restricted to the class of equivariant functions $u(x) = \bar u(r) e^{im\theta}$, yields a space for the radial profile $\bar u$, denoted by $\bX=\dot B^1_{2,1,e}$, with norm 
\begin{equation}
\| \bar u \|_{ \dot B^1_{2,1,e}} = \| u\|_{\dot B^1_{2,1}}   
\end{equation}
Correspondingly for $\psi$ we should  use the space $\LX:=\dot B^0_{2,1,e}$ at fixed time; the notation here is motivated in Section~\ref{s:XLX}. In this setting, as we prove later in Proposition~\ref{uQpsi}, we have the norm relation
\begin{equation}
\| \psi \|_{ \LX} \approx \| u - Q^{\pm 2}_{\alpha,\lambda}\|_{X}.   
\end{equation}
Also to measure the space-time regularity of $\psi$ we will use the $l^1S$ norm, which is the $l^1$ Besov analogue of the earlier $l^2 S$
space. We refer the reader to Section~\ref{defnot} for more details on function spaces.
Our main bounds in the $l^1$ Besov class are as follows:

\begin{t1} \label{tmain-L1}
Let $u$  be a solution to the Schr\"odinger map equation as in Theorem~\ref{tmain-L}. Assume in addition that the initial data $u_0$ satisfies
\begin{equation}
  u_0 - Q^{\pm 2}_{\alpha(0),\lambda(0)} \in  X = \dot B^{1}_{2,1}.   
\end{equation}
 Then we have 
\begin{equation} \label{psiest-m1}
    \| \psi \|_{l^1 S(I_{max})}  \lesssim \|  u_0 - Q^{\pm 2}_{\alpha(0),\lambda(0)}\|_{X}.
    \end{equation}
\end{t1}

We emphasize that the implicit constant is universal, despite the fact that no smallness is assumed for the $l^1$ norm on the right. 
There one may choose any $(\alpha,\lambda)$ as in \eqref{goodal}, but the norm on the left 
potentially depends on our exact choice of the 
modulation parameters.
 \medskip

So far, our results have only considered the evolution of the reduced field $\psi$. The next natural step is to also consider 
the time evolution of the modulation parameters $(\alpha(t),\lambda(t))$, of which $\lambda(t)$ plays the leading role
as it governs the potential blow-up behavior.

Our first result in this direction establishes that there is \emph{no finite time blow-up}, but potentially allows for an infinite time blow-up:

\begin{t1} \label{tmain-G} 
Assume that we have a $\pm 2$-equivariant initial data $u_0 \in \dot H^1$ in the homotopy class of $Q^{\pm 2}$, and with energy as in \eqref{small-energy}.
Then the Schr\"odinger Maps equation \eqref{SM} has a unique global in time solution $u(t)$. Furthermore, 
if we assume the normalization $\lambda(0)=1$, then
the modulation parameter $\lambda(t)$ satisfies 
\begin{equation} \label{lamdagb} 
\la t \ra^{-\frac12} \lesssim   \lambda(t) \les e^{\la t \ra^{C \delta^2}}, \forall t \in \R,
\end{equation}
for some universal constant $C$. 
\end{t1}

Here we recall that, as a consequence of Theorem \ref{tmain-L}, the field associated with this solution satisfies uniform global dispersive bounds as stated in \eqref{psiest-m}.

In \eqref{lamdagb} the bound on the right on $\lambda(t)$ is the interesting one. The one on the left is simply the self-similar scale which follows directly from the local well-posedness result in Theorem~\ref{th2}. The normalization $\lambda(0)$ may be removed by rescaling, in which case the bound \eqref{lamdagb} becomes
\begin{equation} \label{lamdagb+} 
\la t  \lambda(0)^2\ra^{-\frac12} \lesssim  \frac{ \lambda(t)}{\lambda(0)} \les e^{\la t \lambda(0)^2\ra^{C \delta^2}}, \forall t \in \R,
\end{equation}

In particular, the above theorem allows for an infinite time blow-up. In forthcoming work, we aim to show that this can actually happen for well chosen initial data. 

On the other hand, if in addition we assume $l^1$ dyadic summability for the data, then we can prevent this from happening:

\begin{t1} \label{tmain-G1}
Let $u$  be a solution to the Schr\"odinger map equation as in Theorem~\ref{tmain-G}. Assume in addition that the modulation parameters $\alpha(0),\lambda(0)$ corresponding to the initial data $u_0$ satisfy
\begin{equation}\label{X-data}
 \| u_0 - Q^{\pm 2}_{\alpha(0),\lambda(0)}\|_{X} \leq M  < +\infty. 
\end{equation}
Then we have 
\begin{equation} \label{lamdagb1} 
  | \ln \lambda(t) - \ln \lambda(0)| + |\alpha(t) - \alpha(0)| \les M+ M^2, \qquad \forall t \in \R.
\end{equation}
\end{t1}

Remark. The above condition \eqref{X-data} in which we use the specific choice of modulation parameters dictated by \eqref{ldef} can be easily replaced by any good choice of modulation parameters subject to \eqref{goodal}; this simply changes $M$ to $M+C\delta$, where $\delta$ was defined in Theorem \ref{tmain-L}, but by definition $\delta \lesssim M$.

\bigskip

 Our final result below asserts that the soliton $Q^{\pm 2}$ is stable in the $X$ topology and is a direct consequence of the previous Theorem. 
\begin{t1} \label{tstable}  Assume that we have a $\pm 2$-equivariant initial data $u_0 \in \dot H^1$ in the homotopy class of $Q^{\pm 2}$,  whose initial data  satisfies 
\begin{equation}
 \| {u}_0 -  Q^{\pm 2}\|_{X}\leq \gamma \ll 1. 
\label{tdata}\end{equation}
Then there exists a unique global solution $u$ for the Schr\"odinger map equation \eqref{SM}, which satisfies
${u} -  Q^{\pm 2} \in C(\R;X)$ and 
\begin{equation}
 \| u -  Q^{\pm 2}\|_{C(\R;X)} \lesssim \gamma.
\label{tsolution}\end{equation}
\end{t1}
One consequence of the above theorem is that classical perturbations of $Q^{\pm 2}$ that are small in the energy class
yield global solutions which do not blow-up in infinite time; in particular they will obey classical bounds uniformly in time. By classical perturbations we mean perturbations that are smooth and decay at infinity, enough to control the $\dot B^1_{2,1}$ norm (which does not have to be small).  

Another consequence of the above theorem is that any potential infinite time blow-up construction would not be stable in the energy class. Indeed, if an initial data $u_0$, with $\|u_0-Q^{\pm 2}\|_{\dot H^1} \ll 1$, leads to a solution which blows up in an infinite time, then for each $\delta > 0$ one can construct another data $v_0$ with $\|u_0 -v_0\|_{\dot H^1} \leq \delta$ and with $\| v_0 -Q^{\pm 2}\|_{\dot B^1_{2,1}} < + \infty$ by cutting off low and high frequencies of $u_0 -Q^{\pm 2}$.  In the equivariant class of functions this is easily achieved by performing 
a direct Littlewood-Paley truncation
for $u_0-Q^{\pm 2}$, followed by a geometric projection of the new $u_0$ back on the sphere. The instability of the blow-up in the energy class then follows from the fact that the solution with data $v_0$ does not blow-up in infinite time (by Theorem \ref{tstable}).

\subsection{ Acknowledgements} 
The first author was supported by the NSF grant DMS-1900603. The second author was supported by the NSF grant DMS-2103106.
The third author was supported by the NSF grant DMS-2054975 as well as by a Simons Investigator grant from the Simons Foundation.

\section{Definitions and notations.} 
\label{defnot}
In this section we cover some definitions and notations which will be used extensively thoghout the paper. The section does not intend to exhaust all of them, as many objects are defined later as the paper progresses.

While at fixed time our maps into the sphere are functions defined on
$\R^2$, the equivariance condition allows us to reduce our analysis to
functions of a single variable $|x|=r \in [0,\infty)$.  One such
instance is exhibited in \eqref{equiv} where to each equivariant map
$u$ we naturally associate its radial component $\bar u$.  Some other
functions will turn out to be radial by definition, see, for instance,
all the gauge elements in Section \ref{seccoulomb}.  We 
agree to identify such radial functions with the corresponding
one dimensional functions of $r$.  Some of these functions 
are complex valued, and this convention allows us to use the bar 
notation with the standard meaning, i.e. the  complex conjugate.

Even though we work mainly with functions of a single spatial variable
$r$, they originate in two dimensions. Therefore, it is natural to
make the convention that for the one dimensional functions all the
Lebesgue integrals and spaces are with respect to the $rdr$ measure,
unless otherwise specified.

We define $A\les B$, if there is a universal constant $C>0$ such that $A\leq C B$, and $A\ges B$ iff $B\les A$. Further, we define $A\approx B$ iff both $A\les B$ and $B\les A$. Also, we define $A\ll B$ if $A \leq cB$ for some constant $c>0$ that can be chosen very small, depending only on some universal constants. Also, $A\gg B$ iff $B\ll A$. For $y>0$, we write $x=O(y)$ to mean $|x| \leq C y$ for some universal constant $C>0$.

Given a positive parameter $\lambda$, and a function $f:(0,\infty) \rightarrow \C$, we define its dilation by $\lambda$ as follows
\[
f^\lambda(r)= f(\lambda r), \quad \forall r \in (0,\infty). 
\]

On the physical space side, we will use a dyadic partition of $(0,\infty)$  (which has the natural extension  to $\R^2$) into sets  $\{A_j\}_{j \in \Z}$ given by
\[
A_j = \{  2^{j-1} < r < 2^{j+1}\}.
\]
We also use the notation $A_{< k}=\cup_{j < k} A_j$ as well as
$A_{> k}, A_{\geq k}$ which are similarly defined. Throughout the paper we  involve 
smooth approximations of the characteristic functions of these sets. We construct in the usual manner functions $\chi_j,\ j \in \Z$ with the following properties
\[
\chi_{j}(r) = \chi(\frac{r}{2^{j}}), \quad \chi \in C^{\infty}_{c}((\frac{1}{4},2)), 
\quad \chi(r) =1, \quad \frac{1}{2} \leq r \leq 1,
\]
and which satisfy the summation property
\[
\sum_j \chi_j(r) =1, \quad \forall r \in (0,\infty).
\]
Finally we also define 
\[
\chi_{< k} = \sum_{j <k} \chi_j, \quad \chi_{\leq k} = \sum_{j \leq k} \chi_j \quad  \chi_{> k} = \sum_{j > k} \chi_j, \quad \chi_{\geq k} = \sum_{j \geq k} \chi_j.
\]

\bigskip

In Section \ref{HCGS} we encountered the following functions
\[
h_1^m(r)=\frac{2r^m}{r^{2m}+1}, \qquad h_3^m(r)=\frac{r^{2m}-1}{r^{2m} + 1}; 
\]
as detailed in \eqref{FQm} they were the components of the basic ground state $Q^m$ and they will play a crucial role in the analysis in this paper. These functions are continuously differentiable on their domain $[0,\infty)$, have limits at $\infty$, 
$\lim_{r \rightarrow \infty} h_1^m(r)=0, \lim_{r \rightarrow \infty} h_3^m(r)=1$, 
and satisfy the following obvious property:
\[
(h_1^m(r))^2 + (h_3^m(r))^2=1, \quad \forall \ r \in [0,\infty);
\]
their derivatives are given by
\begin{equation} \label{derhf}
\partial_r h_1^m= -m \frac{h_1^m h_3^m}r, \quad \partial_r h_3^m= m \frac{(h_1^m)^2}r.
\end{equation}
We point out that fairly early in the paper we specialize in the case $m=2$ and drop the index $m=2$ from notation, that is we will use $h_1:=h_1^2,\ h_3:=h_3^2$. 

\bigskip 

We will use several modified Fourier transforms, adapted to our setting, and use Greek letters  
such as $\xi,\eta \in \R^+$ to denote frequencies. 
For a positive $\xi$, we use the notation $\ixi = \min\{\xi,1\}$.

In the context of various Fourier transforms, we involve Fourier projectors which are defined in the standard fashion.We start with the functions $m_j, j \in \Z$ satisfying the following properties
\[
m_{j}(\xi) = m(\frac{\xi}{2^{j}}), \quad m \in C^{\infty}_{c}((\frac{1}{4},2)), \quad m(\xi) =1, \quad \frac{1}{2} \leq \xi \leq 1.
\]
and which satisfy the summation property
\[
\sum_j m_j(\xi) =1, \forall \xi \in (0,\infty). 
\]
We also define $m_{<k} , m_{>k}, m_{\geq k}$ in the standard way, that is $m_{<k} : =\sum_{j<k} m_j$ etc.
The reader may notice that these functions are not different from the previously defined $\chi_j$ 
functions; the only reason we use different letters (for the same onject) is that they are used in different context: the $\chi$'s are used on the physical side, while the $m$'s are used on the Fourier side. 

We also make use of the following functions: given $\widetilde{m} \in C^{\infty}_{c}((\frac18,4)), \widetilde{m}(x)=1$ in the support of $m$ then, we let $\widetilde{m}_{j}(x) = \widetilde{m}(\frac{x}{2^{j}})$. $\widetilde{m}_j$ are very similar to $m_j$, except that they do not enjoy the partition of unity property from above; in addition the  following holds true: $\widetilde{m}_j \cdot m_j = m_j$. 

As mentioned earlier we will employ several Fourier transforms and we list them below for convenience. 

\begin{itemize}

\item The standard Fourier transform. This will be briefly used in $\R^2$ in the section just below. In sections \ref{s:mod2} and \ref{s:ode} it will be used in the context of time-dependent functions.

\item The Hankel transform (decomposition into Bessel frames). This will be used in Section \ref{s:XLX}, see in particular Section \ref{SBF}. The use of the Hankel transform is natural in the context of equivariant setup.

\item The Fourier transform in the frames associated to the operators $H$ and $\tilde H$; this theory is fully developed in Section \ref{spectral}. The operators $H$ and $\tilde H$ (see their precise definition in \eqref{H-factor}) occur naturally in our problem once we linearize near the solitons.  

\end{itemize}

We use the following function $\chi_{i=j}=2^{-\frac{|i-j|}2}$ defined for $i,j \in \Z$; this can be seen as a smoother version of the Dirac mass $\delta_{i=j}$.

For $\lambda >0$ and $j \in \Z$ we also use the notation
\[
\chi_{\lambda=2^{j}} = \left(\frac{2^{j/2}}{\sqrt{\lambda}} \mathbbm{1}_{\{\frac{2^{j}}{\lambda} \leq 4\}}+\frac{\lambda}{2^{j}} \mathbbm{1}_{\{\frac{2^{j}}{\lambda} \geq 1\}}\right)
\]
which is essentially $\chi_{i=j}$, except with $i$ formally replaced by $\log_{2}(\lambda)$, and has slightly better decay as $\frac{2^{j}}{\lambda} \rightarrow \infty$.
\subsection{ Sobolev and Besov spaces} \label{SBS-def}
Since equivariant functions are
easily reduced to their one-dimensional companion via \eqref{equiv},
we introduce the one dimensional equivariant version of $\dot H^1$,
\begin{equation}\label{defhe}
\| f \|_{\dHe}^2 = \| \partial_r f\|_{L^2(rdr)}^2 + 4 \| r^{-1} f \|_{L^2(rdr)}^2, \quad 
\| f \|_{\He}^2=\| f \|_{\dHe}^2 + \| f \|_{L^2(rdr)}^2.
\end{equation}
This is natural since for functions $u: \R^2 \to \R^2$ with $u(r,\theta) = e^{2\theta R}
\bar u(r)$ we have
\[
\| u \|_{\dot H^1}= (2\pi)^\frac12 \| \bar u \|_{\dHe},\qquad 
\| u \|_{H^1}= (2\pi)^\frac12 \| \bar u \|_{\He}.
\]
In a similar fashion we define $\dHde$ and $\Hde$ by the norms
\[
\| f \|^2_{\dHde}= \| \partial_r^2 f \|^2_{L^2} +\| r^{-1}\partial_r f \|_{L^2}^2+ \| r^{-2}f \|_{L^2}^2,
\qquad \| f \|^2_{\Hde}=\| f \|^2_{\dHde} + \| f \|_{L^2}^2
\] 
as the natural substitute for $\dot{H}^2$ and $H^2$.

We will also use the dual space $\dot H^{-1}_e$, defined by
\[
\| f \|_{\dot H^{-1}_e} =\sup_{\phi \in \dHe; \|\phi\|_{\dHe}=1} \la f, \phi \ra. 
\]

 We can think of the spaces $\dot H^1_e$, $\dot H^2_e$ as 
the "$H^1$", respectively the  "$H^2$" space associated to the equivariant Laplacian 
\[
-\Delta_m = - \partial_r^2 - \frac{1}{r} \partial_r 
+ \frac{m^2}{r^2},
\]
precisely 
\[
\dot H^1_e = D((-\Delta_m)^\frac12)):=\{ f: (-\Delta_m)^\frac12 f \in L^2\}. 
\]
and
\[
\dot H^2_e = D(-\Delta_m):=\{ f: -\Delta_m f \in L^2\}. 
\]
In our problem the case $m = 2$ is relevant, but to define the Sobolev spaces the choice of $|m| \geq 2$ is not important. The above Sobolev spaces can be characterized using spectral projectors. If $P_k^e$ are the spectral projectors associated 
to $-\Delta_m$, (which can be defined
by restricting the corresponding projectors associated to $-\Delta$ in $\R^2$ to the $m$-equivariant class, or equivalently by using  Bessel frames, as discussed later in Section~\ref{SBF}), then we have 
\[
\| f\|^2_{\dot H^s_e} \approx \sum_{k \in \Z} 2^{2ks} \|P_k^e f\|_{L^2}^2,
\]
where $s$ can be any of the values used above.

In some of the analysis carried in the paper, using these spectral projectors and sharp Fourier decompositions turns out to be counterproductive; instead we can use more robust decompositions. Precisely we have the following
\begin{equation}\label{He-relax}
\| f\|^2_{\dHe}
\approx \inf_{f = \sum f_k} \sum_k (2^k \|f_k\|_{L^2}
+ 2^{-k} \| f_k\|_{\dot H^2_e})^2
\end{equation}
and
\begin{equation}\label{L2-relax}
\| f\|_{L^2}^2
\approx \inf_{f = \sum f_k} \sum_k (2^k \|f_k\|_{\dot H^{-1}_e}
+ 2^{-k} \| f_k\|_{\dot H^1_e})^2.
\end{equation}
It is a straightforward exercise to establish \eqref{He-relax} and \eqref{L2-relax}.

Correspondingly we define associated Besov spaces. The ones we use in this paper are $\dot B^{1}_{2,1}$, $\dot B^{0}_{2,1,e}$ and $\dot B^{0}_{2,1,e}$ with 
norms 
\[
\| u \|_{\dot B^{1}_{2,1}} = \sum_{k}
2^k \|P_{k} f\|_{L^2}, \qquad 
\| f \|_{\dot B^{0}_{2,1,e}} = \sum_{k} \|P^e_{k} f\|_{L^2}, \qquad  
\| f \|_{\dot B^{1}_{2,1,e}} = \sum_{k}
2^k \|P^e_{k} f\|_{L^2}. 
\]
In the first norm $P_k$ represent the standard 
Littlewood-Paley projectors in $\R^2$, while in the latter norms, just as above, $P_k^e$ are the spectral projectors associated 
to $-\Delta_m$. 
For a complete definition, we let $\dot B^{1}_{2,1,e}$ be the subspace of functions $f \in \dHe$ which have the property that $\|f\|_{\dot B^{1}_{2,1,e}} < \infty$; similarly $\dot B^{0}_{2,1,e}$ is the subspace of  functions $f \in L^2$ with the property that $\|f\|_{\dot B^{0}_{2,1,e}} < \infty$. 
The key upgrade that the spaces $\dot B^{1}_{2,1,e}$, and $\dot B^{0}_{2,1,e}$,  have over their counterparts $\dHe$, respectively $L^2$, is that they bring an $l^1 L^2$ structure versus the standard $l^2 L^2$ one; here $l^1$ and $l^2$ indicate the norm used for the sequence $(\|P^e_{k} f\|_{L^2})_{k \in \Z}$. 

For brevity we also introduce some shorter notations, which will be suggestive in different contexts:
\[
X: = \dot B^{1}_{2,1}, \qquad LX: = \dot B^{0}_{2,1},
\]
which will be used for functions defined on $\R^2$, and, the corresponding for one dimensional functions,
\[
 \bX:= \dot B^{1}_{2,1,e},
\qquad \LX :=  \dot B^{0}_{2,1,e},
\]
where we restrict $X$ to functions $u$ of the form  
$u(r,\theta) = e^{2\theta R}\bar u(r)$, so that we
have the algebraic and topological identification $u \in X$ iff $\bar u \in \bX$.

Just as we did above for the Sobolev spaces, we can alternatively define these Besov spaces as follows:

\begin{itemize}
\item for $j= 1$ we can set
\begin{equation}\label{X-relax}
\| u\|_{\dot B^{1}_{2,1,e}} 
= \inf_{u = \sum u_k} \sum_k 2^k \|u_k\|_{L^2}
+ 2^{-k} \| u_k\|_{\dot H^2_e};
\end{equation}
\item  for $j = 0$ we have
\begin{equation}\label{LX-relax}
\| u\|_{\dot B^{0}_{2,1,e}} 
= \inf_{u = \sum u_k} \sum_k 2^k \|u_k\|_{\dot H^{-1}_e}
+ 2^{-k} \| u_k\|_{\dot H^1_e}.
\end{equation}
\end{itemize}

\subsection{Integration operators on radial functions}
Two operators which are often used on radial functions are 
$[\partial_r]^{-1}$ and $[r \partial_r]^{-1}$ defined as
\begin{equation} \label{drdef}
[\partial_r]^{-1} f(r) = - \int_{r}^\infty f(s) ds, \qquad 
[r \partial_r]^{-1}f(r) = - \int_{r}^\infty \frac{1}s f(s) ds.
\end{equation}
A direct argument shows that we have the Hardy type inequality
\begin{equation} \label{rdrm}
\| [r\partial_r]^{-1}f \|_{L^p} \lesssim_p \| f \|_{L^p}, \qquad 1 \leq p < \infty.
\end{equation}
We also have a weighted version 
\begin{equation} \label{rdrmw}
\| w [r\partial_r]^{-1}f \|_{L^p} \lesssim_p \| w f \|_{L^p}, \qquad 1 \leq p < \infty,
\end{equation}
under the assumption that $g(r)=w(r) r^{\frac2{p}}$ is an increasing function satisfying
\[
g(r) \leq (1-\epsilon) g(2r),
\]
for some $\e > 0$. The proof is straightforward.

Given a function $f \in L^2$ we note that the following holds true:
\[
\|f\|_{L^2}^2=\sum_{m \in \Z} \|f\|^2_{L^2(A_m)}. 
\]
Thus we have an $l^2_m$ structure for the sequence $(\|f\|^2_{L^2(A_m)})_{m \in Z}$. The following slightly stronger norm (improving the $l^2$ summation to an $l^1$ one)
\[
\|f\|_{l^1 L^2}=\sum_{m \in \Z} \|f\|_{L^2(A_m)},
\]
will play an important role in our analysis; one can simply define $l^1 L^2$ as the subspace of $L^2$ for which the above norm is finite. 

Just as above, the operator $[\partial_r]^{-1}$ is used to define the space $[\partial_r]^{-1}l^1L^2$, as the completion of  of $H^1_{comp}([0,\infty))$ with respect to the
following norm
\[
\| f \|_{[\partial_r]^{-1}l^1L^2} = \| \partial_r f \|_{l^1L^2} =\sum_{m} \| \partial_r f \|_{L^2(A_m)}.  
\]
Since $\| \partial_r f \|_{L^1(dr)} \lesssim \| f \|_{[\partial_r]^{-1}l^1L^2}$,  it follows that $f$
has limits both at $0$ and $\infty$; and since it is approximated
by functions in $H^1_{comp}([0,\infty))$, it follows that $\lim_{r \rightarrow \infty} f(r)=0$.
We also have the following inequality
\begin{equation} \label{l1l2dr}
\| \partial_r f \|_{L^2} + \| f \|_{L^\infty} \lesssim \| f \|_{[\partial_r]^{-1}l^1L^2}.
\end{equation}
We will also work with space $\dHe + [\partial_r]^{-1}l^1L^2$ and note that $\partial_r (\dHe + [\partial_r]^{-1}l^1L^2) = \partial_r \dHe + l^1L^2$.

We seek to understand a little better the structure of the space $[\partial_r]^{-1}(l^1L^2)$. Since $\partial_r f \in l^1 L^2$, it suffices to consider the case when $\partial_r f = \chi_{m} g \in L^2, g \in L^2$. Then we let
\[
f = \chi_{< m-4} c_m + \chi_{ \geq m-4} \int_r^\infty \chi_{m}  g(s) ds, \quad c_m = \int \chi_{m}  g(s) ds;
\]
 The first observation is that 
\[
|c_m| \les \|\chi_{m}  g\|_{L^2}. 
\]
The second component is obviously in $\dHe$, thus we have obtained that a generic function in $[\partial_r]^{-1}(l^1L^2)$ decomposes as follows:
\begin{equation} \label{drdech1}
f=f_1 + \sum_{m} c_m \chi_{\leq m}, \quad \|f_1\|_{\dHe} + \|(c_m)\|_{l^1} \les \| f\|_{[\partial_r]^{-1}(l^1L^2)}.  
\end{equation}

\section{An outline of the paper}

Due to the complexity of the paper, an overview of the ideas 
and the organization of the paper is helpful before an in-depth reading.

\subsection{ The frame method and the Coulomb gauge}
At first sight the Schr\"odinger Map equation has little to do with
the Schr\"odinger equation. A good way to bring in the Schr\"odinger
structure is by using the frame method. Precisely, at each point
$(x,t)\in \R^{2+1}$ one introduces an orthonormal frame $(v,w)$ in
$T_{u(x,t)} \S^2$. This frame is used to measure the derivatives of
$u$, and reexpress them as the complex valued radial {\em
  differentiated fields}
\[
\psi_1 =  \partial_r u \cdot v + i \partial_r u \cdot w, 
\qquad \psi_2 = \partial_\theta u \cdot v + i \partial_\theta u \cdot w.
\]
Here the use of polar coordinates is motivated by the equivariance
condition.  Thus instead of working with the equation for $u$, one
writes the evolution equations for the differentiated fields.
The frame $(v,w)$ does not appear directly there, but only via the
real valued radial connection coefficients
\[
A_1 = \partial_r v \cdot w, \qquad A_2 =  \partial_\theta v \cdot w,
\qquad A_0 =  \partial_t  v \cdot w.
\]

A-priori the frame is not uniquely determined. To fix it one first
asks that the frame be equivariant, and then that it satisfies an
appropriate condition.  Here it is convenient to use the {\em
  Coulomb gauge}; due to the equivariance this takes a
very simple form, $A_1 =0$.  The construction of the Coulomb gauge is
the first goal in the next section.  In Proposition~\ref{p:gauge} we
prove that for $\dot H^1$ equivariant maps into $\S^2$ close to
$Q$ there exists an unique Coulomb frame $(v,w)$ which
satisfies appropriate boundary conditions at infinity, see
\eqref{bcvw}. In addition, this frame has a $C^1$ dependence on the
map $u$ in a suitable topology. 

In the Coulomb gauge the other spatial connection
coefficient $A_2$, while nonzero, has a very simple form $A_2 = u_3$.
We also compute $A_0$ in terms of $\psi_1$, $\psi_2$ and 
$A_2$,
\begin{equation}
 A_0 = - \frac12 \left( |\psi_1|^2 - \frac{1}{r^2}|\psi_2|^2\right)
+  [r \partial_r^{-1}] \left( |\psi_1|^2 - \frac{1}{r^2}|\psi_2|^2\right).
\label{intro-a0bis}\end{equation}

\subsection{The reduced field \texorpdfstring{$\psi$}{}  }

 The two fields $\psi_1$ and $\psi_2$ are
not independent. Hence it is convenient to work with a single field
\[
\psi = \psi_1 - i r^{-1} \psi_2,
\]
which we will call {\em the reduced field}. The relevance of the
variable $\psi$ comes from the following reinterpretation.  If $\W$ is
defined as the vector
\[
\W = \partial_r u - \frac{m}{r} u \times R u \in T_u(\S^2),
\] 
then $\psi$ is the representation of $\W$ with respect to the frame $(v,w)$. 
On the other hand, a direct computation, see for instance \cite{GuKaTs-2}, leads to
\[
E(u) = \pi \int_0^\infty \left( |\partial_r \bar{u}|^2 + \frac{m^2}{r^2} |\bar{u} \times
R \bar{u}|^2 \right) rdr = \pi \| \bar \W \|_{L^2(rdr)}^2 + 4\pi |m|
\]
where we recall that $u(r,\theta)= e^{m \theta R} \bar{u}(r)$ and, similarly, 
let $\W(r,\theta)= e^{m \theta R} \bar \W(r)$.
Therefore $\psi=0$ is a complete characterization of $u$ being a harmonic map.
Moreover the mass of $\psi$ is directly related to the energy of $u$ via
\begin{equation} \label{basicpsi}
\|\psi\|_{L^2}^2 = \| \bar \W \|_{L^2}^2 = \frac{E(u)-4\pi |m|}{\pi}.
\end{equation}

A second goal of the next section is to derive an equation for
the time evolution of  $\psi$. This is governed by  a  cubic NLS type
equation, 
\begin{equation}
(i \partial_t + \Delta -\frac{2}{r^2} ) \psi =
\left( A_0 - 2 \frac{A_2}{r^2}-\frac{1}r \Im({\psi}_2 \bar{\psi})\right)\psi.
\label{intro-psieq}\end{equation}
In addition, we show that $\psi$ is connected back to
$(\psi_2,A_2)$  via the ODE system
\begin{equation}
\partial_r A_2= \Im{(\psi \bar{\psi}_2)}+\frac{1}r |\psi_2|^2, \qquad \partial_r
\psi_2 = i A_2 \psi - \frac{1}r A_2 \psi_2
\label{intro-comp}\end{equation}
with the conservation law $A_2^2 + |\psi_2|^2 = 1$. 

\subsection{The modulation parameters \texorpdfstring{$(\alpha,\lambda)$}{}}

The system \eqref{intro-psieq} is not self-contained,
as the reduced field $\psi$ does not uniquely determine $(\psi_2,A_2)$ via the system \eqref{intro-comp}. What is  missing is a suitable initial condition.  

The initial condition for \eqref{intro-comp} is closely related to a suitable choice of a "closest soliton" to the map $u$. This is denoted by $Q_{\alpha,\lambda}$, where the soliton parameters 
$(\alpha,\lambda)$ are viewed as \emph{modulation parameters} for our Schr\"odinger map. The modulation parameters are uniquely identified with a 
solution to the compatibility ode \eqref{intro-comp}
via a suitable orthogonality condition, see \eqref{ldef}. The one-to-one correspondence between
the modulation parameters $(\lambda,\alpha)$ 
and the solutions to \eqref{intro-comp} is established in Proposition~\ref{p:mod}.

\subsection{ Linearizations and the operators \texorpdfstring{$H$, $\tilde H$}{}} \label{outlinH}

A key role in our analysis is played by the linearization of the
Schr\"odinger Map equation around the soliton $Q$. A solution to the
linearized flow is a function
\[
u_{lin} : \R^{2+1} \to T_Q S^2.
\]
The Coulomb frame associated to $Q$ has the form 
\[
 v_Q(\theta,r) = e^{m \theta R} \bar v_Q(r), \qquad w_Q(\theta,r) = e^{m \theta R}
\bar w_Q(r)
\]
with 
\[
\bar v_Q(r) = 
\left( \begin{array}{ccc} h_3^m(r) \\ 0 \\ - h_1^m(r) \end{array} \right),
\qquad 
\bar w_Q(r) = 
\left( \begin{array}{ccc} 0 \\ 1 \\ 0 \end{array} \right).
\qquad 
\]
Expressing $u_{lin}$ in this frame,
\[
 \phi_{lin} = \la u_{lin},v_Q\ra + i \la u_{lin}, w_Q\ra
\]
one obtains the Schr\"odinger type equation
\begin{equation}
(i \partial_t - H) \phi_{lin} = 0 
\label{philin}\end{equation}
where the operator $H$ acting on radial functions has the form
\[
H = -\Delta + V, \qquad V(r) =\frac{m^2}{r^2}(1-2(h_1^m)^2).  
\]

On the other hand linearizing the equation \eqref{intro-psieq} around the
soliton $Q$, we obtain a linear Schr\"odinger equation of the form
\begin{equation}
(i \partial_t - \tilde H) \psi_{lin} = 0 
\label{psilin}\end{equation}
where the operator $\tilde H$ acting on radial functions has the
form
\[
\tilde H = -\Delta + \tilde V, \qquad \tilde V(r) = \frac{1+m^2-2mh_3^m}{r^2}. 
\]
The operators $H$ and $\tilde H$ are conjugate operators
and admit the factorizations
\begin{equation}\label{H-factor}
H =  L^* L, \qquad \tilde  H = L  L^* 
\end{equation}
where
\[
L = h_1^m \partial_r \frac{1}{h_1^m} = \partial_r + \frac{m}{r} h_3^m, \qquad
L^{*}=- \frac{1}{h_1^m} \partial_r h_1^m -\frac1{r}= -\partial_r +
\frac{m h_3^m-1}{r}.
\]
The linearized variables 
$\phi_{lin}$ and $\psi_{lin}$ are also conjugated variables, 
\begin{equation} \label{linL}
 \psi_{lin} = L\phi_{lin}.
\end{equation}

The operator $H$ is nonnegative and bounded from $\dot H^1$ to $\dot H^{-1}$,
but it is not positive definite; instead it has a zero mode
$\phi_0$, solving $L \phi_0=0$, namely 
\[
 \phi_0(r) = h_1^m.
\]
This corresponds to the solution $\phi_{lin}$ for \eqref{philin}
obtained by differentiating
the soliton family with respect to either parameter. A consequence 
of this is that the linear Schr\"odinger evolution \eqref{philin}
does not have good dispersive properties, a fact which  is at the heart of our
instability result.

If $|m|=1$ then the zero mode is a resonance, while if $|m| \geq 2$ then the zero resonance 
is replaced by a zero eigenvalue. If $|m| \geq 3$ then this eigenvalue belongs to $\dot H^{-1}$,
which allows for a clean splitting of the $\dot H^1$ space into 
an eigenvalue mode, which is stationary, and an orthogonal component,
which has good dispersive properties. This leads to the stability 
results in \cite{GuKaTs-2}, \cite{GuNaTs}. As already mentioned in the introduction, if $m=1$, it was shown in \cite{BeTa-1}
that instability occurs, and construction of solutions which blow-up in finite time were provided in \cite{MeRaRo} and \cite{Pe}. 
In this paper we address the case $|m|=2$ 
and show that instability still occurs. Hence, despite the fact that $h_1^2$
is an eigenvalue, the result for $|m|=2$ is closer to the $|m|=1$ case.

If $Q$ is replaced by $Q_{\alpha, \l}$ then $H$ and $\tilde H$ are
replaced by their rescaled versions $H_\lambda$ and $\tilde H_\lambda$
where $V$ and $\tilde V$ are replaced by
\begin{equation}\label{Vlambda}
V_\lambda = \lambda^2 V(\lambda r), \qquad \tilde V_\lambda =
\lambda^2 \tilde V(\lambda r).
\end{equation}
Correspondingly, the operator $L$ 
in the factorization \eqref{H-factor}
is replaced by 
\begin{equation}\label{Llambda}
L_\lambda = \partial_r + \frac{m}{r} h_3^m(\lambda r).
\end{equation}

From this point on our analysis becomes specialized to the case $m=2$; the case $m=-2$ is identical. Hence for the reminder of this section we assume that $m=2$.  

\subsection{Spectral theory for \texorpdfstring{$H$}{} and \texorpdfstring{$\tilde H$}{}.}
A first objective of  Section~\ref{spectral} is to describe the spectral theory
for the linear operators $H$ and $\tilde H$. The analysis in the case
of $H$ has already been done in \cite{KrScTa}, and it is easily obtained
via the $L$ conjugation in the case of $\tilde H$. The normalized
generalized eigenfunctions for $H$ and $\tilde H$ are denoted by
$\phi_\xi$, respectively $\psi_\xi$, and satisfy
\[
 H \phi_\xi = \xi^2 \phi_\xi, \qquad  \tilde H \psi_\xi = \xi^2 \psi_\xi,
\qquad L \phi_\xi = \xi \psi_\xi.
\]
Correspondingly we have a generalized Fourier transform $\F_H$
associated to $H$ and a generalized Fourier transform $\F_{\tilde H}$
associated to $\tilde H$.

This quickly leads to generalized eigenfunctions for the rescaled
operators $H_\lambda$ and $\tilde H_\lambda$. A considerable effort is
devoted to the study of the transition from one $\tilde H_\lambda$
frame to another. This is closely related to the transference operator
introduced in \cite{KrScTa}.

One reason we prefer to work with the $\psi$ variable is that 
the operator $\tilde H$ has a good spectral behavior at frequency zero,
therefore we have favorable dispersive decay estimates for the corresponding 
linear Schr\"odinger evolution \eqref{psilin}.

\subsection{ Stronger \texorpdfstring{$\ell^1$}{} Besov topologies: the spaces \texorpdfstring{ $X$, $\LX$}{}.} The bulk of our analysis is done for solutions  $u$ in the energy space $\dot H^1$, which corresponds to $\psi$ in $L^2$. But in these topologies the solitons 
turn out to be unstable, which motivates us to also
seek stronger topologies where the solitons are instead stable. In the $2$-equivariant setting studied  in this article, these spaces, denoted by $X$ for the map $u$, respectively by $\LX$ for the 
for the reduced field $\psi$, turn out to be simply 
the corresponding $\ell^1$ Besov type spaces, namely 
$X =  \dot B^{1}_{2,1}$, respectively $\LX =\dot B^{0}_{2,1,e}$ (see Section \ref{SBS-def} for precise defintions).   

These spaces and their properties are discussed in detail in Section~\ref{s:XLX}, where we establish the 
one-to-one elliptic correspondence between $u \in X$ 
and $\psi \in \LX$, see Proposition~\ref{uQpsi}.
In the same section we also establish equivalent
characterizations of $\bX$ and $\LX$ relative to the 
Littlewood-Paley decompositions associated to the operators $H$, respectively $\tH$. This in particular justifies our notations, by showing that
\[
L: \bX \to \LX,
\]
as a surjective map with a one dimensional kernel.

\subsection{ The linear \texorpdfstring{$\tH$}{} flow} \label{LtHf}
 This represents the main component in the evolution equation for the reduced field $\psi$.  In Section \ref{s:linear} we 
 study this linear flow, with the aim of proving local energy decay and Strichartz bounds.

 One added difficulty is that we need 
 to allow the scaling parameter $\lambda$ 
 to vary as a function of time. But an arbitrary time dependence of $\lambda$ cannot be allowed; instead, we assume that we have a bound on the scale invariant quantity
 \[
\Mf = \left\| \frac{\lambda'}{\lambda^2}\right\|_{L^2_t}.
 \]
This turns out to be the crucial quantity,
whose finiteness guarantees global, universal bounds
for the linear $\tH$ flow.

\subsection{ The nonlinear source term \texorpdfstring{$N(\psi)$}{}}
Once the modulation parameters $(\alpha,\lambda)$
have been selected using our orthogonality condition, the $\psi$ equation \eqref{intro-psieq}
can be recast in the form 
\begin{equation}
(i \partial_t -\tH) \psi = N(\psi).
\label{intro-psi-tH}\end{equation}
where $N(\psi)$ is at least quadratic and contains all nonlinear contributions.

In Section \ref{s:nolinear} we prove that $N(\psi)$
can be estimated perturbatively within the framework of the linear $\tilde H$ equation. This implies global, uniform Strichartz an local energy bounds for the reduced field $\psi$ associated to a Schr\"odinger map $u$ in terms of the initial data
size in $L^2$, under the sole assumption that we control the quantity $\Mf $ above.
In the same section we also show that for initial data $\psi(0)$ in the smaller space $\LX$, the $\ell^1$ Besov structure carries over to the solution
$\psi$, relative to the (time-dependent) Littlewood-Paley decomposition associated to $\tH_\lambda$.

\subsection{ The modulation equations}
Following our study of the $\psi$ equation, the second main step in the proof of our results is the 
study of the modulation equations for the modulation parameters $(\lambda(t),\alpha(t))$. 
This analysis begins in Section~\ref{s:mod},
where we derive the modulation equations as a nonlinear ode system with source terms which depend on $\psi$. These source terms can be divided as follows:

\begin{itemize}
    \item linear in $\psi$; these are the most troublesome ones, which at best can be estimated 
    in $L^2_t$ using the local energy decay bounds.

    \item quadratic and higher, which can be estimated in $L^1_t$ and thus play a perturbative role.
\end{itemize}
Nevertheless, the $L^2$ bounds turn out to be sufficient in order to close the bound for $\Mf $ under the sole assumption that $\psi(0)$ is small in $L^2$. This is a critical step in our analysis, as it implies that the bounds for $\psi$ are universal on the existence time for the solutions, 
even though by itself it does not preclude 
finite time blow-up.

A second, more refined step is undertaken 
in Section~\ref{s:mod2}, where we show that the 
linear in $\psi$ source terms can be in effect placed also  in the space $\dot H_t^{-\frac12}$,
for which perturbative analysis fails but only in a borderline fashion.

\subsection{ ODE analysis for the modulation parameters}
The final steps in our study of the modulation equation are carried out in Section~\ref{s:ode} and Section~\ref{s:final}. 

In a first step we begin by assuming a stronger Besov bound $\dot B^{-\frac12}_{2,1}$ 
for the source terms. In that case the source terms 
can be treated perturbatively, but the difficulty
is that we cannot assume smallness so a careful control of the constants is needed. This step suffices for the case when the initial data $\psi(0)$
is in the smaller space $\LX$.

In the second and final step, we consider
source terms which are small in $\dot H_t^{-\frac12}$
and also with small $\Mf $, which gives an $L^2$ bound at sufficiently high frequencies.
Over compact time intervals, such source terms 
may be placed in the stronger Besov space 
$\dot B^{-\frac12}_{2,1}$, but with a large norm. 
From here, a careful balancing of the time 
interval $[0,T]$ and of the scale of $\lambda$ in $[0,T]$ leads to the proof of the bound
\[
\lambda(t) \lesssim T^{C \delta^2}, \qquad \delta = \|\psi(0)\|_{L^2}
\]
which in turn implies our main result, asserting that no finite time blow-up is possible.

\section{The Coulomb gauge representation of the equation} 
\label{seccoulomb}

Our first goal in this section is to introduce the Coulomb gauge, which allows us to rewrite the Schr\"odinger map equation for equivariant solutions as a semilinear Schr\"odinger system 
for a differentiated field $\psi$, with a nonlinearity that contains no derivatives; this is an important feature that makes the nonlinear analysis much simpler.  However, the recovery of a Schr\"odinger Map state
from its Coulomb gauge representation $\psi$  is not unique, instead it retains two degrees of freedom. We fix these two degrees of freedom via an orthogonality 
condition relative to a  nearby reference soliton $Q_{\alpha,\lambda}$.
The soliton parameters $(\alpha,\lambda)$ will vary as a function 
a time, and will be referred to as \emph{modulation parameters}.
Later in the paper, we will use the orthogonality condition in order
to derive a set of ordinary differential equations, called the \emph{modulation equations}, for the two modulation parameters.
This will allow us to view the Schr\"odinger map evolution as a coupled system with two components:
\begin{enumerate}[label=(\roman*)]
    \item A differentiated field $\psi$ which solves a semilinear Schr\"odinger type equation, and 
    \item Two modulation parameters $(\alpha,\lambda)$ which solve 
    an appropriate ODE system.
\end{enumerate}
Given that we need to introduce and develop an array of notations along the way, it is difficult to summarize the full conclusion of this section here; instead,  we will do so at the end, see Section \ref{Sec3wrap}.

For reference, we note that the use of the Coulomb gauge in the context of the Schr\"odinger map equation originates in the work of Chang, Shatah, Uhlenbeck~\cite{ChShUh}. This strategy was particularly successful in later works on the Schr\"odinger map equation with data (hence solutions) that have radial or, more general, equivariant symmetries; see for instance 
\cite{BIKT-2, BeTa-1,GuKaTs-1,GuKaTs-2,GuNaTs}. The Coulomb gauge has been also used for general data (without symmetries) in \cite{bik} in high dimensions $n \geq 4$, but fails to be the efficient gauge in low dimensions $n=2,3$.

The gauge representation theory in this section is entirely based on the setup developed in \cite{BeTa-1} by the first and last authors, see Chapter 3 there. The orthogonality condition is also similar to the one we have used in \cite{BeTa-1}, but we should mention that it first appeared in the work of Gustafson, Nakanishi, Tsai~\cite{GuNaTs} in a different form, namely at the level of maps rather than gauge components. Indirectly  we are  using some of the analysis in \cite{GuNaTs}.

\subsection{The differentiated maps}

 We let the differentiation operators  $\partial_0,\partial_1,\partial_2$
stand for $\partial_t, \partial_r, \partial_\theta$ respectively. Our strategy 
will be to replace the equation for the Schr\"odinger map $u$ with 
equations for its derivatives $\partial_1 u$, $\partial_2 u$ 
expressed in an orthonormal frame $v,w \in T_u \S^2$. To fix the 
sign in the choice of $w$, we assume that
\[
u \times v = w.
\]
 Since $u$ is $m$-equivariant it is natural to work with $m$-equivariant frames, i.e. 
 \[
v = e^{m \theta R} \bar{v}(r), \qquad w = e^{m \theta R} \bar{w}(r).
\]
Given such a frame  we introduce
the differentiated fields $\psi_k$ and the connection coefficients
$A_k$ by
\begin{equation}
\label{connection}
\begin{split}
\psi_k = \partial_k u \cdot v + i \partial_k u \cdot w, \qquad
A_k = \partial_k v \cdot w, \quad k=0,1,2.
\end{split}
\end{equation}
Due to the equivariance of $(u,v,w)$ it follows that both $\psi_k$ and
$A_k$ are radially symmetric (therefore subject to the 
conventions made in Section \ref{defnot}). Conversely, given $\psi_k$ and $A_k$, we can return to the frame $(u,v,w)$ via the ODE system:
\begin{equation}
\label{return}
\left\{ \begin{array}{l}
\partial_k u = (\Re {\psi_k}) v  + (\Im{\psi_k}) w
\cr
\partial_k v = - (\Re{\psi_k}) u + A_k w
\cr
\partial_k w = - (\Im{\psi_k}) u - A_k v.
\end{array} \right.
\end{equation}

If we introduce the covariant differentiation
\[
D_k = \partial_k + i A_k, \ \ k \in \{0,1,2\},
\]
then it is a straightforward computation to check the compatibility conditions:
\begin{equation} \label{compat}
D_l \psi_k = D_k \psi_l, \ \ \  l,k=0,1,2.
\end{equation}
The curvature of this connection is given by 
\begin{equation}\label{curb}
  D_l D_k - D_k D_l = i(\partial_l A_k - \partial_k
A_l) = i \Im{(\psi_l \overline{\psi}_k)}, \ \ \  l,k=0,1,2.
\end{equation}
An important geometric feature is that $\psi_2, A_2$ are closely
related to the original map. Precisely, for $A_2$ we have:
\begin{equation}
A_2 = m (k \times v) \cdot w= m k \cdot (v \times w)= m k \cdot u = m u_3
\label{a2u3}\end{equation}
and, in a similar manner,
\begin{equation}
\psi_2 = m(w_3-iv_3).
\label{psi2vw3}\end{equation}
Since the $(u,v,w)$ frame is orthonormal, the following relations
also follow: 
\begin{equation}\label{psiA2}
|\psi_2|^2 = m^2(u_1^2 + u_2^2), \qquad
|\psi_2|^2 + A_2^2 = m^2.
\end{equation}

\subsection{The Coulomb gauge} \label{CG}
Now we turn our attention to the choice of the $(v, w)$ at $\theta = 0$, that is to $(\bar v,\bar w)$. Here we have the gauge freedom of an arbitrary rotation
depending on $t$ and $r$. Our gauge choice aims to remove this freedom.
In this article we will use the Coulomb gauge, which for general maps $u$ has the form 
\[
\text{div } A = 0.
\]
In polar coordinates this is written as 
\[
\frac{1}{r}\partial_{1}(r A_{1}) + \frac{1}{r^2}\partial_{2}A_{2} = 0.
\]
However, in the equivariant case $A_2$ is radial,
so we are left with a simpler formulation $A_1 = 0$, or equivalently
\begin{equation}
\partial_r \bar v \cdot \bar w=0
\label{coulomb}\end{equation}
which can be rearranged into a convenient ODE as follows
\begin{equation} \label{cgeq}
\partial_r  \bar v = (\bar v \cdot \bar u) \partial_r \bar u
- (\bar v \cdot \partial_r \bar u) \bar u.
\end{equation}
The first term on the right vanishes and could be omitted,
but it is convenient to add it so that the above linear ODE is solved 
not only by $v$ and $w$, but also by $u$. Then we can write an equation 
for the matrix $ \calO = (\bar v, \bar w,\bar u)$:
\begin{equation} \label{cgeq-m}
\partial_r \calO = M \calO, \qquad M = \partial_r \bar u \wedge \bar u : = 
\partial_r \bar u \otimes \bar u - \bar u \otimes \partial_r \bar u
\end{equation}
with an  antisymmetric matrix $M$.

The ODE \eqref{cgeq} needs to be initialized at some point.
 A change in the initialization leads to a multiplication
of all of the $\psi_k$ by a unit sized complex number. This is 
irrelevant at fixed time, but as the time varies we need to be careful 
and choose this initialization uniformly with respect to $t$, in 
order to avoid introducing a constant time dependent potential 
into the equations via $A_0$.  Since in our results we start with
data which converges asymptotically to $\vec{k}$ as $r \to \infty$,
and the solutions continue to have this property, it is natural to
fix the choice of $\bar{v}$ and $\bar{w}$ at infinity,
\begin{equation}
\label{bcvw}
 \lim_{r \to \infty} \bar v(r,t) = \vec{i} , \qquad \lim_{r \to \infty} \bar w(r,t) =
\vec{j}.
\end{equation}
The justification of the fact that we can impose these conditions at $\infty$ and find the gauge, or equivalently  uniquely solve \eqref{cgeq} with this boundary condition for $\bar{v}$, is provided in the proof of Theorem 3.2 in \cite{BeTa-1}. 

Before considering the general case we begin with the solitons.
The simplest case is when $u = Q^{m}$, when the triplet 
$(\bar v,\bar w,\bar u)$ is given by
\begin{equation}\label{vwq}
\left(  \bar V^{m}, \bar W^{m},\bar Q^{m}\right)
= 
  \left(\begin{array}{ccc}
   h_3^m(r) & 0 & h_1^m(r) \cr 0 & 1 & 0 \cr
  -h_1^m(r) & 0  & h_3^m(r)
\end{array}\right).
\end{equation}
More generally, if $ u = Q^{m}_{\alpha,\lambda}$ then from the above formula,
by rescaling and rotation, we obtain the corresponding triplet $
\left(\bar V^{m}_{\alpha,\lambda}, \bar W^{m}_{\alpha,\lambda},\bar
  Q^{m}_{\alpha,\lambda}\right)$ of the form
\[
 \!\left(\!\!\! \begin{array}{ccc}
\! h_3^m(\lambda r)\cos^2 m\alpha + \sin^2 m\alpha 
 & \! (h_3^m(\lambda r)-1) \sin m\alpha \cos m \alpha
& h_1^m(\lambda r) \cos m\alpha  \cr
 \! (h_3^m(\lambda r)-1) \sin m \alpha \cos m \alpha \!
& \!  h_3^m(\lambda r) \sin^2 m\alpha + \cos^2 m\alpha 
& h_1^m(\lambda r) \sin m\alpha  \cr
\! -h_1^m(\lambda r) \cos m\alpha & \!
- h_1^m(\lambda r) \sin m \alpha & h_3^m(\lambda r) \end{array} \!\!\! \right). \!
\]
For later reference we also note the values of $\psi_1$, $\psi_2$ and $A_2$ 
in this case:
\begin{equation} \label{solref}
\begin{split}
\psi_{\alpha,\lambda,1}^m = - m r^{-1} h_1^m(\lambda r) e^{im\alpha}, &\quad \psi_{\alpha,\lambda,2}^m= i m h_1^m(\lambda r) e^{im\alpha},
\\
 A_{\alpha,\lambda,2}^m =& m h_3^m(\lambda r).  
\end{split}
\end{equation}

In this article we work with maps $u:\R^2 \to \S^2$
which are near a soliton $Q^m_{\alpha,\lambda}$ in the sense that
\begin{equation}
\|u - Q_{\alpha,\lambda}^m\|_{\dot H^1} \leq \delta \ll 1.
\label{l1a}\end{equation}
The following result had been established in \cite{BeTa-1}.
\begin{l1}
 Let $m \geq 1$ and $u: \R^2 \to \S^2$ be an $m$-equivariant map
which satisfies \eqref{l1a}. Then
\begin{equation}
\lim_{r \to 0} u(r,\theta) = -\vec{k}, \qquad \lim_{r \to \infty} u(r,\theta) = \vec{k}
\label{l1b}\end{equation}
and 
\begin{equation}
\| r^{-1} ( u - Q_{\alpha,\lambda}^m)\|_{L^2}+  \| u - Q_{\alpha,\lambda}^m\|_{L^\infty} \lesssim \delta.
\label{l1cc}\end{equation}
\label{l:ns}\end{l1}

To measure the regularity of the frame $(v,w)$ we use the 
Sobolev type space $\dot H^1_C$ of functions $f: \R^2 \to \R^3$,
with norm
\[
\| f\|_{\dot H^1_C} = \| \partial_{r}  \bar f \|_{L^2} + 
\|\bar f\|_{L^\infty} + \|r^{-1} \bar f_3\|_{L^2}, \qquad f(r,\theta)=e^{m\theta R}\bar f(r)
\]
The next proposition, which was also established in \cite{BeTa-1}, shows that the initialization \eqref{bcvw}  is  well-defined for arbitrary maps $u$ close to the soliton family:
\begin{p1} \label{lc}
a) For each  $m$-equivariant map $u: \R^2 \to \S^2$ 
satisfying \eqref{l1a} there exists an unique $m$-equivariant orthonormal frame $(v,w)$ which satisfies the Coulomb gauge condition \eqref{coulomb} and the boundary condition \eqref{bcvw}. This frame satisfies the bounds
\begin{equation}
 \| v - V^m_{\alpha,\lambda}\|_{\dot H^1_C} + \| w- W^m_{\alpha,\lambda}\|_{ \dot H^1_C}\lesssim \delta.
\label{l1e}\end{equation}
b) Furthermore, the maps $u \to v,w$ are $C^1$ from $\dot{H^1}$ 
into $ \dot H^1_C$ as well as from $L^2 \to L^2$. 
\label{p:gauge}\end{p1}

As a direct consequence of part (a) of the above proposition,
we can describe the regularity and properties of the differentiated
fields $\psi_1$, $\psi_2$ and the connection coefficient $A_2$ at 
fixed time:

\begin{c1}
  Let  $u: \R^2 \to \S^2$ be an $m$-equivariant map
as in \eqref{l1a}. Then $\psi_1$, $\psi_2$ and $A_2$ satisfy
\eqref{compat}, \eqref{curb} for $k,l=1,2$
as well as the bounds
\[
 \| \psi_1- \psi_{\alpha,\lambda,1}^m\|_{L^2} + 
\|\psi_2 - \psi_{\alpha,\lambda,2}^m\|_{\dHe} +
\|A_2 - A_{\alpha,\lambda,2}^m\|_{\dHe} \lesssim \delta. 
\]
In addition, the map $u \to (\psi_1, \psi_2, A_2)$ from 
$\dot H^1$ into the above spaces is $C^1$.
\end{c1}

A second step is to consider Schr\"odinger maps with more regularity; this is particularly useful in order to justify formal computations. As a consequence of part (b) of Proposition \ref{lc} we 
have:

\begin{c1}\label{c:coulomb}
  Let $I$ be a compact interval, and  $u: I\times \R^2 \to \S^2$ be an $m$-equivariant map satisfying \eqref{l1a} uniformly in $I$ and which 
has the additional regularity
\[
u \in C(I;\dot H^2), \qquad \partial_t u \in C(I; L^2).
\]
Then $\psi_0$, $\psi_1$, $\psi_2$ and  $A_0$, $A_2$ satisfy
the relations \eqref{compat}, \eqref{curb} for $k,l=0,1,2$
and have the additional regularity
\begin{equation} \label{Cextra}
\begin{split}
&  \psi_0, A_0 \in C(I;L^2), \psi_1- \psi_{\alpha,\lambda,1}^m \in C(I;\dHe), \\
& \psi_2 - \psi^{m}_{\alpha,\lambda,2},
A_2 - A^{m}_{\alpha,\lambda,2} \in C(I;\dot \Hde).
\end{split}
\end{equation}
\end{c1}

This last result is almost identical to the corresponding one in \cite{BeTa-1}, except that we remove the (low frequency) hypothesis $u-Q^{m}_{\alpha,\lambda} \in L^2$ and adjust the conclusion in \eqref{Cextra} by removing the associated (low frequency) information $\psi_2 - \psi^{m}_{\alpha,\lambda,2}, A_2 - A^{m}_{\alpha,\lambda,2} \in C(I;L^2)$. 

In practice, in order to ensure the additional regularity in the corollary above (namely $u \in C(I;\dot H^2), \partial_t u \in C(I; L^2)$) for solutions to the Schr\"odinger map equation, it suffices to assume that $u \in C(I;\dot H^1 \cap  \dot H^2)$, since the second part follows from the equation and the $L^\infty$ embedding. 

\subsection{ Schr\"odinger maps in the Coulomb gauge}
\label{SMCG}
At this point we are ready  to write the evolution equations for the
differentiated fields $\psi_1$ and $\psi_2$ in \eqref{connection}
computed with respect to the Coulomb gauge.

Writing the Laplacian in polar coordinates, a direct computation
using the formulas \eqref{connection} shows that we can rewrite the
Schr\"odinger Map equation \eqref{SM} in the form
\[
\psi_0 = i \left(D_1 \psi_1 + \frac{1}{r} \psi_1 + \frac{1}{r^2} D_2 \psi_2\right).
\]
Applying the operators $D_1$ and $D_2$ to both sides of this equation and using the relations \eqref{compat} and \eqref{curb}, we can derive the evolution equations for $\psi_m$, $ m=1,2$:
\begin{equation}\label{eqn:psi12}
\begin{split}
\partial_t \psi_1 + i A_0 \psi_1 = & \ i \Delta \psi_1  - 2  A_1 \partial_1
\psi_1  - \partial_1 A_1 \psi_1 - \frac{1}{r} A_1 \psi_1 \\
& \ - i A_1^2 \psi_1 - i \frac{1}{r^2} A_2^2 \psi_1  - i \frac1{r^2} \psi_1
+ \frac2{r^3} A_2 \psi_2 - \frac{1}{r^2} \Im{(\psi_1 \bar{\psi}_2)} \psi_2, \\
\partial_t \psi_2 + i A_0 \psi_2 = & \ i \Delta \psi_2  - 2 A_1 \partial_1
\psi_2  - \partial_1 A_1 \psi_2 - \frac{1}{r} A_1 \psi_2 \\
& \ - i A_1^2 \psi_2 - i \frac{1}{r^2} A_2^2 \psi_2 - \Im{(\psi_2 \bar{\psi}_1)}
\psi_1.
\end{split}
\end{equation}
Under the Coulomb gauge $A_1 = 0$ these  equations become
\begin{equation*} \label{smg}
\begin{split}
\partial_t \psi_1 + i A_0 \psi_1 = & i \Delta \psi_1  - i \frac{1}{r^2} A_2^2
\psi_1  -i  \frac1{r^2} \psi_1
+  \frac2{r^3} A_2 \psi_2 - \frac{1}{r^2} \Im{(\psi_1 \bar{\psi}_2)} \psi_2, \\
\partial_t \psi_2 + i A_0 \psi_2 = & i \Delta \psi_2  - i \frac{1}{r^2} A_2^2
\psi_2 - \Im{(\psi_2 \bar{\psi}_1)} \psi_1.
\end{split}
\end{equation*}
while the relations  \eqref{compat} and \eqref{curb} become
an ode system for $(A_2,\psi_2)$, namely
\begin{equation} \label{comp}
\partial_r A_2= \Im{(\psi_1 \bar{\psi}_2)}, \qquad \partial_r \psi_2 = i A_2 \cdot
\psi_1.
\end{equation}
On the other hand from the compatibility relations involving $A_0$ we obtain
\begin{equation} \label{a0}
\partial_r A_0= - \frac1{2r^2} \partial_r (r^2 |\psi_1|^2 - |\psi_2|^2),
\end{equation}
which determines $A_0$ modulo constants.
This allows us to derive an expression for $A_0$,
\begin{equation}
 A_0 = - \frac12 \left( |\psi_1|^2 - \frac{1}{r^2}|\psi_2|^2\right)
+ [r\partial_r]^{-1} \left( |\psi_1|^2 - \frac{1}{r^2}|\psi_2|^2\right) 
\label{a0bis}\end{equation}
where we recall the definition of  $[r \partial_r]^{-1}f(r)$
from \eqref{drdef}. To justify \eqref{a0bis}, we are using the fact that $[r \partial_r]^{-1}$ is bounded on $L^1$ and that $|\psi_1|^2 - \frac{1}{r^2}|\psi_2|^2 \in L^1$; this guarantees that the right-hand side of \eqref{a0bis} belongs to $L^1$ (and also to 
$r^{-1} L^\infty$). To dispense with the additional possible 
constant in $A_0$ we simply note that an additional regularity 
assumption on the map $u$ in this problem gives us that $A_0 \in L^2$, see Corollary \ref{c:coulomb}. 
Since a constant is not integrable in any sense, it follows
that the only solution of \eqref{a0} which has some decay at $\infty$ must be the one in \eqref{a0bis}.

There is quite a bit of redundancy in the equations for $\psi_1$ and $\psi_2$; we eliminate this by introducing a single \emph{differentiated field} $\psi$ by
\begin{equation}\label{def-psi}
 \psi=\psi_1 -i \frac{\psi_2}{r}.
\end{equation}
The size of $\psi$  is closely related to the energy of the original map $u$, precisely 
\begin{equation}
\pi \|\psi\|_{L^2}^2 =   E(u) - 4\pi m,  
\end{equation}
and in particular $\psi$ vanishes iff $u$ is a soliton.

A direct computation yields the Schr\"odinger type equation for $\psi$:
\[
i \partial_t \psi  +  \Delta \psi = A_0 \psi - 2 \frac{A_2}{r^2} \psi +
\frac{1}{r^2} \psi + \frac{A_2^2}{r^2} \psi - \frac1{r} \Im{(\psi_2
\bar{\psi}_1)} \psi.
\]
By replacing $\psi_1 = \psi + i r^{-1} \psi_2$ and using $A_2^2 + |\psi_2|^2=m^2$,
 we obtain the key evolution equation we work with in this paper,
 namely
\begin{equation} \label{psieq}
i \partial_t \psi  +  \Delta \psi - \frac{1+m^2}{r^2} \psi = A_0 \psi - 2
\frac{A_2}{r^2} \psi   - \frac1{r}\Im{(\psi_2 \bar{\psi})} \psi.
\end{equation}
Our strategy will be to use this equation in order to obtain estimates for
$\psi$. The functions $A_2$ and $\psi_2$ are related to $\psi$
via the system of ODE's
\begin{equation}
 \label{comp1}
 \left\{ \begin{array}{l}
\partial_r A_2= \Im{(\psi \bar{\psi}_2)}+\dfrac{1}r |\psi_2|^2, \cr \cr
\partial_r \psi_2 = i A_2 \psi - \dfrac{1}r A_2 \psi_2
\end{array} \right.
\end{equation}
which is derived from \eqref{comp}, together with the compatibility condition from \eqref{psiA2}, 
\begin{equation}
A_2^2 + |\psi_2|^2 = m^2.
\label{sphere}\end{equation}
 However, they are not uniquely determined by $\psi$, instead we have two additional degree of freedom,  which correspond to prescribing initial data in the above ode subject to the constraint provided by \eqref{sphere}. This 
will lead us to the modulation parameters described in the next subsection.

Once  $(A_2,\psi_2)$ are given, the $A_0$
connection coefficient is given by \eqref{a0bis} which now becomes
\begin{equation} \label{aoef}
 A_0(r) =  -\frac12 |\psi|^2  + \frac{1}r \Im (\psi_2 \bar \psi)
+ [r\partial_r]^{-1}( |\psi|^2 - \frac{2}{r} \Im (\psi_2 \bar \psi)).
\end{equation}
Finally, given $\psi$, $A_2$ and 
$\psi_2$, we can return to the  Schr\"odinger map $u$
via the system \eqref{return} with the boundary condition at infinity
given by \eqref{bcvw}.

\subsection{The modulation parameters \texorpdfstring{$\lambda(t), \alpha(t)$}{} } 
\label{s:modulation-param}

Our paper is concerned with the behaviour of maps near the $Q^{\pm 2}$ soliton (this is meant to include the more general objects $Q^{\pm 2}_{\alpha, \lambda}$).
The cases $m = 2$ and $m=-2$ can be identified via 
a reflection, therefore, beginning with this subsection, we simply set $m=2$.

This allows us to drop 
the upper script $m$ from $h_1^m$ and $h_3^m$ and simply use $h_1,h_3$. This allows us to introduce
another upper script convention
\[
h_1^\l(r)=h_1(\l r),\qquad h_3^\l(r)=h_3(\l r),
\] 
which is very useful due to the key role the scaling parameter $\l$ plays in
our analysis.

Given a $2$-equivariant map $u:\R^2 \to \S^2$ in the homotopy class of
$\mathcal Q^2_e$ and satisfying  $0 \leq E(u) - E(\mathcal Q^2) \ll 1$, we recall that, by \eqref{stab}, $u$ must be close in $H^1$ to one of the solitons $Q^2_{\alpha,\lambda}$:
\[
\text{dist}_{\dot H^1} (u,\mathcal Q^2_e)^2= \inf_{\alpha,\lambda}
\| Q^2_{\alpha,\lambda} - u\|_{\dot H^1}^2 \lesssim E(u) - E(\mathcal Q^2)= E(u) - 8\pi.
\]
 For such  $u$, it is important to identify a specific soliton $Q^2_{\alpha,\lambda}$ to be the "closest" soliton to $u$.
 We will think of the chosen $(\alpha,\lambda)$ as the modulation parameters associated to $u$, and later on we will study their evolution as $u$ evolves along the Schr\"odinger map flow. At the same time,
 we want to use $(\alpha,\lambda)$ to remove the two remaining degrees 
 of freedom in the choice of $(\psi_2,A_2)$ in the previous section.

 While the soliton $Q^2_{\alpha,\lambda}$ does not need to be the minimizer of the above distance, a natural condition to impose is that 
\[
\| Q^2_{\alpha,\lambda} - u\|_{\dot H^1}^2 \lesssim E(u) - 8\pi.
\]
This still leaves us with infinitely many choices for the parameters $\alpha$ and $\lambda$ and 
it is important to seek an efficient rule that selects a unique choice for the two parameters. We will select our parameters based on the following orthogonality condition
\begin{equation} \label{ldef}
\la \psi_2-2ie^{2 i \alpha(t)} h_1(\lambda r), \varkappa(\lambda r) \ra = 0,
\end{equation}
where $\varkappa$ can be subject to various choices. Next we discuss how these choices were made in prior works on this problem, and then what is our choice here.

One natural choice would have been to choose $\varkappa=h_1$ (or better $h_1^2$ in the more general context of $m$-equivariant setup which we discuss here), which is the eigenvalue of the linearized operator $H_\lambda$. Such a choice (that is, setting $\varkappa=h_1^m$) works well for $m \geq 4$, see \cite{GuKaTs-2}.
However the finite energy bound for $u$ corresponds to 
$r^{-1}(\psi_2-2ie^{2 i \alpha(t)} h_1(\lambda r)) \in L^2$, therefore in order to make sense of \eqref{ldef} we would need $r \varkappa \in L^2$. But $r h_1^m \notin L^2$ when $m=2$, therefore the choice $\varkappa=h_1^2$ is unsuitable for us. 

Another choice would be a point-type condition by simply setting $\varkappa=\delta_{r=1}$. This has been used in the case $m=1$ by the first and third author, see \cite{BeTa-1}. However there are several elements of the analysis in \cite{BeTa-1} that are not reproducible in the case $m=2$, and this is why, although the point condition will be partially used in our analysis, it will not be the one we ultimately choose to determine our parameters $\alpha$ and $\lambda$. 

Our choice here (that is in our context with $m=2$) is similar to the one used by Gustafson, Nakanishi, Tsai~\cite{GuNaTs}: $\varkappa$ is a smooth function, compactly supported in $[ \frac12 ,2 ]$  and subject to the following two nondegeneracy conditions:
\begin{equation}\label{chi-nondeg}
    \la \varkappa, h_1 \ra \ne 0, \qquad \la \varkappa, h_1 h_3 \ra \ne 0.
\end{equation}
Thus there is a lot of freedom in choosing $\varkappa$. However we note that the actual orthogonality condition used in \cite{GuNaTs} is different than ours - in \cite{GuNaTs} the orthogonality involves a linearization of the actual map, while our condition \eqref{ldef} involves the linearization of the gauge components. 

To keep formulas shorter in the remaining of this section and the rest of the paper, we introduce the notation
\[
\delta^{\lambda,\alpha} \psi_2= \psi_2 - 2ie^{2 i \alpha} h_1^\lambda, \quad \delta^\lambda A_2 = A_2 - 2h_3^\lambda. \]

The first issue that needs to be addressed is the existence 
of modulation parameters $\alpha$ and $\lambda$ satisfying the condition \eqref{ldef}. 
For this purpose we prove the following

\begin{p1}[Modulation parameters]
\label{p:mod}
Assume that $\varkappa$ satisfies the nondegeneracy conditions \eqref{chi-nondeg}. Then the following properties hold for $\delta > 0$ small enough:

i) Given any 2-equivariant map $u: \R^2 \rightarrow \S^2$, $u \in \dot H^1$ satisfying \eqref{eqvl} and
such that  $E(u)-8\pi \leq \delta^2$,  there exist a unique pair of modulation parameters 
$(\alpha,\lambda) \in [0,2\pi) \times \R^+$  such that 
\begin{equation} \label{lell}
\la \delta^{\lambda,\alpha} \psi_2, \varkappa^\lambda \ra = 0
\end{equation}
and
\begin{equation} \label{lellf}
\|\delta^{\lambda,\alpha} \psi_2\|_{L^\infty} \les \delta.
\end{equation} 
Further, $(\alpha,\lambda)$ have a Lipschitz dependence  on $\psi_2$ in the following sense:
\begin{equation} \label{Liplp}
|\ln \lambda(\psi_2) - \ln \lambda(\tilde \psi_2)| + |\alpha(\psi_2)-\alpha(\tilde \psi_2)| \les \min \{\lambda_{max}^{1-s} \|\psi_2 - \tilde \psi_2 \|_{\dot H^s_e}; s=1,0,-1\},
\end{equation}
where $\lambda_{max}=\max(\lambda(\psi_2), \lambda(\tilde \psi_2))$. 

ii) Assume that the modulation parameters $(\alpha,\lambda)$ are chosen as in \eqref{lell} and \eqref{lellf}. Then the following holds true
\begin{equation} 
\| \delta^{\lambda,\alpha} \psi_2 \|_{\dHe} + \| \delta^\lambda A_2\|_{\dHe} \lesssim  \| \psi \|_{L^2}.  
\label{compsolh1-rr}\end{equation}
As a consequence these modulation parameters are a good choice in the sense of \eqref{goodal}, that is
\begin{equation} \label{goodal2}
\| u - Q_{\alpha,\lambda}^2\|_{\dot H^1} \lesssim \|\psi\|_{L^2}.
\end{equation}
Furthermore, under additional hypothesis, we have the following estimates

\begin{itemize}
\item  if $\psi \in L^4$ (not necessarily small) then the following holds true
\begin{equation} 
\| \frac{\delta^{\lambda,\alpha} \psi_2}{r} \|_{L^4} + \| \frac{\delta^\lambda A_2}{r}\|_{L^4} \lesssim  \| \psi \|_{L^4}.  
\label{psiL4}\end{equation}

\item  if $\dfrac{\psi}r \in L^2$ (not necessarily small) then the following holds true
\begin{equation} 
\| \frac{\delta^{\lambda,\alpha} \psi_2}{r^2} \|_{L^2} + \| \frac{\delta^\lambda A_2}{r^2}\|_{L^2} \lesssim  \| \frac{\psi}r \|_{L^2}.  
\label{psiLE2}\end{equation}
\end{itemize}
\end{p1}

\begin{r1}
The pointwise smallness condition in \eqref{lellf} is needed in order to guarantee that we choose the correct pair $(\alpha,\lambda)$, in particular they satisfy \eqref{goodal2}. 
\end{r1}

 We also have a natural converse for the proposition above:

\begin{p1}[Modulation parameters - part 2]
\label{p:modb}
Assume that $\varkappa$ satisfies the nondegeneracy conditions \eqref{chi-nondeg} and $\delta > 0$ is small enough. Then given any radial map $\psi \in L^2$ 
such that  $\|\psi\|_{L^2}  \leq \delta$,  and any pair of modulation parameters 
$(\alpha,\lambda) \in [0,2\pi) \times \R^+$,
there exists a unique solution $(\psi_2,A_2)$ for the system \eqref{comp} so that both the
orthogonality relation \eqref{lell} and the 
bound \eqref{lellf} hold. In addition there exists an 
unique $2$-equivariant map $u: \R^2 \rightarrow \S^2$, $u \in \dot H^1$ satisfying \eqref{eqvl} and $E(u)-8\pi \leq \delta^2$, and with the property that $\psi, \psi_2, A_2$ are its Coulomb gauge components as constructed in Section \ref{SMCG} and $\alpha,\lambda$ are its correct modulation parameters as described in Proposition \ref{p:mod}.  
\end{p1}

\begin{proof} 
[Proof of Proposition~\ref{p:mod}]
i) Following \cite{GuNaTs}, the strategy here is to actually establish the existence of the parameters subject to the point-type condition (when $\varkappa=\delta_{r=1}$) and then use the inverse function theorem to find the parameters for the smooth choice of $\varkappa$. In addition to the intrinsic technical advantages of this approach, this will also provide us with an easy way to prove the estimates claimed in ii). 

Given the above strategy, the first step is to observe that we can find parameters $\alpha_0,\lambda_0$ so that
\[
\psi_2(\lambda_0^{-1}) =2ie^{2i\alpha_0} h_1(1)=2ie^{2i\alpha_0}.
\]
To see this, we note that $\lim_{r \rightarrow 0} A_2(r)=-2$ and $\lim_{r \rightarrow \infty} A_2(r)=2$ and the continuity of $A_2$ imply that there exists $\lambda_0$ with $A_2(\lambda_0^{-1})=0$. Then the  compatibility relation \eqref{psiA2} implies that 
$|\psi_2(\lambda_0^{-1})|=2$, from which it follows that there exists $\alpha_0$ such that $\psi_2(\lambda_0^{-1}) =2ie^{2i\alpha_0}$. 

The next step is to examine $\delta^{\lambda, \alpha} \psi_2$ and $\delta^{ \lambda} A_2$ corresponding to this pointwise choice of $\alpha$ and $\lambda$.  For further use, we phrase the result in a slightly more general fashion:

\begin{l1} \label{psito2} Consider the system \eqref{comp1} with $\psi \in L^2$, small, and initial data
\begin{equation} \label{ldef-r}
\delta^{\lambda,\alpha} \psi_2(r_0)=0, \qquad \delta^{\lambda} A_2(r_0)=0.
\end{equation}
for some $r_0 \approx \lambda^{-1}$. Then the system admits a unique solution $(\psi_2,A_2)$ with regularity $(\delta^{\lambda,\alpha} \psi_2,  \delta^\lambda A_2) \in \dHe \times \dHe$  and which satisfies \eqref{sphere}. 
Furthermore, this solution satisfies the following $\dHe$ bound:
\begin{equation} 
\| \delta^{\lambda,\alpha} \psi_2 \|_{\dHe} + \| \delta^\lambda A_2\|_{\dHe} \lesssim  \| \psi \|_{L^2}.  
\label{compsolh1-r}\end{equation}

If in addition $\psi \in L^4$ (not necessarily small) then the following holds true:
\begin{equation} 
\| \frac{\delta^{\lambda,\alpha} \psi_2}{r} \|_{L^4} + \| \frac{\delta^\lambda A_2}{r}\|_{L^4} \lesssim  \| \psi \|_{L^4}.  
\label{psiL44}\end{equation}

If in addition $\frac{\psi}r \in L^2$ (not necessarily small) then the following holds true:
\begin{equation} 
\| \frac{\delta^{\lambda,\alpha} \psi_2}{r^2} \|_{L^2} + \| \frac{\delta^\lambda A_2}{r^2}\|_{L^2} \lesssim  \| \frac{\psi}r \|_{L^2}.  
\label{psiLE22}\end{equation}

\end{l1}

\begin{proof} We first establish \eqref{compsolh1-r}. Using the compatibility condition \eqref{sphere}, the equations \eqref{comp1} become 
\begin{equation}
 \label{sysAp2}
\left\{ \begin{aligned}
L_\lambda \delta^{\lambda,\alpha} \psi_2  = & \ 2i h_3^\lambda \psi + \delta^\lambda A_2 \psi - \frac{1}r \delta^\lambda A_2 (2ie^{2 i \alpha}h_1^\lambda + \delta^{\lambda,\alpha} \psi_2), \\
L_\lambda \delta^\lambda A_2  = & \ - 2 h_1^\lambda \Re{(e^{2i\alpha} \psi)} + \Im(\psi \overline{\delta^{\lambda,\alpha} \psi}_2) - \frac1r (\delta^\lambda A_2)^2.
\end{aligned} \right.
\end{equation}
with $L_\lambda$ as in \eqref{Llambda}.
The solution to the homogeneous $L_\lambda$ equation is given by $h_1^\lambda$, so the equations above can be rewritten in the integral form:
\[
\begin{split}
& \delta^{\lambda,\alpha} \psi_2(r) = h^\lambda_1(r) \int_{r_0}^r \frac1{h_1^\lambda}\left( 2i h_3^\lambda \psi + \delta^{\lambda} A_2 \psi - \frac{1}s \delta^{\lambda} A_2 (2ie^{2i\alpha}h_1^\lambda + \delta^{\lambda,\alpha} \psi_2) \right) ds \\
& \delta^\lambda A_2 =  h^\lambda_1(r) \int_{r_0}^r \frac1{h_1^\lambda}\left( - 2 h^\lambda_1 \Re{ (e^{2i\alpha} \psi)} + \Im(\psi \overline{\delta^{\lambda,\alpha} \psi}_2) - \frac{(\delta^\lambda A_2)^2}s \right) ds. 
\end{split}
\]
As a tool to obtain bounds on the solutions, we record the following simple inequality,
\begin{equation} \label{Linvest}
\| \frac1r h_1^\lambda \int_{r_0}^r \frac1{h_1^\lambda} f ds \|_{L^p} \les \|f\|_{L^p}, \qquad  1 \leq p \leq \infty,
\end{equation}
which in the particular case $p = 2$ can be augmented to 
\begin{equation} \label{Linvest+}
\| \frac1r h_1^\lambda \int_{r_0}^r \frac1{h_1^\lambda} f ds \|_{L^2} 
+ \| h_1^\lambda \int_{r_0}^r \frac1{h_1^\lambda} f ds \|_{L^\infty}
\les \|f\|_{L^2}, \qquad  1 \leq p \leq \infty.
\end{equation}
The only role that $r_0$ plays here is that is sits at the "peak" of $h_1^\lambda$, that is $h_1^\lambda(r_0) \approx 1$ and $h_1^\lambda $ decays away from a neighborhood of size $\lambda^{-1}$ of $r_0$.

Under the assumption that $\|\psi\|_{L^2} \ll 1$,
the above inequality allows us to obtain the solutions $ (\delta^{\lambda,\alpha} \psi_2,\delta^\lambda A_2)$ to the above integral system using the 
contraction principle in the space $Y$ with norm
\[
\|f\|_{Y} = \| f \|_{L^\infty} + \| r^{-1} f\|_{L^2}.
\]
Returning to \eqref{sysAp2}, the $L^2$ bound for 
$\partial_r (\delta^{\lambda,\alpha} \psi_2,\delta^\lambda A^2)$ also follows, completing the proof of the $\dot H^1_e$ bound in the lemma. The uniqueness of the (small) solution $(\delta^{\lambda,\alpha} \psi_2 \dHe \times \dHe, \delta^\lambda A_2)$ in $\dHe \times \dHe$ follows from the above fixed point argument.

The argument for \eqref{psiL44} is entirely similar: it uses \eqref{Linvest} with $p=4$ and the already established smallness of $\|\delta^{\lambda,\alpha} \psi_2\|_{L^\infty}+\|\delta^\lambda A_2\|_{L^\infty}$. 

To obtain \eqref{psiLE22}, we rely instead on the following estimate in the special case $p = 2$,
\begin{equation} \label{Linvest1r}
\| \frac1{r^2} h_1^\lambda \int_{r_0}^r \frac1{h_1^\lambda} f ds \|_{L^p} \les \|\frac{f}r\|_{L^p}, \qquad 1 \leq p < \infty,
\end{equation}
and proceed just as above.

\end{proof}

Now we return to our initial guess $(\alpha_0,\lambda_0)$. 
From Lemma~\ref{psito2} above it follows that the functions $(\psi_2,A_2)$ 
satisfy
\begin{equation}\label{dpsi0}
\| \psi_2 -  2ie^{2 i \alpha_0} h_1^{\lambda_0} \|_{\dHe} + \| A_2 - h_3^{\lambda_0}\|_{\dHe} \les \| \psi \|_{L^2}.
\end{equation}
We use this as as starting point in order to establish the existence of 
modulation parameters satisfying the orthogonality condition \eqref{lell} with the smooth choice of $\varkappa$. Here it is convenient to work on an exponential scale for $\lambda$, so we set $\lambda = e^{\gamma} \lambda_0$.
Then it is natural to consider the following function
\[
F(\psi_2,\alpha,\gamma)= \lambda^2 \la \psi_2-2ie^{2 i \alpha} h_1(\lambda r), \varkappa(\lambda r) \ra.
\]
and seek a zero of $F$ near $(\alpha_0,0)$. From \eqref{dpsi0} it follows that 
\[
\| \frac{\psi_2-2ie^{2 i \alpha_0} h_1(\lambda_0 r)} r\|_{L^2} \les \|\psi\|_{L^2},
\]
and this implies that 
\[
|F(\psi_2,\alpha_0,0)| \les \| \psi\|_{L^2} \les \delta.
\]
Further, $F$ is uniformly $C^2$ with respect to $\alpha,\gamma$ and 
\[
D_{\alpha,\gamma}F|_{\alpha_0,0} = 4 e^{2i\alpha_0} ( \la h_1, \chi \ra d \alpha - i  \la h_1 h_3, \varkappa \ra d \gamma) + O(\delta d\gamma). 
\]
Assuming that both $\la h_1, \varkappa \ra \ne 0$ and $\la h_1 h_3, \varkappa \ra \ne 0$ and that $\delta$ is small enough, by the inverse function theorem 
the function $F$ has a unique zero in a small neighbourhood of $(\alpha_0,0)$, so that 
\[
|\alpha - \alpha_0| + |\sigma| \les \| \psi \|_{L^2} \les \delta.
\]
This in turn implies both the pointwise smallness condition
\[
\|\psi_2-2ie^{2 i \alpha} h_1(\lambda r)\|_{L^\infty} \les \| \psi\|_{L^2} \les \delta,
\]
and the bound \eqref{compsolh1-rr}
as a consequence of the similar bound for $(\alpha_0,\lambda_0)$. 

To obtain the slightly more general uniqueness result in the proposition,
we observe that if we had another solution $(\alpha_1,\lambda_1)$ which 
also satisfies the second bound in \eqref{lell} then we must have 
\[
\| e^{2i\alpha} h_1^\lambda - e^{2i\alpha_1} h_1^{\lambda_1}\|_{L^\infty}
\lesssim \delta,
\]
which in turn implies that 
\[
|\alpha - \alpha_1| + |\ln \lambda - \ln \lambda_1| \lesssim \delta,
\]
which places $(\alpha_1,\lambda_1)$ within the range of applicability 
of the previous uniqueness statement.

\bigskip

Next we consider the dependence of $(\alpha,\gamma)$ with respect to $\psi_2$. For this we note the following dependence of $F$ with respect to $\psi_2$, 
\[
|F(\psi_2,\alpha,\gamma)- F(\tilde \psi_2,\alpha,\gamma)| \les \| \psi_2 - \tilde \psi_2\|_{\dHe}, 
\] 
This implies that we have Lipschitz dependence of the parameters $\alpha$ and $\gamma$ with respect to the $\dot H^1_e$ norm of $\psi_2$, just as claimed in \eqref{Liplp}. In addition we can establish a similar Lipschitz dependence with respect to rougher norms of $\psi_2$ - this is possible because we test with a smooth $\chi$. A straightforward computation shows that
\[
|F(\psi_2,\alpha,\gamma)- F(\tilde \psi_2,\alpha,\gamma)| \les \lambda \| \psi_2 - \tilde \psi_2\|_{L^2}, 
\]
and this establishes the second part of the claim in \eqref{Liplp}.

Finally the last part of \eqref{Liplp} follows from the straightforward inequality 
\[
|F(\psi_2,\alpha,\gamma)- F(\tilde \psi_2,\alpha,\gamma)|   \les \lambda^2 \| \psi_2 - \tilde \psi_2\|_{\dot H^{-1}_e}. 
\]

ii) Lemma \ref{psito2} above provides the desired estimates when the choice of parameters is made with the point-type condition versus the smooth one that we work with. Thus our goal here is to prove that the correction coming from the change in parameters can still be controlled in the same way. By definition,  
\[
F(\psi_2,\alpha_0,\lambda_0) =  \lambda_0^2 \la \psi_2-2ie^{2 i \alpha_0} h_1(\lambda_0 r), \varkappa(\lambda_0 r) \ra. 
\]
Here $(\psi_2-2ie^{2 i \alpha_0} h_1^{\lambda_0})(\lambda_0^{-1})=0$
while the inner product depends only on $\psi_2$ in the region $\{ r \approx \lambda_0^{-1}\}$. To determine $\psi_2$ in this region from $\psi$ via the ode \eqref{comp1}, it suffices to know $\psi$ in the same region. By a local ode stability analysis it follows that
\[
\| \psi_2-2ie^{2 i \alpha_0} h_1(\lambda_0 r)\|_{L^\infty(r \approx \lambda_0^{-1})} \lesssim \| \psi \|_{L^2(r \approx \lambda^{-1})}.
\]
This implies that
\[
|F(\psi_2,\alpha_0,\lambda_0)| \les \| \psi \|_{L^2(r \approx \lambda^{-1})}. 
\]
from which the inverse function theorem used above gives the improved estimate
\[
|\alpha-\alpha_0| + |\ln \lambda - \ln \lambda_0| \les  \| \psi \|_{L^2(r \approx \lambda^{-1})}. 
\]
From this last bound it follows that
\[
\| \frac{\psi_{\alpha,\lambda,2} -\psi_{\alpha_0,\lambda_0,2}}{r} \|_{L^2} 
+ \| \frac{A_{\alpha,\lambda,2} -A_{\alpha_0,\lambda_0,2}}{r} \|_{L^2} \les \| \psi \|_{L^2},
\]
\[
\| \frac{\psi_{\alpha,\lambda,2} -\psi_{\alpha_0,\lambda_0,2}}{r} \|_{L^4} 
+ \| \frac{A_{\alpha,\lambda,2} -A_{\alpha_0,\lambda_0,2}}{r} \|_{L^4} \les \| \psi \|_{L^4}
\]
and 
\[
\| \frac{\psi_{\alpha,\lambda,2} -\psi_{\alpha_0,\lambda_0,2}}{r^2} \|_{L^2} 
+ \| \frac{A_{\alpha,\lambda,2} -A_{\alpha_0,\lambda_0,2}}{r^2} \|_{L^2} \les \| \frac{\psi}r \|_{L^2}.
\]
Combined with \eqref{psiL44} and \eqref{psiLE22}, this finishes the proof of the two claims \eqref{psiL4} and \eqref{psiLE2}.

The first inequality combined with \eqref{compsolh1-r} gives
\[
\| \frac{\delta^{\lambda,\alpha} \psi_2}r \|_{L^2} + \| \frac{\delta^\lambda A_2}r\|_{L^2} \lesssim  \| \psi \|_{L^2}. 
\]
Then we use this estimate and the system \eqref{sysAp2} to conclude that
\[
\| \partial_r \delta^{\lambda,\alpha} \psi_2 \|_{L^2} + \| \partial_r \delta^\lambda A_2 \|_{L^2} \lesssim  \| \psi \|_{L^2};
\]
this concludes the proof of \eqref{compsolh1-rr}. 

\bigskip

The proof of Proposition~\ref{p:mod} is concluded once we establish \eqref{goodal2}. For this purpose we need the following result. 

\begin{l1} \label{lZL2}
Consider the system of ODE
\begin{equation} \label{mpL2}
\partial_r Z = N Z + F, \qquad \lim_{r \rightarrow \infty} Z(r) =  0.
\end{equation}
If $N,F$ are in in $ \partial_r \dHe + l^1L^2$, with $\|N\|_{\partial_r \dHe + l^1L^2}$ small, then the above equation has a unique solution $Z \in \dHe + [\partial_r]^{-1} l^1L^2$ satisfying
\begin{equation}
\| Z \|_{\dHe + [\partial_r]^{-1} l^1L^2} \lesssim  \| F \|_{\partial_r \dHe + l^1L^2}. 
\end{equation}
Furthermore, the map from $N,F \in \partial_r \dHe + l^1L^2$ to $Z \in \dHe + [\partial_r]^{-1} l^1L^2$ is analytic.
\end{l1}

Some remarks are in order here:

\begin{itemize}

\item We did not specify the size of our system, but we will be mainly interested in the following two cases: i) $Z,N,F$ are $3 \times 3$ matrices and ii) $Z,F$ are $3 \times 1$, $N$ is $3 \times 3$.

\item A similar result holds true for systems of type $\partial_r Z =  Z N + F$. 
\end{itemize}

\begin{proof}  We claim the following basic inequality:
\begin{equation} \label{bildr}
\|f \cdot g \|_{\partial_r \dHe + l^1L^2} \les 
\| f \|_{\partial_r \dHe + l^1L^2} \|  g \|_{\dHe + [\partial_r]^{-1} l^1L^2}
\end{equation}
We write $f=\partial_r f_1+f_2, f_1 \in \dHe, f_2 \in l^1L^2$ and $g=g_1+g_2, g_1 \in \dHe, \partial_r g_2 \in l^1L^2$ with $g_2(\infty)=0$. Then from the simple bound
\begin{equation} \label{l2li}
\|  h \|_{l^2L^\infty} \lesssim \| h \|_{\dHe},
\end{equation}
it follows that $\| \partial_r f_1 \cdot g_1\|_{l^1 L^2} \les \|f_1\|_{\dHe} \|g_1\|_{\dHe} $. Next we write $\partial_r f_1 \cdot g_2 = \partial_r (f_1 g_2)- f_1 \partial_r g_2$; $f_1 \partial_r g_2 \in l^1L^2$ is obvious since $\|f_1\|_{L^\infty} \les \|f_1\|_{\dHe}$. We claim that $f_1 g_2 \in \dHe$; indeed this follows from the representation \eqref{drdech1} applied for $g_2$, the algebra property of $\dHe$ and the fact that $\dHe$ is stable under multiplication by $\chi_{A_{\leq m}}$. We also have the trivial inequality $\| f_2 g_1\|_{l^1L^2} \les \| f_2 \|_{l^1L^2}\| g_1\|_{L^\infty} \les \| f_2 \|_{l^1L^2}\| g_1\|_{\dHe}$.   Finally $f_2 g_2 \in l^1 L^2$ given the bound \eqref{l1l2dr} that places $g_2$ in $L^\infty$.

Back to our problem, the solution $Z$ is obtained via a Picard iteration in the space 
$ \dHe + [\partial_r]^{-1} l^1L^2$ as follows: from \eqref{bildr} we obtain   
\[
\| N Z + F \|_{\partial_r \dHe + l^1L^2}  \lesssim 
 \| N \|_{\partial_r \dHe + l^1L^2} \|  Z \|_{\dHe + [\partial_r]^{-1} l^1L^2} + \| F \|_{\partial_r \dHe + l^1L^2},
\]
and $[\partial_r]^{-1}$ (the operator that gives the solution to the linear inhomogeneous ODE) is bounded from $\partial_r \dHe + l^1L^2$ to $ \dHe + [\partial_r]^{-1} l^1L^2$; finally the convergence of the iterations is insured by the smallness
of $\| N \|_{\partial_r \dHe + l^1L^2}$.  
\end{proof}

We return to the proof of \eqref{goodal2}. We use the system  \eqref{return}, which we recast in a matrix form as an equation for $\calO=(\bar v, \bar w, \bar u)$ as follows
\begin{equation} \label{sysM2}
\partial_r \calO = \calO R(\psi), \quad  \calO(\infty) = I_3. 
\end{equation}
with 
\[
R(\psi)=\left( 
\begin{array}{lll}
0 & 0 &   \Re \psi_1 \\
0 & 0 &  \Im \psi_1 \\
- \Re \psi_1 &  - \Im \psi_1 &  0
\end{array}
\right)
\]
If $\psi = 0$ then $\psi_2 = 2 i e^{2i\alpha} h_1^\lambda$, which yields $\psi_1= - 2 e^{2i\alpha} \frac{h_1^\lambda}{r}$, hence
\begin{equation}
R(0) =
2\frac{h_1^\lambda}r \left( 
\begin{array}{lll}
0 & 0 &  - \cos{2\alpha}  \\
 0 & 0 & -\sin{2\alpha} \\
\cos{2\alpha} & \sin{2\alpha} & 0
\end{array}
\right). \label{EE}
\end{equation}
The solution is given by (see the generalization of \eqref{vwq})
\[
\calO_0 = \!\left(\!\!\! \begin{array}{ccc}
\! h_3^\lambda \cos^2 2\alpha + \sin^2 2\alpha 
 & \! (h_3^\lambda -1) \sin 2\alpha \cos 2 \alpha
& h_1^\lambda  \cos 2\alpha  \cr
 \! (h_3^\lambda -1) \sin 2 \alpha \cos 2 \alpha \!
& \!  h_3^\lambda  \sin^2 2\alpha + \cos^2 2\alpha 
& h_1^\lambda  \sin 2\alpha  \cr
\! -h_1^\lambda  \cos 2\alpha & \!
- h_1^\lambda  \sin 2 \alpha & h_3^\lambda  \end{array}\!\!\! \right).\!
\]
We note that $\calO_0^{-1}= \calO_0^t$. We will prove that 
\begin{equation}\label{Rest2}
\| R(\psi) - R(0)\|_{\partial_r \dHe + l^1L^2} \lesssim \|\psi\|_{L^2}.
\end{equation}
Suppose this is done. Then we write the solution to \eqref{sysM2} is of the form
\begin{equation} \label{Yexpr2}
\calO(r) = (I+Y(r))\calO_0(r)
\end{equation}
where $Y$  solves  the differential equation
\begin{equation} \label{Ysysup2}
\partial_r Y =  Y N + G, \qquad Y(\infty) = 0 
 \qquad N=G=\calO_0 (R(\psi) - R(0))  \calO_0^{-1}.
\end{equation}
We apply Lemma \ref{lZL2} to solve this system and conclude that
\[
\|\bar u - \bar Q^2_{\alpha,\lambda}\|_{\dHe+[\partial_r]^{-1} l^1L^2} 
+  \| \bar v - \bar V^2_{\alpha,\lambda}\|_{\dHe+[\partial_r]^{-1} l^1L^2} 
+ \| \bar w- \bar W^2_{\alpha,\lambda}\|_{ \dHe+[\partial_r]^{-1} l^1L^2}\lesssim \|\psi\|_{L^2}.
\]
To finish our claim we need an additional bound
for $\|r^{-1}(\bar u- \bar Q_{\alpha,\lambda})\|_{L^2}$.  We first remark that the last row of
$\calO$ is a-priori known, namely $(\bar v_3, \bar w_3, \bar u_3) = \frac12 (-\Im \psi_2,\Re
\psi_2,A_2)$; this already shows that 
\[
\| r^{-1}(\bar v_3 +  h_1^\lambda \cos{2\alpha} )\|_{L^2} + \| r^{-1} (\bar w_3+ h_1^\lambda \sin{2\alpha}) \|_{L^2} + \| r^{-1}(\bar u_3 - h_3^\lambda)\|_{L^2}  
\lesssim \|\psi\|_{L^2}. 
\]
To transfer this information to $\bar u_1$ and $\bar u_2$ we use again the orthogonality
of $\calO$. To keep the computations below compact we assume that $\alpha=0$; this does not restrict the generality of the argument. For $\bar u_1$ we have 
\[
\bar u_1 = \bar v_2 \bar w_3 - \bar v_3 \bar w_2= \bar v_2 \bar w_3  - (\bar v_3+h_1^\lambda) \bar w_2 + h_1^\lambda  \bar w_2, 
\]
from which it follows that
\[
\frac{\bar u_1 - h_1^\lambda}r =  \bar v_2 \frac{\bar w_3}r  - \frac{\bar v_3+h_1^\lambda}r \bar w_2 + \frac{h_1^\lambda}r  (\bar w_2-1);  
\]
From the above estimates it follows that $\|\frac{\bar u_1 - h_1^\lambda}r\|_{L^2} \les \|\psi\|_{L^2}$. A similar argument shows that $\|\frac{\bar u_2}r\|_{L^2} \les \|\psi\|_{L^2}$. This concludes the proof of \eqref{goodal2}.

It remains to prove the bound \eqref{Rest2}. Using the second relation in 
\eqref{comp1} we have
\[
\begin{split}
\psi_1 = &  \psi + i\frac{\psi_2}{r} = \frac14 (-i  A_2 \partial_r \psi_2 
+  |\psi_2|^2  \psi + i  \frac{\psi_2}{r} |\psi_2|^2)
\\
= & -2e^{2i\alpha} \frac{h_1}{r} -i A_2 \partial_r \delta^{\lambda,\alpha} \psi_2 + \delta^\lambda A_2 \cdot 2 e^{2i\alpha} \partial_r h_1^\lambda
+  |\psi_2|^2  \psi + \frac1{4r} (  i  \psi_2 |\psi_2|^2 + 8 e^{2i\alpha} (h_1^\lambda)^3)
\end{split}
\]
From this we obtain:
\begin{equation} \label{dpsi1}
\psi_1 + 2e^{2i\alpha} \frac{h_1}{r} = - 2i h_3^\lambda \partial_r \delta^{\lambda,\alpha} \psi_2 + B,
\end{equation}
where
\[
B= - i \delta^\lambda A_2 \partial_r \delta^{\lambda,\alpha} \psi_2 + \delta^\lambda A_2 \cdot 2 e^{2i\alpha} \partial_r h_1^\lambda
+  |\psi_2|^2  \psi + \frac1{4r} (  i  \psi_2 |\psi_2|^2 + 8 e^{2i\alpha} (h_1^\lambda)^3).
\]
Using \eqref{compsolh1-rr}, \eqref{l2li} and the simple observation that $\|h_1^\lambda\|_{l^2L^\infty} \les 1$, it is a straightforward exercise to check that
\[
\|B\|_{l^1L^2} \les \|\psi\|_{L^2}. 
\]
Finally, we write
\[
h_3^\lambda \partial_r \delta^{\lambda,\alpha} \psi_2 = \partial_r  (h_3^\lambda \delta^{\lambda,\alpha} \psi_2) + \partial_r h_3^\lambda \cdot  \delta^{\lambda,\alpha} \psi_2,
\]
and note that $\| h_3^\lambda \delta^{\lambda,\alpha} \psi_2 \|_{\dHe} \les \| \delta^{\lambda,\alpha} \psi_2 \|_{\dHe} \les \|\psi\|_{L^2}$, while 
\[
\| \partial_r h_3^\lambda \cdot  \delta^{\lambda,\alpha} \psi_2 \|_{l^1L^2} \les \| \frac{\delta^{\lambda,\alpha} \psi_2}r\|_{L^2} \les \|\psi\|_{L^2}
\]
(just as one does for the terms in $B$). This concludes the proof of \eqref{Rest2}, and, in turn, the proof of our Proposition.

\end{proof}

\begin{proof} 
[Proof of Proposition~\ref{p:modb}]
The argument is similar to the proof of Proposition~\ref{p:mod}. To emphasize this similarity and avoid cluttering the notations we redenote $(\alpha,\lambda)$ in the 
proposition by $(\alpha_0,\lambda_0)$.

We fix some $r_0 \approx \lambda_0^{-1}$. Then for $(\alpha,\lambda = \e^{\gamma} \lambda_0)$ close to $(\alpha_0,\lambda_0)$ 
we use Lemma~\ref{psito2} to solve the system \eqref{comp1} with initial data as in \eqref{ldef-r}. Denoting 
\[
G(\alpha,\gamma) = \la \delta^{\lambda_0,\alpha_0} \psi_2, \varkappa^{\lambda_0} \ra,
\]
we look for $(\alpha,\gamma)$ so that $G(\alpha,\gamma)=0$. Here $G$ is uniformly smooth, and by   \eqref{compsolh1-rr} we have 
\[
|G(\alpha_0,\lambda_0)| \lesssim \|\psi\|_{L^2}
\]
Further, a direct computation shows that 
the differential $DG(\alpha_0,\lambda_0)$ is nondegenerate at $(\alpha_0,0)$. Then the inverse function theorem yields $(\alpha,\gamma)$ so that 
$G(\alpha,\gamma)=0$, and 
\[
|\alpha-\alpha_0| + |\gamma| \lesssim \|\psi\|_{L^2}.
\]
This in turn allows us to switch from $(\alpha,\lambda)$ to $(\alpha_0,\lambda_0)$ thus completing the argument for the first part of the Proposition, concerning the recovery of $\psi_2,A_2$ and the verification of \eqref{lell} and \eqref{lellf}.

Concerning the second part, which requires the reconstruction of an actual map $u$, we note that this has been essentially done in the argument for \eqref{goodal2}. Indeed, there we start with the analysis of the system \eqref{sysM2} which has $\psi_1=\psi+i\frac{\psi_2}r$ as its  main input (see $R(\psi)$), and where we have already recovered $\psi_2$ in the first part of the Proposition. The analysis of the  the system \eqref{sysM2} essentially reconstructs $u$ (along with its gauge components $v,w$). At this point we need to establish the following:

\begin{enumerate}[label=(\roman*)]
\item  the gauge elements corresponding to the map $u$ are indeed $\psi_1,\psi_2, A_1=0, A_2$, and the gauge field is $\psi$;

\item the modulation parameters subject to \eqref{chi-nondeg} are indeed $(\alpha,\lambda)$ from above. 
\end{enumerate}

The fact that $\bar u, \bar v, \bar w$ satisfy the system \eqref{sysM2} shows that $(v,w)$
is the correct gauge corresponding to $A_1=0$; this is better seen in the formulation \eqref{return} of the system where we take $k=1$. It also follows that $\psi_1$ is the correct representation of $\partial_r u$ in the frame $(v,w)$.

We need to do the same for the representation of $\partial_\theta u$, the angular derivative of $u$. At this point we have the pair $(A_2,\psi_2)$ which was constructed from $\psi$ in the first part of the Proposition, and the modulation parameters in the first part of the Proposition, and we have also the pair $(A_2^{true},\psi_2^{true})$ which are the "true" gauge components of the map $u$ in the Coulomb gauge as described in Section \ref{CG}. Our goal is to establish that they are the same, that is 
$(A_2,\psi_2)=(A_2^{true},\psi_2^{true})$. if we achieve this then we are done since we have already established that $\psi_1=\psi_1^{true}$; thus $\psi=\psi_1-i\frac{\psi_2}r$ is the correct gauge field, and the modulation parameters satisfy the correct orthogonality condition \eqref{chi-nondeg}, so they are also the correct modulation parameters for $u$. 

To establish the equality $(A_2,\psi_2)=(A_2^{true},\psi_2^{true})$ we note the following:

\begin{enumerate}[label=(\roman*)]
\item The vector $\textbf{v} =(v_3,w_3, u_3)$ satisfies the system 
\[
\partial_r \textbf{v} = \textbf{v} R(\psi), \quad  \textbf{v}(\infty) = (0,0,1),
\]
which is a consequence of \eqref{sysM2};

\item The vector $\textbf{v'} =\frac12(- \Im \psi_2, \Re \psi_2, A_2)$ satisfies the same exact system
\[
\partial_r \textbf{v'} = \textbf{v'} R(\psi), \quad  \textbf{v'}(\infty) = (0,0,1).
\]
\end{enumerate}

Indeed $\psi_2,A_2$ obey the system \eqref{sysAp2} which, using the compatibility condition \eqref{sphere}, implies that they obey \eqref{comp1} and in turn \eqref{comp}; the later system \eqref{comp} implies the claim above that $\textbf{v'}$ obeys the system, while the condition at $\infty$
follows from the fact that $\delta^\lambda A_2, \delta^{\lambda,\alpha} \psi_2$ are both in $\dHe$, hence have zero limit at $\infty$, and thus the limit is given by the limit of $\frac12 (-\Im (2ie^{2i\alpha} h_1^\lambda), \Re (2ie^{2i\alpha} h_1^\lambda), 2h_3^\lambda)$ which is precisely $(0,0,1)$. Just as in the analysis of \eqref{sysM2}, we recast this so that it fits into the framework of Lemma \ref{lZL2}; precisely, we write the system for $\textbf{v} - (0,0,1)$ and $\textbf{v'} - (0,0,1)$, so that we have the zero data at infinity. To invoke the uniqueness part from Lemma \ref{lZL2}, we need establish that 
$\textbf{v} - (0,0,1)$ and $\textbf{v'} - (0,0,1)$ belong to $\dHe + [\partial_r]^{-1} l^1L^2$; this is clear for $\textbf{v'} - (0,0,1)$ since this is where we solved the system \eqref{sysM2}. As for 
$\textbf{v} - (0,0,1)= \frac12(- \Im \psi_2, \Re \psi_2, A_2-2)$, we note that from the above arguments we already know that
\[
\| \delta^{\lambda,\alpha} \psi_2 \|_{\dHe} + \| \delta^\lambda A_2 \|_{\dHe} \lesssim  \| \psi \|_{L^2}.
\]
This is enough to place the first two components of $\textbf{v} - (0,0,1)$ in $\dHe$. Next we write 
$A_2-2 = \delta^\lambda A_2 + 2(h_3^\lambda-1)$ and note that $\delta^\lambda A_2 \in \dHe$, while 
 $h_3^\lambda-1 \in [\partial_r]^{-1} l^1L^2$ (this follows from $\partial_r (h_3-1)= 2 \frac{(h_1)^2}r \in l^1 L^2$ when $\lambda=1$). Therefore we conclude that $\textbf{v} = \textbf{v'}$, which completes the proof of Proposition \ref{p:modb}.

\end{proof}

As we consider solutions for the Schr\"odinger map equation, the modulation parameters will vary as functions of time. Later we will study this dependence in much greater detail, but for now we are   content with proving that they are $C^1$ functions of time:

\begin{c1} \label{cdif}

The modulation parameters $\lambda(t)$ and $\alpha(t)$ associated to a map $u$ based on \eqref{lell} and \eqref{lellf} are continuously differentiable on the maximal interval of existence $I_{max}$ and obey the pointwise bounds
\[
|\alpha'|+|\frac{\lambda'}{\lambda}| \les \lambda^2. 
\]
\end{c1}

\begin{proof}
We simply invoke \eqref{Liplp} and seek an estimate for $\partial_t \psi_2$ in appropriate norms. Recall from \eqref{eqn:psi12} with $A_1=0$ that 
\[
\partial_t \psi_2=  - i A_0 \psi_2 +  i \Delta \psi_2  - i \frac{1}{r^2} A_2^2
\psi_2 - \Im{(\psi_2 \bar{\psi}_1)} \psi_1.
\]
Based on \eqref{Liplp}, it suffices to obtain a uniform bound in $\dot H^{-1}_e$ for the right-hand side. The main observation here is that we have the uniform in time bound
\begin{equation} \label{psi2dhe}
\| \psi_2(t) \|_{\dHe} \les 1.
\end{equation}
Indeed, from the relation $|\psi_2|^2 = m^2(u_1^2 + u_2^2)$ and the expression of the energy in \eqref{energy}, it follows that $\|\frac{\psi_2}r\|_{L^2} \les 1$. From the compatibility condition \eqref{comp}, we have $\partial_r \psi_2 = i A_2 \cdot \psi_1$, thus
\[
\|\partial_r \psi_2\|_{L^2} \les \| \psi_1\|_{L^2} \les \| \partial_r u \|_{L^2} \les E(u) \les 1. 
\]
Based on \eqref{psi2dhe}, we estimate as follows:
\[
|\la \Delta \psi_2, \phi \ra| \les \| \partial_r \psi_2 \|_{L^2} \cdot \| \partial_r \phi \|_{L^2}\les  \|\phi\|_{\dHe},
\]
\[
|\la \frac{1}{r^2} A_2 \psi_2, \phi \ra| \les \| \frac{\psi_2}r \|_{L^2} \cdot \| \frac{\phi}r \|_{L^2}\les \|\phi\|_{\dHe},
\]
\[
|\la \Im{(\psi_2 \bar{\psi}_1)} \psi_1, \phi \ra| \les \| \psi_2 \|_{L^\infty} \cdot \|\psi_1\|_{L^2}^2 \cdot \| \phi \|_{L^\infty} \les \|\psi_2\|_{\dHe} \|\phi\|_{\dHe} \les \|\phi\|_{\dHe},
\]
\[
|\la A_0 \psi_2, \phi \ra| \les \| \psi_2 \|_{L^\infty} \cdot \|A_0\|_{L^1} \cdot \| \phi \|_{L^\infty} \les 
\|\psi_2\|_{\dHe} \|\phi\|_{\dHe} \|\psi\|_{L^2} \|\frac{\psi_2}r\|_{L^2} \les \|\phi\|_{\dHe}. 
\]
The above estimates suffice to conclude the proof. 

\end{proof}

\subsection{The complete setup of the problem} \label{Sec3wrap}

Here we summarize the conclusions of this section. Each $2$-equivariant solution to the Schr\"odinger map equation with homotopy degree $2$ and energy 
just above the ground state energy is uniquely represented by 
two components, which vary as functions of time:

\begin{enumerate}
    \item the differentiated field $\psi \in L^2$;
    \item the modulation parameters $(\alpha,\lambda)$.
\end{enumerate}

In order to identify the  reference soliton associated to a map and thus fix the \emph{modulation parameters}, we use the orthogonality conditions  
\eqref{lell} and \eqref{lellf}, which we recall here for convenience:
\begin{equation}\label{which-lambda}
\la \delta^{\lambda,\alpha} \psi_2, \varkappa^\lambda \ra = 0, \qquad 
\|\delta^{\lambda,\alpha} \psi_2\|_{L^\infty} \les \delta.
\end{equation}
Having made this choice, the PDE driving the evolution of $\psi$, namely the equation \eqref{psieq}, becomes
\begin{equation} \label{psieq2}
i \partial_t \psi  - \tilde H_\lambda \psi = N(\psi), 
\end{equation}
where we have used the notations
\[
\tilde H_\lambda : = -  \Delta + \frac{5-4h_3^\lambda}{r^2}, 
\quad N(\psi):=W_\lambda \cdot  \psi, \quad 
W_\lambda : = A_0 - 2 \frac{\delta^\lambda A_2}{r^2}   - \frac1{r}\Im{(\psi_2 \bar{\psi})}.
\]
The initial data of \eqref{psieq2} satisfies $\|\psi(0)\|_{L^2} \les \delta$. 
Our analysis shows that the auxiliary functions $(\psi_2,A_2)$ and then $A_0$ are uniquely determined by $\psi$ together with $(\alpha,\lambda)$.

 It is easy to see that all terms in $N$ are at least quadratic in $\psi$, or can be estimated by quantities that are at least quadratic in $\psi$. We place it in  right-hand side of \eqref{psieq2} in order to emphasize our goal to treat $N$ as a nonlinear perturbation. 

The linear part of \eqref{psieq2} shows that the following linear PDE
\begin{equation} \label{Linpsi}
(i \partial_t - \tilde H_\lambda) \Psi =0 
\end{equation}
should play crucial role in our analysis. It is important to note here that $\lambda=\lambda(t)$, and thus we are dealing with a variable coefficient linear PDE. In Section \ref{s:linear} we carry a full analysis of \eqref{Linpsi} (including its inhomogeneous counterpart), and in Section~\ref{s:nolinear}
we prove an array of estimates fo $N(\psi)$.

 The analysis of the time evolution of the modulation  parameters $(\alpha,\lambda)$, on the other hand, begins in  Section \ref{s:mod},
 where we use the orthogonality condition \eqref{which-lambda} in order to derive the modulation  equations governing the evolution of $\lambda$ and $\alpha$ as a function of time. This analysis is further refined in Section \ref{s:mod2}.

At this point we are also able to explain the role of the parameter $\delta$ involved in this paper. Recall that we consider  $2$-equivariant maps $u$ 
satisfying \eqref{eqvl} and 
whose energy satisfies $E(u) \leq 8\pi + \delta^2$; this translates into the mass constraint $\|\psi(0)\|_{L^2_r} \leq \delta$, which is propagated along the flow on the maximal interval of existence. The motivation for  the smallness condition $\delta \ll 1$ comes from two different sources:

i) the fixed time theory developed in this section, more precisely the constraint imposed in Proposition \ref{p:mod};

ii) the dynamic (in time) analysis of the PDE \eqref{psieq2} described above, as well as the dynamics (in time) of the ODE system describing the evolution of the modulation parameters $\alpha(t)$ and $\lambda(t)$ (to be detailed in Section \ref{s:mod}).

\section{Spectral analysis for  the operators 
\texorpdfstring{$H$, $\tilde H$}{}}
\label{spectral}

As discussed in the previous section, see in particular \eqref{Linpsi}, the linearization of the Schr\"odinger map problem in its gauge representation near  the gauge components of the soliton $Q^2_{\alpha,\lambda}$ reveals that the operator $\tilde H_\lambda$ plays an important role in our problem. If we were instead to directly linearize 
the Schr\"odinger map problem near the the soliton $Q^2_{\alpha,\lambda}$, then the corresponding linearization in \eqref{Linpsi} would replace $\tilde H_\lambda$ with $H_\lambda$, see \cite{GuNaTs} for details. This second linearization plays no role in our analysis, so apriori there is no reason to study the spectral properties of $H_\lambda$. Indeed, our main goal in this section is to develop the spectral theory for $\tilde H_\lambda$. But for technical reasons, it is convenient to do so first for $H_\lambda$ and then use this to derive the one for $\tilde H_\lambda$; in particular we heavily rely on the analysis developed by Krieger, Schlag and Tataru in Section 4 of \cite{kstym}, where the spectral theory is developed for an operator related to $H_\lambda$ rather than $\tilde H_\lambda$. 

The spectral theory for $\tilde H_\lambda$ will play a key role in the study of the dispersive properties of the linear Schr\"odinger equation \eqref{Linpsi}, which in turn will be critical in the study of \eqref{psieq2}, the main PDE governing the evolution of the gauge field $\psi$.

Let us recall from Section~\ref{outlinH}  the definition of the two operators
\[
\begin{aligned}
H_\lambda = & \  -\Delta + V_\lambda, \quad V_\lambda(r) =\frac{4}{r^2}(1-2(h_1^\lambda)^2),
\\
\tilde H_\lambda =& \ -\Delta + \tilde V_\lambda, \quad \tilde V_\lambda(r) = \frac{5-4h_3^\lambda}{r^2}. 
\end{aligned}
\]
The operators $H_\lambda$ and $\tilde H_\lambda$ are conjugate operators
and admit the factorizations
\[
H_\lambda =  L^*_\lambda L_\lambda, \qquad \tilde  H_\lambda = L_\lambda  L^*_\lambda,
\]
where
\begin{equation} \label{Llambdadef}
L_\lambda = h_1^\lambda \partial_r \frac{1}{h_1^\lambda} = \partial_r + \frac{2}{r} h_3^\lambda, \qquad
L^{*}_\lambda=- \frac{1}{h_1^\lambda} \partial_r h_1^\lambda -\frac1{r}= -\partial_r +
\frac{2 h_3^\lambda-1}{r}.
\end{equation}
 In these expressions, $\lambda$ plays the role of a scaling parameter,  therefore for the spectral theory it suffices  to consider the case $\lambda=1$. For convenience we denote the
corresponding operators by $H=H_1, \tilde H=\tilde H_1, L=L_1, L^*=L^*_1$. Note that we also have the scaling relations:
\[
L_{\lambda}(v)(r) = \lambda L(v(\frac{\cdot}{\lambda}))(r\lambda), \quad 
L^*_{\lambda}(v)(r) = \lambda L^*(v(\frac{\cdot}{\lambda}))(r\lambda).
\]
In the first part of the section we develop the spectral theory for the operators $H$ and $\tilde H$, and at the end we use it to derive the spectral theory for the rescaled operators $H_\lambda$ and $\tilde H_\lambda$, considering in particular the \emph{transference operator}
which describes the $\lambda$ dependence in the spectral theory.

We note that the spectral theory for $H$ in the case $m=1$ was studied in detail by Krieger-Schlag-Tataru in \cite{KrScTa}. Then it was adapted to a setup similar setup to ours (in particular extending the theory to $\tilde H$) in \cite{BeTa-1}.

In the current paper we need the spectral theory for $H$ and $\tilde H$ in the case $m=2$. Just as in prior works, see \cite{KrScTa,kstym}, one shows that $H$ fits into the theory developed by Gesztesy-Zinchenko in \cite{GeZi}, which provides the abstract framework for the spectral theory. We first characterize the generalized eigenfunctions for $H$; once the spectral theory for $H$ is complete, we use it together with the conjugation described above in order to derive the one for $\tilde H$. 

To show that $H$ fits into the theory developed by Gesztesy-Zinchenko
in \cite{GeZi}, we proceed as in Section 4 of \cite{kstym}: if we define the operator $\mathcal{L}$ by conjugating $H$ with respect to the weight
$\sqrt{r}$, so that it is selfadjoint in $L^2(dr)$,
\begin{equation}\label{lcaldef}
\mathcal{L} = \sqrt{r} H r^{-\frac12} = -\partial_r^2 + \left(\frac{15}{4r^{2}}-32\frac{r^{2}}{(1+r^{4})^{2}}\right):=-\partial_r^2+V(r).
\end{equation}
Then, by Example 3.10 in \cite{GeZi}, the potential $V$ satisfies hypothesis 3.1 in \cite{GeZi}; the interested reader may find some more details in \cite{kstym}, which studies the spectral theory for the operator $\mathcal{L}_{KST}$ defined by
\[
\mathcal{L}_{KST}:= -\partial_r^2 + \left(\frac{15}{4r^{2}}-\frac{24 }{\left(r^2+1\right)^2}\right).
\]
This provides us with the abstract spectral theory described below. \\
\\
Precisely, we consider $H$ acting as an unbounded selfadjoint operator in
$L^2(rdr)$. Then $H$ is nonnegative, and its spectrum $[0,\infty)$ is  absolutely continuous. 
$H$ has a zero eigenvalue, namely $\phi_0=h_1$,
\[
 H h_1 = 0, \qquad
h_1(r)=\frac{2r^2}{r^{4}+1}.
\]
For each $\xi > 0$ one can choose a normalized generalized eigenfunction
$\phi_\xi$, smooth at $r = 0$,  
\begin{equation}\label{phi-xi}
 H \phi_\xi = \xi^2 \phi_\xi.
\end{equation}
These generalized eigenfunctions are unique up to a $\xi$ dependent multiplicative factor,
which is chosen as described below. To these one associates a generalized Fourier transform $\FH$. 

Because of the zero eigenvalue, we regard $\FH$ as a two-component vector, defined 
by 
\begin{equation}\label{htransvecnotation}
\FH:L^2(rdr) \to \R \times L^2(d\xi), \qquad 
\FH(f)(\xi) = \begin{bmatrix} a\\
g(\xi)\end{bmatrix}\end{equation} 
and
\[
a=\langle f, \phi_{0}\rangle_{L^{2}(r dr)}\|\phi_{0}\|^{-1}_{L^{2}(r dr)},
\] 
\[
g(\xi) = \int_0^\infty \phi_\xi(r) f(r)  rdr, \quad \xi >0,
\]
where the integral above is considered in the singular sense. This is the same notation as in \cite{kstym}. On the image we 
will use  the norm
\[
\|\begin{bmatrix} a\\
	g\end{bmatrix}\|_{L^{2,0}}^{2}:= |a|^{2}+\int_{0}^{\infty} |g(\xi)|^{2} d\xi.
\]
This makes the generalized Fourier transform an $L^2$ isometry,
and we have the inversion formula
\begin{equation} \label{invf}
\FH^{-1}(\begin{bmatrix} a\\
g(\cdot)\end{bmatrix})(r) = a \frac{\phi_{0}(r)}{\|\phi_{0}\|_{L^{2}(r dr)}}+ \int_{0}^{\infty} g(\xi) \phi_{\xi}(r) d\xi.
\end{equation}

Compared to the setup of \cite{kstym}, there are several differences in our setting here, stemming from the different
choice of the spectral parameter in \eqref{phi-xi}, namely 
$\xi^2$ rather than $\xi$, as well as our use of the $rdr$ measure rather than $dr$. For quick reference, we provide
a brief dictionary. 

The counterparts, in our setting, of $\phi(r,\xi)$ from Proposition 4.5 of  \cite{kstym} are $\phi_{KST}(r,\xi)$, which we define, for $\xi >0$, so that
\begin{equation}\label{phikstintro}
\sqrt{\rho(\xi)} \phi_{KST}(r,\xi) = \frac{\sqrt{r}}{\sqrt{2}\xi^\frac14}  \phi_{\sqrt{\xi}}(r)\end{equation}
where $\rho$ is the density of the continuous part of the spectral measure of $\mathcal{L}$, which exactly corresponds to the function $\rho$ from Theorem 4.3 of \cite{kstym}. For the zero mode, we have
\[
\phi_{KST}(r,0)=\frac{\sqrt{r}}{2} \phi_{0}(r).
\]
The counterpart of the function $\psi^{+}$ from Theorem 4.3 of \cite{kstym}
is 
\[
\psi_{KST}^{+}(r,\xi) := \frac{\sqrt{r}}{\xi^{1/4}} \phi_{\sqrt{\xi}}^{+}(r).
\]
where $\phi_{\xi}^{+}$ is defined in \eqref{repphi+}.

\medskip

The functions $\phi_\xi$ are smooth with respect to both $r$ and $\xi$.
To describe them one considers two distinct regions, $r \xi \lesssim 1$ 
and $r \xi \gtrsim  1$. The main properties of $\phi_{\xi}$  are summarized in the next theorem.
\begin{t1} \label{htransthm}
a) 	In the inner region $\{r\xi \lesssim 1\}$ the functions $\phi_\xi$
	admit a power series expansion of the form
	\begin{equation} \label{repphi}
		\phi_\xi (r)= q(\xi) \left( \phi_0(r) + \sum_{j=1}^\infty (r\xi)^{2j}
		\phi_j(r^2)\right)
	\end{equation}
	where $\phi_0=h_1$ and the functions $\phi_j$ are analytic and satisfy 
	\begin{equation}\label{phijest}|\phi_{j}(u)| \leq \frac{C_{1}}{((j-1)!)^{2}} \left(\frac{3}{2^{3/2}}\right)^{j-1} \frac{u}{\langle u \rangle}, \quad u>0, \ j \geq 1,\end{equation}
as well as\footnote{ Here $\phi_j$ are also analytic at $0$; we only wrote the bounds in this manner in order to gain uniformity in $j$.} 
	\begin{equation}\label{drphijest}
	|\phi_{j}^{(n)}(u)| \leq \frac{M_{n}  j^{2n-1}}{((j-1)!)^{2}} \left(\frac{3}{2^{3/2}}\right)^{j-1} \frac{u}{u^{n} \langle u\rangle}, \quad u>0, \ j \geq 1, \ n \geq 1,
	\end{equation}
for some $M_{n}>0$.
	The smooth positive weight $q$ satisfies 
	\begin{equation}\label{qest}
		q(\xi) \approx \left\{ \begin{array}{ll}
			\displaystyle \xi^\frac12,  &  \xi \ll 1 \cr\cr
			\xi^\frac52,   &  \xi \gg 1
		\end{array} \right., \qquad
		|(\xi \partial_\xi)^\alpha q| \lesssim_\alpha q.
	\end{equation}
In addition, if 
\begin{equation} \label{defmk1}
	\omega_k^1(r)=
	\left\{
	\begin{array}{ll}
		\min\{1, 2^{2k} r^2, 2^{2k} r^4 \} ,  & \ k < 0 \\ \\
		\min \{1, 2^{4k} r^4 \},  & \ k \geq 0
	\end{array}
	\right.
\end{equation}
then, we have
\begin{equation} \label{pointphilow}
	|(\xi \partial_\xi)^\alpha (r \partial_r)^\beta \left( \phi_\xi(r) - q(\xi) \phi_0(r)\right)|
	\lesssim_{\alpha\beta} 2^{\frac{k}2} \omega_k^1(r), \qquad \xi \approx 2^k,\  r\xi \lesssim 1. 
\end{equation}
\medskip

b) In the outer region, $\{r \xi \gtrsim  1\}$,  we have the representation
\begin{equation} 
	\phi_{\xi}(r)=a(\xi) \phi^{+}_\xi(r) + \overline{a(\xi) \phi^{+}_\xi(r)}
\end{equation}
where the complex valued weight $a$ satisfies

\begin{equation} \label{thmabound}
	|a(\xi)| = \frac{1}{\sqrt{2 \pi}}, \qquad | (\xi \partial_\xi)^\alpha a(\xi)| \lesssim_\alpha 1
\end{equation}
and 
\begin{equation} 
	\phi^{+}_\xi(r)= r^{-\frac12} e^{ir\xi} \sigma(r\xi,r), 
	\qquad r\xi \gtrsim  1.
\end{equation}
Moreover, $\sigma$ above admits the  asymptotic expansion
\begin{equation}\label{sigma-exp}
\sigma(y,r) \approx \sum_{j=0}^\infty y^{-j} \phi^{+}_j(r), 
\qquad \phi_0^{+}=1
\end{equation}
with smooth functions $\phi_j^+$ satisfying 
symbol type bounds
\[
\sup_{r > 0} |(r\partial_r)^k \phi^{+}_j| < \infty
\]
and $\phi_1^+$ with asymptotic behavior
\begin{equation}
\begin{split}
\phi_{1}^{+}(r) =\frac{15 i}{8}-\frac{16 i}{5 r^4} + O\left(\frac{1}{r^{8}}\right), \quad r \rightarrow \infty.
\end{split}
\end{equation}
\end{t1}
Here the expansion \eqref{sigma-exp} holds in 
in the sense that for all $n,m \geq 0$, there exists $N_{n,m}>0$ such that, for all $N \geq N_{n,m}$, there exists $C_{N,n,m}$ such that
\begin{equation} \label{sigmaasymp}
\begin{split}
|(r\partial_{r})^{n}(y\partial_{y})^{m} \left(\sigma(y,r)-\sum_{k=0}^{N} y^{-k}\phi_{k}^{+}(r)\right)|\leq \frac{C_{N,n,m}}{y^{N+1}}, \quad y>1, r>0. 
\end{split}
\end{equation}
	We remark that the theorem in particular implies that
\begin{equation}\label{phiximinusqphi0}
|\phi_{\xi}(r) -q(\xi)\phi_{0}(r)| \leq C q(\xi) r^{2}\xi^{2} \frac{r^{2}}{\langle r^{2}\rangle}, \quad  r \xi \lesssim 1,
\end{equation}
\begin{equation}\label{philgr}|\phi_{\xi}(r)| \leq \frac{C}{\sqrt{r}}, \quad 1 \lesssim r\xi.
\end{equation}

\bigskip

The leading role in this paper the spectral theory for the operator $\tilde H$. This is derived from the spectral theory 
for $H$ due to the conjugate representations
\[
H = L^* L, \qquad \tilde H = L L^*.
\]
This allows us to define generalized eigenfunctions $\psi_\xi$ for
$\tilde H$ using the generalized eigenfunctions $\phi_\xi$ for
$ H$,
\begin{equation}\label{psifromphi}
 \psi_\xi = \xi^{-1} L \phi_\xi, \qquad L^* \psi_\xi = \xi \phi_\xi.
\end{equation}
It is easy to see that $\psi_\xi$ are real, smooth, vanish at $r = 0$
and solve
\[
\tilde H \psi_\xi = \xi^2 \psi_\xi.
\]
With respect to this frame we can define the generalized Fourier transform
adapted to $\tilde H$ by 
\[
\FtH{f}(\xi)=\int_0^\infty \psi_\xi(r) f(r)  rdr
\]
where the integral above is again considered in the singular sense. This is again an $L^2$ isometry, and we have the inversion formula
\begin{equation} \label{FTL0}
	f(r) = \int_0^\infty \psi_\xi(r) \FtH(f)(\xi)  d\xi.
\end{equation}
To see this we compute, for a Schwartz function $f$:  
\begin{equation} \label{HtH}
	\begin{split}
		\FtH(Lf)(\xi) & =\! \int_0^\infty  \psi_\xi(r) L f(r)  rdr
		= \!\int_0^\infty L^* \psi_\xi(r)  f(r)  rdr
		\\ &= \!\int_0^\infty \xi \phi_\xi(r)  f(r)  rdr = \xi \FH(f)(\xi).
	\end{split}
\end{equation}
Hence
\[
\| \FtH(Lf)\|_{L^2}^2 = \| \xi \FH{f}(\xi)\|_{L^2}^2 = \langle H f,f\rangle_{L^2(rdr)}
= \|Lf\|_{L^2}^2
\]
which suffices since $Lf$ spans a dense subset of $L^2$.\\
\\
The representation of $\psi_\xi$ in the two regions $r\xi \lesssim 1$
and $r\xi \gtrsim 1$ is obtained from the similar representation of $\phi_\xi$. The main properties of $\psi_{\xi}$ which we use are summarized below.
\begin{t1}\label{httransthm}
a) 	In the region $r\xi \lesssim 1$ the functions $\psi_\xi$
	admit a power series expansion of the form
	\[
	\psi_\xi = \xi q(\xi) \left(\psi_0(r) +  r \sum_{j \geq 1} (r\xi)^{2j} 
	{\psi}_j(r^2)\right),
	\]
	where the coefficients $\psi_j$ are given by
	\[
	{\psi}_j(r)= ( 2 h_3(\sqrt{r})+2 +2j) \phi_{j+1}(r) + 2 r \partial_r \phi_{j+1}(r),
	\]
	and satisfy
	\[
	|(r \partial_r)^\alpha \psi_j| \lesssim_\alpha \frac{C^j}{(j-1)!} \frac{r^2}{\langle r \rangle^2},
	\]
	\[
	\psi_{0}(r) =L(r^{2}\phi_{1}(r^{2}))= \frac{r^{2}-(1+r^{4})\arctan(r^{2})}{2r^{3}}.
	\]
  In particular, 
  \[
  \psi_{0}(r) = -\frac{r^3}{3} +\frac{r^7}{15}+ O(r^{11}), \quad r \rightarrow 0
  \]
 and
 \[
 \psi_{0}(r) = -\frac{\pi  r}{4}+\frac{1}{r} +O(r^{-3}), \quad r \rightarrow \infty.
 \]
 
	In addition, if \begin{equation}\label{defmk}
		\omega_k(r)=
		\left\{
		\begin{array}{ll}
			& \min(1,\frac{r^2}{\langle r^2 \rangle} 2^k r), \ \ \mbox{if} \ k < 0 \\ \\
			& \min(1,2^{3k} r^3) , \ \ \mbox{if} \ k \geq 0
		\end{array}
		\right.
	\end{equation}
	then, we have 
	\begin{equation} \label{pointtp}
		| (r \partial_r)^\alpha (\xi \partial_\xi)^\beta \psi_{\xi}(r) | \lesssim_{\alpha\beta}
		2^{\frac{k}2} \omega_k(r) , \qquad \xi \approx 2^k,\  r\xi \lesssim 1.
	\end{equation}
\medskip
b) In the region $r \xi \gtrsim  1$ we have
the representation
\begin{equation} \label{psirep}
	\psi_{\xi}(r)=a(\xi) \psi^{+}_\xi(r) + \overline{a(\xi) \psi^{+}_\xi(r)}
\end{equation}
where
\begin{equation} \label{reppsi}
	\psi^{+}_\xi(r)= r^{-\frac12} e^{ir\xi} \tilde\sigma(r\xi,r), 
	\qquad r\xi \gtrsim  1
\end{equation}
and $\tilde \sigma$ is given by 
\[
\tilde\sigma(y,r)  = i \sigma(y,r) -\frac12 y^{-1} \sigma(y,r)
+ \frac{\partial}{\partial y} \sigma(y,r)+ \frac{r}{y} L \sigma(y,r) 
\]
and has exactly the same properties as $\sigma$. In particular,
for fixed $\xi$ we have
\begin{equation}
	\tilde{\sigma}(r\xi,r) = i - \frac{3}{8} r^{-1}\xi^{-1} + O(r^{-2}).
\end{equation}

\end{t1}

The last two theorems describe the spectral theory for the operators $H_\lambda$ and $\tilde H_\lambda$ in the case $\lambda=1$. For general $\lambda$, it is easy to see that the following functions
\[
\phi_\xi^\lambda(r)= \lambda^\frac12 \phi_{\lambda^{-1} \xi} (\lambda r), \qquad
\psi_\xi^\lambda(r)= \lambda^\frac12 \psi_{\lambda^{-1} \xi} (\lambda r).
\]
are the eigenvalues of the operators $H_\lambda$, respectively $\tilde H_\lambda$; obviously for $H_\lambda$ we also have the zero-eigenvalue $\phi_0^\lambda=h_1^\lambda$. Using these new eigenfunctions, the spectral theory for $H_\lambda$ and $\tilde H_\lambda$ is similar to the corresponding one for $\lambda=1$. 

For instance, we have
\begin{equation}
\label{Flambda-def}    
\FtH_\lambda (u)(\xi) = \int_{0}^{\infty} \psi^\lambda_{\xi}(r) u(r) r dr = \frac{1}{\lambda^{3/2}} \FtH(u(\frac{\cdot}{\lambda}))(\frac{\xi}{\lambda}).
\end{equation}
Then, the Fourier inversion formula for $\FtH$ gives
\begin{equation}\label{fhtildelambdainverse}\FtH^{-1}_\lambda (u)(r) = \int_{0}^{\infty} \psi^\lambda_{\xi}(r) u(\xi) d\xi.
\end{equation}
Both integrals above are considered in the singular sense.

\subsection{The generalized eigenfunctions of \texorpdfstring{$H$}{}: Proof of Theorem \ref{htransthm}
}

We will carry out the proof in three steps:

\begin{description}
\item[STEP 1] We construct $\phi_\xi$ in the region $r \xi \lesssim 1$, modulo the determination of the normalization coefficient $q(\xi)$. 
We also complete $\phi_\xi$ to a basis of solutions for \eqref{phi-xi} by constructing a second solution $\theta_{\xi}$ satisfying a Wronskian condition.

\item[STEP II]  We construct $\phi_\xi^+$ in the region $r \xi \gtrsim 1$, modulo the determination of the normalization coefficient $a(\xi)$. 

\item[STEP III] We determine the choice of the normalization coefficients $q$ and $a$ so that they match, and so that we have 
the $L^2$ isometry property.
\end{description}
\bigskip

{\bf STEP I:} \emph{Generalized eigenfunctions in the region $r \xi \lesssim 1$.}
We consider a basis $(\phi_0,\theta_{0})$ in the null space of $H$,
\[
H\phi_{0}=H\theta_{0}=0,
\]
where $\theta_{0}$ is unbounded at $r= 0$, and given by
\begin{equation}\label{theta0}
\theta_{0}(r) = \frac{-r^8+r^4-8 r^4 \log (r)+1}{8 \left(r^6+r^2\right)} = \frac{1}{8r^{2}} + \frac{r^{2}}{8} \left(-1+\frac{1}{1+r^{4}}\right) - \frac{\log(r)}{2} \phi_{0}(r).
\end{equation}
With this choice of $\theta_{0}$, we have the Wronskian normalization
\[
r\left(\partial_r \phi_{0}(r)\theta_{0}(r)-\phi_{0}(r)\partial_{r}\theta_{0}(r) \right) =1.
\]
We use these two functions to define a solution operator for the inhomogeneous problem 
\[
H \psi  =  g.
\]
Precisely, for $g$ vanishing quadratically at $0$,  
the unique solution $\psi$ which also vanishes of order four at $r=0$ is given by
the variation of parameters formula
\begin{equation}\label{partsoln}
\psi(r) = \cT(f)(r):= -\phi_{0}(r) \int_{0}^{r} F(s) \theta_{0}(s) s ds + \theta_{0}(r) \int_{0}^{r} F(s) \phi_{0}(s) s ds.
\end{equation}
We rewrite this in the form 
\begin{equation}
\cT(f)(r) = \int_{0}^r  T(r,s) f(s) s ds,
\end{equation}
where the kernel $T$ is given by 
\begin{equation}\label{T-kernel}
 T(r,s) =    \left(\theta_{0}(r)\phi_{0}(s)-\phi_{0}(r)\theta_{0}(s)\right) = \frac{(s^{4}-r^{4})(1+r^{4}s^{4})+8 r^4 s^4 (\log (\frac{s}{r}))}{4 r^2 s^2 \left(r^4+1\right) \left(s^4+1\right)},
\end{equation}
and has the property
\[
 T(r,s) \leq 0, \quad 0 < s\leq r. 
\]

 We can use this to iteratively compute the functions 
$\phi_j$ in \eqref{repphi}. Given the equation 
\eqref{phi-xi}, these functions
must satisfy the recurrence relations
\begin{equation}\label{induction-phij}
\begin{cases} H(r^{2}\phi_{1}(r^{2}))=\phi_{0}(r)\\
H(r^{2j}\phi_{j}(r^{2})) = r^{2(j-1)} \phi_{j-1}(r^{2}), \quad j \geq 2.
\end{cases}
\end{equation}
This implies that we must have
\begin{equation}
\begin{cases} r^{2}\phi_{1}(r^{2})=\cT(\phi_{0})(r)\\
r^{2j}\phi_{j}(r^{2}) = \cT((\cdot)^{2(j-1)} \phi_{j-1}((\cdot)^{2}))(r), \quad j \geq 2,
\end{cases}
\end{equation}
which we take as the inductive definition of the functions $\phi_j$. 
Our next objective is to show that these functions satisfy the bounds \eqref{phijest}
and \eqref{drphijest} in the theorem.
Once this is done, we can conclude that 
$\phi_\xi$ given by \eqref{repphi} solves indeed the eigenfunction equation \eqref{phi-xi}.

\bigskip
For $j = 1$ we have
\begin{equation}\label{phi1formula}
\phi_{1}(u)=-\frac{2 u^2 \int _0^u\frac{2 \tan ^{-1}(x)}{x}dx+\left(u^4-1\right) \tan ^{-1}(u)-3 u^3+u}{8 \left(u^4+u^2\right)} \ .
\end{equation}
For $j \geq 2$, we have
\[
r^{2j}\phi_{j}(r^{2}) = \int_{0}^{r} s^{2(j-1)} \phi_{j-1}(s^{2})  T(r,s) s ds,
\]
which gives
\[
u^{j}\phi_{j}(u) = \int_{0}^{u} y^{j-1}\phi_{j-1}(y)\frac{T(\sqrt{u},\sqrt{y})}{2} dy.
\]
By \eqref{phi1formula}, there exists $C_{1}>0$ such that
\[
|\phi_{1}(u)| \leq C_{1}\frac{u}{\langle u \rangle}, \quad u \geq 0,
\]
which gives \eqref{phijest} for $j = 1$. For $j\geq 1$, we 
let $C_{j} =\frac{C_{1}}{((j-1)!)^{2}} \left(\frac{3}{2^{3/2}}\right)^{j-1}$
and prove \eqref{phijest} inductively.
If we suppose, for some $j \geq 1$, that 
\[
|\phi_{j}(u)| \leq C_{j} \frac{u}{\langle u \rangle}, \qquad u>0,
\]
then
\[
u^{j+1}\phi_{j+1}(u) = \int_{0}^{u} y^{j} \phi_{j}(y) \left(\frac{(y^{2}-u^{2})(1+u^{2}y^{2})+4 u^{2}y^{2} \log(\frac{y}{u})}{8 u y (u^{2}+1)(y^{2}+1)}\right)dy.
\]
If $u \leq 1$, we estimate this by
\begin{equation}
\begin{split}
|u^{j+1}\phi_{j+1}(u)| &\leq \int_{0}^{u} y^{j} C_{j} y \left(\frac{u^{2}-y^{2}}{8 u y} + \frac{u y \log(\frac{u}{y})}{2}\right) dy \\
&\leq C_{j}\left(\frac{1}{8 u} \frac{2 u^{j+3}}{(3+4j+j^{2})} + u^{j+4} \int_{0}^{1} \frac{x^{j+2}}{2} \log(\frac{1}{x}) dx\right)\\
&\leq \frac{3}{4j^{2}} C_{j}u^{j+2}, 
\end{split}
\end{equation}
from which it follows that
\[
|\phi_{j+1}(u)| \leq \frac{3 \sqrt{2}}{4 j^{2}} C_{j} \frac{u}{\langle u \rangle}, 
\quad u \leq 1.
\]
On the other hand, if $u \geq 1$, then we have
\begin{equation}
\begin{split} 
|u^{j+1}\phi_{j+1}(u)| &\leq  \int_{0}^{u} y^{j} C_{j} \left(\frac{(u^{2}-y^{2})}{8 u y} +\frac{1}{2u y} \log(\frac{u}{y})\right) dy\\
&\leq \frac{C_{j} u^{j+1}}{4 j^{2}}+\frac{C_{j} u^{j-1}}{2j^{2}}, \quad u  \geq 1. 
\end{split}
\end{equation}
Combining the two cases, we get
\[
|\phi_{j+1}(u)| \leq  \frac{3 C_{j} \sqrt{2}}{4 j^{2}} \frac{u}{\langle u \rangle}=C_{j+1} \frac{u}{\langle u \rangle}, \quad u >0,
\]
which completes the inductive proof of \eqref{phijest}. 
\medskip

Next, we consider higher order derivatives of $\phi_j$ in order to prove \eqref{drphijest}.
From \eqref{phi1formula} we have
\[
|\phi_{1}'(u)| \leq \frac{C}{\langle u \rangle}, \quad u>0
\]
For $j\geq 2$ we compute
\begin{equation}
\partial_{u}\left(u^{j}\phi_{j}(u)\right) 
= -\int_{0}^{u} y^{j-1}\phi_{j-1}(y) h^1(y,u)\,dy 
\end{equation}
where 
\[
h^1(y,u) = \frac{y^{2}+u^{2}(1-u^{2}+(7+7u^{2}+u^{4})y^{2}+(-1+u^{2})y^{4})+4u^{2}(-1+u^{2})y^{2}\log(\frac{y}{u})}{8 u^{2}(1+u^{2})^{2}y(y^{2}+1)}.
\]
Then, we use \eqref{phijest} to  directly estimate the integral. This yields
\[
|\partial_{u}\left(u^{j}\phi_{j}(u)\right)| \leq \frac{C C_{1}}{((j-2)!)^{2}(j-1)} \left(\frac{3}{2^{3/2}}\right)^{j-2}\frac{u^{j}}{\langle u \rangle}.
\]
Hence, for all $j \geq 1$ we have
\[
|\phi_{j}'(u)| \leq \frac{C C_{1}}{((j-1)!)^{2}} \left(\frac{3}{2^{3/2}}\right)^{j-1}\frac{j}{\langle u \rangle}
\]
which is the case $n=1$ of \eqref{drphijest}.

For larger $n$ we denote
\[
f_{j}(r) = r^{2j} \phi_{j}(r^{2}).
\]
Then by \eqref{induction-phij} we have
\[
f_{j}''(r) = \frac{-1}{r}f_{j}'(r) + \frac{4(1-6r^{4}+r^{8})f_{j}(r)}{r^{2}(1+r^{4})^{2}} - f_{j-1}(r),
\]
which implies
\[
|f_{j}''(r)| \leq \frac{C C_{1} j^{2}}{((j-1)!)^{2}} \left(\frac{3}{2^{3/2}}\right)^{j-1} \frac{r^{2j}}{\langle r^{2}\rangle},
\]
thereby completing the proof of \eqref{drphijest} for $n=2$.

For larger $n$, if we inductively assume that
\[
|f_{j}^{(n)}(r)| \leq \frac{M_{n} j^{n}}{((j-1)!)^{2}} \left(\frac{3}{2^{3/2}}\right)^{j-1} \frac{r^{2j+2}}{\langle r^{2}\rangle r^{n}}
\]
for all $0 \leq n \leq k+1$ for some $k \geq 1$, then, using
\[
f_{j}^{(k+2)}(r) = -\sum_{l=0}^{k} \binom{k}{ l} \partial_{r}^{k-l} \left(\frac{1}{r}\right) f_{j}^{(l+1)}(r) + \sum_{l=0}^{k} \binom{k}{l} \partial_{r}^{k-l}\left(\frac{4(1-6r^{4}+r^{8})}{r^{2}(1+r^{4})^{2}}\right) f_{j}^{(l)}(r) - f_{j-1}^{(k)}(r)
\]
we obtain
\[
|f_{j}^{(k+2)}(r)| \leq \frac{M_{k+2} j^{k+2}}{r^{k+2} ((j-1)!)^{2}} \left(\frac{3}{2^{3/2}}\right)^{j-1} \frac{r^{2j+2}}{\langle r^{2}\rangle},
\]
and this completes the proof of \eqref{drphijest} by induction.

This finishes the analysis in part (a) of the theorem, except for the choice of the normalization coefficient $q(\xi)$. 
This choice will only be considered in the last step of the proof, where we connect the solutions of \eqref{phi-xi} near $r=0$ with the solutions
near $r = \infty$. In order to accomplish this, we supplement the above solution $\phi_\xi$ 
with a second solution $\theta_{\xi}$,

\[
H(\theta_{\xi}) = \xi^{2}\theta_{\xi}
\]
with the Wronskian normalization
\begin{equation}\label{W-phi-theta}
r\left(\partial_r \phi_{\xi}(r)\theta_{\xi}(r)-\phi_{\xi}(r)\partial_{r}\theta_{\xi}(r) \right) =1
\end{equation}
which determines $\theta_{\xi}$ uniquely up to a multiple of $\phi_\xi$:
\begin{l1} In the region $r \xi \lesssim 1$ one could choose
$\theta_{\xi}$ of the form
$$\theta_{\xi}(r) = \frac{1}{8 q(\xi) r^{2}}+\frac{\xi^{2}}{32 q(\xi)} + \frac{1}{8 q(\xi)r^{2}} \sum_{k=0}^{\infty} \widetilde{\theta_{k}}(r) r^{2k} \xi^{2k} - \frac{1}{q(\xi)^{2}} \left(\frac{1}{2}+\frac{\xi^{4}}{256}\right)\log(r)\phi_{\xi}(r)$$
where we have
\[\widetilde{\theta_{0}}(r) = \frac{-r^{8}}{1+r^{4}}\] and
the bounds
\[
|\widetilde{\theta_{k}}(r)| \leq \begin{cases} C r^{4}, \quad k=1,2\\
\frac{T^{k}}{k!}(1+r^{4}), \quad k \geq 3\end{cases}.
\]

\end{l1}

Note that our definition of $\theta_{0}$ is such that
\[
\lim_{\xi \rightarrow 0^{+}}(\theta_{\xi}(r) q(\xi)) =\theta_{0}(r)
\]
The $\log r$ factor serves to provide a constant term in $H \theta_\xi$; all other contributions 
should be analytic in $r$ at $r = 0$.  This exactly 
corresponds to a similar term in $\theta_0$, see \eqref{theta0}. 
\begin{proof} 
It is natural to try the ansatz
\[
\theta_{\xi}(r) = \frac{c_{0}(\xi)}{r^{2}}+c_{1}(\xi) + c_{0}(\xi) r^{2}\sum_{k=0}^{\infty} d_{k}(r) r^{2k}\xi^{2k}+c_{2}(\xi) \log(r) \phi_{\xi}(r)
\]
where $d_k$ are analytic in $r^2$ and all contributions from the sum are at least $O(r^4)$ at $r = 0$.
We have
\begin{equation}\begin{split} H(\theta_{\xi}(r))&=\frac{-32}{(1+r^{4})^{2}} c_{0}(\xi) + c_{1}(\xi) \cdot \frac{4(1-6r^{4}+r^{8})}{r^{2}(1+r^{4})^{2}} + H\left(r^{2}\sum_{k=0}^{\infty} d_{k}(r) r^{2k}\xi^{2k}\right) c_{0}(\xi)\\
&+c_{2}(\xi) \left(\log(r)\xi^{2}\phi_{\xi}(r)-\frac{2}{r}\partial_{r}\phi_{\xi}(r)\right)\end{split}\end{equation}
and
\begin{equation}
\xi^{2}\theta_{\xi}(r) = \frac{\xi^{2}c_{0}(\xi)}{r^{2}}+\xi^{2}c_{1}(\xi)+r^{2}\xi^{2}\sum_{k=0}^{\infty} d_{k}(r) r^{2k}\xi^{2k} c_{0}(\xi) + \xi^{2}c_{2}(\xi) \log(r)\phi_{\xi}(r)
\end{equation}
Then the equation $H\theta_{\xi} =\xi^{2}\theta_{\xi}$ implies
\begin{equation}
\begin{split}
&-\frac{32c_0(\xi)}{(1+r^{4})^{2}} +c_{1}(\xi)\left(\frac{4}{r^{2}}-\frac{32 r^{2}}{(1+r^{4})^{2}}\right)+ c_0(\xi) H\left(r^{2}\sum_{k=0}^{\infty} d_{k}(r) r^{2k}\xi^{2k}\right) 
\\
&=\frac{\xi^{2}c_{0}(\xi)}{r^{2}}+\xi^{2}c_{1}(\xi)+r^{2}\xi^{2}\sum_{k=0}^{\infty} d_{k}(r) r^{2k}\xi^{2k}c_{0}(\xi) 
+ c_{2}(\xi)\frac{2}{r}\partial_{r}\phi_{\xi}(r)
\end{split}
\end{equation}
If we look for a solution for which 
\[
\lim_{r \rightarrow 0}H\left(r^{2}\sum_{k=0}^{\infty} d_{k}(r) r^{2k}\xi^{2k}\right)=0
\]
then, by examining the $\dfrac{1}{r^{2}}$ terms on each side, we obtain
\[
c_{1}(\xi) = \frac{\xi^{2}c_{0}(\xi)}{4}.
\]
By considering the limit as $r \rightarrow 0^{+}$ of both sides, we get 
\[
-32 c_{0}(\xi) -2\lim_{r\rightarrow 0} \left(c_{2}(\xi) q(\xi) \frac{\phi_{0}'(r)}{r}\right) = \xi^{2}c_{1}(\xi)
\]
and this gives
\[
c_{2}(\xi) = -\left(4+\frac{\xi^{4}}{32}\right)\frac{c_{0}(\xi)}{q(\xi)}.
\]
Finally, the coefficient $c_0$ is determined by examining
the Wronskian relation \eqref{W-phi-theta} in the limit as $r$
approaches $0$. This gives
\[
c_{0}(\xi) = \frac{1}{8 q(\xi)}.
\]

Then we can factor out all the $c_0$ factors, and the equation $H\theta_{\xi} = \xi^{2}\theta_{\xi}$ reduces to

\begin{equation}\begin{split}&\frac{32  r^{4}(2+r^{4})}{(1+r^{4})^{2}} - \frac{8 \xi^{2} r^{2}}{4 (1+r^{4})^{2}} + H\left(r^{2}\sum_{k=0}^{\infty} d_{k}(r) r^{2k}\xi^{2k}\right) \\
& = r^{2}\xi^{2} \sum_{k=0}^{\infty} d_{k}(r) r^{2k}\xi^{2k} 
- \left(4+\frac{\xi^{4}}{32}\right)\left( \frac{2}{r q(\xi) }\partial_{r  }\phi_{\xi}(r) - 8 \right).
\end{split}\end{equation}
At this point we use the expansion \eqref{repphi} for $\phi_\xi$,
to obtain
\begin{equation}
\begin{split} 
 H\left(r^{2} \sum_{k=0}^{\infty} d_{k}(r) r^{2k}\xi^{2k}\right)
 = &\ \frac{32 r^{4}}{(1+r^{4})^{2}} + \frac{8 \xi^{2} r^{2}}{(1+r^{4})^{2}} +\frac{\xi^4 r^{4}(3+r^{4})}{4(1+r^{4})^{2}} \\
& \hspace{-1in} + r^{2}\xi^{2} \sum_{k=0}^{\infty} d_{k}(r) r^{2k}\xi^{2k}
- \frac{2}{r}\left(4+\frac{\xi^{4}}{32}\right) \sum_{k=1}^{\infty}\partial_{r}\left( r^{2k} \phi_{k}(r^{2})\right)\xi^{2k}.
\end{split}
\end{equation}
Now we simply match the coefficients of $\xi^{2j}$. By considering the $\xi^{0}$ terms, we get
\[
H\left(r^{2}d_{0}(r)\right) =  \frac{32 r^{4}(3+r^{4})}{(1+r^{4})^{2}} -  \frac{32 r^{4}(2+r^{4})}{(1+r^{4})^{2}} 
\]
which gives
\[
d_{0}(r) = -1+\frac{1}{1+r^{4}}.
\]
By examining the $\xi^{2}$ terms, we obtain
\[
  H\left(r^{4}d_{1}(r)\right) =\frac{8 r^{2}}{(1+r^{4})^{2}} - \frac{8}{r} \partial_{r}\left(r^{2} \phi_{1}(r^{2})\right) + r^{2}d_{0}(r).
\]
The right hand side has size $O(r^2)$ at zero, so the
above equation has an unique solution of size $O(r^4)$, namely
the one given by \eqref{partsoln},
\[
r^{4}d_{1}(r) = \int_{0}^{r} \left(\frac{-s^{2}(-8+s^{4}+s^{8})}{(1+s^{4})^{2}}-\frac{8}{s} \partial_{s}\left(s^{2} \phi_{1}(s^{2})\right)\right) T(r,s) s ds.
\]
The integrand is an analytic, odd function of $s$, except for the log factor.  But the log factor can be eliminated integrating by parts, for instance
\[
\int_0^r s^{2k+1} \log(s/r) ds = - \frac{1}{2k+2} \int_0^r s^{2k+1}  ds =  - \frac{1}{(2k+2)^2} r^{2k+2}.
\]
So we obtain an analytic, even function of $r$. Using the estimates on $\phi_{1}$ we get
\[
|d_{1}(r)| \lesssim 1, \quad r>0.
\]
Similarly, we consider the $\xi^{4}$ terms, to get
\[
H\left(r^{6}d_{2}(r)\right) = r^{4}d_{1}(r) - \frac{8}{r}\partial_{r}\left( r^{4}\phi_{2}(r^{2})\right) + \frac{r^{4}}{4} \frac{(3+r^{4})}{(1+r^{4})^{2}}
\]
which is solved again using \eqref{partsoln} with the kernel $T$ as in \eqref{T-kernel} to obtain
\[
r^{6}d_{2}(r) =  \int_{0}^{r} \left(s^{4}d_{1}(s) - \frac{8}{s} \partial_{s}\left(s^{4}\phi_{2}(s^{2})\right) + \frac{s^{4}}{4} \frac{(3+s^{4})}{(1+s^{4})^{2}}\right) T(r,s) s ds.
\]
We obtain a function $d_2$ which is even, analytic, and
using the estimates on $\phi_{j}$ and $d_{1}$ above, we get
\[
|d_{2}(r)| \lesssim 1, \quad r>0.
\]
Finally, we let
\[
g_{k}(r) = r^{2k+1}d_{k}(r), \quad k \geq 2,
\]
and note that $H\theta = \xi^{2}\theta$ is completely solved if, for all $k \geq 3$,
\begin{equation}\label{gkeqn}
\begin{split}
r g_{k}(r)
= \int_{0}^{r}\left(s g_{k-1}(s) -\frac{1}{16s}\partial_{s}\left( s^{2k-4} \phi_{k-2}(s^{2})\right)-\frac{8}{s} \partial_{s}\left(s^{2k}\phi_{k}(s^{2})\right)\right) T(r,s) s ds.
\end{split}
\end{equation}
Using the $\phi_{j}$ estimates and  the $d_{2}$ estimate, we obtain
\[
|g_{3}(r)| \leq \widetilde{C_{3}} r^{3}(1+r^{4}).
\]
If we suppose, for some $k-1 \geq 3$, that 
\[
|g_{k-1}(r)| \leq \widetilde{C_{k-1}} r^{2k-5}(1+r^{4}).
\]
then, \eqref{gkeqn} gives the following. There exists $C>0$ so that
\[
|g_{k}(r)| \leq C r^{2k-3}(1+r^{4}) \left(\frac{C^{k}}{k!}+\frac{\widetilde{c_{k-1}}}{k}\right).
\]
Therefore, there exists $T$ sufficiently large so that, by induction, for all $k \geq 3$, and $r>0$,
\[
|g_{k}(r)| \leq \frac{T^{k}}{k!} r^{2k-3} (1+r^{4}).
\]
Setting $\widetilde{\theta_{k}}(r) = d_{k}(r) r^{4}$
yields
\[\widetilde{\theta_{0}}(r) = \frac{-r^{8}}{1+r^{4}}, \quad |\widetilde{\theta_{k}}(r)| \leq \begin{cases} C r^{4}, \quad k=1,2\\
\frac{T^{k}}{k!}(1+r^{4}), \quad k \geq 3\end{cases}\]
which finishes the proof of the claim. 

\end{proof}

This completes the first step of the proof, i.e. the analysis of the $r \xi \lesssim 1$ region.

\bigskip

{\bf STEP II:} \emph{Generalized eigenfunctions in the region $r \xi \gtrsim 1$.}  Here our aim is to construct the  generalized eigenfunctions
\begin{equation} \label{repphi+}
\phi^{+}_\xi(r)= r^{-\frac12} e^{ir\xi} \sigma(r\xi,r), 
\qquad r\xi \gtrsim  1
\end{equation}
solving 
\[
H \phi^{+}_\xi = \xi^2 \phi^+_\xi.
\]
Substituting, this yields a second order
equation for the amplitude $\sigma$:
\begin{equation}\label{sigmaeqn}
\begin{aligned}
0 = & \ e^{-i r \xi} \sqrt{r}\left(H-\xi^{2}\right)\left(\frac{e^{i r \xi} \sigma(r\xi,r)}{\sqrt{r}}\right) 
\\
= 
& \ - \partial_{r}^{2}\left(\sigma(r\xi,r)\right)   -2 i \xi \partial_{r}\left(\sigma(r\xi,r)\right)+ V(r) \sigma(r\xi,r), 
\end{aligned}
\end{equation}
where we recall, see \eqref{lcaldef}, that
\[
V(r)=\frac{15}{4r^{2}}-32\frac{r^{2}}{(1+r^{4})^{2}}.
\]
We first look for a WKB type asymptotic series expansion of $\sigma(r\xi,r)$ in the region of large $r\xi$. This is derived by considering the system of equations arising from formally substituting 
\[
\sigma(r\xi,r)= \sum_{j=0}^{\infty} f_{j}(r) \xi^{-j}, \qquad f_{0}(r)=1
\]
into \eqref{sigmaeqn}.  This yields
\[
f_{k+1}'(r) =\frac{1}{2i}\left(V(r)f_{k}(r)-f_{k}''(r)\right), \quad k \geq 0,
\]
where we look for the unique solutions $f_{k}$ satisfying a boundary condition at infinity, 
\[
\lim_{r \rightarrow \infty}f_{k}(r) =0, \quad k \geq 1.
\]
Setting 
\begin{equation}\label{phiplusintermsoff}\phi_{k}^{+}(r) = r^{k}f_{k}(r)\end{equation}
we compute directly
\[
\phi_{1}^{+}(r) = \frac{-4 i}{r^4+1}-\frac{i r \left(\log \left(\frac{r^{2}+\sqrt{2}r+1}{r^2-\sqrt{2} r+1}\right)+2(\pi-\tan ^{-1}\left(\sqrt{2} r-1\right)- \tan ^{-1}\left(\sqrt{2} r+1\right)) \right)}{\sqrt{2}}+\frac{47 i}{8}.
\]
In particular we note that
\[
\phi_{1}^{+}(r) = \frac{15 i}{8}-\frac{16 i}{5 r^4} + O\left(\frac{1}{r^{8}}\right), \quad r \rightarrow \infty.
\]
By inspection one easily sees that there exist constants $C_{m}>0$ such that 
\[
|(r \partial_{r})^{m}\phi_{1}^+(r)| \leq C_{m}, \quad r>0.
\]
Then, by induction, it easily follows for all $m \geq 0$, $j \geq 0$ there exist constants $C_{m,j}>0$ such that 
\[
|(r \partial_{r})^{m}\phi_{j}^{+}(r)| \leq C_{m,j}, \quad r>0.
\]
Overall,  we thus have the asymptotic expansion
\[
\sigma(y,r) \sim \sum_{j=0}^{\infty} y^{-j}\phi_{j}^{+}(r).
\]
The difficulty at this point is that we do not know that this series converges. To rectify this, as usual in WKB approximations,
we construct an approximate sum by 
truncating the terms in the series at well chosen points. Precisely, we define $\sigma_{ap}$ by
\[
\sigma_{ap}(y,r) = \sum_{k=0}^{\infty} y^{-k} \chi_{\leq 1}(\frac{1}{y \delta_{k}}) \phi_{k}^{+}(r),
\]
where 
\[
\chi_{\leq 1}(x) = \begin{cases} 1, \quad x \leq \frac{1}{2}\\
0, \quad x \geq 1\end{cases}, \quad \chi_{\leq 1}\in C^{\infty}(\mathbb{R}),
\]
and $\delta_{k} \rightarrow 0$ rapidly enough 
in order to ensure convergence.

Comparing the approximate sum with the exact partial sums, we have that for  each $N,m,n \geq 0$, 
\begin{equation}\label{asympsum}
\left|(y \partial_{y})^{m}\left(r \partial_{r}\right)^{n}\left(\sum_{k=0}^{\infty} y^{-k} \chi_{\leq 1}(\frac{1}{y \delta_{k}}) \phi_{k}^{+}(r) - \sum_{k=0}^{N} y^{-k}\phi_{k}^{+}(r)\right)\right| \leq \frac{C_{N,n,m}}{y^{N+1}}, \quad r>0, \ y \geq 1.
\end{equation}
This implies that $\sigma_{ap}$ is a good approximate solution for \eqref{sigmaeqn}
near infinity, in the sense that its associated source term
\[
e(r\xi,r) = -\left(\partial_{r}^{2}+2 i \xi \partial_{r}-V(r)\right) \left(\sigma_{ap}(r\xi,r)\right)
\]
is rapidly decreasing at infinity, i.e.
for all $N,m,n \geq 0$,
\begin{equation}\label{eest}
\left|(y \partial_{y})^{m} (r\partial_{r})^{n}  e(y,r)\right| \leq \frac{C_{N,n,m}}{r^{2}y^{N}}, \quad r>0, \ y \geq 1.
\end{equation}
We then need to complete $\sigma_{ap}$ to an exact solution to \eqref{sigmaeqn}. Writing
\[
\sigma(r\xi,r) = \sigma_{ap}(r\xi,r) + f(r\xi,r),
\]
the correction $f$ should solve 
\begin{equation} \label{f-e}
\left(-\partial_{r}^{2}-2 i \xi \partial_{r}+V(r)\right) f(r\xi,r)
= e(r\xi,r).
\end{equation}

We will show that, given $e$ as in \eqref{eest}, there exists a unique solution $f$ satisfying 
the bounds 
\begin{equation}\label{fest}
\left|(y \partial_{y})^{m} (r\partial_{r})^{n}  f(y,r)\right| \leq \frac{C_{N,n,m}}{y^{N}}, \quad r>0, \ y \geq 1.
\end{equation}
To achieve this we interpret the ordinary differential equation \eqref{f-e} as a Cauchy problem with zero 
Cauchy data at infinity. Here it is convenient to use the variable $\xi$ as a parameter, instead of $y$ 
which also contains some $r$ dependence.  By the chain rule we can interchange the operators $(r\partial_r, y \partial_y)$ with $(r \partial_r,\xi \partial_\xi)$.  Then, changing notations to view 
$f$ and $e$ as functions of $r$ and $\xi$, it remains to prove the following

\begin{p1}
Consider the equation
\begin{equation} \label{f-e+}
\left(-\partial_{r}^{2}-2 i \xi \partial_{r}+V(r)\right) f(\xi,r)
= e(\xi,r)
\end{equation}
Assume that $e$ satisfies the bounds
\begin{equation}\label{eest+}
\left|(\xi \partial_{\xi})^{m} (r\partial_{r})^{n}  e(\xi,r)\right| \leq \frac{C_{N,n,m}}{r^{2}(r\xi)^{N}}, \quad r>0, \ r\xi \geq 1.
\end{equation}
Then there exists a unique solution satisfying the bound
\begin{equation}\label{fest+}
\left|(\xi \partial_{\xi})^{m} (r\partial_{r})^{n}  f(\xi,r)\right| \leq \frac{C_{N,n,m}}{(r\xi)^{N}}, \quad r>0, \ r\xi \geq 1.
\end{equation}
\end{p1}

\begin{proof}

We first drop the vector fields $(r \partial_r,\xi \partial_\xi)$ and consider our problem for fixed $\xi$.
We summarize the result in the following Lemma:

\begin{l1}
Suppose $k > k_0$ and
\[
|e(r)| \leq r^{-k-2}, \qquad r\xi > 1.
\]
Then the equation \eqref{f-e+} admits a unique 
solution $f$ which satisfies the bound 
\begin{equation}\label{what-f}
|f(r)| \lesssim_k r^{-k}, \qquad |\partial_r f(r)| \lesssim_k r^{-k-1}, \qquad |\partial_r^2 f(r)| \lesssim_k r^{-k-2}
\end{equation}
\end{l1}

\begin{proof}
The heart of the proof is an energy estimate 
which holds for any local solution to \eqref{f-e+}. This has the form
\begin{equation}
	\label{fenest0}\partial_{r}\left((|f|^{2}+r^{2}|\partial_{r}f|^{2})r^{C}\right) \geq -C r^{3+C} |e(r)|^{2}
\end{equation}
with a universal $C$.  This directly gives uniqueness of solutions as in \eqref{what-f}. 
Existence is obtained in a standard manner, by solving a truncated problem
\begin{equation} \label{f-e-cut}
\left(-\partial_{r}^{2}-2 i \xi \partial_{r}+V(r)\right) f_{r_0}(r)
= 1_{\{r < r_0\}}e(r)
\end{equation}
The desired solution $f$ is then obtained as the limit of $f_{_0}$ as $r_0 \to \infty$.

For later use we remark that we could also replace the sharp truncation with a smooth truncation on the 
dyadic scale. This has the advantage that also preserves the vector field bounds \eqref{eest+}.
\end{proof}

Now we use the lemma to conclude the proof of the proposition. Since we restrict ourselves to the range
$r \xi \geq  1$, it suffices to work with $N$ large enough. This provides sufficient decay for $e$ 
to allow us to apply the lemma.
Then the uniqueness of $e$ is immediate. Further, the discussion above, it suffices 
to prove it under the qualitative assumption that $e$ 
has compact support. Then the smoothness of 
$f$ as a function of both $r$ and $\xi$ directly follows, and it remains to establish the bounds
\eqref{fest+}.

We prove these bounds by induction on $m$. For $m = 0$
we start with the bounds in the lemma and get the 
higher $r$ derivatives directly from the equation. 
Suppose now, by induction, that we have the bounds
\eqref{fest+} for $m=k$, and prove them for $m = k+1$.
Denoting 
\[
f_m = (\xi \partial_\xi)^m f, \qquad 
e_m = (\xi \partial_\xi)^m e
\]
we write an equation for $f_{k+1}$:
\[
\left(-\partial_{r}^{2}-2 i \xi \partial_{r}+V(r)\right) f_{k+1}(\xi,r)
= e_{k+1}(\xi,r) + 2i(k+1)\xi \partial_r f_k
\]
By the induction hypothesis we can include
the second term on the right into $e_{k+1}$.
This allows us to conclude by directly applying
the $m = 0$  result.
\end{proof}

{\bf STEP 3:} The normalization coefficients $a$ and $q$.
Now, we compute the asymptotic behavior of $a$ and $q$. For this, it is convenient to follow the procedure of \cite{kstym}. For comparison purposes we note that the theory developed in \cite{kstym}, which uses also parts from \cite{KrScTa}, has three normalization features that differ from our setup: 
\begin{itemize}
\item the corresponding eigenfunctions $\phi_\xi$ there are used for the eigenvalues $\xi$ versus $\xi^2$;

\item the calculus on the physical side is with respect to the measure $dr$ versus $rdr$ in our case; this corresponds to  the conjugation of $H$ by $\sqrt{r}$ in \cite{kstym}, which yields the operator $\mathcal{L}_{KST}$.
\item  the calculus on the Fourier side is with respect to a measure $\rho(\xi) d\xi$ versus $d\xi$
in our case; here we compensate this by using the additional coefficient $q(\xi)$ in the expression of our eigenfunctions $\phi_\xi(r)$.
\end{itemize}
Below we modify the elements of our spectral theory so as to fit the setup developed in \cite{kstym} and \cite{KrScTa} (which in part was based on the setup in \cite{GeZi}). This allows us to use that theory and then by reversing the process we can recover $a(\xi)$ and $q(\xi)$ in our setup.

Therefore, we define the following functions, which are counterparts of the eigenfunctions in Theorem 4.3 of \cite{kstym}. In other words, repeating the identical procedure of \cite{kstym} for our equation yields the functions $\phi_{KST}(r,\xi), \theta_{KST}(r,\xi), \psi^{+}_{KST}(r,\xi)$ given below, which correspond to the functions $\phi(r,\xi), \theta(r,\xi), \psi^{+}(r,\xi)$ from Theorem 4.3 of \cite{kstym}. We let
\[
\phi_{KST}(r,\xi) = \phi_{0,KST}(r)+r^{-\frac{3}{2}} \sum_{j=1}^{\infty}(r^{2}\xi)^{j} \phi_{j,KST}(r^{2})
\]
with
\[
\phi_{0, KST}(r) = \frac{\sqrt{r} \phi_{0}(r)}{2} = \frac{r^{5/2}}{1+r^{4}}.
\]
We also set
\[\phi_{j,KST}(u) = \frac{u}{2} \phi_{j}(u), \quad \theta_{KST}(r,\xi) = 2 q(\sqrt{\xi}) \sqrt{r}\theta_{\sqrt{\xi}}(r)\]
and
\[
\psi_{KST}^{+}(r,\xi) = \frac{\sqrt{r}}{\xi^{1/4}} \phi_{\sqrt{\xi}}^{+}(r) = \frac{e^{i r \sqrt{\xi}}}{\xi^{1/4}} \sigma(r\sqrt{\xi},r).
\]
With these notations we have
\[ \mathcal{L}(\phi_{KST}(\cdot,\xi))(r) = \xi \phi_{KST}(r,\xi) \]
\[ \mathcal{L}(\theta_{KST}(\cdot,\xi)) = \xi \theta_{KST}(r,\xi) \]
\[ \mathcal{L}(\psi^{+}_{KST}(\cdot,\xi)) = \xi \psi_{KST}^{+}(r,\xi), \]
where $\mathcal{L}$ is given in  \eqref{lcaldef}.
From our functions $\theta$ and $\phi$ we inherit the Wronskian normalization
\[
W(\theta_{KST}(\cdot,\xi),\phi_{KST}(\cdot,\xi))=1.
\]
On the other hand, using the above formulae for $\psi^{+}_{KST}$ and $\psi^{-}_{KST}=\overline{\psi^{+}_{KST}}$, we obtain
\[
W(\psi^{+}_{KST},\psi^{-}_{KST}) = -2i.
\]
Finally, since
\[
\phi_{KST}(r,\xi) = \frac{W(\psi^{-}_{KST},\phi_{KST})}{W(\psi^{-}_{KST},\psi^{+}_{KST})} \psi^{+}_{KST}(r,\xi) + \frac{W(\phi_{KST},\psi^{+}_{KST})}{W(\psi^{-}_{KST},\psi^{+}_{KST})} \psi_{KST}^{-}(r,\xi)
\]
and $\phi_{KST}$ is real-valued, while $\psi^{-}_{KST}=\overline{\psi^{+}_{KST}}$, we obtain
\[
\phi_{KST}(r,\xi) = a_{KST}(\xi)\psi^{+}_{KST}(r,\xi) + \overline{a_{KST}(\xi) \psi^{+}_{KST}(r,\xi)}
\]
with
\[
a_{KST}(\xi) = \frac{i W(\phi_{KST},\psi^{-}_{KST})}{2}.
\]
By using, for instance, Faa di Bruno's formula and the symbol-type estimates on $\sigma$, we obtain that, for all $\alpha \geq 0$, and $r \sqrt{\xi} \sim 1$,
\[
|\partial_{\xi}^{\alpha} \phi_{KST}| \leq \frac{C_{\alpha} r^{5/2}}{\xi^{\alpha}\langle r^{2}\rangle}, \quad |\partial_{\xi}^{\alpha} \psi^{-}_{KST}(r,\xi)| \leq \frac{C_{\alpha}}{\xi^{\alpha+\frac{1}{4}}}, \quad |\partial_{\xi}^{\alpha}\partial_{r}\phi_{KST}| \leq \frac{C_{\alpha} r^{3/2}}{\langle r^{2}\rangle \xi^{\alpha}}, \quad |\partial_{\xi}^{\alpha}\partial_{r} \psi^{-}_{KST}| \leq \frac{C_{\alpha}}{\xi^{\alpha-\frac{1}{4}}}
\]
and these give, for all $\alpha \geq 0$
\[
|a_{KST}^{(\alpha)}(\xi)| \leq \frac{C_{\alpha}}{\xi^{\alpha+1} \langle \frac{1}{\xi}\rangle} \leq \frac{C_{\alpha}}{\xi^{\alpha}} \begin{cases} 1, \quad \xi \lesssim 1\\
\frac{1}{\xi}, \quad \xi \gtrsim 1\end{cases}.
\]
On the other hand, from \eqref{repphi}, we have
\[
\phi_{KST}(r,\xi) = \frac{r^{5/2}}{1+r^{4}} + \frac{r^{2}\xi}{2} \sqrt{r} \phi_{1}(r^{2}) + \sum_{j=2}^{\infty} \frac{(r^{2}\xi)^{j}}{2} \sqrt{r} \phi_{j}(r^{2}).
\]
Therefore,
\begin{equation}
\begin{split} 
\partial_{r}\phi_{KST}(r,\xi) &= \frac{r^{3/2}(5-3r^{4})}{2(1+r^{4})^{2}} + \frac{5}{4} r^{3/2} \xi \phi_{1}(r^{2}) + r^{7/2} \xi \phi_{1}'(r^{2}) \\
&+\sum_{j=2}^{\infty} \frac{(2j+\frac{1}{2})r^{2j-\frac{1}{2}}}{2} \xi^{j} \phi_{j}(r^{2}) + \sum_{j=2}^{\infty} r^{2j+\frac{1}{2}} \xi^{j} r \phi_{j}'(r^{2}).
\end{split}
\end{equation}
Now, we use the same argument as the proof of Lemma 4.7 of \cite{kstym}. In particular, using the formula for $\phi_{1}$, namely \eqref{phi1formula}, we get, for some small, fixed $\delta>0$,
\[
|\partial_{1}\phi_{KST}(\frac{\delta}{\sqrt{\xi}},\xi)| \geq \begin{cases} \frac{C}{\xi^{3/4}}, \quad \xi \gtrsim 1\\
C \xi^{1/4} ,\quad \xi \lesssim 1\end{cases}.
\]
Although \eqref{sigmaasymp} was stated for $y >1$, we still have that
\[
|\partial_{1}\psi^{-}_{KST}(\frac{\delta}{\sqrt{\xi}},\xi)| \leq C \xi^{1/4}.
\]
Then,
\begin{equation}
\begin{split}
|a_{KST}(\xi)| \geq C \frac{|\partial_{r}\phi_{KST}|}{|\partial_{r}\psi^{+}_{KST}|}\geq C \begin{cases} \frac{1}{\xi}, \quad \xi \gtrsim 1\\
1, \quad \xi \lesssim 1\end{cases}.
\end{split}
\end{equation}
In particular,
\begin{equation}\label{asymbolint}
\xi^{\alpha} |a_{KST}^{(\alpha)}(\xi)| \leq C_{\alpha}|a_{KST}(\xi)|, \quad \alpha \geq 0.
\end{equation}
As in Theorem 4.3 of \cite{kstym}, the density of the continuous part of the spectral measure is given by
\[
\rho(\xi) = \frac{1}{\pi} Im(m(\xi))\Bigr|_{\xi >0} = \frac{1}{\pi} Im\left(\frac{W(\theta,\psi^{+})}{W(\psi^{+},\phi)}\right), \quad \xi >0.
\]
As in the proof of Proposition 5.7 of \cite{KrScTa}, we get 
\begin{equation}\label{rhoexpr}
\rho(\xi) = \frac{1}{4 \pi |a_{KST}(\xi)|^{2}} \sim \begin{cases} \xi^{2}, 
\quad 1 \lesssim \xi\\
1, \quad \xi \lesssim 1\end{cases}.
\end{equation}
Also, we have
\begin{equation}\label{qexpr}q(u) = \sqrt{\frac{u \rho(u^{2})}{2}}\end{equation}
and thus
\begin{equation}
q(\xi) \sim \begin{cases} \xi^{5/2}, \quad 1 \lesssim \xi\\
\sqrt{\xi}, \quad \xi \lesssim 1\end{cases}.
\end{equation}
Since $|a_{KST}(\xi)|^{2} >0$ (which, for instance follows from the above lower bounds on $|a_{KST}|$), we get that, for all $\alpha \geq 0$,
\begin{equation}\label{reciprocalmodasymbol}|\partial_{\xi}^{\alpha}\left(\frac{1}{|a_{KST}(\xi)|}\right)| = |\partial_{\xi}^{\alpha}\left(\frac{1}{\sqrt{|a_{KST}(\xi)|^{2}}}\right)| \leq \frac{C_{\alpha}}{\xi^{\alpha}} \begin{cases} \xi, \quad 1 \lesssim \xi\\
	1, \quad \xi \lesssim 1. 
 \end{cases}\end{equation}
Using Faa di Bruno's formula, we obtain that
\[
|\partial_{\xi}^{\alpha} q(\xi)| \leq \frac{C_{\alpha}}{\xi^{\alpha}} q(\xi), 
\quad \alpha \geq 0.
\]
Then we have the representation
\begin{equation} \label{phipsi}
\phi_{\xi}(r)=a(\xi) \phi^{+}_\xi(r) + \overline{a(\xi) \phi^{+}_\xi(r)}.
\end{equation}
where the complex valued function $a$ satisfies
\begin{equation}\label{aexpr}
a(\xi) = \frac{2 q(\xi) a_{KST}(\xi^{2})}{\sqrt{\xi}}.
\end{equation}
Using \eqref{qexpr} and \eqref{rhoexpr}, we get that
\[
a(\xi) = \frac{a_{KST}(\xi^{2})}{\sqrt{2 \pi} |a_{KST}(\xi^{2})|}.
\]
Then, using \eqref{asymbolint} and \eqref{reciprocalmodasymbol}, we get that
\begin{equation} \label{abound}
|a(\xi)| =\frac{1}{\sqrt{2\pi}}, \qquad |\partial_{\xi}^{\alpha}(a(\xi))| \leq \frac{C_{\alpha}}{\xi^{\alpha}}, \quad \alpha \geq 0.
\end{equation}

\subsection{The generalized eigenfunctions of \texorpdfstring{$\tH$}{}: Proof of Theorem \ref{httransthm}}
We begin with the region $r\xi \lesssim 1$. Here, we use \eqref{psifromphi} to get
\[
\psi_\xi = \xi q(\xi) \left(\psi_0(r) +  r \sum_{j \geq 1} (r\xi)^{2j} 
{\psi}_j(r^2)\right)
\]
where
\[
{\psi}_j(r)= ( 2 h_3(\sqrt{r})+2 +2j) \phi_{j+1}(r) + 2 r \partial_r \phi_{j+1}(r).
\]
Next, \begin{equation}\label{drpsijests}
 |(r \partial_r)^\alpha \psi_j| \lesssim_\alpha \frac{C^j}{(j-1)!} \frac{r^2}{\la r \ra^2}
\end{equation} follows from \eqref{drphijest}
In addition, $\psi_0$ solves $L^*\psi_0 = \phi_0$  therefore
a direct computation shows that

\[
\psi_{0}(r) =L(r^{2}\phi_{1}(r^{2}))= \frac{r^{2}-(1+r^{4})\arctan(r^{2})}{2r^{3}}.
\]
 The estimates \eqref{pointtp} follow from a direct estimation using \eqref{qest} and \eqref{drpsijests}.\\
 On the other hand in the region $r \xi \gtrsim  1$ we define 
\[
 \psi_{\xi}^+ = \xi^{-1} L\phi_{\xi}^+ 
\]
and therefore get \eqref{psirep}. The theorem now follows from the properties of $\phi^{+}_{\xi}$ from Theorem \ref{htransthm}.

\subsection{The transference Identity}
The operators $H$ and $\tilde H$ have variable coefficients 
so they do not admit a scaling symmetry. However, the behavior of 
the generalized Fourier transform under scaling is important,
as the scale parameter $\lambda$ varies significantly along our 
Schr\"odinger map flow. The commutator of the Fourier transform
with scaling is captured by the so-called  "transference operator" $\mathcal K$, following \cite{kstym}. The aim of this section 
will be to define and study the transference operator,
which provides a convenient expression for $\mathcal{F} (r \partial_{r})$ (see \eqref{htransdef}) called the "transference identity". \\
\\
We use the same notation as in \cite{kstym}, regarding the elements in the range of $\mathcal{F}$ as a two-component vector:

\begin{d1}
The transference operator $\K$ is defined (a priori, for  $u \in C^{\infty}_{c}((0,\infty))$) by
\begin{equation}\label{htransdef}\mathcal{F}(r \partial_{r}u)(\xi) = \begin{bmatrix} -1 & 0 \\ 0 &
- \xi \partial_{\xi} - \frac32 \end{bmatrix} \mathcal{F}(u) + \mathcal{K}(\mathcal{F}(u))(\xi). \end{equation}
\end{d1}
We remark that for the standard Fourier transform applied to radial functions in $\R^2$ one has $\K = -\frac32 I$. Precisely,
the $3/2$ factor arises due to the different measures used in the physical space ($rdr$) and in the Fourier space ($d\xi$).

Since the generalized Fourier transform is an isometry, it is easily verified that $\K$ is a skew-adjoint operator. 
Hence, in light of \eqref{htransvecnotation}, we can regard $\mathcal{K}$ as a matrix: 
\[
\mathcal{K} = \left(
\begin{array}{cc}
 0 & \mathcal{K}_{12} \\
 \mathcal{K}_{21} & \mathcal{K}_{22} \\
\end{array}
\right)
\]
where the entries satisfy 
\[
\K_{12}^* = -\K_{21}. \qquad \K_{22}^* = -\K_{22}
\]
and can be interpreted as integral operators with kernels as follows:

\begin{p1}\label{p:K}
a) The entries $\K_{ij}$ of the transference operator are integral operators of the form
\begin{equation}
 \begin{aligned}
& \K_{12} f = -\int_{0}^{\infty} K(\xi) f(\xi) \, d \xi,   
 \\
 & \K_{21} (\xi) = K(\xi),
 \\
& \K_{22} f(\xi) = \int_{0}^{\infty} K_{22}(\xi,\eta)  f(\eta) \, d\eta,
 \end{aligned}   
\end{equation}
where the kernels $K$ and $K_{22}$ are given by
\begin{equation}\label{K}
K(\xi) =  \|\phi_{0}\|_{L^{2}(r dr)}^{-1}\int_{0}^{\infty} r^{2}\phi_{0}'(r) \phi_{\xi}(r) dr,
\end{equation}
\begin{equation}\label{FH}
K_{22}(\xi,\eta) =  \text{p.v.} \frac{F(\eta,\xi)}{\xi^{2}-\eta^{2}}, \qquad 
F(\xi,\eta) = \int_{0}^{\infty}  \phi_{\eta}(r)  \left(\frac{-128 r^{2}(r^{4}-1)}{(r^{4}+1)^{3}}\right)\phi_{\xi}(r)r dr.
\end{equation}

b) These kernels satisfy the bounds 
\begin{equation}\label{kest}
  |K(\xi)| \lesssim \begin{cases}  \sqrt{\xi} , \quad \xi\leq 1\\
 	\xi^{-2}, \quad \xi \geq 1\end{cases},
 \end{equation}
respectively, with $\ixi = \min\{\xi,1\}$, 
\begin{equation}\label{fhests}|F(\xi,\eta)| \lesssim  
(\ixi^2+ \ieta^2)\ixi^\frac12 \ieta^\frac12
\end{equation}
and
\begin{equation}\label{dfhests}
\ixi |\partial_{\xi}F(\xi,\eta)| \lesssim  (\ixi^2+ \ieta^2)\ixi^\frac12 \ieta^\frac12
\end{equation}

c) In particular. we have the $L^2$ bound
\begin{equation}
\|\K\|_{L^{2,0} \to L^{2,0}} \lesssim 1.
\end{equation}
\end{p1}

We remark that the bounds \eqref{fhests} and \eqref{dfhests}
above are not sharp as written. One could gain off-diagonal decay, 
i.e. $\la \xi-\eta\ra^{-1}$ factors. We do not pursue this here 
because it is not needed, but we will obtain such an improvement 
later for the operator $\tcK$ associated to $\tH$.

\begin{proof}
We start by describing each $\mathcal{K}_{jk}$. For $\mathcal{K}_{12}$ we have
\begin{equation*}\label{k12def}
\begin{split}
\mathcal{K}_{12}f &= \int_{0}^{\infty}  r \frac{\phi_{0}(r)}{\|\phi_{0}\|_{L^{2}(r dr)}} r \partial_{r}\left(\int_{0}^{\infty}  \phi_{\xi}(r) f(\xi) d\xi\right) dr\\
&= -2 \int_{0}^{\infty}  \frac{\phi_{0}(r)}{\|\phi_{0}\|_{L^{2}(r dr)}} \int_{0}^{\infty}  \phi_{\xi}(r) f(\xi) d\xi \, rdr - \int_{0}^{\infty} \int_{0}^{\infty}  \phi_{\xi}(r) f(\xi) \frac{\phi_{0}'(r)}{\|\phi_{0}\|_{L^{2}(r dr)}} r^2 d \xi dr \\
&=
- \int_{0}^{\infty} \int_{0}^{\infty}  \phi_{\xi}(r) f(\xi) \frac{\phi_{0}'(r)}{\|\phi_{0}\|_{L^{2}(r dr)}} r^2 d \xi dr.
\end{split}
\end{equation*}
We obtain
\[
\mathcal{K}_{12}f  = -\int_{0}^{\infty} K(\xi)  f(\xi)d\xi,
\]
 where $K$ is as in \eqref{K}.
 
 Next we prove \eqref{kest}, which in particular implies that $K \in L^2$. We separate two cases.  When $\xi \geq 1$, we write $H(\phi_{\xi}) = \xi^{2} \phi_{\xi}$ and integrate by parts. On the other hand, when $\xi \leq 1$, we directly estimate $K$, as follows.
 \begin{equation}\nonumber
 \begin{split}
 |K(\xi)| &\lesssim \int_{0}^{\frac{1}{\xi}} r dr |\phi_{\xi}(r)-q(\xi)\phi_{0}(r)| |\phi_{0}'(r)| r + \int_{0}^{\frac{1}{\xi}} r dr q(\xi) \phi_{0}(r) |\phi_{0}'(r)| r\\
 &+ \int_{\frac{1}{\xi}}^{\infty} r dr |\phi_{\xi}(r)| |\phi_{0}'(r)| r\\
 &\lesssim  \sqrt{\xi},
 \end{split}
 \end{equation}
 where we have used \eqref{phiximinusqphi0}, \eqref{philgr}, \eqref{qest}.  Then, Cauchy-Schwartz implies the $L^2$ bound
 \[
 |\mathcal{K}_{12}(f)| \lesssim  \|f\|_{L^{2}}.
 \]
 Having $\K_{12}$, we also directly obtain its adjoint 
 \[
 \mathcal{K}_{21}(\xi) =  K(\xi).
 \]

Finally we consider $\mathcal{K}_{22}f$, 
for which we compute integrating by parts
\begin{equation*}
\begin{split}
&\mathcal{K}_{22}f(\xi)\\
&= \int_{0}^{\infty} r dr \phi_{\xi}(r) \int_{0}^{\infty} r \phi_{\eta}'(r) f(\eta) d\eta + \int_{0}^{\infty} r dr \phi_{\xi}(r) \int_{0}^{\infty} \eta f'(\eta) \phi_{\eta}(r) d\eta + \frac{3}{2}f(\xi)\\
&= \int_{0}^{\infty} r dr \phi_{\xi}(r) \int_{0}^{\infty} r \phi_{\eta}'(r) f(\eta) d\eta - \int_{0}^{\infty} r dr \phi_{\xi}(r) \int_{0}^{\infty} d\eta f(\eta) (\phi_{\eta}(r) + \eta \partial_{\eta}\phi_{\eta}(r))+\frac{3}{2}f(\xi)\\
&= \int_{0}^{\infty} r dr \phi_{\xi}(r) \int_{0}^{\infty} d\eta (r \phi_{\eta}'(r) - \eta \partial_{\eta}\phi_{\eta}(r)) f(\eta)+\frac{1}{2}f(\xi)\\
&:= \mathcal{K}_{220}f(\xi) +\frac{1}{2}f(\xi).
\end{split}
\end{equation*}
In order to capture the off-diagonal behavior of $\K_{220}$ we commute the $H$ operator inside this representation, using the relation
\[
(H-\eta^2) \partial_\eta \phi_\eta = 2 \eta \phi_\eta.
\]
This yields
\begin{equation*}
\begin{split}
\xi^{2} \mathcal{K}_{220}f(\xi) &= \int_{0}^{\infty} r dr H(\phi_{\xi})(r) \int_{0}^{\infty} d\eta (r \phi_{\eta}'(r) - \eta \partial_{\eta}\phi_{\eta}(r)) f(\eta)\\
&=\int_{0}^{\infty} dr \phi_{\xi}(r) r \int_{0}^{\infty} d\eta H(r\phi_{\eta}'(r) - \eta \partial_{\eta}\phi_{\eta}(r)) f(\eta)\\
&=\int_{0}^{\infty} dr \phi_{\xi}(r) r \int_{0}^{\infty} d\eta \eta^{2}(r \partial_{r}-\eta \partial_{\eta})(\phi_{\eta}(r)) f(\eta) \\
&+ \int_{0}^{\infty} dr \phi_{\xi}(r) r \int_{0}^{\infty} d\eta \left(\frac{-128 r^{2}(r^{4}-1)\phi_{\eta}(r)}{(r^{4}+1)^{3}}f(\eta)\right).
\end{split}
\end{equation*}
The kernel of the last term is a smooth function 
of $\xi$ and $\eta$, which we denote by $F(\xi,\eta)$ as in \eqref{FH}. Hence we have proved that
\begin{equation}
\xi^{2}\mathcal{K}_{22}f(\xi) - \mathcal{K}_{22}((\cdot)^{2} f)(\xi) = \int_{0}^{\infty} F(\xi,\eta) f(\eta) d\eta.
\end{equation}
This shows that away from the diagonal, the kernel of $\K_{22}$ is given by 
\begin{equation}\label{K22}
K_{22}(\xi,\eta) = \frac{F(\xi,\eta)}{\xi^2-\eta^2}.
\end{equation}
It remains to determine the diagonal behavior of $\mathcal{K}_{22}$. In the region $1 \lesssim r \xi$ we have
\[
\phi_{\xi}(r) = 2 \text{Re}\left(\frac{a(\xi)}{\sqrt{r}}e^{i r \xi} (1+\frac{15 i}{8 r \xi})\right) + O\left(\frac{1}{r^{11/2} \xi}\right)+O\left(\frac{1}{r^{5/2}\xi^{2}}\right),
\]
which yields the asymptotic
\[
\left(r \partial_{r}-\xi \partial_{\xi}\right) \phi_{\xi}(r) = -\frac{1}{2} \phi_{\xi}(r) +  O\left(\frac{1}{r^{11/2} \xi}\right)+O\left(\frac{1}{r^{5/2}\xi^{2}}\right) -2 \text{Re}\left(\frac{\xi a'(\xi) e^{i r \xi}}{\sqrt{r}} (1+\frac{15 i}{8 r\xi})\right),
\]
where the decay of the tail limits the 
singular behavior on the diagonal. 
This is similar to the analysis in \cite{kstym},
so we do not repeat it here. Using also 
the antisymmetry of $\K_{22}$, we obtain that
its kernel is given by the principal value expression in \eqref{K22}, which concludes the proof of \eqref{FH}.

\bigskip

To conclude the proof of the proposition we need to study the kernel
\[
F(\xi,\eta) = \int_{0}^{\infty}  \left(\frac{-128 r^{2}(r^{4}-1)}{(r^{4}+1)^{3}}\right)\phi_{\xi}(r)\phi_{\eta}(r) r dr.
\]
For this, we recall \eqref{phiximinusqphi0}, \eqref{philgr}, and  note that
\[
|\partial_{\xi}\phi_{\xi}(r)-q'(\xi)\phi_{0}(r)| \lesssim \frac{ q(\xi) r^{4} \xi}{\langle r^{2}\rangle}, \quad r \xi \lesssim 1
\]
\[
|\partial_{\xi}\phi_{\xi}(r)| \lesssim \sqrt{r}, \quad 1 \lesssim r \xi.
\]
In addition, we observe the cancellation
\[
\int_{0}^{\infty}  \phi_{0}(r)  \left(\frac{-128 r^{2}(r^{4}-1)}{(r^{4}+1)^{3}}\right)\phi_{0}(r)\, r dr =0.
\]
This can be seen by direct computation, or one could note that
\[
\int_{0}^{\infty}  \phi_{0}(r)  \left(\frac{-128 r^{2}(r^{4}-1)}{(r^{4}+1)^{3}}\right)\phi_{0}(r) \, rdr= \int_{0}^{\infty}  \phi_{0}(r)  \left([H,r \partial_{r}]\phi_{0} - 2 H(\phi_{0})\right)\, rdr
\]
and use $H(\phi_{0})=0$ and the fact that $H$ is self-adjoint on $L^{2}(r dr)$. Combining these properties, it is a direct computation 
to prove the bounds \eqref{fhests} and \eqref{dfhests}.

It remains to prove the $L^2$ bound for the operator 
$\mathcal{K}_{22}$, which is expressed in terms of $F$ as
in \eqref{FH}. For this we decompose the kernel into a leading near diagonal part, and two milder terms
\begin{equation}\nonumber
\begin{split}
\text{p.v.} \frac{F_{H}(\xi ,\eta )}{\eta ^2-\xi ^2}= & \ 
\text{p.v.}\frac{\mathbbm{1}_{|\xi-\eta| \leq 1}}{\eta-\xi}
\cdot \frac{F_{H}(\eta ,\eta )}{2 \eta}+
\frac{\mathbbm{1}_{|\xi-\eta| \leq 1}}{\eta -\xi } \left(\frac{F_{H}(\xi ,\eta )}{\eta +\xi }-\frac{F_{H}(\eta ,\eta )}{2 \eta }\right)+\frac{\mathbbm{1}_{|\xi-\eta| > 1} F_{H}(\xi ,\eta )}{(\eta -\xi ) (\eta +\xi )}
\\ = & \ T \frac{F_{H}(\eta,\eta)}{2\eta} + 
K_{22}^1 + K_{22}^2
\end{split}
\end{equation}
where we have factored the first term on the right into a multiplication with a bounded function and a frequency localized Hilbert transform,
\begin{equation}
\label{topdef}T(g)(\eta) = \text{p.v.}\int_{\mathbb{R}} \frac{\mathbbm{1}_{|\xi-\eta| \leq 1}}{\eta-\xi} g(\xi)d\xi. 
\end{equation}
This is a Fourier multiplier with bounded symbol, so is bounded on $L^{2}(\mathbb{R})$. The remaining terms are nonsingular, and we can simply estimate their kernels in $L^2$. This yields the $L^2$ bound
\begin{equation}\label{pvest}
\begin{split}
\| \K_{22} f\|_{L^{2}}
&\lesssim  
\|f\|_{L^{2}(d\xi)} \left(\left\|\frac{F_{H}(\eta,\eta)}{2\eta}\right\|_{L^{\infty}_{\eta}}+
\|K_{22}^1\|_{L^2} + \|K_{22}^2\|_{L^2} \right).
\end{split}
\end{equation}
It remains to bound the norms in the last sum.

\medskip

Using \eqref{fhests}, we directly see that $\dfrac{F_{H}(\eta,\eta)}{2\eta} \in L^{\infty}_{\eta}$.

\medskip

Next, by the fundamental theorem of calculus, we get
\begin{equation}
\begin{split}
&\left|\frac{F_{H}(\xi,\eta)}{(\eta+\xi)} - \frac{F_{H}(\eta,\eta)}{2\eta}\right|\lesssim \int_{\min\{\eta,\xi\}}^{\max\{\eta,\xi\}} \left(\frac{|\partial_{1}F_{H}(\omega,\eta)|}{(\eta+\omega)} + \frac{|F_{H}(\omega,\eta)|}{(\eta+\omega)^{2}}\right)d\omega
\end{split}
\end{equation}  
Using \eqref{fhests} and \eqref{dfhests} we can bound this by
\begin{equation}
\begin{split}
&\left|\frac{F_{H}(\xi,\eta)}{(\eta+\xi)} - \frac{F_{H}(\eta,\eta)}{2\eta}\right|\lesssim 
|\xi-\eta| \frac{\ixi^2+\ieta^2}{\xi+\eta},
\end{split}
\end{equation}
which shows that
\[
K_{22}^1(\xi,\eta) = \frac{1}{\eta-\xi} \mathbbm{1}_{|\xi-\eta| \leq 1} \left(\frac{F_{H}(\xi,\eta)}{(\eta+\xi)} - \frac{F_{H}(\eta,\eta)}{2\eta}\right) \in L^{2}(d\xi)L^{2}(d\eta).
\]
Finally, \eqref{fhests} gives that $|F_{H}(\xi,\eta)| \lesssim 1$ in the region $\xi+\eta \gtrsim  1$. This suffices in order to 
show that
\[
K_{22}^2(\xi,\eta) = 
\frac{\mathbbm{1}_{|\xi-\eta| \geq 1}}{(\eta-\xi)} \frac{F_{H}(\xi,\eta)}{(\eta+\xi)} \in L^{2}(d\xi)L^{2}(d\eta).
\]
In light of \eqref{pvest}, this completes the proof of the $L^{2}$ boundedness of $\K_{22}$,
and thus the proof of the proposition.

\end{proof}

In this article we will also use 
even more the transference operator associated the the  $\widetilde{H}$ operator, which we will denote by $\tcK$. Compared with $\mathcal{K}$ this is simpler as we no longer 
need to separate the zero mode:
\begin{d1}\label{d:tK}
The transference operator $\tcK$ associated to $\tH$ is defined (a priori, for  $u \in C^{\infty}_{c}((0,\infty))$) by
\begin{equation}\label{thtransdef}
\FtH(r \partial_{r}u)(\xi) = 
 (-\xi \partial_\xi - \frac32) \FtH(u) + \tcK(\FtH(u))(\xi), 
 \end{equation}
\end{d1}
Here the $\frac32$ correction is again chosen
in order to ensure that $\tcK$ is antisymmetric.

One may think of $\tcK$ as a scaling derivative 
of the Fourier transform $\FtH_\lambda$, in the sense that
\begin{equation}\label{dlambda-F1}
{\frac{d}{d \lambda} {\FtH_\lambda}}_{|\lambda = 1} =
- \tcK \FtH.
\end{equation}
With this interpretation the antisymmetry requirement is clear, since the Fourier transforms are isometries. 
This is immediately seen by differentiating
with respect to $\lambda$ in \eqref{Flambda-def},
and applying the chain rule. For the inverse Fourier 
transform this directly gives
\begin{equation}\label{dlambda-F1-}
{\frac{d}{d \lambda} {\FtH_\lambda^{-1}}}_{|\lambda = 1} =
\FtH^{-1} \tcK. 
\end{equation}

We summarize the properties of this operator in the next proposition, which is the direct counterpart of Proposition~\ref{p:K}:

\begin{p1}\label{p:cK}
a) The transference operator $\tcK$ is an integral operator of the form
\begin{equation}
 \begin{aligned}
\tcK f(\xi) = \int_{0}^{\infty} \tK(\xi,\eta)  f(\eta) \, d\eta,
 \end{aligned}   
\end{equation}
where the kernel $\tK$ is given by
\begin{equation}\label{tFH}
\tK(\xi,\eta) =  \text{p.v.}\, \frac{\tF(\eta,\xi)}{\xi^{2}-\eta^{2}}, \qquad 
\tF(\xi,\eta) = -
\int_{0}^{\infty}   \frac{\psi_{\eta}(r)\psi_{\xi}(r) 32 r^{2}}{(1+r^{4})^{2}} rdr.
\end{equation}

b) This kernel satisfies the bounds 
\begin{equation}\label{fhests-th}|\tF(\xi,\eta)| \lesssim  
\ixi^\frac32 \ieta^\frac32 
\frac{\min\{\la\xi\ra,\la\eta\ra\}}{\la \xi+\eta\ra \la \xi-\eta \ra},
\end{equation}
and
\begin{equation}\label{dfhests-th}
\ixi |\partial_{\xi}\tF(\xi,\eta)| \lesssim  \ixi^\frac32 \ieta^\frac32.
\end{equation}

c) In particular. we have the $L^2$ bound
\begin{equation}
\|\tcK\|_{L^{2} \to L^{2}} \lesssim 1.
\end{equation}
\end{p1}

\begin{proof}
a) From the definition of $\tcK$ we have
\[
\tcK f(\xi) = \FtH(r \partial_{r} \FtH^{-1}(f))(\xi) + \xi f'(\xi) + \frac32 f(\xi), \quad f \in C^{\infty}_{c}((0,\infty)).
\]
Proceeding as for $\mathcal{K}$, we rewrite this in the form
\[
\tcK f(\xi) = \int_{0}^{\infty} r dr \psi_{\xi}(r) \int_{0}^{\infty} f(\eta) \left(r \partial_{r}\psi_{\eta}(r) - \eta \partial_{\eta} \psi_{\eta}(r)\right) \, d\eta +\frac12 f(\xi).
\]
Commuting the operator $\tH$ inside this representation
we arrive at 
\[
\xi^{2}\tcK(f)(\xi) - \tcK((\cdot)^{2}f(\cdot))(\eta) = \int \tF(\xi,\eta) f(\eta) d\eta,
\]
where the smooth symmetric function $\tF$
is given by 
\begin{equation}\label{tF}
\tF(\xi,\eta) =-  \int_{0}^{\infty} \frac{\psi_{\eta}(r)\psi_{\xi}(r) 32 r^{2}}{(1+r^{4})^{2}} rdr.
\end{equation}
This yields the off-diagonal behavior for the kernel $\tK$ as in \eqref{tFH}.
For its diagonal behavior, we start with
\[
\psi_{\xi}(r) = 2 \text{Re}\left(\frac{a(\xi)}{\sqrt{r}} e^{i r \xi} \widetilde{\sigma}(r\xi,r)\right), \quad r \xi \geq C
\]
and
\[
(r \partial_{r}-\xi \partial_{\xi})(\psi_{\xi}(r)) = -\frac{1}{2}\psi_{\xi}(r) + O\left(\frac{1}{r^{3/2} \xi}\right) - 2 \text{Re}\left(\frac{\xi}{\sqrt{r}} a'(\xi) e^{i r \xi} \widetilde{\sigma}(r\xi,r)\right)
\]
where, for fixed $\xi$,
\[
\tilde{\sigma}(r\xi,r) = i - \frac{3}{8} r^{-1}\xi^{-1} + O(r^{-2}).
\]
Arguing as in \cite{kstym} and using the antisymmetry, we obtain the p.v. diagonal kernel behavior as in \eqref{tFH}.
\medskip

b)	The estimates on $\psi_{\xi}$ from Theorem \ref{httransthm} directly give
\begin{equation}\label{Fkernelests}
|\tF(\xi,\eta)| \lesssim   \ixi^{3/2}\ieta^{3/2}
\end{equation}
and
\begin{equation}\label{dxiFkernelests}
\ixi|\partial_{\xi}\tF(\xi,\eta)| \lesssim  \ixi^{3/2}\ieta^{3/2}
\end{equation}
Here \eqref{dxiFkernelests} is what we need, but \eqref{Fkernelests} only suffices in the region $|\xi-\eta| \lesssim 1$.
To prove the better bound \eqref{fhests-th} we need
an off-diagonal refinement of \eqref{Fkernelests} when  $\xi,\eta \gtrsim 1$. For this, we recall that
\[
\tF(\xi,\eta) = -32 \int_{0}^{\infty} \frac{r^{2}}{(1+r^{4})^{2}} \psi_{\xi}(r) \psi_{\eta}(r) \, rdr.
\]
Then,
\begin{equation*}
	\begin{split}
			\xi^{2}\tF(\xi,\eta) &= -32 \int_{0}^{\infty} \frac{r^{3} dr}{(1+r^{4})^{2}} \widetilde{H}(\psi_{\xi})(r) \psi_{\eta}(r) \\
			&= -32 \langle \frac{r^{2}}{(1+r^{4})^{2}} \eta^{2} \psi_{\eta}+[\widetilde{H},\frac{r^{2}}{(1+r^{4})^{2}}]\psi_{\eta},\psi_{\xi}\rangle_{L^{2}(r dr)}\\
			&= \tF(\xi,\eta) \eta^{2} + 32 \int_{0}^{\infty}\psi_{\xi}(r) [\widetilde{H},\frac{r^{2}}{(1+r^{4})^{2}}]\psi_{\eta} \, rdr.
\end{split}
\end{equation*}
Computing the commutator, for $\xi \neq \eta$ we obtain 
\[
\tF(\xi,\eta) = \frac{32}{(\xi^{2}-\eta^{2})} \int_{0}^{\infty} \frac{\psi_{\xi}(r) \cdot 4 r}{(1+r^{4})^{4}} \left(\left(1-14 r^4+9 r^8\right) \psi_{\eta}(r)+r \left(1-2r^{4}-3r^{8}\right) \psi_{\eta}'(r)\right)dr.
\]
Hence, using the bounds on $\psi_\xi$ from from Theorem \ref{httransthm}, and the symmetry of $\tF$, we obtain  
\begin{equation*}\label{Fkernellg} 
|\tF(\xi,\eta)| \lesssim \frac{ \text{min}\{\xi,\eta\}}{|\xi^{2}-\eta^{2}|}, \quad \xi,\eta \gtrsim 1
\end{equation*}
\[|\tF(\xi,\eta)| \leq \frac{C \min\{\xi,\eta\}^{3/2}}{|\xi^{2}-\eta^{2}|}, \quad 0 < \eta \lesssim 1 \lesssim \xi \text{ or } 0 < \xi \lesssim 1 \lesssim \eta,\]
which completes the proof of \eqref{fhests-th}.

c) Now, we study the operator
\[
f \mapsto \text{p.v.}\left(\int_{0}^{\infty} \frac{f(\xi) \tF(\xi,\eta) d\xi}{\xi^{2}-\eta^{2}}\right).
\]
As for $\mathcal{K}$, 
we start with the decomposition
\begin{equation}\label{fkernaldecomp}\frac{\tF(\xi ,\eta )}{\xi ^2-\eta ^2}=\frac{\mathbbm{1}_{|\xi-\eta| > 1} \tF(\xi ,\eta )}{(\xi -\eta ) (\eta +\xi )}+\frac{\mathbbm{1}_{|\xi-\eta| \leq 1} \left(\frac{\tF(\xi ,\eta )}{\eta +\xi }-\frac{\tF(\eta ,\eta )}{2 \eta }\right)}{\xi -\eta }+\frac{\mathbbm{1}_{|\xi-\eta| \leq 1} \tF(\eta ,\eta )}{2 \eta  (\xi -\eta )}\end{equation}
and estimate the contributions of each of the three kernels separately. We use \eqref{Fkernelests} and \eqref{fhests-th} to get
\[
\mathbbm{1}_{|\xi-\eta| >1} |\tF(\xi,\eta)| \leq C \begin{cases} \xi + \eta, \quad \eta,\xi \lesssim 1 \\
1, \quad \text{otherwise.}\end{cases}
\]
This implies that
\[
\frac{\tF(\xi,\eta) \mathbbm{1}_{|\xi-\eta|>1}}{(\xi-\eta)(\eta+\xi)} \in L^{2}(d\xi) L^{2}(d\eta),
\]
which suffices for the first term in \eqref{fkernaldecomp}.
Next, we note that
\[
\left|\frac{\tF(\xi,\eta)}{\eta+\xi} - \frac{\tF(\eta,\eta)}{2 \eta}\right| \leq C \int_{\min\{\eta,\xi\}}^{\max\{\eta,\xi\}} \left(\frac{|\partial_{1}\tF(\omega,\eta)|}{\eta+\omega} + \frac{|\tF(\omega,\eta)|}{(\eta+\omega)^{2}}\right) d\omega.
\]
From \eqref{Fkernelests},
\[
\frac{|\tF(\omega,\eta)|}{(\eta+\omega)^{2}} \leq C \left(\mathbbm{1}_{\{4 \geq \eta\}} + \frac{1}{\eta} \mathbbm{1}_{\{\eta \ge 4\}}\right).
\]
Also, from \eqref{dxiFkernelests}, we get
\[
\frac{|\partial_{1}\tF(\omega,\eta)|}{\eta+\omega} \leq C \left(\mathbbm{1}_{\{4 \geq \eta\}} + \frac{1}{\eta} \mathbbm{1}_{\{\eta \ge 4\}}\right).
\]
Therefore,
\[
\left|\frac{\tF(\xi,\eta)}{\eta+\xi} - \frac{\tF(\eta,\eta)}{2 \eta}\right| \leq C |\eta-\xi|  \left(\mathbbm{1}_{\{4 \geq \eta\}} + \frac{1}{\eta} \mathbbm{1}_{\{\eta \ge 4\}}\right)
\]
which gives
\[
\frac{\mathbbm{1}_{|\xi-\eta| \leq 1}}{|\xi-\eta|}\left(\frac{\tF(\xi,\eta)}{(\eta+\xi)}-\frac{\tF(\eta,\eta)}{2 \eta}\right) \in L^{2}(d\xi) L^{2}(d\eta),
\]
as needed for the second term in \eqref{fkernaldecomp}.

Next, as in \eqref{pvest}, we have
\begin{equation*}
\begin{split}
&\left\|   \frac{\tF(\eta,\eta)}{2\eta} \text{p.v.} \int_{0}^{\infty} \frac{d\xi \mathbbm{1}_{|\xi-\eta| \leq 1} \tF(\eta,\eta)}{\xi-\eta} d\xi\right\|_{L^{2}(d\eta)}\leq C \|T\|_{\mathcal{L}(L^{2}(\mathbb{R}), L^{2}(\mathbb{R}))} \|f\|_{L^{2}(d\xi)} \|\frac{\tF(\eta,\eta)}{2 \eta}\|_{L^{\infty}_{\eta}}
\end{split}
\end{equation*}
where we recall the definition of $T$ in \eqref{topdef}. Finally, \eqref{Fkernelests} shows that
\[
\|\frac{\tF(\eta,\eta)}{\eta}\|_{L^{\infty}_{\eta}} \leq C
\]
which suffices for the third term in \eqref{fkernaldecomp}, and finally shows that the operator $\tcK$
is bounded on $L^{2}((0,\infty))$.
\end{proof}

The transference operators $\cK$ and $\tcK$ were defined above in the context of the operators
$H$ and $\tH$. However, we will also need them in the 
rescaled setting, associated to $H_\lambda$ and $\tH_\lambda$.  Their rescaled versions are denoted
by $\cK_\lambda$ and $\tcK_\lambda$, and are still defined by \eqref{htransdef}, respectively \eqref{thtransdef}, but with $\FH$ and $\FtH$ replaced
by $\FH_\lambda$ , respectively $\FtH_\lambda$. 
These are obtained from $\cK$ and $\tcK$ by rescaling, for instance
\begin{equation}\label{klambdadef}
\tcK_{\lambda}f(\xi) = \frac{1}{\lambda}\tcK(f(\cdot\lambda))(\frac{\xi}{\lambda}).
\end{equation}
The representation \eqref{tFH} remains valid, but with $\tF$ replaced by its rescaled version $\tF_\lambda$ 
given by 
\begin{equation}\label{tflambdadef}
\tF_\lambda(\xi,\eta)  = \tF(\frac{\xi}{\lambda},\frac{\eta}{\lambda}).
\end{equation}
Finally, the relations \eqref{dlambda-F1} and \eqref{dlambda-F1-} also carry through, and give
\begin{equation}\label{dlambda-F}
{\frac{d}{d \lambda} {\FtH_\lambda}} =
-\tcK_\lambda \FtH_\lambda,\qquad
\frac{d}{d \lambda} {\FtH_\lambda^{-1}} =
\FtH^{-1}_\lambda \tcK_\lambda. 
\end{equation}

\subsection{The operator \texorpdfstring{$r \FtH^{-1}$}{}}
Our goal here is to study the operator of multiplication by $r$ on the Fourier side. This will allow us later on to obtain more refined, spatially localized bounds for the solutions to the 
linear $\tH_\lambda$ flow.

\begin{l1}\label{multrfourier} For all $f \in C^{\infty}_{c}((0,\infty))$ we have 
$$\|r \FtH^{-1}f\|_{L^{2}(r dr)} \lesssim  \|\frac{f(\xi)}{\xi}\|_{L^{2}(d\xi)} +  \|f'(\xi)\|_{L^{2}(d\xi)}.$$
\end{l1}
In particular, the above inequality is true for all $f$ in the closure of $C^{\infty}_{c}((0,\infty))$ under the norm
\[
\|f\|^2:=\|f(\xi)\|_{L^{2}(d\xi)}^{2} + \|\frac{f(\xi)}{\xi}\|_{L^{2}(d\xi)}^{2} + \|f'(\xi)\|_{L^{2}(d\xi)}^{2}
\]

\begin{proof}
We start with an analog of
\begin{equation}\label{j2j1rel}r J_1(r \xi)=\partial_{\xi} J_2( r\xi)+\frac{2 J_2(r \xi)}{\xi} \end{equation}
in our setting. (Recall that $J_{1}(r\xi)$ are eigenfunctions of the large $r$ part of $\widetilde{H}$, and $J_{2}(r \xi)$ are eigenfunctions of the large $r$ part of $H$). More precisely,
our starting point is the identity
\[
\psi_{\xi}(r)= \xi^{-1}L\phi_{\xi},
\]
where we recall that 
\[
Lf= \partial_{r}f+\frac{2}{r}h_{3}(r) f, \quad h_{3}(r) = \frac{r^4-1}{r^4+1}.
\]
This allows us to obtain the identity
\begin{equation}
\begin{split} 
r \psi_{\xi}(r) &=\frac{1}{\xi}\left(r \partial_{r}\phi_{\xi}+2 h_{3}\phi_{\xi}\right)\\
&=\frac{1}{\xi}(r \partial_{r}-\xi\partial_{\xi})\phi_{\xi} + \partial_{\xi}\phi_{\xi} + \frac{2 \phi_{\xi}}{\xi} + \frac{2 (h_{3}(r)-1)}{\xi}\phi_{\xi}.
\end{split}
\end{equation}
This is the analog of \eqref{j2j1rel}, with the only error terms either decaying quickly at infinity (e.g. $|h_{3}(r)-1| \leq \frac{C}{r^{4}}$) or involving an operator which annihilates functions of the form $f(r\xi)$.\\
\\
Therefore, for $f \in C^{\infty}_{c}((0,\infty))$ we have 
\begin{equation}\label{rinversehtilde}
\begin{split}r \FtH^{-1}(f)(r) &= r \int_{0}^{\infty} \psi_{\xi}(r) f(\xi) d\xi
\\
&=\int_{0}^{\infty} (r\partial_{r}-\xi\partial_{\xi} )(\phi_{\xi}(r)) \frac{f(\xi)}{\xi} d\xi + \int_{0}^{\infty} \partial_{\xi}\phi_{\xi}(r) f(\xi) d\xi
\\
&+\int_{0}^{\infty} \frac{2\phi_{\xi}}{\xi} f(\xi) d\xi + \int_{0}^{\infty} \frac{2(h_{3}(r)-1)}{\xi} \phi_{\xi}(r) f(\xi) d\xi.\end{split}\end{equation}

We use the transference operator $\K$ to rewrite the first term on the right-hand side of \eqref{rinversehtilde} in the form
\begin{equation}\begin{split} 
\int_{0}^{\infty} (r\partial_{r}-\xi\partial_{\xi})\phi_{\xi}(r) \frac{f(\xi)}{\xi}d\xi &= r \partial_{r}\FH^{-1}\begin{bmatrix} 0\\
f/\xi \end{bmatrix} + \FH^{-1}\begin{bmatrix} 0\\
f'(\cdot)\end{bmatrix} 
\\
&=\FH^{-1}\begin{bmatrix} 0\\
- \xi \partial_{\xi} \left(f/\xi\right)\end{bmatrix} + \FH^{-1}(\mathcal{K}-\frac32) \begin{bmatrix} 0\\
f/\xi\end{bmatrix} + \FH^{-1}\begin{bmatrix} 0\\
f'(\cdot)\end{bmatrix} 
\\
&= \FH^{-1}(\mathcal{K}-\frac12)\begin{bmatrix} 0\\
f/\xi\end{bmatrix}.
\end{split}
\end{equation}
Overall, \eqref{rinversehtilde} becomes
\begin{equation}
\begin{split} 
r \FtH^{-1}f &= 
  \FH^{-1}\left(\mathcal{K}+\frac32\right)\begin{bmatrix} 0\\
f/\xi\end{bmatrix}- \FH^{-1}\begin{bmatrix} 0\\
f'(\cdot)\end{bmatrix}+ 2(h_{3}(r)-1) \FH^{-1}\begin{bmatrix} 0\\
f/\xi \end{bmatrix}
\end{split}
\end{equation}
which, by the $L^2$ boundedness of the Fourier transform and of the transference operator (see Proposition~\ref{p:K}) gives
\[
\|r \FtH^{-1}(f)\|_{L^{2}(r dr)} \lesssim  \|\frac{f(\xi)}{\xi}\|_{L^{2}(d\xi)} +  \|f'(\xi)\|_{L^{2}(d\xi)},
\]
as needed.
\end{proof}

\subsection{Littlewood-Paley projectors in the \texorpdfstring{$\tilde H_\lambda$}{} frame} 
\label{L-PP}
In this section we seek to understand the properties of the Littlewood-Paley projectors in the $\tilde H_\lambda$ frame. Here one could proceed as we did in the previous subsections, with $\lambda = 1$, and then rescale. But for reference purposes we preferred 
to state the results with the parameter $\lambda$ included.

We recall from Section \ref{defnot} the functions $m_j, j \in \Z$ with the properties
\[
m_{j}(x) = m(\frac{x}{2^{j}}), \quad m \in C^{\infty}_{c}((\frac{1}{4},2)), \quad m(x) =1, \quad \frac{1}{2} \leq x \leq 1.
\]
and 
\[
\sum_j m_j =1.
\]
Based on these functions we define our projectors
\begin{equation}\label{pjlambdadef}
P^\lambda_{j}  = \FtH_{\lambda}^{-1} m_{j} \FtH_{\lambda}.
\end{equation}

We note that the notation $P_j^\lambda$ does not carry the $\tilde{\cdot}$ symbol which we used for the operator $\tilde H_\lambda$. The reasons we do so is because throughout most of this paper we work with projectors only in the  $\tilde H_\lambda$ calculus, and not with the ones in the $H_\lambda$ calculus; the only exception is in Section \ref{s:XLX} where we will carefully differentiate between the notation used for the projectors in the two different calculi. This allows us to reserve the notation $\tilde P_j^\lambda$ for the purpose described below (which is standard in the literature). 

We also need the projectors $\tilde P_j^\lambda$ which also localize at frequency $2^j$ and enjoy the property $ \tilde P_j^\lambda P_j^\lambda =  P_j^\lambda$; they are constructed using the functions $\widetilde{m}_{j}$ by
\[
\tP^\lambda_{j}  = \FtH_{\lambda}^{-1} \tm_{j} \FtH_{\lambda}.
\]
We recall from Section \ref{defnot} that $\widetilde{m}_{j}(x) = \widetilde{m}(\frac{x}{2^{j}})$ and that $\widetilde{m} \in C^{\infty}_{c}((0,\infty)), \widetilde{m}(x)=1$ in the support of $m$.

If $u\in C^{\infty}_{c}((0,\infty))$, then we  can represent  $P_{j}^\lambda u$
using an integral kernel,
\[
P_{j}^\lambda(u)(r) =\int_{0}^{\infty}  K_{j}^\lambda(r,s) u(s) s ds,
\]
where the kernel of the projector is given by
\[
K_{j}^\lambda (r,s) = \int_{0}^{\infty}  \psi^\lambda_{\xi}(r) m_{j}(\xi) \psi^\lambda_{\xi}(s) d\xi = \lambda\int_{0}^{\infty}  \psi_{\lambda^{-1}\xi}(r\lambda) m_{j}(\xi) \psi_{\lambda^{-1} \xi}(s\lambda) d\xi.
\]
This kernel and its properties will play an important role in this paper.
\medskip

There is another kernel whose characterization we need, which arises when we seek to express
a function $\psi$ in terms of $L^* \psi$ in an elliptic fashion, for a frequency localized function $\psi$. This is done based on the identity
\[
P^\lambda_{j}\psi = \widetilde{P_{j}^\lambda}(P_{j}^\lambda\psi) = \widetilde{P_{j}^\lambda} \widetilde{H}_{\lambda}^{-1} L_{\lambda} L^{*}_{\lambda}(P_{j}^\lambda\psi).
\]
which shows that the operator we need to consider is $\widetilde{P_{j}^\lambda} \widetilde{H}_{\lambda}^{-1} L_{\lambda}$. We also represent this operator via an integral 
kernel, writing for $u \in C^{\infty}_{c}((0,\infty))$
\begin{equation}\begin{split}
\widetilde{P_{j}}^{\lambda}(\widetilde{H}_{\lambda}^{-1} L_{\lambda}(u))(r) =  \int_{0}^{\infty} \psi^\lambda_{\xi}(r) \widetilde{m}_{j}(\xi) \FtH_\lambda (\widetilde{H}_{\lambda}^{-1} L_{\lambda} u)(\xi) d\xi=\int_{0}^{\infty}   K_{j}^{1,\lambda}(s,r) u(s) s ds.\end{split}\end{equation}
Here we recall that $\xi \psi^\lambda_\xi = L_\lambda \phi_\xi$. This allows us to use 
the $H_\lambda$ associated Fourier representation of $u$ in order to write
\[
K_{j}^{1,\lambda}(s,r) =  \int_{0}^{\infty}  \psi^\lambda_{\xi}(r) \widetilde{m}_{j}(\xi) \frac{\phi^\lambda_{\xi}(s)}{\xi} d\xi=\lambda \int_{0}^{\infty}  \psi_{\frac{\xi}{\lambda}}(r\lambda) \widetilde{m}_{j}(\xi) \frac{\phi_{\frac{\xi}{\lambda}}(s\lambda)}{\xi} d\xi.
\] 
Another related operator that will be used in this paper is 
$L^{-1}_{\lambda} \widetilde{P_j^\lambda}$, which can be defined as
\[
 L^{-1}_{\lambda} \widetilde{P_j^\lambda} := L^*_{\lambda} \widetilde{H}_{\lambda}^{-1} \widetilde{P_j^{\lambda}}.
\]
We remark that its adjoint is given by
\[
( L^{-1}_{\lambda} \widetilde{P_j^{\lambda}})^* = \widetilde{P_j^\lambda} (L_{\lambda}^*)^{-1} := \widetilde{P_j^\lambda} \tilde H_{\lambda}^{-1} L_{\lambda},
\]
therefore the kernel of $L^{-1}_{\lambda} \widetilde{P_j^\lambda}$ is simply 
$K_j^{1,\lambda}(r,s)$.

Our aim in this section is to characterize the two kernels introduced above, $K_j^\lambda(r,s)$ (the kernel of the projector $P_j^\lambda$) and $K_j^{1,\lambda}(r,s)$. For this purpose we introduce two auxiliary weight functions
\begin{equation} \label{omdef}
\omega_{j,\lambda}(r) = \begin{cases} \text{min}\{1,r^{3}2^{3j}\}, \quad 2^{j} \geq \lambda \\
\text{min}\{1,2^{j} r \dfrac{r^{2}\lambda^{2}}{1+r^{2}\lambda^{2}}\}, \quad 2^{j} \leq \lambda.\end{cases}
\end{equation}
and
\begin{equation} \label{tomdef}
\tilde \omega_{j,\lambda}(r) = \begin{cases}\text{min}\{1,\dfrac{r^{2}\lambda^{2}}{\langle r^{2} \lambda^{2}\rangle^{2}} + \dfrac{r^{2}\lambda^{2} (2^{j}r)^{2}}{\langle r^{2}\lambda^{2} \rangle}\}, \quad 2^{j}\leq \lambda\\
\text{min}\{1, r^{2}2^{2j}\}, \quad 2^{j}\geq \lambda\end{cases}.
\end{equation}
With this notation in place, we state our main result in this section. 

\begin{p1} \label{p-lp} The kernels $K^\lambda_{j}$
and $K_{j}^{1,\lambda}$ satisfy the following bounds:

\begin{equation}\label{ker1}
|K^\lambda_{j}(r,s)| \leq \frac{C_{N}2^{2j} \omega_{j,\lambda}(r)\omega_{j,\lambda}(s)}{(1+2^{j}(s+r))(1+2^{j}|r-s|)^{N}},
\end{equation}

\begin{equation}\label{ker2}
|\partial_{r}K^\lambda_{j}(r,s)| \leq \frac{C_{N} 2^{3j}}{(1+2^{j}(r+s))} \frac{\omega_{j,\lambda}(s)}{(1+2^{j}|r-s|)^{N}} \left(\frac{\omega_{j,\lambda}(r)}{r 2^{j}}+
\tilde \omega_{j,\lambda}(r) \right),
\end{equation} 

\begin{equation}\label{ker3}
|K_{j}^{1,\lambda}(r,s)| \leq \frac{C_{N} 2^{j}}{(1+2^{j}(s+r))} \frac{\omega_{j,\lambda}(s)}{(1+2^{j}|s-r|)^{N}} \tilde \omega_{j,\lambda}(r). 
\end{equation}

In addition, if  $2^{j} r \leq 1$ and $2^{j} \leq \lambda$, then $K_{j}^{1,\lambda}$ can be decomposed as
\[
K_{j}^{1,\lambda}(r,s) = K_{j,res}^{1,\lambda}(r,s)+K_{j,reg}^{1,\lambda}(r,s)
\]
where
\[
K_{j,res}^{1,\lambda}(r,s) = \lambda \phi_{0}(r\lambda) \int_{0}^{\infty}  \frac{q(\xi) \widetilde{m_{j}}(\xi)}{\xi} \psi_{\lambda^{-1}\xi}(s\lambda) d\xi
\]
and we have
\begin{equation}\label{ker4}
\begin{split}
& |K_{j,res}^{1,\lambda}(r,s)| \leq C_{N}\phi_{0}(r \lambda) \mathbbm{1}_{\{r 2^{j} \leq 1\}} 2^{j} \frac{\min \{1,\dfrac{s 2^{j} s^{2}\lambda^{2}}{\langle s^{2}\lambda^{2}\rangle}\}}{(s 2^{j}+1)^{N}},  \\
& |K_{j,reg}^{1,\lambda}(r,s)| \leq \frac{C_{N} 2^{j}}{((r+s) 2^{j} +1)} \frac{\omega_{j,\lambda}(r)}{(1+2^{j}|r-s|)^{N}} \omega_{j,\lambda}(s).
\end{split}
\end{equation}

\end{p1}
We remark that in these statements $\lambda$ is just a scale parameter, and all the bounds are equivalent to the ones for $\lambda = 1$. The same applies to the proof below.\\
\\
We also note that the estimates on $K_{j}^{\lambda}$ in the proposition above are also true for the kernel of $\widetilde{P_{j}^{\lambda}}$.
\begin{proof}
We begin with \eqref{ker1}, for which we consider several cases:
\medskip

(i)  $2^{j}r \leq 1$ and $2^{j}s \leq 1$.  From the formulae for $\psi_{\xi}(r)$, we obtain the bound
\[
|\psi_{\frac{\xi}{\lambda}}(r\lambda)| \leq \frac{C 2^{j/2}}{\sqrt{\lambda}} \omega_{j,\lambda}(r), \quad \xi \sim 2^{j}, \quad r\xi \lesssim 1.
\]
Then  we directly estimate using the definition of $K_{j}^{\lambda}$ to get
\[
|K_{j}^{\lambda}(r,s)| \leq C 2^{2j} \omega_{j,\lambda}(r)\omega_{j,\lambda}(s), \quad r,s \leq 2^{-j}.
\]
\medskip
(ii)  $r \leq 2^{-j} \leq s$ (and the symmetric case)
Here we have
\[
K_{j}^{\lambda}(r,s) = \lambda \frac{2}{\sqrt{s \lambda}} \text{Re}\left(\int_{0}^{\infty}  m_{j}(\xi) \psi_{\frac{\xi}{\lambda}}(r\lambda)a(\frac{\xi}{\lambda}) e^{i s \xi} \widetilde{\sigma}(s\xi,s) d\xi \right).
\]
Then we integrate by parts, using the symbol-type estimates on $\widetilde{\sigma}$ and $a$, to obtain
\[
|K_{j}^{\lambda}(r,s)| \leq \frac{C_{N} 2^{2j}}{(2^{j}s)^{N}} \omega_{j,\lambda}(r).
\]
This suffices, since $2^{j} s \geq 1$, so $\omega_{j,\lambda}(s) \sim 1$. 

\medskip
(iii) $2^{-j} \leq r \leq s$. Then  we have
\begin{equation}
K_{j}^{\lambda}(r,s) = \int_{0}^{\infty} d\xi \frac{m_{j}(\xi)}{\sqrt{r}\sqrt{s}} \begin{aligned}[t]&\left(e^{-i(r+s)\xi}(\overline{a}(\frac{\xi}{\lambda}))^{2} \overline{\widetilde{\sigma}}(r\xi,r\lambda) \overline{\widetilde{\sigma}}(s\xi,s\lambda)\right.\\
&\left.+e^{-i(r-s)\xi}|a(\frac{\xi}{\lambda})|^{2} \overline{\widetilde{\sigma}}(r\xi,r\lambda) \widetilde{\sigma}(s\xi,s\lambda)\right.\\
&\left.+e^{-i(s-r)\xi}|a(\frac{\xi}{\lambda})|^{2} \widetilde{\sigma}(r\xi,r\lambda) \overline{\widetilde{\sigma}}(s\xi,s\lambda)\right.\\
&\left.+e^{i(r+s)\xi}(a(\frac{\xi}{\lambda}))^{2} \widetilde{\sigma}(r\xi,r\lambda) \widetilde{\sigma}(s\xi,s\lambda)\right).
\end{aligned}
\end{equation}
If  $2^{j}|r-s| \lesssim  1$, then we estimate the above integral directly to get
\[
|K_{j}^{\lambda}(r,s)| \leq \frac{C 2^{2j}}{\sqrt{r}\sqrt{s} 2^{j}}  \leq \frac{C 2^{2j} \omega_{j,\lambda}(r)\omega_{j,\lambda}(s)}{(1+2^{j}(r+s))}.
\]
If $\frac{s}{2} \leq r \leq s$ and $2^{j}|r-s| \geq 1$ then we integrate by parts to obtain
\[
|K_{j}^{\lambda}(r,s)| \leq \frac{C 2^{j}}{\sqrt{rs}} \frac{1}{2^{jN}}\frac{1}{|r-s|^{N}}  \leq \frac{C 2^{2j} \omega_{j,\lambda}(r)\omega_{j,\lambda}(s)}{(1+2^{j}(r+s))(1+2^{j}|r-s|)^{N}}.
\]
Finally, if  $r \leq \frac{s}{2}$ then we integrate by parts, to obtain instead
\[
|K_{j}^{\lambda}(r,s)| \leq \frac{C 2^{2j}}{(2^{j} s)^{N}} \omega_{j,\lambda}(r)\omega_{j,\lambda}(s).
\]
The estimates for $K^{1,\lambda}_j$ are obtained using the same procedure. More precisely, 
the expression for $K_{j}^{1,\lambda}$ is exactly that of $K_{j}^{\lambda}$, except that $\psi_{\frac{\xi}{\lambda(t)}}(s \lambda(t))$ is replaced by $\xi^{-1} \phi_{\frac{\xi}{\lambda(t)}}(s\lambda(t))$. This explains the similar form of the corresponding estimates, and the extra factor of $2^{-j}$ in the $K_{j}^{1,\lambda}$ estimate, compared to that of $K_{j}^{\lambda}$ . To estimate $\partial_r K_j^\lambda$, we use the fact that
\[
L^{*}_\lambda= -\partial_r + \frac{2 h_3^\lambda-1}{r}
\]
and estimate $L^{*}_{\lambda} K^\lambda_{j}(r,s)$. We have
\[
L^{*}_{\lambda} K_{j}^{\lambda}(r,s) = \lambda \int_{0}^{\infty}  \phi_{\frac{\xi}{\lambda}}(r\lambda) \xi m_{j}(\xi) \psi_{\frac{\xi}{\lambda}}(s\lambda)  d\xi.
\]
which has a similar form as $K_{j}^{1,\lambda}(r,s)$, except for an extra factor of $\xi^{2}$ in the integrand.

For \eqref{ker4} we write 
\[
\phi_{\frac{\xi}{\lambda}}(r\lambda) =\phi_{0}(r\lambda) q(\frac{\xi}{\lambda}) +  \phi_{\frac{\xi}{\lambda}}(r\lambda) - \phi_{0}(r\lambda) q(\frac{\xi}{\lambda}),
\]
and note that $K_{j,res}^{1,\lambda}$ is obtained by inserting the first term of the above decomposition into the definition of $K_{j}^{1,\lambda}$.
Then we have
 \[
 |K_{j,res}^{1,\lambda}(r,s)| \leq C \lambda \phi_{0}(r \lambda) \mathbbm{1}_{\{r 2^{j} \leq 1\}} \int_{0}^{\infty} \frac{q(\frac{\xi}{\lambda}) |\widetilde{m_{j}}(\xi)|}{\xi} |\psi_{\frac{\xi}{\lambda}}(s\lambda)| \, d\xi.
 \]
  If $2^{j} s \lesssim 1$ then we directly estimate the integral, while for $2^{j} s \gtrsim 1$ we write
  \[
  K_{j,res}^{1,\lambda} = \lambda \phi_{0}(r\lambda) \mathbbm{1}_{\{r 2^{j} \leq 1\}} \int_{0}^{\infty} d\xi \frac{q(\frac{\xi}{\lambda}) \widetilde{m_{j}}(\xi)}{\xi} \cdot 2 \text{Re}\left(a(\frac{\xi}{\lambda}) \psi_{\frac{\xi}{\lambda}}^{+}(s\lambda)\right)
  \]
and then we integrate by parts.  Next, we have
 \[
 K_{j,reg}^{1,\lambda} = \lambda \mathbbm{1}_{\{r 2^{j} \leq 1\}} \sum_{l=1}^{\infty} \int_{0}^{\infty} d\xi q(\frac{\xi}{\lambda}) (r\xi)^{2l} \psi_{l}(r^{2}\lambda^{2}) \frac{\widetilde{m_{j}}(\xi)}{\xi} \psi_{\frac{\xi}{\lambda}}(s\lambda).
 \]
 We again directly estimate the integral in the case $2^{j} s \lesssim 1$ and integrate by parts as above for the case $2^{j} s \gtrsim 1$.
\end{proof}

 Using the above kernel type estimate we derive the following Bernstein-type estimate.
\begin{l1}\label{l:BernsteinG} 
 Assume that $f \in L^{2}(r dr)$ is localized at frequency $2^k$ in the $\tilde H_\lambda$ calculus; then for any $1 \leq q \leq p \leq \infty$ the following holds true:
\begin{equation} \label{BtH}
\| f\|_{L^{p}(r dr)} \les  2^{2k(\frac1q-\frac1p)} \|f\|_{L^{q}(r dr)}.
\end{equation}

\end{l1}

\begin{proof}
 Without restricting the generality of the argument we can assume that $f=\tilde P_k^\lambda f$,
thus we have
\[
f(r) = \int \tilde K_k^\lambda (r,s) f(s) sds.
\]
The bounds on the kernel $K_k^\lambda (r,s)$ are provided in \eqref{ker1}. We write
\[
f= \chi_{A_{\leq -k}} f + \sum_{j \geq 1} \chi_{[2^{-k}j,2^{-k}(j+1)]} f: = f_0 + \sum_{j \geq 1} f_j. 
\]
For the first term we simply use
\[
|\int \tilde K_k^\lambda (r,s) f_0(s) sds| \les \int 2^{2k} (1+2^k r)^{-N} |f_0(s)| sds,
\]
from which we obtain
\[
\| \int \tilde K_k^\lambda (r,s) f_0(s) sds \|_{L^p} \les 2^{2k(\frac1q-\frac1p)} \|f_0\|_{L^{q}(r dr)}
\]
with improved decay away from $A_{-k}$, that is
\[
\| \int \tilde K_k^\lambda (r,s) f_0(s) sds \|_{L^p[2^{-k}j,2^{-k}(j+1)]} \les \la j \ra^{-N} 2^{2k(\frac1q-\frac1p)} \|f_0\|_{L^{q}(r dr)}, \quad \forall j \geq 1.
\]
For the terms in the sum, we use 
\[
|\int \tilde K_k^\lambda (r,s) f_j(s) sds| \les \int 2^{2k} (1+2^k(r+s)) (1+2^k |r-s|)^{-N} f(s) sds,
\]
and then rely on Schur's test to obtain 
\[
\|\int \tilde K_k^\lambda (r,s) f_j(s) sds \|_{L^p} \les 2^{2k(\frac1q-\frac1p)} \|f_j\|_{L^{q}(r dr)}.
\]
Improved decay is away from $[2^{-k}j,2^{-k}(j+1)]$ is also available, and this allows us to retrieve the estimate \eqref{BtH} from the estimates above.

\end{proof}

\subsection{The time dependence of the Littlewood-Paley projectors}
Throughout this section, we have studied the spectral properties of the operators $H_\lambda$ and $\tilde H_{\lambda}$ where $\lambda \in (0,\infty)$ is a scaling parameter. 
As noted in Section~\ref{Sec3wrap}, this scale parameter is needed because the linear operator $\tilde H_\lambda=-\Delta + \tilde V_\lambda$, that naturally occurs in the PDE governing the dynamics of the field $\psi$, has in fact a time dependent potential coming precisely from the fact that $\lambda=\lambda(t)$. 

In Section \ref{s:linear} we will seek to establish estimates for the linear equation \eqref{Linpsi} 
\[
(i \partial_t - \tilde H_\lambda) \Psi =0,
\]
when $\lambda=\lambda(t)$. In doing so, we will use the appropriate form of a Littlewood-Paley decomposition, at which stage we encounter the operator $[\partial_{t},P_{k}^\lambda]$, whose various operator norms need to be estimated. 

For convenience, we found it more efficient to prepare the study of this operator in this section. Thus, in what follows, we assume $\lambda: I \rightarrow \R$, where $I \subset \R$ can be any interval. In terms of regularity, we will eventually work with $\lambda \in \dot H^\frac12 \cap L^\infty$, although for all practical purposes one can assume that $\lambda$ is continuously differentiable, by Corollary \ref{cdif}.

Recall that the Fourier transform associated to $\tilde H_\lambda$ 
is obtained by rescaling from the Fourier transform associated to $\tilde H$, see \eqref{Flambda-def}.

We also recall, from \eqref{pjlambdadef}, the spectral projectors $P^\lambda_{j}  $ associated to $\tilde H_\lambda$.
We are now interested in the time derivative of the spectral projectors.
The first step will be to express it in terms of the transference operator.

\begin{l1} \label{dtP}
The commutator of $\partial_{t}$ with $P_{k}^{\lambda(t)}$ is given by the following:
\begin{equation}\label{dt-P}
[\partial_{t},P_{k}^{\lambda(t)}]u=\lambda'(t) \FtH^{-1}_{\lambda(t)} [\tcK_{\lambda(t)},m_{k}]\FtH_{\lambda(t)} u
\end{equation}
\end{l1}
where $\tcK_{\lambda}$ is as in \eqref{klambdadef}.

\begin{proof}
 The proof is a direct computation using the relations \eqref{dlambda-F}:
\[
\partial_{t} P_{k}^{\lambda(t)} = 
\lambda'(t) \partial_\lambda (\FtH^{-1}_\lambda m_{k} \FtH_\lambda) = \lambda'(t) (\FtH^{-1}_\lambda \tcK_\lambda m_{k} \FtH_\lambda - \FtH^{-1}_\lambda m_{k}  \tcK_\lambda \FtH_\lambda) 
\]
and \eqref{dt-P} follows. 

\end{proof}

Given the above commutator formula, it is natural
to seek estimates for the commutator $[\tcK_\lambda,m_j]$, where we recall the notation \eqref{klambdadef}. Heuristically, we expect this commutator 
to be largest when $\lambda$ and $j$ are matched,
i.e. $\lambda \approx 2^j$. To measure the decay away from this region we will use the weight
\[
\chi_{\lambda=2^{j}} = \left(\frac{2^{j/2}}{\sqrt{\lambda}} \mathbbm{1}_{\{\frac{2^{j}}{\lambda} \leq 4\}}+\frac{\lambda}{2^{j}} \mathbbm{1}_{\{\frac{2^{j}}{\lambda} \geq 1\}}\right).
\]
While $m_j$ provided frequency localization, 
$\tcK_\lambda$ does not, so neither does the commutator.
However, we also expect to have decay at frequencies 
$2^k$ when $k$ is away from $j$. For this we will use 
the weight 
\[
\chi_{k=j} = 2^{-\frac{|k-j|}2}.
\]
Note that if $|k-j| > 2$ then we have
\[
m_k [ \tcK_\lambda,m_j] =
m_k \tcK_\lambda m_j  = [m_k, \tcK_\lambda] m_j,
\]
so it does not matter if we place the $m_k$ localization on the left or on the right. For symmetry we will include both. Our main bounds in this section are collected in the following Lemma.

\begin{l1}\label{l:mj-com}
The following fixed time estimates hold true:
\begin{equation}\label{com-est}
 \begin{split}
 \|m_l [\tcK_\lambda,m_{j}]m_{k}u\|_{L^{2}(d\xi)} \lesssim \frac{2^j}{\lambda^2} \chi_{\lambda=2^{k}} \chi_{j=k} \chi_{j=l} \|u\|_{L^{2}(d\xi)}
 \end{split}
 \end{equation}
respectively
\begin{equation}\label{rcom-est}
\begin{split}
&\| r \FtH_\lambda^{-1} m_l [\tcK_{\lambda},m_{j}] m_k u\|_{L^{2}(r dr)} \lesssim \frac{1}{\lambda^{2}} \chi_{\lambda=2^{j}} \chi_{j=k} \chi_{j=l}  \|u\|_{L^{2}(d\xi)}.
\end{split}
\end{equation}

\end{l1}
\begin{proof}
 We first remark that $\lambda$ is a scaling parameter
in this Lemma. Hence, without any loss in generality 
we can simply set $\lambda = 1$ in the proof.
Using the  representation in \eqref{tFH} for the kernel of $\tcK$, it is easily seen that the operator
$m_l[m_j,\tcK]m_k$ can be expressed in the form
\begin{equation}\label{commutatorKformula}\begin{split} &m_l[m_j,\tcK]m_kf(\xi):=- 
		\int_{0}^{\infty}  \frac{m_j(\xi)-m_j(\eta)}{\xi^{2}-\eta^{2}} m_l(\xi) m_k(\eta) \tF(\xi,\eta) f(\eta)  d\eta
\end{split}\end{equation}
with the integral kernel
\[
K_{jkl}(\xi,\eta) =  \frac{m_{j}(\xi)-m_{j}(\eta)} {\xi^{2}-\eta^{2}} m_l(\xi) m_k(\eta) \tF(\xi,\eta),
\]
where we note that one consequence of taking the commutator is that the kernel no longer has a singularity on the diagonal. This kernel vanishes 
unless we have either $|l-j| < 4$ or $|k-j|< 4$.
Assuming this holds, we note that $K_{jkl}(\xi,\eta)$
is supported in the region $\{ \xi \approx 2^l, \
\eta \approx 2^k \}$, and satisfies 
\[
|K_{jkl}(\xi,\eta)| \lesssim \frac{|\tF(\xi,\eta)|}{(\xi+\eta)^2}, \qquad |\partial_\xi K_{jkl}(\xi,\eta)|
+ \xi^{-1}|K_{jkl}(\xi,\eta)|
\lesssim  \frac{|\tF(\xi,\eta)|}{\xi (\xi+\eta)^2}
+  \frac{|\partial_\xi \tF(\xi,\eta)|}{(\xi+\eta)^2}.
\]
To prove \eqref{com-est} it suffices to have the
\[
\|K_{jkl}(\xi,\eta)\|_{L^2_{\xi,\eta}}
\lesssim 2^j \chi_{k=0} \chi_{j=k} \chi_{j=l}
\]
which is easily verified using the above kernel bound and \eqref{fhests-th}. 

For \eqref{rcom-est} we use Lemma~\ref{multrfourier},
which reduces it to the bounds
\[
\|\partial_\xi K_{jkl}(\xi,\eta)\|_{L^2_{\xi,\eta}}
+ \|\xi^{-1} K_{jkl}(\xi,\eta)\|_{L^2_{\xi,\eta}}
\lesssim \chi_{k=0} \chi_{j=k} \chi_{j=l}
\]
which is again easily verified using the above kernel bound together with both \eqref{fhests-th} and \eqref{dfhests-th}. 
\end{proof}

\section{Elliptic analysis and \texorpdfstring{$\ell^1$}{} Besov structures}
\label{s:XLX}

Most of the analysis in this paper is concentrated around the gauge field $\psi \in L^2$ and the modulation parameters $\alpha$ and $\lambda$. In Section \ref{seccoulomb} we have shown why this is
essential in order to understand  the dynamics of the original map $u(t)$, which is uniquely determined by $(\psi,\alpha,\lambda)$.
In particular one may think of  $Q_{\alpha(t),\lambda(t)}$ as the closest soliton to $u$, and the orthogonality condition \eqref{ldef} ensures that our choice of $\alpha,\lambda$ satisfies the condition   (see   \eqref{goodal2})
\[
\|u-Q^2_{\alpha,\lambda}\|_{\dot H^1} \approx \|\psi\|_{L^2}.
\]

Our results assert that the smallness of $\psi$ in $L^2$ 
suffices in order to guarantee global well-posedness for the Schr\"odinger map flow,
but likely this does not preclude blow-up at infinity.
Then a natural question becomes whether there exists some slightly stronger topology for the initial data $u(0)$, and correspondingly for $\psi(0)$,  where we have soliton stability, and  in particular no blow-up can happen.

This turns out to be indeed the case, and the appropriate
spaces for $u(0)$, respectively $\psi(0)$  are Besov spaces
with $\ell^1$ summation, precisely 
\begin{equation}
 X = \dot B^{1}_{2,1}, \qquad \bX = \dot B^{1}_{2,1,e}, \qquad \LX = \dot B^{0}_{2,1,e},   
\end{equation}
see Theorem~\ref{tmain-L1}, Theorem~\ref{tmain-G1} and
Theorem~\ref{tstable}. The definitions of these spaces, along with alternative characterizations, are provided in Section \ref{defnot}.
The notations $X$ and $\LX$ are inherited from \cite{BeTa-1}, and justified by the equivalent characterization of these spaces in the next subsection.

Our analysis is done primarily at the level of $(\psi, \lambda,\alpha)$, while the original problem is stated at the level of the map $u$. The goal of this section is to describe the transition between the two settings, and in particular to establish the norm equivalence
\begin{equation}
    \|u-Q^{2}_{\alpha,\lambda}\|_{X} \approx \|\psi\|_{\LX}
\end{equation}
The main result of this section is the following:

\begin{p1} \label{uQpsi}
We assume the setup from Proposition \ref{p:mod} part i). Then

i)  If $\|u-Q^2_{\alpha,\lambda}\|_{X} < \infty$, then the associated field obeys $\|\psi\|_{\LX} \les \|u-Q^2_{\alpha,\lambda}\|_{X}$.

ii) Vice versa, if the field $\psi \in \LX$, then $\|u-Q^2_{\alpha,\lambda}\|_{X} \les \|\psi\|_{\LX}$.
\end{p1}
The rest of this section is devoted to the proof of the above 
proposition, and is organized as follows. First we cover some basic properties of these spaces including the mapping properties of the operator $L$ in this context and some algebra properties. Then we introduce a companion space $\tilde X$ which, while morally at the level of $\bX$, allows for slightly more general functions which arise in our analysis. This space is well suited for characterizing the solutions of an ODE which plays a key role in the analysis of  the transfer information between $u-Q^2_{\alpha,\lambda}$ and $\psi$; in fact this system was analyzed earlier in the context of energy setup, see Lemma \ref{lZL2}. With all these at hand, we finish the section with the proof of the Proposition above. 

\subsection{An equivalent characterization of the \texorpdfstring{$\bX$}{} and \texorpdfstring{$\LX$}{} spaces}

Above (and in more detail in Section \ref{defnot}) we have defined $\bX$ and $\LX$ as classical Besov spaces restricted to equivariant functions. On the other hand in our 
earlier work on the $1$-equivariant case in \cite{BeTa-1} we have defined 
the counterparts of these spaces based on 
the spectral decomposition associated to the  $H_\lambda$ and $\tilde H_\lambda$ operators. It is then a natural question whether
the same can be done here, which we shall answer in the affirmative. 

To both justify the notation $\LX$ and as a starting point for the subsequent discussion,
we first show that the spaces $\bX$ and $\LX$
can be connected via the operator $L$:

\begin{l1} \label{LLwf}
a) The operator $L$ maps $\bX$ to $\LX$, and
\begin{equation}\label{XtoLX}
\| Lu \|_{\LX} \lesssim \|u\|_{\bX}    
\end{equation}

b) Conversely, given $f \in \LX$, there exists a solution $w$ 
to $Lw = f$ with 
\begin{equation}\label{LXtoX}
\|w \|_{\bX} \lesssim \| f\|_{\LX}  
\end{equation}
This holds in particular if $w$ is chosen either so that 
$w(1) = 0$, or if $w$ is chosen to be orthogonal to $h_1$. 
\end{l1}

In this result the spaces $\bX$ and $\LX$ are homogeneous, while the operator $L_\lambda$ satisfies the scaling relation
\[
L_\lambda u^\lambda = \lambda (Lu)^\lambda .
\]
Hence by rescaling one may easily replace $L$ by $L_\lambda$, without affecting the implicit constants in the estimates in the lemma.

\begin{proof}

a) It is easily seen that
\[
L: \dot H^2_e \to \dot H^1_e, 
\qquad L: L^2 \to \dot H^{-1}_e,
\]
where the second can also be seen as the dual property to 
\[
 L^*: \dot H^1_e \to L^2.
\]
Then the property \eqref{XtoLX} directly follows
from the equivalent characterization of the two norms in \eqref{X-relax}, \eqref{LX-relax}.

\bigskip

b) Given $f \in  \LX$ we split it 
into $f = f_1+f_2$ where the two terms are supported in  $r < 2$, respectively $r > 1$. In the two regions  we replace the equation 
\begin{equation}\label{L-inhom}
Lw = f,
\end{equation}
where
\[
L = \partial_r + \frac{2}r h_3 = \partial_r + \frac{2}r \frac{r^{4}-1}{r^4+1},
\]
with
\[
(\partial_r - \frac{2}{r} ) w_1 = f_1, \qquad w_1([2,\infty)) = 0,
\]
respectively
\[
(\partial_r + \frac{2}{r}) w_2 = f_2. \qquad w_2[(0,1])= 0.
\]
Here the solutions $w_1,w_2$ have a similar support,
and are easily seen to satisfy good scale invariant Besov bounds
\[
\|w_1\|_{\dot B^1_{2,1,e}} \lesssim \|f_1\|_{\dot B^0_{2,1,e}}, \qquad \|w_2\|_{\dot B^1_{2,1,e}} \lesssim \|f_2\|_{\dot B^0_{2,1,e}}.
\]
These are proved by interpolating between the corresponding $\dot H^1_e \to \dot H^2_e$ and $\dot H^{-1}_e \to L^2_e$ bounds, which in turn reduce 
to one dimensional Hardy type inequalities.

 Then we get an approximate solution
\[
w_{app} = w_1 + w_2
\]
for \eqref{L-inhom}, so that 
\[
Lw_{app} = f+ g, 
\]
where the extra source term $g$ has the form
\[
 g \approx r^3 w_1  + r^{-5} w_2,
\]
which has regularity $g \in H^1_e$ and is also localized near $r=1$.

Now it remains to solve $Lz = g$, where it is convenient to start with the initial condition $g(1) = 0$, using the fundamental solution for $L$. This  yields a solution $z \in H^2_e$, which is much better than needed.

The solution $w = w_{app} - z$ which we have constructed satisfies the bound \eqref{LXtoX}, which in particular implies the bounds
\[
|w(1)| \lesssim \| f\|_{LX}, \qquad 
\la w,h_1\ra \lesssim \| f\|_{LX}
\]
This allows us to correct $w$ with a well chosen multiple of $h_1$ to insure that either $w(1) = 0$ or $\la w,h_1\ra  = 0$.

\end{proof}

An immediate consequence of the above Lemma and its proof is the following:
\begin{equation} \label{drLX}
\|\partial_r f \|_{\LX} + \|\frac{f}r\|_{\LX} \les \|f\|_{\bX}, \quad \forall f \in \bX. 
\end{equation}
This can be easily seen from the arguments used in the proof of part a) of the Lemma. 
 
\bigskip

To set the stage for what follows, we note that throughout this paper we use the notation $P_k^\lambda$ for the Fourier projectors in the $\tilde H_\lambda$ frame; this is because most of the analysis is carried out at the level of the gauge field $\psi$ whose dynamics uses the $\tilde H_\lambda$ operator. This is the only section in the paper where we need to use projectors in both frames $H_\lambda$ and $\tilde H_\lambda$, and we need to differentiate between them at the notation level. Thus we use the notation $P_k^{H_\lambda}$ for the Fourier projectors in the $H_\lambda$ frame, 
and $P_k^{\tilde H_\lambda}$ for the Fourier projectors in the $\tilde H_\lambda$. 
If $\lambda=1$, then we drop the index $\lambda$ from the operators and simply use 
$P_k^{H}$ and $P_k^{\tilde H}$ instead. 

Now we are ready to state our equivalent characterizations of $\bX$ and $\LX$:

\begin{p1}
a) The space $\LX$ can be equivalently characterized as the space of functions
$f \in L^2$ for which the following sum is finite, with equivalent norms:
\begin{equation}\label{LX-equiv}
\|f \|_{\LX} \approx \sum_{k \in \Z} \| P^{\tH_\lambda}_k f\|_{L^2}.
\end{equation}

b) The space $\bX$ can be equivalently characterized as the space of functions
$u \in \dot H^1_e$ for which the following sum is finite, with equivalent norms:
\begin{equation} \label{X-equiv}
\|u \|_{\bX} \approx \|u\|_{\dot H^1_e}+ \sum_{k \in \Z} 2^k \| P^{H_\lambda}_k u\|_{L^2}.
\end{equation}
\end{p1}

\begin{proof}
To start with, we note that the left hand side 
in both \eqref{LX-equiv} and \eqref{X-equiv} does not depend on $\lambda$, but the right hand side apriori does. Hence $\lambda$ only plays the role
of a scaling parameter, which we  can harmlessly  set  to $\lambda = 1$.

We begin with the relation \eqref{LX-equiv},
for which we use the equivalent $\LX$ norm given by 
\eqref{LX-relax}. On the other hand for the right hand side we have a similar equivalent norm but 
using the $\dot H^1$ and $\dot H^{-1}$ spaces 
associated to the operator $\tH$.

But then \eqref{LX-equiv} is  straightforward due to the fact that the standard $\dot H^1_e$ norm  and the one associated to $\tH$ are equivalent (which is to say, $\tH$ is coercive in $\dot H^1_e$), and correspondingly the associated $\dot H^{-1}$ norms are equivalent.

The bound \eqref{X-equiv} is more interesting. Using the relation
\[
(\F^{\tilde H} P_k^{\tilde H} Lf)(\xi)= m_k(\xi) \xi (\F^{H} f)(\xi) = \xi (\F^{H} P_k^H f)(\xi)
\]
and \eqref{LX-equiv}
we immediately see that the right hand side 
of \eqref{X-equiv} may be equivalently 
written as 
\[
\|u\|_{\dot H^1_e} + \|Lu\|_{\LX}. 
\]
By part (a) of Lemma~\ref{LLwf} we immediately get the bound
\[
\|u\|_{\dot H^1_e} + \|Lu\|_{\LX} \lesssim \|u\|_{\bX}.
\]
For the opposite inequality let $w \in \bX$ be, as in Lemma~\ref{LLwf}, the unique solution to $Lw = Lu$ which is orthogonal to $h_1$.
Then we have $L(u-w) = 0$ which implies that 
$u-w$ is a multiple  of $h_1$.
Hence, using part (b) of Lemma~\ref{LLwf},
\[
\| u\|_{\bX} \lesssim \| u-w\|_{\bX} + \| w\|_{\bX}
\approx \| u-w\|_{\dot H^1} + \| w\|_{\bX}
\lesssim \| u\|_{\dot H^1} +\|Lw\|_{\LX},
\]
as desired.
\end{proof}

In connection with \eqref{X-equiv}, we remark that 
$P_k^H u$ is orthogonal to $h_1$, therefore 
by Lemma~\ref{LLwf} we have the equivalence
\[
\| P_k^H u \|_{\bX} \approx \| P_k^{\tH} Lu\|_{L^2}
\approx 2^k \| P_k^H u \|_{L^2}.
\]
This implies that for $u \in \bX$ the sum 
\[
\sum_{k \in \Z} P_k^H u 
\]
converges in $\bX$ and is orthogonal to $h_1$. This 
in turn yields the representation
\begin{equation}\label{X-rep}
u = c h_1 + \sum_{k \in \Z} P_k^H u,
\end{equation}
where \eqref{LXtoX} can be equivalently interpreted as 
\begin{equation}\label{LX-equiv+}
\|u \|_{\bX} \approx |c| + \sum_{k \in \Z} 2^k \| P^{H_\lambda}_k f\|_{L^2}.
\end{equation}

The above analysis highlights a key property that functions in the space $\bX$ enjoy, and which functions in $\dHe$ do not. In the above we have established that for any $u \in \bX$ we are able to meaningfully project $u$ on $h_1$. This cannot be done for $f \in \dHe$ because in the classical sense we would need to make sense of the quantity $\la f, h_1 \ra= \la \frac{f}r, r h_1 \ra$; the problem is that while $\|\frac{f}r\|_{\dHe} \les \|f\|_{\dHe}$, $rh_1 \notin L^2$ and there is no other way to fix this, that is to make sense of  $\la f, h_1 \ra$.

\subsection{Spectral analysis in the Bessel frame} \label{SBF}
This is needed in order to work with classical Littlewood-Paley projectors
restricted to the equivariant class of functions. We recall the definition of the Bessel function of the first kind (for instance, from page 511 of \cite{sneddon1995})  
\begin{equation}\label{j2series} J_{m}(x) = \sum_{k=0}^{\infty} \frac{(-1)^{k} \left(\frac{x}{2}\right)^{m+2k}}{k! (k+m)!}, \quad m \in \mathbb{N}.\end{equation} 
From page 52-53 of \cite{sneddon1995}, we recall the definition of the Hankel transform of order $m$ of $f$:
\[
\mathcal{F}_{m}(f)(\xi) = \int_{0}^{\infty} f(r) J_{m}(r \xi) r dr.
\]
The inversion formula is
\[
f(r) = \int_{0}^{\infty} \mathcal{F}_{m}(f)(\xi) J_{m}(r \xi) \xi d\xi.
\]
Page 60 of \cite{sneddon1995} gives the $L^{2}$ isometry property of $\mathcal{F}_{m}$:
\[
\int_{0}^{\infty}  f(r) g(r) r dr = \int_{0}^{\infty}  \mathcal{F}_{m}(f)(\xi) \mathcal{F}_{m}(g)(\xi) \xi d\xi.
\]
We will also use the following asymptotic expansion of Bessel functions of the first kind (for instance, when $m \in \mathbb{N}$)
from page 199 of \cite{watson1922} 
\begin{equation}\label{j2asymp}J_{m}(x) \sim \left(\frac{2}{\pi x}\right)^{1/2} \left(\cos(x-\frac{\pi}{2} m  - \frac{\pi}{4}) \sum_{k=0}^{\infty}  \frac{a_{k}(m)}{x^{2k}} - \sin(x-\frac{\pi}{2} m  - \frac{\pi}{4}) \sum_{k=0}^{\infty} \frac{b_{k}(m)}{x^{2k+1}}\right).\end{equation}

If $f_k= P_k f = \mathcal{F}_2^{-1} ( m_k \mathcal{F}_2 f)$ (this is the projector at frequency $\approx 2^k$ in the $J_2$ frames), then we record the following pointwise bound:\begin{l1} If $ \mathcal{F}_{2}f$ is localized at frequencies $\approx 2^{k}$, then,
\begin{equation} \label{fjpd}
|f_{k}(r)|  \leq C \|f_{k}\|_{L^{2}} 2^{k} \begin{cases} (2^{k} r)^{-1/2}, \quad   r 2^{k+2} \geq 1\\
   ( 2^{k} r)^{2}, \qquad  r 2^{k+2} < 1.\end{cases}
\end{equation}\end{l1}

This follows from directly estimating the definition of $f_{k}$, using the following pointwise estimates on $J_{2}$ (which follow from \eqref{j2series} and \eqref{j2asymp}): $$|J_{2}(x)| \leq C \begin{cases} x^{2}, \quad 0<x \leq 1\\
\frac{1}{\sqrt{x}}, \quad x >1\end{cases}$$
We will also make use of the following lemma, which describes the action of $\partial_{r}$ (conjugated by a power of $r$) on frequency localized functions.
\begin{l1} \label{gg1dec}
Assume $\mathcal{F}_{2} g \in L^2$ is localized at $\xi \approx 2^k$. Then we can write
\[
g= (\partial_r + \frac{3}r) g^1, 
\]
where $\mathcal{F}_{3} g^1$ is localized at frequency $\xi \approx 2^k$. In addition the following holds true
\begin{equation}\label{reverse-g}
 2^k \| g^1 \|_{L^2} \les \|g\|_{L^2}, \quad 
\|g^1\|_{\dHe} \les  \|g\|_{L^2}.
\end{equation}
\end{l1}

\begin{proof} The conclusion of the lemma  can be obtained by replacing $m$ by 2 everywhere in the proof below (which is true for any $m \geq 1$). Using the Hankel transform of order $m$, we write:
\begin{equation}
\begin{split}
\Delta \Delta^{-1}g &= \left(\partial_{r} + \frac{m + 1}{r}\right) \left(\partial_{r} - \frac{m}{r}\right) \left(\int_{0}^{\infty}  J_{m}(r \xi) \frac{1}{-\xi^{2}}\mathcal{F}_{m}(g)(\xi)\xi d\xi\right)\\
&=\left(\partial_{r}+\frac{m+1}{r}\right)\left(\int_{0}^{\infty}  J_{m+1}( r\xi) \frac{1}{\xi} \mathcal{F}_{m}(g)(\xi) \xi d\xi\right)\\
&:=\left(\partial_{r}+\frac{m+1}{r}\right)g^{1}(r)
\end{split}
\end{equation}
Taking $m=2$ leads to our specific choice of $g^1$. It is also clear from the above that $g^1$ has the support of its Hankel transform (of order $3$) compactly supported in the region $\xi \approx 2^k$. Next, using the above formulas, we have
\[
\|g^1\|_{L^2} \approx \| \frac{\mathcal{F}_{m}(g)}{\xi} \|_{L^2} \approx 2^{-k} \|g\|_{L^2}. 
\]
We obtain the second bound in \eqref{reverse-g} by noting that $g^{1}$ is naturally identified with an $m+1$-equivariant function on $\mathbb{R}^{2}$, so that we have
$$\|g^{1}\|_{\dot{H}^{1}_{e}} = \|\xi \mathcal{F}_{m+1}(g^{1})\|_{L^{2}(\xi d\xi)} = \|\xi \cdot \frac{1}{\xi} \mathcal{F}_{m}(g)(\xi)\|_{L^{2}(\xi d\xi)} = \|\mathcal{F}_{m}(g)(\xi)\|_{L^{2}(\xi d\xi)} = \|g\|_{L^{2}(r dr)}$$

\end{proof}

\subsection{ Algebra properties for \texorpdfstring{$\bX$}{} and \texorpdfstring{$\LX$}{}.}

Here we use the Besov characterization of $\bX$ and $\LX$ in order to study multiplicative properties in these spaces:

\begin{l1} \label{XLXb}
The spaces $\bX,\LX$ satisfy the following properties:

i) $\bX$ is an algebra. Further, 
if $f,g \in \bX$ then following refined bilinear estimate holds true:
\begin{equation} \label{Xalgref}
\|fg\|_{\bX} \les \|f\|_{\dHe} \|g\|_{\bX} + \|f\|_{\bX} \|g\|_{\dHe};
\end{equation}

ii) $\LX$ is stable under multiplication by functions in $\bX$; moreover, the following refined bilinear estimate holds true:
\begin{equation} \label{XLXLX}
\|fg\|_{\LX} \les \|f\|_{\bX} \|g\|_{L^2} + \|f\|_{\dHe} \|g\|_{\LX}, 
\qquad \forall f \in \bX, \ g \in \LX.
\end{equation}

\end{l1}

\begin{proof}
We first  prove 
the $\bX$ algebra property \eqref{Xalgref}.
This is most easily done using equivariant extensions in $\R^2$, where we can use classical dyadic frequency localization. If for $f,g$ below we use 2-equivariant extensions then $fg$ is 4-equivariant, but this is not a problem. 

We have
\[
\begin{aligned}
\| f g\|_{X} \lesssim & \ \sum_k 2^k  \|P_k (fg)\|_{L^2}
\\
\lesssim & \ \sum_k 2^k ( \|P_{< k} f\|_{L^\infty}\| P_k g\|_{L^2} +  \|P_{< k} g\|_{L^\infty}\| P_k f\|_{L^2}
+ 2^k \sum_{j \geq k} \|P_j f\|_{L^2} \| P_j g\|_{L^2})
\\
\lesssim & \ \| f\|_{\dot H^1} \|g\|_{X} + \|f\|_X \|g\|_{\dot H^1}.
\end{aligned}
\]
Next we prove the  product estimate in \eqref{XLXLX}, recalling the notation $LX$ for the 
$2$-equivariant two dimensional lift of $\LX$ functions. 
We have
\[
\begin{aligned}
\| f g\|_{LX} \lesssim & \ \sum_k \|P_k (fg)\|_{L^2}
\\
\lesssim & \ \sum_k  (\|P_{< k} f\|_{L^\infty}\| P_k g\|_{L^2} +  \|P_{< k} g\|_{L^\infty}\| P_k f\|_{L^2}
+ \sum_{j \geq k} 2^{k} \|P_j f\|_{L^2} \| P_j g\|_{L^2})
\\
\lesssim & \   \|f\|_{\dot H^1_e}\|g\|_{LX} +  \|g\|_{L^2}\| f\|_{X}.
\end{aligned}
\]

\end{proof}

\subsection{The companion space \texorpdfstring{$\tX$}{} and an ODE result}

The spaces $\bX$ and $\LX$ are used to measure $u-Q \in X$, respectively $\psi \in \LX$, but unfortunately they cannot be used in the context of the gauge elements $v,w$. This issue was already highlighted in the analysis in Section \ref{seccoulomb}, see the statement of Proposition \ref{lc}. In that context, we could not use the native energy space $\dot H^1$ for the frame, but instead we used $\dot H^1_C$. 

In our context, we measure the frame elements (or better their one-dimensional reduced version) in $\tilde X$ - this a Sobolev type companion $\tilde X$ for $\bX$; thus it is meant to measure functions defined on $(0,+\infty)$. We introduce the following atoms. Given $k \in \Z$ we say that the function 
$\varphi: [0,+\infty) \rightarrow \C, \varphi \in C^2$ (that is $\varphi$ is twice continuously differentiable) is an $\A_k$ atom provided that it satisfies the following properties

1) $\varphi$ is supported in the interval $[0,2^{-k+10})$;

2) $\partial_r \varphi(0)=0$;

3) it obeys the following $\|\varphi\|_{\A_k} = \sum_{\alpha=0}^2 \| 2^{-k \alpha} \partial_r^\alpha \varphi \|_{L^\infty} < +\infty$.

It is obvious that $\chi_{< -k}$ is an $\A_k$ atom; in fact the $\A_k$ atoms are meant to act as generalization of this basic function $\chi_{< -k}$. We record the following basic inequality for the atoms in $\A_k$:
\begin{equation} \label{Akbasic}
2^k  \|\partial_r  \varphi \|_{\dot H_e^{-1}}+2^{-k} \|\partial_r  \varphi \|_{\dHe} \les \|\varphi\|_{\A_k};
\end{equation}
its proof is straightforward and left as an exercise. 

 We define $\A$ to be the space of functions  $\varphi: [0,+\infty) \rightarrow \C$ which admit an atomic decomposition 
\[
\varphi = \sum_{k \in \Z} \varphi_k, \quad \varphi_k \in \A_k \ \mbox{and} \ \sum_{k \in \Z} \| \varphi_k\|_{\A_k} < \infty;  
\]
we endow this space with a norm in the standard fashion, that is
\[
\| \varphi \|_{\A} = \inf_{rep} \{ \sum_{k \in \Z} \| \varphi_k\|_{\A_k}; \varphi = \sum_{k \in \Z} \varphi_k, \varphi_k \in \A_k\},
\]
where the infimum is taken over all possible representations of $\varphi$ as a sum of atoms. 

We record here a basic inequality:
\begin{equation} \label{atomlinf}
\|\varphi\|_{L^\infty} \les \|\varphi\|_{\A},
\end{equation}
which follows from the atomic structure of $\varphi$ and the $L^\infty$ bounds available for its atoms.

The space $\tilde X$ is the sum of  $\bX$ and  $\A$; its precise definition is as follows. A function $f$ is in $\tilde X$ if it admits a decomposition
\begin{equation} \label{tXrep}
f = f_1 + f_2, \qquad f_1  \in \bX \ \mbox{and} \ f_2 \in \A;
\end{equation}
the norm in $\tilde X$ is defined as follows 
\[
\|f\|_{\tilde X} = \inf \{ \|f_1\|_{\bX} + \|f_2\|_{\A}; f=f_1+f_2, f_1  \in \bX \ \mbox{and} \ f_2 \in \A \}.
\]

A natural question to ask is how different are these two components. An element in $\varphi \in \A$ with an atomic decomposition $\varphi = \sum_k \varphi_k, \varphi_k \in \A_k$ contains an $l^1_k$ structure just as $\bX$. The characterization of elements in $\A_k$, shows that $\varphi_k$ is morally at frequency $2^k$ and \eqref{Akbasic}
establishes good bounds for $\| \partial_r \varphi_{k} \|_{L^2}$; what is missing is the equivalent bound for $\| \frac{\varphi_k}r \|_{L^2}$ and this is why we cannot place the atoms in $\dHe$. This failure is simply related to the fact that the atoms $\varphi_k$ are not required to satisfy the condition $\varphi_k(0)=0$; indeed note that we can easily estimate $\|\frac{\varphi_k-\varphi_k(0)}r \cdot \chi_{\leq -k+10}\|_{L^2}$ and this would place $\varphi_k-\varphi_k(0) \cdot \chi_{\leq -k+10}$  in $\dHe$. The need to augment the structure $\bX$ with these atoms, and arrive at the $\tilde X$ structure, stems from the analysis below where integrals of type $\int_r^\infty$ fail the basic cancellation property that $\int_0^\infty =0$, thus producing functions which do not have zero limit at $r=0$, but otherwise have properties that are very similar to those of elements in $\bX$. The result in \eqref{cLX} in Lemma \ref{cdot}  highlights again how close these elements are to $\bX$, as we establish that $\partial_r \varphi \in \LX$ for any $\varphi \in \A$.

We will also need to bound the two components on different scales, for which we use the subset  $B^{\tilde X}_{M,m}$ of $\tilde X$ defined as the functions which admit at least a representation of type \eqref{tXrep} and for which
\begin{equation} \label{BMn}
\|f_1\|_{\bX} \leq M, \quad \|f_1\|_{\dHe} + \| f_2 \|_{\A} \leq m.
\end{equation}

Our first result seeks to understand some more of the properties of functions in $\A$, and in particular how  they interact with functions in the more classical spaces introduced earlier, such as $L^2, \dHe, \LX,\bX$. 

\begin{l1} \label{cdot}
If  $c \in \A$ then the following hold true:
\begin{equation} \label{cLX}
\| \partial_r c \|_{\LX} \les \| c \|_{\A}.
\end{equation}
In addition the following multiplicative estimates hold true
\begin{equation} \label{cH1e}
\|c \cdot f\|_{\dHe} \les  \| c \|_{\A} \|f\|_{\dHe}, 
\quad \|c \cdot g\|_{L^2} \les \| c \|_{\A} \|g\|_{L^2},
\end{equation}
for any $f \in \dHe$ and $g \in L^2$, and
\begin{equation} \label{cXX}
\|c \cdot f\|_{\bX} \les  \| c \|_{\A} \|f\|_{\bX}, 
\quad \|c \cdot g\|_{\LX} \les  \| c \|_{\A} \|g\|_{\LX},
\end{equation}
for any $f \in \bX$ and $g \in \LX$. 
\end{l1}

\begin{proof} Throughout this argument we work with a decomposition of $c$ as follows
\[
c=\sum_k c_k, \quad \sum_k \|c_k\|_{\A_k} \les \|c\|_{\A}.
\]
The estimate \eqref{cLX} is a direct consequence of the definition of $\LX$ and the estimates on $\partial_r c_k$ in $\dot H^{-1}_e$ and $\dHe$ from \eqref{Akbasic}. 

The proof of \eqref{cH1e} is straightforward; using \eqref{cLX} (which implies an $L^2$ bound for $\partial_r c$) and \eqref{atomlinf} we obtain the first estimate in \eqref{cH1e} as follows
\[
\|\partial_r (cf)\|_{L^2} + \|\frac{c f}r\|_{L^2} \les \|\partial_r c\|_{L^2} \|f\|_{L^\infty} + \|\partial_r f\|_{L^2} \|c\|_{L^\infty}+ \|c\|_{L^\infty} \|\frac{f}r\|_{L^2} \les \| c \|_{\A} \|f\|_{\dHe}.
\]
The second estimate in \eqref{cH1e} follows from \eqref{atomlinf}. 

The estimate \eqref{cXX} requires a bit more work. For the estimate in $\bX$ it suffices to consider a single component $c_k f_j$ where $c_k \in \A_k$ and $f_j$ is localized at frequency $2^j$ in the sense that its 2-equivariant extension to 
$\R^2$ is as described in Section \ref{defnot}; the most important thing is that we have control on the quantity
\[
2^{j} \|f_j\|_{L^2} + 2^{-j} \|f_j\|_{\dHde},  
\]
as highlighted in \eqref{X-relax}).  The term $c_k f_j$ is essentially treated as being at frequency $\max(2^k,2^j)$.  We consider two cases.
\medskip 

a) If $k \leq j$, then we estimate
\[
2^j \|c_k f_j\|_{L^2} \les \|c_k\|_{L^\infty} 2^j \|f_j\|_{L^2} \les 
 \|c_k\|_{\A_k} 2^j \|f_j\|_{L^2}, 
\]
and 
\[
\begin{split}
2^{-j} \|c_k f_j\|_{\dHde} & \les \|c_k \|_{L^\infty} 2^{-j}\|f_j\|_{\dHde}
+ 2^{-2j}\| \partial_r^2 c_k \|_{L^\infty} 2^j \|f_j\|_{L^2} \\
& + 2^{-j}\| \partial_r c_k \|_{L^\infty} \|\partial_r f_j\|_{L^2} 
+ 2^{-j}\| \partial_r c_k \|_{L^\infty} \|\frac{ f_j}r\|_{L^2} \\
& \les \|c_k \|_{\A_k}  (2^{j} \|f_j\|_{L^2} + 2^{-j} \|f_j\|_{\dHde}). 
\end{split}
\]
This we conclude with the bound
\[
2^j \|c_k  f_j\|_{L^2} + 2^{-j} \|c_k  f_j\|_{\dHde} \les 
\|c_k \|_{\A_k}  (2^{j} \|f_j\|_{L^2} + 2^{-j} \|f_j\|_{\dHde}), 
\]
which suffices in this case. 

\medskip

b) If $k \geq j$, on the other hand, then we estimate as follows
\[
2^k \|c_k f_j\|_{L^2} \les 2^k \|c_k \|_{L^2} \|f_j\|_{L^\infty} \les 
 \|c_k\|_{\A_k} \|f_j\|_{\dHe}, 
\]
and 
\[
\begin{split}
2^{-k}\|c_k f_j\|_{\dHde} & \les 2^{-k} \|c_k \|_{L^\infty} \|f_j\|_{\dHde}
+ 2^{-k}\| \partial_r^2 c_k \|_{L^2} \|f_j\|_{L^\infty} \\
& + 2^{-k} \| \partial_r c_k \|_{L^\infty} \|\partial_r f_j\|_{L^2} 
+ 2^{-k} \| \partial_r c_k  \|_{L^\infty} \|\frac{ f_j}r\|_{L^2} \\
& \les \|c_k\|_{\A_k} (2^{j} \|f_j\|_{L^2} + 2^{-j} \|f_j\|_{\dHde}). 
\end{split}
\]
From the two estimates above we obtain
\[
2^k \|c_k  f_j\|_{L^2} + 2^{-k} \|c_k  f_j\|_{\dHde} \les 
\|c_k\|_{\A_k} (2^{j} \|f_j\|_{L^2} + 2^{-j} \|f_j\|_{\dHde}). 
\]
The two bounds above suffice in order to conclude the proof of the first estimate in \eqref{cXX}. 

\medskip

Similarly, for the estimate in $\LX$ it suffices to consider a single component $c_k  g_j$ where $c_k$ is an atom in $\A_k$ and $g_j$ is just as $f_j$ above except that it sits at the $L^2$ level and we have control on the following quantity
\[
2^{j} \|g_j\|_{\dot H^{-1}_e} + 2^{-j} \|g_j\|_{\dHe}.   
\]
Just as above, this component is essentially treated as being at frequency $\max(2^k,2^j)$. We consider the same two cases:

\medskip

a) If $k \leq j$, then we test 
\[
2^j |\la c_k  g_j, \varphi \ra| = |\la  g_j, c_k \varphi \ra| \les 
\|g_j\|_{\dot H_e^{-1}} \| c_k \varphi \|_{\dHe}
\les 2^{j} \|g_j\|_{\dot H_e^{-1}} \| c_k \|_{\A_k} 
\| \varphi\|_{\dHe},
\]
where we have used that
\[
\| c_k \varphi \|_{\dHe} \les \| c_k \|_{\A_k} \| \varphi\|_{\dHe}, 
\]
which is essentially contained in \eqref{cH1e}. We also have
\[
\begin{split}
2^{-j} \|c_k g_j\|_{\dHe} & \les \|c_k \|_{L^\infty} 2^{-j} \|g_j\|_{\dHe}
+ 2^{-j} \| \partial_r c_k \|_{L^\infty} \|g_j\|_{L^2} \\
& \les \|c_k\|_{\A_k} (2^{j} \|g_j\|_{\dot H^{-1}_e} + 2^{-j} \|g_j\|_{\dHe})
\end{split}.
\]
From the above estimates we obtain
\[
2^j \|c_k g_j\|_{\dot H_e^{-1} } + 2^{-j} \|c_k  g_j\|_{\dHe} \les 
\|c_k \|_{C_k} (2^{j} \|g_j\|_{\dot H^{-1}_e} + 2^{-j} \|g_j\|_{\dHe}),
\]
which suffices.

\medskip

b) If $k \geq j$, on the other hand, then we estimate
\[
|\la c_k g_j, \varphi \ra| \les \|r c_k \|_{L^\infty} \|g_j\|_{L^2} 
\| \frac{\varphi}r\|_{L^2} \les 
2^{-k} \|c_k \|_{\A_k}  \|g_j\|_{L^2} \|\varphi\|_{\dHe}.
\]
This implies
\[
\| c_k g_j \|_{\dot H_e^{-1}} \les 2^{-k} \|c_k \|_{\A_k}  \|g_j\|_{L^2}.
\]
Also, just like in the previous case, we have
\[
\begin{split}
2^{-k} \|c_k g_j\|_{\dHe} \les 2^{-k} \|c_k\|_{L^{\infty}} \|g_j\|_{\dHe} +  
2^{-k} \|\partial_r c_k\|_{L^\infty} \|g_j\|_{L^2} \les \|c_k\|_{\A_k} 
(2^{j} \|g_j\|_{\dot H^{-1}_e} + 2^{-j} \|g_j\|_{\dHe}).
\end{split}
\]
Thus we conclude with the estimate
\[
2^k \|c_k g_j\|_{\dot H^1_e} + 2^{-k} \|c_k  g_j\|_{\dHe} \les 
\|c_k \|_{\A_k} (2^{j} \|g_j\|_{\dot H^{-1}_e} + 2^{-j} \|g_j\|_{\dHe}). 
\]
The bounds above for the two cases suffice in order to conclude the proof of the second estimate in \eqref{cXX}. 

\end{proof}

As a consequence of the above results we also claim the following bounds:
\begin{equation} \label{h1X}
\|h_1 f\|_{\bX} + \|h_3 f\|_{\bX} \les \|f\|_{\bX}, \quad 
\|h_1 g\|_{\LX} + \|h_3 g\|_{\LX} \les \|g\|_{\LX}.
\end{equation}
It is straightforward to check that $h_1 \in \bX$, so Lemma \ref{XLXb} justifies the  the estimates for the components involving multiplication with $h_1$. It is a straightforward exercise to check that $h_3-1=\frac{-2}{r^4+1} \in \A$ and then Lemma \ref{cdot} justifies the estimates for the components involving multiplication with $h_3$; here we simply split $h_3-1= \chi_{\leq 0} (h_3-1) + (1-\chi_{\leq 0}) (h_3-1)$ and notice that the first term is an atom in $\A_0$, while the second term belongs to $\bX$.

The next Lemma highlights the context in which the $\tilde X$ structure comes in handy; the integrals considered below appear naturally when solving some ODEs in the following section.  

\begin{l1} \label{intfg}
Assume $f \in \bX$ and $g \in \LX$. Then we have the following representation:
\begin{equation}
\int_r^\infty f(s) g(s) ds = l + c, \quad c=\sum_{k \in Z} c_k,
\end{equation}
where $l \in X$ and $c_k$ are supported in $[0,2^{-k+10})$. Moreover the following holds true:
\begin{equation} \label{lestH1X}
 \| l\|_{\bX}  \les \|f\|_{\bX} \|g\|_{L^2} + \|f\|_{\dHe} \|g\|_{\LX}, 
\end{equation}
and
\begin{equation} \label{l1c}
\begin{split}
 \| l\|_{\dHe} + \sum_k  \| c_k \|_{\A_k}  \les  \|f\|_{\dHe} \|g\|_{L^2}.
\end{split}
\end{equation}
In addition $\partial_r c_k(0)=0$, hence $c$ has the atomic structure in $\A$.
\end{l1}

The integral $\int_r^\infty f(s) g(s) ds$ appears in the process of solving an ODE in the following subsection. Our result essentially states that the integral yields a component $l \in \bX$ and a second component $c$ which, while not in $\bX$ nor in $\dHe$,  retains some of their features, in particular $\partial_r c \in \LX$ ( which is a consequence of \eqref{l1c} and \eqref{cLX}).   

\begin{proof}  We define  $l$ by
\begin{equation}
\begin{split}
l & =\sum_k l_k = \sum_k  l^1_{k} + l_{k}^2  \\
& \ :=\sum_{k} \chi_{ \geq -k}(r) \left(\int_{r}^{\infty} f_{<k}(s) g_{k}(s) ds  + \int_{r}^{\infty} f_{k}(s) g_{\leq k}(s) ds\right),
\end{split}
\end{equation}
and 
\[
c_k = (1-\chi_{ \geq -k}(r))(c^1_{k} + c_{k}^2) ; \quad c^1_{k}   := \int_{r}^{\infty} f_{<k}(s) g_{k}(s) ds, \quad c_{k}^2  := \int_{r}^{\infty} f_{k}(s) g_{\leq k}(s) ds. 
\]
This provides a decomposition as follows:
\[
\int_r^\infty f(s) g(s) ds = l + \sum_k c_k.
\]
It remains to  prove that this decomposition satisfies the claims in the Lemma. 

\medskip

Concerning $l$, in order to conclude the proof of \eqref{lestH1X} it suffices to prove the following
\begin{equation} \label{lkb}
2^k \| l^1_k \|_{L^2} + 2^{-k} \| l^1_k \|_{\dHde} \les \|f \|_{\dHe} \|g_k\|_{L^2}, \quad 
2^k \| l^2_k \|_{L^2} + 2^{-k} \| l^2_k \|_{\dHde} \les \|f_k \|_{\dHe} \|g\|_{L^2},
\end{equation}
 We start with a pointwise bound. For $g_k$ we use Lemma \ref{gg1dec} to write $g_k=(\partial_{s}+\frac{3}{s})g_{k}^{1}$ and then estimate
\begin{equation}
\begin{split}
& |\int_{r}^{\infty} f_{<k}(s) g_{k}(s) ds|
=|\int_{r}^{\infty} f_{<k}(s)\left(\partial_{s}+\frac{3}{s}\right)g_{k}^{1}(s) ds|\\
=& \ |-f_{<k}(r)g_{k}^{1}(r) + \int_{r}^{\infty} (-g_{k}^{1}(s) f_{<k}'(s) + f_{<k}(s) \frac{(3)}{s} g_{k}^{1}(s)) ds|\\
\les & \ |f_{<k}(r)| |g_{k}^{1}(r)| +  \int_{0}^{\infty} \mathbbm{1}_{\{s \geq r\}} (|g_{k}^{1}(s)| |f_{<k}'(s)| + \frac{|f_{<k}(s)|}{s} |g_{k}^{1}(s)|) ds.
\end{split}
\end{equation}
Therefore,
\begin{equation}
\begin{split}
&\|\chi_{\geq -k}(r)\int_{r}^{\infty} f_{<k}(s) g_{k}(s) ds\|_{L^{2}(r dr)}\\
\les & \ \|f_{<k}(r) g_{k}^{1}(r)\|_{L^{2}(r dr)} +  \int_{0}^{\infty} \|\mathbbm{1}_{\{s \geq r\}}\|_{L^{2}(r dr)} (|g_{k}^{1}(s)| |f_{<k}'(s)| + \frac{|f_{<k}(s)|}{s} |g_{k}^{1}(s)|) ds\\
\les & \ \|f_{<k}\|_{L^{\infty}} \|g_{k}^{1}\|_{L^{2}} +  \int_{0}^{\infty} \left(\int_{0}^{s} r dr\right)^{1/2} (|g_{k}^{1}(s)| |f_{<k}'(s)| + \frac{|f_{<k}(s)|}{s} |g_{k}^{1}(s)|) ds\\
\les & \ \|f_{<k}\|_{L^{\infty}} \|g_{k}^{1}\|_{L^{2}} +  \int_{0}^{\infty} s (|g_{k}^{1}(s)| |f_{<k}'(s)| + \frac{|f_{<k}(s)|}{s} |g_{k}^{1}(s)|) ds\\
\les & \  2^{-k} \|f_{<k}\|_{\dot{H}^{1}_{e}} \|g_{k}\|_{L^{2}} + \|g_{k}^{1}\|_{L^{2}(s ds)} \|f_{<k}'(s)\|_{L^{2}(s ds)} + \|\frac{f_{<k}(s)}{s}\|_{L^{2}(s ds)} \|g_{k}^{1}\|_{L^{2}(s ds)} \\
\les & \  2^{-k} \|f_{<k}\|_{\dot{H}^{1}_{e}} \|g_{k}\|_{L^{2}},
\end{split}
\end{equation}
where we have used the bound $\|g_k^1\|_{L^2} \les 2^{-k} \|g_k\|_{L^2}$ from Lemma \ref{gg1dec}. 
This provides the desired estimate for $l_k$ in $L^2$, namely
\begin{equation} \label{l1kL2}
\|l^1_k\|_{L^2} \les 2^{-k} \|f_{<k}\|_{\dot{H}^{1}_{e}} \|g_{k}\|_{L^{2}}. 
\end{equation}

We now turn to estimating $l_k$ in $\dHde$. From the above it follows that 
\begin{equation} \label{l1kr2}
\|\frac{l_k^1}{r^{2}} \|_{L^{2}}
\les  2^{2k} \|l_k^1\|_{L^{2}} \les   2^{k} \|f_{<k}\|_{\dot{H}^{1}_{e}} \|g_{k}\|_{L^{2}}. 
\end{equation}
Next, we begin the estimate for the $\partial_r^2 l_k^1$ as follows:
\begin{equation}
\begin{split}
\|\partial_{r}^{2} l_k^1 \|_{L^{2}}
\les & \  \|\frac{\chi_{>-k}''(r)}{2^{-2k}} \int_{r}^{\infty} f_{<k}(s) g_{k}(s) ds\|_{L^{2}} +  \|\frac{\chi_{>-k}'(r)}{2^{-k}}  \cdot f_{<k}(r) g_{k}(r)\|_{L^{2}}\\
 + & \ \|\chi_{>-k}(r) f_{<k}'(r) g_{k}(r)\|_{L^{2}} + \|\chi_{>-k}(r) f_{<k}(r) g_{k}'(r)\|_{L^{2}}. 
\end{split}
\end{equation}
Below we estimate each of the four terms on the right-hand side in the order in which they appear. From the estimates above for $l_k^1$ in $L^2$, it follows that
\[
\|\frac{\chi_{>-k}''(r)}{2^{-2k}} \int_{r}^{\infty} f_{<k}(s) g_{k}(s) ds\|_{L^{2}} \les 
2^{-2k} 2^{k} \|f_{<k}\|_{\dot{H}^{1}_{e}} \|g_{k}\|_{L^{2}} 
= 2^k  \|f_{<k}\|_{\dot{H}^{1}_{e}} \|g_{k}\|_{L^{2}}. 
\]
Next we have the straightforward estimate
\[
\|\frac{\chi_{>-k}'}{2^{-k}}  \cdot f_{<k} g_{k}\|_{L^{2}} \les 2^{-k} \|f_{<k}\|_{L^\infty}
\|g_k\|_{L^2} \les 2^k  \|f_{<k}\|_{\dot{H}^{1}_{e}} \|g_{k}\|_{L^{2}}. 
\]

A consequence of \eqref{L2-relax} is that $\|g_k'\|_{L^2} \les 2^k \|g_k\|_{L^2}$, thus
\[
\|\chi_{>-k}(r) f_{<k}(r) g_{k}'(r)\|_{L^{2}} \les \|f_{< k}\|_{L^\infty} \|g_k'\|_{L^2}
\les 2^k \|f_{<k}\|_{\dot{H}^{1}_{e}} \|g_{k}\|_{L^{2}}. 
\]
Similarly, we obtain
\[
\begin{split}
\|f'_{<k} \|_{L^\infty} & \les \sum_{j < k} \|f'_{< k}\|_{\dHe} \les \sum_{j < k} \|f_j'\|_{\dHe} \les \sum_{j < k} \|f_j\|_{\dHde} \les \sum_{j < k} 2^j \|f_j\|_{\dHe} \\
& \les  2^k \sum_{j < k} 2^{j-k} \|f_j\|_{\dHe} \les 2^k \left( \sum_{j<k} \|f_j\|_{\dHe}^2 \right)^\frac12 \les 2^k \| f_{< k}\|_{\dHe}. 
\end{split}
\]
Using this we can estimate
\[
\|\chi_{>-k}(r) f_{<k}' g_{k} \|_{L^{2}} \les 
\| f_{<k}' \|_{L^\infty} \|g_{k}\|_{L^{2}} \les 
2^k \|f_{<k}\|_{\dot{H}^{1}_{e}} \|g_{k}\|_{L^{2}}.
\]
As a consequence of all the bounds above we obtain
\begin{equation} \label{ddrl1k}
\|\partial_{r}^{2} l_k^1 \|_{L^{2}} \les 2^k \|f_{<k}\|_{\dot{H}^{1}_{e}} \|g_{k}\|_{L^{2}}.
\end{equation}
We also claim that
\begin{equation} \label{drl1k}
\|\frac{\partial_{r} l_k^1}r \|_{L^{2}} \les 2^k \|f_{<k}\|_{\dot{H}^{1}_{e}} \|g_{k}\|_{L^{2}};
\end{equation}
the argument is similar to the one provided for \eqref{ddrl1k} and the details are left as an exercise.

From \eqref{l1kr2}, \eqref{ddrl1k} and \eqref{drl1k} we obtain that
\[
\|l^1_k\|_{\dHde} \les 2^k \|f_{<k}\|_{\dot{H}^{1}_{e}} \|g_{k}\|_{L^2}.
\]
Together with \eqref{l1kL2} this provides the correct contribution to \eqref{lkb} for $l^1_k$. 

The bound for $l_k^2$ is similar but simpler since the integration by parts is not needed. We start with
\begin{equation}
\begin{split}
 \left|\int_{r}^{\infty} f_{k}(s) g_{\leq k}(s) ds\right|
\les   \int_{0}^{\infty} \mathbbm{1}_{\{s \geq r\}} |f_{k}(s)| |g_{\leq k}(s)| ds .
\end{split}
\end{equation}
Based on this we estimate
\[
\|l_k^2\|_{L^2(rdr)} \les \int_{0}^{\infty} s |f_{k}(s)| |g_{\leq k}(s)| ds \les \|f_k\|_{L^2} \|g_{< k} \|_{L^2} \les 2^{-k} \|f_k\|_{\dHe} \|g_{\leq  k} \|_{L^2},
\]
as needed for the  $L^2$ bound for $l_k^2$ in \eqref{lkb}.  From this it follows that 
\begin{equation} \label{l3kr2}
\|\frac{l_k^2}{r^{2}} \|_{L^{2}}
\les  2^{2k} \|l_k^2\|_{L^{2}} \les   2^{k} \|f_{k}\|_{\dot{H}^{1}_{e}} \|g_{\leq k}\|_{L^{2}}. 
\end{equation}
Finally, the estimate 
\begin{equation} \label{ddrl3k}
\|\partial_{r}^{2} l_k^2 \|_{L^{2}} + \| \frac{\partial_r l^2_k}r \| \les 2^k \|f_{k}\|_{\dot{H}^{1}_{e}} \|g_{\leq k}\|_{L^{2}},
\end{equation}
follows in a similar manner, using the same steps as in the corresponding estimate for $l^1_k$. 

\medskip

We now turn our attention to the $c_k$ terms and prove \eqref{l1c}. 
For $c_k^1$ term  we use the same integration by parts used earlier for $l_k^1$,
\begin{equation}
\begin{split}
c_k^1=\sum_{j < k} \int_{r}^{\infty} f_{j}(s) g_{k}(s) ds
=  \sum_{j<k} \left(-f_{j}(r) g_{k}^{1}(r) - \int_{r}^{\infty}\left(g_{k}^{1}(s) f_{j}'(s) - f_{j}(s) \frac{3}{s} g_{k}^{1}(s)\right) ds\right). 
\end{split}
\end{equation}
From  \eqref{fjpd} we have
\[
|f_{j}(r)| \les (2^j r)^2 \|f_{j}\|_{\dot{H}^{1}_{e}};
\]
which allows us to bound the first term by 
\begin{equation}
\begin{split}
|f_{j}(r) g_{k}^{1}(r)| \les  (2^j r)^2 \|f_{j}\|_{\dot{H}^{1}_{e}} \|g_{n}^{1}\|_{\dot{H}^{1}_{e}}
\les 2^{2(j-k)} \|f_{j}\|_{\dot{H}^{1}_{e}} \|g_k\|_{L^{2}}, \quad r \les  2^{-k}.
\end{split}
\end{equation}
We recall from Section \ref{defnot} that
\begin{equation}
\| \partial_s^2 f_j\|_{L^2}+\|\frac{\partial_s f_{j}(s)}{s}\|_{L^{2}} + \|\frac{f_{j}(s)}{s^{2}}\|_{L^{2}} \les 2^{j}\|f_{j}\|_{\dot{H}^{1}_{e}}.
\end{equation}
Using the two estimates above we obtain the following for $r \les  2^{-k}$,
\begin{equation}
\begin{split}
|c_k^1(r)|
&\les  \sum_{j < k}\left(|f_{j}(r)| |g_{k}^{1}(r)| + \int_{r}^{\infty} |g_{k}^{1}(s)| \left(\frac{|f_{j}'(s)|}{s}+\frac{|f_{j}(s)|}{s^{2}}\right) s ds\right)\\
&\les  \sum_{j < k}\left(|f_{j}(r)| |g_{k}^{1}(r)| + \|g_{k}^{1}(s)\|_{L^{2}} \left(\|\frac{f_{j}'(s)}{s}\|_{L^{2}} + \|\frac{f_{j}(s)}{s^{2}}\|_{L^{2}}\right)\right)\\
&\les  \sum_{j < k} 2^{j-k} \|f_{j}\|_{\dot{H}^{1}_{e}} \|g_{k}\|_{L^{2}}.
\end{split}
\end{equation}
For $c_k^2$ we estimate as follows
\begin{equation}
\begin{split}
|c_k^2(r)| &\les \sum_{j \leq k} \int_{r}^{\infty} |f_{k}(s)| \cdot |g_{j}(s)| ds
\les  \sum_{j \leq k} \int_{r}^{\infty} |f_{k}(s) \frac{g_{j}(s)}{s}| s ds\\
&\les  \sum_{j \leq k} \|f_{k}\|_{L^{2}} \|\frac{g_{j}}{s}\|_{L^{2}}
\les  \sum_{j \leq k} 2^{j-k} \|f_{k}\|_{\dot{H}^{1}_{e}} \|g_{j}\|_{L^{2}},
\end{split}
\end{equation}
where we have used the bounds
\[
\|f_{k}\|_{L^{2}} \les  2^{-k} \|f_{k}\|_{\dot{H}^{1}_{e}}, \quad 
\|\frac{g_{j}(s)}{s}\|_{L^{2}} \les \|g_j\|_{\dHe} \les  2^j \|g_{j}\|_{L^{2}}
\]
from Section \ref{defnot}. 
We combine the estimates above to conclude that
\begin{equation}
\begin{split}
\sum_{k} \|c_{k} \|_{L^\infty(r \leq 2^{-k+10})} & \les   \sum_{k} \|c_{k}^1 \|_{L^\infty} + \|c_{k}^3 \|_{L^\infty}
\les \sum_{k} \sum_{j\leq k} 2^{j-k} (\|f_{k}\|_{\dot{H}^{1}_{e}} \|g_{j}\|_{L^{2}}   +  \|f_{j}\|_{\dot{H}^{1}_{e}} \|g_{k}\|_{L^{2}}) \\
&\les \left( \sum_k \|f_{k}\|_{\dot{H}^{1}_{e}}^2 \right)^\frac12 \left( \sum_k \|g_{k}\|^2_{L^{2}} \right)^\frac12 
\les  \|f\|_{\dot{H}^{1}_{e}} \|g\|_{L^{2}};
\end{split}
\end{equation}
Next, we compute
\[
\partial_r c_k^1 = -f_{<k} g_k = -\sum_{j < k} f_j g_k,
\]
and estimate, for $r \les 2^{-k}$, 
\[
|\partial_r c_k^1(r)| \les \sum_{j < k} |f_j(r)| |g_k(r)| \les \sum_{j < k} (2^j r)^2 \|f_j \|_{\dHe} \|g_k \|_{\dHe}
\les \sum_{j < k} 2^{2(j-k)} \|f_j \|_{\dHe} 2^k \|g_k \|_{L^2}.
\]
where we have used again \eqref{fjpd}. In a similar way
\[
\partial_r c_k^2 = -f_{k} g_{\leq k} = -\sum_{j \leq k} f_k g_j,
\]
and
\[
|\partial_r c_k^2(r)| \les \sum_{j \leq k} |f_k(r)| |g_j(r)| \les \sum_{j \leq k} 2^j (2^j r)^2 \| g_j\|_{L^2} \|f_k\|_{\dHe}
\les 2^k \sum_{j \leq k} 2^{2(j-k)} \| g_j\|_{L^2} \|f_k\|_{\dHe}. 
\]
Based on the two estimates above we conclude with the following estimates
\[
\sum_k 2^{-k} \|\partial_r c_k\|_{L^\infty(r \leq 2^{-k+10})} \les \|f\|_{\dHe} \|g\|_{L^2};
\]
we note that the bound above are stable under multiplication by 
$1-\chi_{ \geq -k}$ (recall that this is how $c_k$ is obtained from $c_k^1+c_k^2$). 

Finally we can differentiate $c_k^i,i=1,2$ twice, argue as above and conclude with
\[
\sum_k 2^{-2k} \|\partial^2_r c_k\|_{L^\infty((r \leq 2^{-k+10})} \les \|f\|_{\dHe} \|g\|_{L^2}.  
\]
We also have the obvious properties that $c_k$ is supported in $[0,2^{-k+10})$
and $\partial_r c_k(0)=0$ is obvious from its definition and the fact that $f_k(0)=f_{< k}(0)=0$ since they are elements in $\dHe$. 

This completes the proof of the claims made for $c$, and in turn it concludes the proof of the Lemma.

\end{proof}

Using these spaces  we  now revisit the result in Lemma \ref{lZL2} and provide a similar result in the context of the refined structures $\tilde X$. 

\begin{l1} \label{lZ}
Consider the vector valued ODE
\begin{equation} \label{mp}
\partial_r Z = N Z + N, \qquad \lim_{r \rightarrow \infty} Z(r) =  0,
\end{equation}
where the coefficients $N$ have the form 
\[
N=\partial_r H + F, \qquad H \in \bX, \quad F_{ij}= \sum_l^{finite} f_{l,ij} g_{l,ij}
\]
with entries satisfying the bounds
\[
\|H\|_{\bX} + \sum_l \|g_{l,ij}\|_{\LX} \leq M, \quad \sum_l \|f_{l,ij}\|_{\bX} \leq C+M,
\]
\[
\|H\|_{\dHe}  + \sum_l \|g_{l,ij}\|_{L^2} \leq m, \quad \sum_l \|f_{l,ij}\|_{\dHe} \leq C+m,
\]
where $C$ is a universal constant. Then, assuming that $m$ is small enough\footnote{The smallness of $m$ is universal and independent on $M$}, the above equation has a unique solution in $B^{\tX}_{2(C+M),2m}$ (see Definition \ref{BMn}).
Furthermore, the map from the coefficients $N$ in the norms above  to $Z \in \tilde X$ is analytic.
\end{l1}

We make several remarks concerning this result:

\begin{enumerate}[label=(\roman*)]
\item The statement above can be easily adapted to systems of type $\partial_r Z =  Z N + N$. 

\item The statement above (including the one made in Remark 1 above) can be generalized to systems of the form $\partial_r Z = N_1 Z + N_2$ (or $\partial_r Z =  Z N_1 + N_2$), where $N_1,N_2$ satisfy similar bounds as $N$. 

\item The universal constant $C$ may be very well thought of as $C=1$; in practice though it depends on $h_1$ and $h_3$ and it is larger than $1$, but it is independent on $M$ used in the statement. 
\end{enumerate}

\begin{proof} The system \eqref{mp}is iterated as follows:
\[
\partial_r Z_1=N, \quad Z_1(\infty)=0
\]
and
\[
\partial_r Z_{n+1}= N Z_n, \quad Z_{n+1}(\infty)=0, 
\]
where the final solution is $Z=\sum_{n=1}^\infty Z_n$, provided we establish the convergence of this series in $\tX$.
We note that the first iteration is given by
\[
Z_1 = H(r) - \int_r^\infty F(s) ds= H + L_1 + C_1,  
\]
where $L_1$ and $C_1$ are given by Lemma \ref{intfg}. In particular
\[
\|L_1\|_{\dHe} \les m, \quad \|L_1\|_{\bX} \les (1+m)(1+M), \quad 
\|C_1\|_{\A} \les m(1+m), 
\]
where the $1+m$ factors can be omitted as $m \ll 1$.
We note that although the estimate for $\|L_1\|_{\bX}$ does not yet contain a factor of $m$ (this is the small parameter), the bound for $\|L_1\|_{\dot{H}^{1}_{e}}$ does so; the point is that all further iterates will add factors of $m$ in the estimates in both $\dot{H}^{1}_{e}$ and $\bX$. 

Next, we solve for the second iterate
\[
\begin{split}
Z_2(r) & =  -\int_r^\infty (F(s)  + \partial_s H(s)) (H(s) + L_1(s) + C_1(s)) ds \\
& = -\int_r^\infty (F(s) + \partial_s H(s)) (H(s) + L_1(s))  ds \\
& - \int_r^\infty F(s)  C_1(s) ds + H(r) C_1(r) + \int_r^\infty H(s) \partial_s C_1(s) ds. 
\end{split}
\]
Based on the bounds for the components of $Z_1$, by the results in Lemma \ref{cdot} and Lemma \ref{XLXb} we obtain the following
representation:
\[
Z_2 = L_2 + C_2,
\]
\[
\|L_2\|_{\dHe} \les m^2, \quad 
\|L_2\|_{\bX} \les m(1+m)(1+M), \quad
\|C_{2}\|_{\A} \les m^2(1+m). 
\]
One notices that at this point even the estimate for $\|L_2\|_{\bX}$ contains a factor of $m$. For general $n$ we obtain a similar representation $Z_n=L_n + C_n$  obeying the following bounds
\[
\|L_n\|_{\dHe} \leq D^n m^n, \quad 
\|L_n\|_{\bX} \leq D^n m^{n-1} (1+M), \quad
\|C_n\|_{\A} \leq D^n m^n. 
\]
Here $D$ is the universal constant that appears in the use of $\les$ above. It is clear that if $m$ is sufficiently small, depending on $D$, but not on $M$, then the above iteration scheme converges.

\end{proof}

At this point we are ready to start proving Proposition \ref{uQpsi} (i).
To keep things streamlined, it is preferable to work with the reduced map $\bar u$ and the corresponding reduced gauge $\bar v, \bar w$ (recall that $u = e^{m \theta R} \bar{u}(r)$ and similarly for $v,w$), since these are functions of $r$ only.  

In this context, it suffices to  show that the  differentiated field $\psi$
satisfies the bound
\begin{equation}\label{u-to=psi}
\| \psi\|_{\LX} \lesssim \| \bar u -\bar Q^2_{\alpha,\lambda}\|_{\bX}.   
\end{equation}
This is done in two stages: first we transfer information from $\bar u$ to the gauge elements $\bar v, \bar w$ and 
then we transfer all the information we have (including the one on the gauge elements) to $\psi$. By scaling and rotation we can  assume the parameter choice $\alpha=0,\lambda=1$, and drop the indices for $Q^2,V^2,W^2$
in the  arguments that follow. 

\subsection{The transition from \texorpdfstring{$u$}{} to \texorpdfstring{$(v,w)$}{}} 
\label{utovw}
Our aim in this first step is to show that we have 
\begin{equation}\label{u-to-vw}
\| \bar v -\bar V\|_{\tX} +\| \bar w -\bar W  \|_{\tX} \lesssim \|\bar u -\bar Q\|_{\bX},
\end{equation}
as well as
\begin{equation}\label{u-to-dpsi2A2}
\| \delta \psi_2\|_{\tX} +\| \delta A_2 \|_{\tX} \lesssim \|\bar u -\bar Q\|_{\bX}.
\end{equation}

We use the equation \eqref{cgeq-m} for the matrix $\calO =  ( \bar v, \bar w,\bar u)$,
namely
\begin{equation} \label{sysM}
\partial_r \calO = M(\bar u) \calO, \qquad \calO(\infty) = I_3, \qquad M(\bar u) = \partial_r \bar u \wedge \bar u.
\end{equation}
If $u = Q$ then $M(\bar u)$ has the form
\begin{equation}
M(\bar Q) =2
\left( 
\begin{array}{lll}
0 & 0 &  -\frac{h_1}r  \\
 0 & 0 & 0 \\
\frac{h_1}r & 0 & 0
\end{array}
\right). \label{E}
\end{equation}

We start with the following identity:
\[
\begin{split}
M(\bar u) - M(\bar Q) & = \partial_r ( \bar u -\bar Q) \wedge (\bar u
-\bar Q) + 2 \partial_r ( \bar u -\bar Q) \wedge \bar Q + \partial_r
\left( \bar Q \wedge (\bar u - \bar Q)\right) \\
& = \partial_r ( \bar u -\bar Q) \wedge (\bar u
-\bar Q) + \partial_r
\left( \bar Q \wedge (\bar u - \bar Q)\right) - 2( \bar u -\bar Q) \wedge \partial_r \bar Q \\
& = \partial_r ( \bar u -\bar Q +2 \bar Q) \wedge (\bar u
-\bar Q) + \partial_r \left( \bar Q \wedge (\bar u - \bar Q)\right). 
\end{split}
\]
We notice that we can write
\begin{equation} \label{Murep}
M(\bar u) - M(\bar Q) =  \partial_r H + F,
\end{equation}
where 
\[
\|H\|_{\bX} \les \|\bar u - \bar Q\|_{\bX}, \quad \| H \|_{\dHe}  \les \|\bar u - \bar Q\|_{\dHe},
\]
and the entries in $F$ are of the form $\sum_{l} f_{l,ij} g_{l,ij}, i,j=1,3$ ($l=2$ in this case) and, in view of \eqref{h1X}, obey the estimates
\[
\sum_l \| f_{l,ij} \|_{\dHe}   \les 1+ \|\bar u - \bar Q\|_{\dHe}, 
\quad \sum_l \| f_{l,ij} \|_{\bX}  \les 1+ \|\bar u - \bar Q\|_{\bX}, 
\]
\[
\sum_l \| g_{l,ij} \|_{L^2}   \les \|\bar u - \bar Q\|_{\dHe}, 
\quad \sum_l \| g_{l,ij} \|_{\bX}  \les  \|\bar u - \bar Q\|_{\bX}. 
\]
Returning to \eqref{sysM}, we start with the solution $\calO_0$  for the case 
$u = Q$, which is given by 
\begin{equation}\label{Fsol}
\calO_0=
\left( 
\begin{array}{ccc}
h_3 & 0  & h_1 \\
0  & 1 & 0 \\
-h_1 & 0 & h_3
\end{array}
\right), \qquad
\calO_0^{-1}= \calO_0^t.
\end{equation} 
Then we express the solution to \eqref{sysM} in the form
\begin{equation} \label{Yexpr}
\calO(r) = \calO_0(r)(I+Y(r)),
\end{equation}
where $Y$  solves  the differential equation
\begin{equation} \label{Ysysup}
\partial_r Y = N Y + N, \qquad Y(\infty) = 0 
 \qquad N=\calO_0^{-1} (M(\bar u) - M(\bar Q))  \calO_0.
\end{equation}
By invoking \eqref{h1X}, it is clear that $N$ has a similar representation as in Lemma \ref{lZ} with $m \approx \|\bar u - \bar Q\|_{\dHe}$ and $M \approx \|\bar u - \bar Q\|_{\bX}$.  Thus  we obtain a solution $Y$ for this system satisfying $Y= L + C$
with $\|L\|_{\dHe} \les \|\bar u - \bar Q\|_{\dHe}, \|L\|_{\bX} \les \|\bar u - \bar Q\|_{\bX}, \|C\|_{\A} \les \|\bar u - \bar Q\|_{\dHe}$. This information is then easily transferred to $ \bar v -\bar V, \bar w -\bar W$, and \eqref{u-to-vw}
follows.

It is useful to note that these are precisely the same type of bounds that have been used in the proof of \eqref{goodal2}, see the analysis of \eqref{sysM2} and the particular setup  in \eqref{Ysysup2}. Here we essentially upgrade that theory from the $\dHe$ framework to the $\bX$ framework. 

We can further improve the bounds for $\bar v_3, \bar w_3$; for instance  
\[
\bar w_3= \bar u_1 \bar v_2 - \bar u_2 \bar v_1= (\bar u_1 -h_1) \bar v_2 + h_1 \bar v_2 - \bar u_2 \bar v_1. 
\]
We have just proved that $\bar v_2=  f +c$ where $\|f\|_{\bX} \les \|\bar u -\bar Q\|_{\bX}$ and $\|f\|_{\dHe} + \|c\| \les \|\bar u -\bar Q\|_{\dHe}$; also we have the trivial estimate $\|\bar u_1-h_1\|_{\bX} \les \|\bar u -\bar Q\|_{\bX}$ and $\|u_1-h_1\|_{\dHe} \les \|\bar u -\bar Q\|_{\dHe}$. Then using \eqref{Xalgref} we obtain
\[
\| (\bar u_1 -h_1) f \|_{\bX} \les \|\bar u -\bar Q\|_{\bX} \|\bar u -\bar Q\|_{\dHe}, 
\]
while using \eqref{cXX} gives
\[
\| (\bar u_1 -h_1) c \|_{\bX} \les \|\bar u -\bar Q\|_{\bX} \|\bar u -\bar Q\|_{\dHe}. 
\]
In a similar manner, using the fact that $h_1 \in X$, we obtain 
\[
\| h_1 \bar v_2 \|_{\bX} \les \|\bar u -\bar Q\|_{\bX}.
\]
The term $\bar u_2 \bar v_1$ is estimated in the same fashion and the final conclusion of this analysis is $\| \bar w_3 \|_{X} \les \|\bar u -\bar Q\|_{X}$. A similar analysis shows that
$\|\bar v_3 +h_1 \|_{\bX} \les \|\bar u -\bar Q\|_{\bX}$; in particular this bound implies that $\|\delta \psi_2\|_{\bX} \les \|\bar u -\bar Q\|_{\bX}$; since $A_2 = 2u_3$,  we also  automatically obtain $\|\delta A_2\|_{\bX} \les  \|\bar u -\bar Q\|_{\bX}$,
concluding the proof of \eqref{u-to-dpsi2A2}.

\subsection{The transition from \texorpdfstring{$(u,v,w)$}{} to  \texorpdfstring{$\psi$}{}}
Here we consider $\psi$, which is  represented as 
\[
\psi = \overline{\W} \cdot \bar v+ i \overline{\W} \cdot \bar w, \qquad \overline{\W}=\partial_r \bar u - \frac{1}{r} \bar u \times R \bar u.
\]
The structure for $(\bar v, \bar w)$, including the representation in the space $\tilde X$, and the multiplicative estimates  for $\LX$, that is \eqref{XLXLX}, \eqref{h1X} and \eqref{cXX}, show that, in order to complete the proof of the result in part i) of Proposition \ref{uQpsi}, precisely the bound \eqref{u-to=psi}, it suffices to establish that  
\begin{equation}
\| \overline{\W} \|_{\LX} \lesssim \|\bar u - \bar Q\|_{\bX}. 
\end{equation}
 Since $\overline{\W}$ vanishes if $u = Q$, we can write
\[
\begin{split}
\overline{\W} \!= & \partial_r (\bar u - \bar Q) - \frac{1}{r} (\bar u -\bar Q) \times \! R (\bar u - \bar Q) 
-  \frac{1}{r} \bar Q \times R (\bar u - \bar Q) 
-  \frac{1}{r} (\bar u -\bar Q) \times R \bar Q
\\
= & L (\bar u - \bar Q) - \frac{1}{r} (\bar u -\bar Q) \times R (\bar u - \bar Q) 
+ \tilde \W.
\end{split}
\]
The first term is in $\LX$ by definition and the second belongs to $\LX$ by using \eqref{drLX}
and \eqref{XLXLX}. It remains to consider 
the last component
\[
\tilde \W = - \frac{h_3}{r} (\bar u - \bar Q) -  \frac{1}{r} \bar Q \times R (\bar u - \bar Q) 
-  \frac{1}{r} (\bar u -\bar Q) \times R \bar Q
\]
This is easily estimated in $\LX$ using \eqref{drLX} and \eqref{h1X}. 
This concludes the proof of \eqref{u-to=psi}, completing the argument for part i) of Proposition \ref{uQpsi}.

\subsection{The transition from  \texorpdfstring{$\psi$}{} to \texorpdfstring{$u$}{}}
In this subsection we establish part ii) of Proposition \ref{uQpsi}, that is  recovering information on the map $u$ from information on the differentiated field $\psi$. Since the value of $\alpha$ plays no particular role in the arguments and $\lambda$ is just a scaling parameter, we do not restrict the generality of the argument by assuming $\alpha=0$ and $\lambda=1$; also, to keep notation compact, we simply write $\delta \psi_2=\delta^{1,0} \psi_2$ and $\delta A_2 = \delta^1 A_2$. 

Our goal here is to establish the following
\begin{equation}\label{psi-to-u}
\| \bar u -\bar Q\|_{\bX}  \lesssim    \| \psi\|_{\LX}.
\end{equation}
This is also done in two steps:  first we transfer the $\LX$ information from $\psi$ to information at the level of $\bX$ for $\delta \psi_2, \delta A_2$, and then we transfer the information to $\bar u$.
 
 For the first step, our main claim here is the following:
\begin{equation} \label{psitodpA}
\| \delta \psi_2 \|_{\bX} + \| \delta A_2 \|_{\bX} \les \|\psi\|_{\LX}.
\end{equation}
This is done as in the Lemma \ref{psito2}, where we transferred $L^2$ information on $\psi$ to $\dot H^1_e$ on $\delta \psi_2, \delta A_2$. We rely on  the system \eqref{sysAp2} which we recall here:
\begin{equation}
 \label{sysAp3}
\left\{ \begin{aligned}
L \delta \psi_2  = & \ 2i h_3 \psi + \delta A_2 \psi - \frac{1}r \delta A_2 (2ih_1 + \delta \psi_2), \\
L \delta A_2  = & \ - 2 h_1 \Re{\psi} + \Im(\psi \overline{\delta \psi}_2) - \frac1r (\delta A_2)^2.
\end{aligned} \right.
\end{equation}
The strategy is to produce an iteration scheme in $\bX$ for the above system with an appropriate initialization. The easiest initialization that we can impose is in the form of an initial data, just as we did in Lemma \ref{psito2}; to be more precise, we let $\lambda_0,\alpha_0$ be the parameters obtained by imposing the pointwise orthogonality condition 
$\delta^{\lambda_0,\alpha_0} \psi_2(r_0)=0, \delta^{\lambda_0} A_2(r_0)=0$ for some $r_0 \approx 1$; recall that our true parameters were set to $\lambda=1,\alpha=0$. The proof of Proposition \ref{p:mod} establishes that 
\[
|\ln \lambda- \ln \lambda_0|+ |\alpha-\alpha_0| \les \|\psi\|_{L^2},
\]
which in turn gives 
\[
\begin{split}
\|\delta \psi_2 - \delta^{\lambda_0,\alpha_0} \psi_2 \|_{\bX} 
& \approx \| h_1 - e^{2i\alpha_0} h_1^{\lambda_0}\|_{\bX} \\
& \les |\ln \lambda- \ln \lambda_0|+ |\alpha-\alpha_0|  \les \|\psi\|_{L^2},
\end{split}
\]
together with
\[
\begin{split}
\|\delta A_2 - \delta^{\lambda_0} A_2 \|_{\bX} 
 \approx \| h_3 - h_3^{\lambda_0}\|_{\bX} 
 \les |\ln \lambda- \ln \lambda_0| \les \|\psi\|_{L^2}.
\end{split}
\]
Thus it suffices to establish $X$ bounds for $\delta^{\lambda_0,\alpha_0} \psi_2$ and $\delta^{\lambda_0} A_2$. In Lemma \ref{psito2} we have established that these were the unique solutions to \eqref{sysAp3} with the above initial conditions, which were obtained via an iteration scheme. Our strategy here is simple: assuming the additional structure 
$\psi \in \LX$, we show that the same iteration scheme improves the structure of the solution to the desired $X$ norms
\begin{equation}
\| \delta^{\lambda_0,\alpha_0} \psi_2 \|_{\bX} 
+ \| \delta^{\lambda_0} A_2 \|_{\bX}
\les \|\psi\|_{\LX}. 
\end{equation}

The improvement in the iteration scheme relies on two basic ingredients:

\begin{enumerate}
    \item 
the bilinear estimates in Lemma \ref{XLXb}, namely the estimates in \eqref{h1X} and \eqref{drLX}; it is important for the convergence of the iteration scheme  that the bilinear estimates in \eqref{XLXb} allow us to gain smallness from the quantities $\|\delta^{\lambda_0,\alpha_0}  \psi_2\|_{\dHe}$ and $\|\delta^{\lambda_0}  A_2\|_{\dHe}$,  as 
 we do not have any smallness in $\bX$.

\item the result in Lemma \ref{LLwf} (b), which allows us to solve the linear $L$ equation in $\bX$. 
\end{enumerate}

This completes the proof of \eqref{psitodpA}. 

\vspace{.2in}

The next step is to transfer the $LX$ bounds from $\psi$ together with the $X$ bounds on $\delta \psi_2$ and $\delta A_2$ to $\bar u- \bar Q$ and conclude the proof of \eqref{psi-to-u}. The strategy here follows the same steps as in the proof of  \eqref{goodal2}. We recall the system \eqref{sysM2}  for $\calO=(\bar v, \bar w, \bar u)$ 
\begin{equation} \label{sysM3}
\partial_r \calO = \calO R(\psi), \quad  \calO(\infty) = I_3. 
\end{equation}
with 
\[
R(\psi)=\left( 
\begin{array}{lll}
0 & 0 &   \Re \psi_1 \\
0 & 0 &  \Im \psi_1 \\
- \Re \psi_1 &  - \Im \psi_1 &  0
\end{array},
\right)
\]
as well as the fact that if $\psi = 0$ then $\psi_2 = 2 i h_1$, which yields $\psi_1= - 2  \frac{h_1}{r}$, hence
\begin{equation}
R(0) =
2\frac{h_1}r \left( 
\begin{array}{lll}
0 & 0 &  - 1 \\
 0 & 0 & 0 \\
1 & 0 & 0
\end{array}
\right). \label{EE3}
\end{equation}
The solution is given by (see the generalization of \eqref{vwq})
\[
\calO_0 = \!\left(\!\!\! \begin{array}{ccc}
\! h_3 
 & \! 0
& h_1  \cr
 \! 0 \!
& \!  1
& 0  \cr
\! -h_1 & \!
0 & h_3  \end{array}\!\!\! \right).\!
\]
We note that $\calO_0^{-1}= \calO_0^t$. We will prove that 
\begin{equation}\label{Rest3}
R(\psi) - R(0)= (\sum_l f_{l,ij} g_{l,ij})_{i,j=1,3} + \partial_r H,
\end{equation}
where
\[
\|H\|_{\dHe} + \sum_l \|g_{l,ij}\|_{L^2} \les \|\psi\|_{L^2}, 
\quad \|H\|_{\bX} + \sum_l \|g_{l,ij}\|_{\LX} \les \|\psi\|_{\LX}
\]
and
\[
\sum_l \|f_{l,ij}\|_{\dHe} \les 1 + \|\psi\|_{L^2}, 
\quad \sum_l \|f_{l,ij}\|_{\bX} \les 1 + \|\psi\|_{\LX}. 
\]
Suppose this is done. Then we write the solution to \eqref{sysM2} in the form
\begin{equation} \label{Yexpr3}
\calO(r) = (I+Y(r))\calO_0(r)
\end{equation}
where $Y$  solves  the differential equation
\begin{equation} \label{Ysysup2+}
\partial_r Y =  Y N + N, \qquad Y(\infty) = 0 
 \qquad N=\calO_0 (R(\psi) - R(0))  \calO_0^{-1}.
\end{equation}
The matrix $N$ inherits a similar representation to the one in \eqref{Rest3} and with similar bounds.
Thus we can apply Lemma \ref{lZ} to solve this system and conclude that
\[
Y=L+C
\]
with
\begin{equation} \label{YLCb}
\|L\|_{\bX} \les \|\psi\|_{\bX}, \quad \|L\|_{\dHe} \les \|\psi\|_{L^2},  \quad \|C\| \les \|\psi\|_{L^2}. 
\end{equation}
This provides a similar representation for all columns of $Y$, in particular for $\bar u- \bar Q$. 
To finish the proof of our claim \eqref{Rest3} we need to upgrade the information on $\bar u- \bar Q$ to $\bX$. 
We first remark that the last row of
$\calO$ is a-priori known, namely $(\bar v_3, \bar w_3, \bar u_3) = \frac12 (-\Im \psi_2,\Re
\psi_2,A_2)$; this shows that from \eqref{psitodpA} we obtain 
\[
\| \bar v_3 +  h_1 \|_{\bX} + \| \bar w_3 \|_{\bX} + \| \bar u_3 - h_3\|_{\bX}  \lesssim \|\psi\|_{\LX}. 
\]
To transfer this information to $\bar u_1$ and $\bar u_2$ we use again the orthogonality
of $\calO$.  For $\bar u_1$ we have 
\[
\bar u_1 - h_1= \bar v_2 \bar w_3 - \bar v_3 \bar w_2= \bar v_2 \bar w_3  - (\bar v_3+h_1) \bar w_2 + h_1  (\bar w_2-1), 
\]
It is important that in each product above we have one term which is apriori bounded in $\bX$ (in order of appearance, $\bar w_3, \bar v_3 +h_1, h_1$), while the other component of each product is one of the entries in $Y=L+C$ and inherits the bounds in  \eqref{YLCb}. Thus we can we use \eqref{cXX} and \eqref{XLXb} and the bounds in \eqref{YLCb} (or better for entries in those matrices) to conclude that 
\[
\|\bar u_1 -h_1\|_{\bX} \les \|\psi\|_{\LX}. 
\]
A similar argument shows that $\| \bar u_2\|_{\bX} \les \|\psi\|_{\LX}$. This concludes the proof of \eqref{psi-to-u}.

It remains to prove the bounds claimed for the elements in the decomposition \eqref{Rest3}. From \eqref{dpsi1} we have:
\[
\psi_1 + 2 \frac{h_1}{r} = - 2i h_3 \partial_r \delta \psi_2 + B,
\]
where
\[
B= - i \delta A_2 \cdot \partial_r \delta \psi_2 + \delta A_2 \cdot 2  \partial_r h_1
+  |\psi_2|^2  \psi + \frac1{4r} (  i  \psi_2 |\psi_2|^2 + 8 (h_1)^3).
\]
Using \eqref{psitodpA} and the apriori bounds in the energy spaces from Lemma \ref{psito2}, it is a straightforward exercise to check that
\[
B=\sum_{j} f_j g_j, \quad \|f_j \|_{\bX} \les 1+ \|\psi\|_{\LX}, \|f_j \|_{\dHe} \les 1+ \|\psi\|_{L^2}, 
\|g_j \|_{\LX} \les \|\psi\|_{\LX}, \|g_j \|_{L^2} \les \|\psi\|_{L^2}. 
\]
Finally, we write
\[
h_3 \partial_r \delta \psi_2 = \partial_r  (h_3 \delta \psi_2) + \partial_r h_3 \cdot  \delta \psi_2,
\]
and note that $\| h_3 \delta \psi_2 \|_{\bX} \les \| \delta \psi_2 \|_{\bX} \les \|\psi\|_{\LX}$, while $\| \partial_r h_3  \|_{\LX} \les 1$. This completes the proof of the claims regarding the decomposition in \eqref{Rest3}, and in turn of \eqref{psi-to-u}. 

This completes the proof of Proposition \ref{uQpsi}.

\section{The linear \texorpdfstring{$\tilde H$}{} Schr\"odinger equation}
\label{s:linear}

The main goal of this section is the study of estimates for the linear $\tilde H$  evolution,
\begin{equation}\label{wlin-eq}\begin{cases} 
i \partial_{t}\psi  - \tilde H_{\lambda} \psi =f\\
\psi(0)=\psi_{0}\end{cases}\end{equation}
where we recall that
\[
\tilde H_{\lambda} = - \partial_{r}^{2}  - \frac{1}{r}\partial_{r} +\left(\frac{1}{r^{2}}+\frac{8}{r^{2}(1+\lambda(t)^{4}r^{4})}\right).
\]
One of the important features to highlight is that $\lambda$ is allowed to depend on $t$. These linear estimates will be crucial in the study of the nonlinear equation \eqref{psieq2} and the evolution of the modulation parameters $\lambda$ and $\alpha$. 

Our analysis will provide estimates in two categories of spaces. The first one includes the classical energy, Strichartz and local energy decay norms. The second category is a refinement of  the first, considering each of the above elements at the dyadic level in the context  of our time dependent Littlewood-Paley decomposition.

The first spaces are straightforward; with the local energy norm and its dual defined by 
\[
 \| \psi\|_{LE} = \| \frac{\psi}{r} \|_{L^2_{t}L^{2}_{r}(r dr)}, \quad 
  \| \psi\|_{LE^*} = \| r\psi \|_{L^2_{t}L^{2}_{r}(r dr)},
\]
we define the space $S$ for solutions to
\eqref{wlin-eq} and the dual type space $N$ (precisely, $S=N^*$)
 for the inhomogeneous term in \eqref{wlin-eq} as follows:
\[
S = L^\infty_t L^2_r \cap L^4_{tr} \cap LE, \qquad
N = L^1_t L^2_r +  L^\frac43_{tr} +  LE^*.
\]  
For instance the $S$ and $N$ structures are robust enough to close the main result in Section \ref{s:mod},
where we get a first insight into the dispersive behaviour of the problem and the crucial control on $\|\frac{\lambda'}{\lambda^2}\|_{L^2_t}$; the relevance of this quantity will become apparent in Proposition \ref{p:mainS} below. \\
\\
On the other hand, for many of our estimates we need to be more precise and work with a
dyadic Littlewood-Paley decomposition in the $\tilde H_{\lambda(t)}$-frequency, 
$ P_k^{\lambda(t)} \psi$; in fact even our strategy to derive the linear estimates involving the spaces $S$ and $N$ uses these finer structures. To measure frequency $2^k$ waves we define a local energy
space $LE_k$,
\[
\|\psi\|_{LE_{k}} = 2^{k}\|\psi\|_{L^{2}_{t}L^{2}_{r}(A_{<-k})} + \sup_{m > -k} 2^{\frac{k-m}{2}} \|\psi\|_{L^{2}_{t}L^{2}_{r}(A_{m})},
\]
as well as the dual space $LE_k^*$. Here we note that these norms vary slowly with $k$,
\begin{equation}\label{lekpert}
 \|\psi\|_{LE_{k}} \approx \|\psi\|_{LE_{k'}} , \quad 
 \mbox{if} \ |k'-k| \leq 10.
\end{equation}
Verification is straightforward and left as an exercise.

We aggregate these norms in an $\ell^2$-Besov fashion, and set
\[
\| \psi\|_{\ell^2 LE}^2 = \sum_k
\| P_k \psi \|_{L^\infty L^2 \cap LE_k}^2,
\qquad \| f\|_{\ell^2 LE^*}^2 = \sum_k
\| P_k f \|_{L^1 L^2 + LE_k^*}^2.
\]

Following \cite{BeTa-1}, we also define an adapted $L^4_k$ norm, which is allowed due to the radial symmetry:
\[ 
\| \psi\|_{L^4_k} = \sup_{m} \max\{2^{-\frac{m+k}2},2^{\frac{m+k}8}\}
\|\psi\|_{L^{4}_{t}L^{4}_{r}(A_{m})}.
\]
The dual norm is denoted by $L^\frac43_k$.
The frequency adapted versions of the $S$ and $N$
norms are
\[
S_k = L^\infty_t L^2_r \cap L^4_k \cap LE_k, \qquad
N_k = L^1_t L^2_r +  L^\frac43_{k} +  LE_k^*, \qquad S_k = N_k^*.
\]  
Square summing these norms we obtain the spaces $l^2 S$ and $l^2 N$ with norms
\begin{equation}\label{l2S}
\| \psi\|_{l^2 S}^2 = \sum_{k\in \Z} \|P_k^\lambda \psi\|_{ S_k}^2,
\qquad \|f\|_{l^2 N}^2 = \sum_{k\in \Z} \|P_k^\lambda f\|_{N_k}^2.
\end{equation}
Given the nice bound \eqref{ker1} on the kernel of the projectors 
$P_k$, it is easy to see that these are dual spaces, thus justifying 
our notation. We recall from Section \ref{L-PP} that $P_k^\lambda$ are the standard projectors in the $\tilde H_\lambda$ calculus, or, equivalently, in the $\FtH_\lambda$ frame; also, $\tilde P_k^\lambda$ are similar projectors with the additional property that  $ \tilde P_k^\lambda P_k^\lambda = P_k^\lambda$.

We will establish in the following subsection the following relation between the two structure introduced above
\begin{equation} \label{Sl2S}
    l^2 S \subset S, \qquad N \subset l^2 N.
\end{equation}

The main result of this section is the following linear estimate. 

\begin{t1} \label{mainTl2}
Assume the time dependent function $\lambda$ satisfies
\begin{equation}\label{good-lambda}
\| \lambda' \lambda^{-2}\|_{L^2(0,T)} \lesssim 1.    
\end{equation}
Then the evolution \eqref{wlin-eq}
is well-posed in $L^2$, and the following estimate holds in $[0,T]$:
\begin{equation}
 \| \psi\|_{l^2 S} \lesssim \|\psi_0\|_{L^2} + \|f\|_{l^2 N}.
\label{mainS} \end{equation}
\end{t1}

In particular, we obtain the following estimates in the $S$ and $N$ spaces.
\begin{c1}\label{c:le}
If $\psi$ is as in the above Theorem, then the following estimate holds in $[0,T]$:
\[\| \psi\|_{S} \lesssim \|\psi_0\|_{L^2} + \|f\|_{ N}\]
\end{c1}

Theorem~\ref{mainTl2} provides a bound corresponding to $\ell^2$ dyadic summation, which is natural when working in the finite energy setting, i.e. $\psi(0) \in L^2$. However, later in the article we also investigate the more restrictive case of data with $\ell^1$ dyadic summation, namely $\psi(0) \in \LX$. Our result is as follows:

\begin{t1}\label{mainTl1}
Assume that $\lambda$ satisfies \eqref{good-lambda}. 
Then the following bound holds in $[0,T]$ for solutions to the 
equation \eqref{wlin-eq}:
\begin{equation}
\label{lejfinal++}\begin{split}  \|\psi\|_{\ell^1 S} &\lesssim 
 \|\psi(0)\|_{\LX} + 
  \|f\|_{\ell^1 N}.\end{split}\end{equation}
\end{t1}
Here the $\ell^1 S$ and $\ell^1 N$ norms are defined as in \eqref{l2S} but with $\ell^1$ 
summation instead.

We remark that one may also write a frequency envelope version of the above bound. This would assert that if $\psi(0)$ and $f$ can be placed under a slowly varying frequency envelope $c_k$ in either $\ell^1$ or  $\ell^2$, then the $S_k$ norm of $P_k \psi$ can be placed under a similar frequency envelope.

The proof of the two main theorems above 
uses Littlewood-Paley decomposition associated to the time dependent operator $\tH_\lambda$
and begins in the next subsection with some 
elliptic estimates for frequency localized functions. The main building block of our analysis is the proof of local energy decay 
for frequency localized functions. This is carried out in Section~\ref{s:led} under  smallness assumption for our control norm $\| \lambda'/\lambda^2\|_{L^2}$, and then expanded 
in the next subsection to solutions with a control norm which is large but finite. 
The Strichartz component of the $S$ and $N$ norms  is added in Section~\ref{ss-str}. The last step of the analysis is to assemble the dyadic bounds into the full bounds in the theorems, which is achieved by perturbatively estimating the frequency localization errors, which are related to the transference operator studied in Section~\ref{spectral}.

In order to keep notation compact, in what follows, whenever a space-time is involved, the time interval is assumed to be restricted to  the interval $[0,T]$ in the condition \eqref{good-lambda}; the estimates are uniform with respect to $T$.

\subsection{Properties of function spaces} 
In this section we establish the basic relation $l^2 S \subset S$ and its dual $N \subset l^2 N$, along with some other properties of the spaces defined at the beginning of this section. 

We begin with some simple estimates which are helpful for later arguments. For convenience we recall the definition of the function $\omega_{j,\lambda}$ from \eqref{omdef}:
\[
\omega_{j,\lambda}(r) = \begin{cases} \text{min}\{1,r^{3}2^{3j}\}, \quad 2^{j} \geq \lambda \\
\text{min}\{1,2^{j} r \dfrac{r^{2}\lambda^{2}}{1+r^{2}\lambda^{2}}\}, \quad 2^{j} \leq \lambda.\end{cases}
\]
Our first result is the following. 

\begin{l1} For any $k,j \in \Z$ and $\lambda:I \rightarrow (0,+\infty)$, the following holds true:
\begin{equation}\label{point-le}
\|\omega_{k,\lambda(t)}(r)^{-1} P_{k}^{\lambda(t)} \psi(t,r)\|_{L^{2}_{t}L^{\infty}_{r}(A_{j})} 
\lesssim  \|P_{k}^{\lambda(t)} \psi\|_{LE_{k}},
\end{equation}
\begin{equation}\label{point-e}\|(1+2^{k}r)^{1/2} \omega_{k,\lambda(t)}^{-1} P_{k}^{\lambda(t)}\psi\|_{L^{\infty}_{t}L^{\infty}_{r}} \lesssim 2^{k} \|P_{k}^{\lambda(t)}\psi\|_{L^{\infty}_{t}L^{2}_{r}},
\end{equation}
with universal implicit constants, independent of $k,j$
and the function $\lambda$. 
\end{l1}

\begin{proof}
Representing 
\begin{equation*}
P_{k}^{\lambda(t)}\psi(r) =  \widetilde{P_{k}^{\lambda(t)}}(P_{k}^{\lambda(t)}(\psi)) = \int_{0}^{\infty} \widetilde{K_{k}^{\lambda(t)}}(r,s) P_{k}^{\lambda(t)}\psi(s) s ds ,
\end{equation*}
we use the spectral projector kernel bound \eqref{ker1}
with $N = 1$ to estimate at fixed time
\begin{equation*}
\begin{split}
|\omega_{k,\lambda(t)}(r)^{-1}P_{k}^{\lambda(t)}\psi(r)|&\lesssim \int_{0}^{\infty}  \frac{ 2^{2k} \omega_{k,\lambda(t)}(s) |P_{k}^{\lambda(t)}(\psi)(s)|}{(1+2^{k}(s+r))(1+2^{k}|r-s|)}s ds\\
&\lesssim  2^{2k}\int_{0}^{\infty} \frac{ |P_{k}^{\lambda(t)}\psi(s)|}{(1+2^{k}(s+r))(1+2^{k}|r-s|)} s ds.
\end{split}
\end{equation*}
We bound the right hand side using the dyadic spatial decomposition in the $LE_k$ norm and applying the Cauchy-Schwarz inequality.  
For $r \in A_m$ this gives
\begin{equation*}
|\omega_{k,\lambda}(r)^{-1}P_{k}^{\lambda(t)}\psi(r)|
\lesssim 2^{k}\|P_{k}^{\lambda(t)}\psi\|_{L^{2}_{r}(A_{<-k})}  + 2^{2k}\sum_{j=-k+1}^{\infty} \|P_{k}^{\lambda(t)}\psi\|_{L^{2}(A_{j})} c_{kjm},
\end{equation*}
where 
\[
c_{kjm}^2= \sup_{r \in A_m} \int_{A_j} \frac{2^{j} ds}{(1+2^{k+j})^{2} (1+2^{k}|r-s|)^{2}} \approx  \begin{cases} 
2^{-2k} 2^{-j-k}, \quad {j-3} \leq m \leq {j+3} \\ 
2^{-2k} 2^{-2(j+k)} , \quad \text{otherwise}.
\end{cases}
\]
Taking the $L^2_t L^\infty_r$ norm we arrive at
\begin{equation*}
\begin{aligned}
\| \omega_{k,\lambda(t)}^{-1}P_{k}^{\lambda(t)}\psi\|_{L^2_t L^\infty_r(A_m)}
\lesssim & \ 2^{k}\|P_{k}^{\lambda(t)}\psi\|_{L^2_t L^{2}_{r}(A_{<-k})}  + 2^{2k}\sum_{j=-k+1}^{\infty} \|P_{k}^{\lambda(t)}\psi\|_{L^2_t L^2_r(A_{j})} c_{kjm}
\\
\lesssim &\ \| P_{k}^{\lambda(t)}\psi\|_{LE_k}  \sum_{j=-k}^{\infty} 2^{2k} c_{kjm} 2^{-\frac{j+k}2}
\\
\lesssim &\ \| P_{k}^{\lambda(t)}\psi\|_{LE_k}.
\end{aligned}
\end{equation*}

Similarly, we have the fixed time bound
\[
|(1+2^{k} r)^{1/2} \omega_{k,\lambda(t)}^{-1} P_{k}^{\lambda(t)}\psi(r)| \lesssim (1+2^{k}r)^{1/2} \int_{0}^{\infty}  \frac{ 2^{2k} \omega_{k,\lambda(t)}(s) |P_{k}^{\lambda(t)}\psi(s)|}{(1+2^{k}(s+r))(1+2^{k}|r-s|)} s ds,
\]
which gives, by Cauchy-Schwarz,
\[
|(1+2^{k} r)^{1/2} \omega_{k,\lambda(t)}^{-1} P_{k}^{\lambda(t)}\psi| \lesssim \| P_{k}^{\lambda(t)} \psi\|_{L^{2}} c_k(r)^\frac12,
\]
where
\begin{equation*}
\begin{split}
 c_k(r)  = \int_{0}^{\infty} \frac{ (1+2^k r)  2^{4k}}{(1+2^{k}(s+r))^{2}(1+2^{k}|r-s|)^{2}} s ds 
 \lesssim \int_{0}^{\infty} \frac{ 2^{3k} \, ds }{(1+|2^{k}(r-s)|)^{2}} \lesssim 2^{2k}.
\end{split}
\end{equation*}
Hence \eqref{point-e} follows.

\end{proof}

We are now ready to proceed with the arguments for the inclusions $l^2 S \subset S$ and $N \subset l^2 N$. First we prove that the $LE$ norm is controlled by the square sum of the $LE_k$ norms and the corresponding dual estimate. 

\begin{l1}
The following estimates hold true:
\begin{equation} \label{lek-to-le}
    \|\psi\|_{LE}^{2} \lesssim \sum_{k \in \mathbb{Z}} \|P_{k}^{\lambda(t)}\psi\|_{LE_{k}}^{2},   
\end{equation}

\begin{equation}\label{le-to-lek}
\sum_{k} \|P_{k}^{\lambda(t)}f\|_{LE_{k}^{*}}^{2} \lesssim \|f\|_{LE^{*}}^{2}.
\end{equation}
\end{l1}
\begin{proof}
In each spatial dyadic region we decompose in frequency and estimate as follows:
\begin{equation*}
\|\frac{\psi}{r}\|_{L^{2}_{r}(A_{j})} 
\lesssim \frac{1}{2^{j}}\|\psi\|_{L^{2}_{r}(A_{j})} 
\lesssim \frac{1}{2^{j}} \sum_{k \in \mathbb{Z}}\|P_{k}^{\lambda(t)}\psi\|_{L^{2}_{r}(A_{j})}.
\end{equation*}
To estimate the terms in the last sum we consider two cases.
If $k+j \geq 0$ we estimate directly
\begin{equation}\label{high-k}
  2^{-j} \|P_{k}^{\lambda(t)}\psi\|_{L^{2}_{r}(A_{j})} \lesssim  2^{-\frac{j+k}2} \|P_{k}^{\lambda(t)}\psi\|_{LE_{k}}.
\end{equation}
If $k+j \leq 0$, we use the $L^2 L^\infty$ bound 
from \eqref{point-le}, and note that 
in $A_j$ we have $\omega_{k,\lambda(t)} \lesssim 2^{j+k}$;
thus we obtain
\begin{equation}\label{low-k}
 2^{-j} \|P_{k}^{\lambda(t)}\psi\|_{L^{2}_{r}(A_{j})} \lesssim   
     \|P_{k}^{\lambda(t)}\psi\|_{L^2 L^\infty(A_j)}
\lesssim 2^{j+k}   \|P_{k}^{\lambda(t)}\psi\|_{LE_{k}}.
\end{equation}
Combining the bounds \eqref{high-k} and \eqref{low-k}
we arrive at 
\begin{equation}
\|\frac{\psi}{r}\|_{L^{2}_{t}L^{2}_{r}(A_{j})} \lesssim \sum_{k=-\infty}^{\infty} 2^{\frac{-|j+k|}{2}} \|P_{k}^{\lambda(t)}\psi\|_{LE_{k}}.
\end{equation}
Finally, using a discrete convolution estimate, we have 
\[
\|\psi\|_{LE}^{2} \leq \sum_{j \in \mathbb{Z}} \|\frac{\psi}{r}\|_{L^{2}_{t}L^{2}_{r}(A_{j})}^{2} \lesssim \sum_{k \in \mathbb{Z}} \|P_{k}^{\lambda(t)}\psi\|_{LE_{k}}^{2},
\]
as needed. This finishes the proof of \eqref{le-to-lek}; the estimate  \eqref{lek-to-le} follows by duality.
\end{proof}

Next we prove that the $L^4_{tr}$ norm is controlled by the square sum of the $L^4_k$ norms and the corresponding dual estimate. 

\begin{l1} We have the following estimates: 
\begin{equation}\label{l4froml4j}
\|\psi\|_{L^{4}_{t,r}} \lesssim \left(\sum_{m} \|P_{m}^{\lambda(t)}\psi\|_{L^{4}_{m}}^{2}\right)^{1/2},
\end{equation}
\begin{equation}\label{l43froml43j}
\left(\sum_{m} \|P_{m}^{\lambda(t)}f\|_{L_{m}^{4/3}}^{2}\right)^{1/2} \lesssim \|f\|_{L^{4/3}_{t,r}}.
\end{equation}
\end{l1}
\begin{proof} We have
\begin{equation}
\begin{split}
\|\psi\|_{L^{4}_{t,r}(A_{m})} &\leq \sum_{j} \|P_{j}^{\lambda(t)}\psi\|_{L^{4}_{t,r}(A_{m})}\\
&\leq \sum_{j > -m} \|P_{j}^{\lambda(t)}\psi\|_{L^{4}_{t,r}(A_{m})} 
+ \sum_{j \leq -m} \|P_{j}^{\lambda(t)}\psi\|_{L^{4}_{t,r}(A_{m})}\\
&\leq \sum_{j > -m} 2^{\frac{-(m+j)}{8}}2^{\frac{(m+j)}{8}} \|P_{j}^{\lambda(t)}\psi\|_{L^{4}_{t,r}(A_{m})} + 
 \sum_{j \leq -m} 2^{\frac{-(m+j)}{2}}2^{\frac{(m+j)}{2}} \|P_{j}^{\lambda(t)}\psi\|_{L^{4}_{t,r}(A_{m})}\\
&\leq \sum_{j \in \mathbb{Z}} 2^{\frac{-|m+j|}{8}} \|P_{j}^{\lambda(t)}\psi\|_{L^{4}_{j}}.
\end{split}
\end{equation}
Then, we use Young's inequality for discrete convolutions to obtain
\[
\sum_{m \in \mathbb{Z}} \|\psi\|_{L^{4}_{t}L^{4}_{r}(A_{m})}^{2} \lesssim \sum_{m} |\sum_{j \in \mathbb{Z}} 2^{\frac{-|m+j|}{8}} \|P_{j}^{\lambda(t)}\psi\|_{L^{4}_{j}}|^{2} \lesssim  \sum_{n \in \mathbb{Z}} \|P_{n}^{\lambda(t)}\psi\|_{L^{4}_{n}}^{2}.
\]
Based on this we obtain
\[
\|\psi\|_{L^{4}_{t,r}}^{2} \lesssim  
\left(\sum_{m} \|\psi\|_{L^{4}_{t,r}(A_{m})}^{4}\right)^{1/2} 
\lesssim  \sum_{m}(\|\psi\|_{L^{4}_{t,r}(A_{m})}^{4})^{1/2} 
\lesssim \sum_{n \in \mathbb{Z}} \|P_{n}^{\lambda(t)}\psi\|_{L^{4}_{n}}^{2},
\]
which implies \eqref{l4froml4j}. We also have the dual estimate \eqref{l43froml43j}.

\end{proof}

It is clear that \eqref{Sl2S} follows from the previous two Lemmas.

\subsection{Local energy decay for frequency localized solutions.} \label{s:led}

A first step in the analysis of the linear Schr\"odinger equation \eqref{wlin-eq} is to derive local energy decay estimates for frequency localized functions. We recall the equation \eqref{wlin-eq}
here for convenience: 
\[
i \partial_{t}\psi  - \tilde H_{\lambda(t)} \psi =f, 
\qquad
\psi(0)=\psi_{0}.
\]

Throughout the rest of this section we say 
that a time dependent function $\psi$ is 
localized at a dyadic frequency $2^j$ if 
$\FtH_{\lambda(t)} \psi(t)$ is supported in the region $|\xi| \approx 2^j$. Notably,
the scale function $\lambda$ is allowed to depend on time. 

In the following Lemma we establish the basic linear estimates in the spaces $LE_j$ for solutions to \eqref{wlin-eq} which are frequency localized. 

\begin{l1} \label{l:LEj}
Assume that $\psi$ is a solution to \eqref{wlin-eq} which is localized at frequency $2^j$. Then the following holds true:
\begin{equation}\label{bllej+}
\|\psi\|_{LE_{j} \cap L^{\infty}_{t}L^{2}_{r}} + 2^{-j}\|\partial_{r}\psi\|_{LE_{j}} \lesssim \|\psi(0)\|_{L^2} +  \|f\|_{L^{1}_{t}L^{2}_{r}+LE_{j}^{*}}.
\end{equation}
\end{l1}
We remark that the frequency localization of $u$ does not guarantee
a similar frequency localization of $f$. This is due to the time dependence of $\lambda$.

\begin{proof} We first establish a weaker version of \eqref{bllej+}, namely
\begin{equation}\label{bllej}
\|\psi\|_{LE_{j}} + 2^{-j}\|\partial_{r}\psi\|_{LE_{j}} \lesssim \|\psi\|_{L^{\infty}_{t}L^{2}_{r}} +  \|f\|_{L^{1}_{t}L^{2}_{r}+LE_{j}^{*}}.
\end{equation}
Then we show that a direct energy estimate allows us to replace the uniform energy bound by the initial data size and conclude with with \eqref{bllej+}.

Our approach is in the spirit of the one used by the third author in
\cite{Ta-1}, see also \cite{MMT}, using the positive commutator method. First we say that a sequence $\{\alpha_n\}_{n \in \Z}$ of positive numbers is slowly varying if 

\[
|\ln{\alpha_j} - \ln{\alpha_{j-1}}| \leq 2^{-10}, \qquad \forall j \in \Z.
\]
Based on such a sequence we introduce the normed space $X_{k,\alpha}$ and
its dual $X_{k,\alpha}'$ as follows
\[
\begin{split}
\| u \|_{X_{k,\alpha}}^2 & = 2^{2k}\| u \|^2_{L^2(A_{< -k})} + 2^k \sum_{l \geq -k } \alpha_l 2^{-l} \| u \|^2_{L^2(A_l)} \\
\| u \|_{X_{k,\alpha}'}^2 & = 2^{-2k}\| u \|^2_{L^2(A_{< -k})} + 2^{-k}
\sum_{l \geq -k } \alpha_l^{-1} 2^{l} \| u \|^2_{L^2(A_l)}.
\end{split}
\]
For all slowly varying sequences $\{\alpha_n\}_{n \in \Z}$ with $\sum_n
\alpha_n =1$, we claim that
\begin{equation}\label{ldp}
\|\psi\|_{X_{j,\alpha}} + 2^{-j} \|\partial_{r} \psi\|_{X_{j,\alpha}} \lesssim \|\psi\|_{L^{\infty}_{t}L^{2}_{r}} +  \|f\|_{X_{j,\alpha}'+L^1_t L^2_r}.
\end{equation}
If we assume that \eqref{ldp} is true, then, we can let $\beta$ be a slowly varying sequence with $\sum_{n}\beta_{n}=1$. Then, $\sum_{n}(\frac{1}{2}(\alpha_{n}+\beta_{n}))=1$,  and we
apply \eqref{ldp} for $\frac{1}{2}(\alpha_{n}+\beta_{n})$  to
obtain, for instance, that
\begin{equation}
\|\psi\|_{X_{j,\alpha+\beta}} + 2^{-j} \|\partial_{r}\psi\|_{X_{j,\alpha+\beta}} \lesssim \|\psi\|_{L^{\infty}_{t}L^{2}_{r}}
+ \|f\|_{X_{j,\alpha+\beta}'+L^1_t L^2_r},
\end{equation}
from which we derive the weaker estimate
\begin{equation} \label{ldp2}
\|\psi\|_{X_{j,\alpha}} + 2^{-j}\|\partial_{r}\psi\|_{X_{j,\alpha}} \lesssim \|\psi\|_{L^{\infty}_{t}L^{2}_{r}}
+ \|f\|_{X_{j,\beta}'+L^1_t L^2_r}.
\end{equation}
Since any $l^1$ sequence can be dominated by a slowly varying sequence
with a comparable $l^1$ size, we can drop the assumption in
\eqref{ldp2} that $\alpha$ and $\beta$ are slowly varying. By maximizing the
 left-hand side with respect to $\alpha \in l^1$ and by minimizing the right-hand side with respect to $\beta \in l^1$, we obtain \eqref{bllej}.

The remaining part of this step is devoted to the proof of \eqref{ldp}. 
Let 
\[
Q_{j}(u) =  \chi(2^{j}r) r \partial_{r} u + (2+r \partial_{r})(\chi(2^{j}r) u),
\]
where $\chi$ will be chosen to be a smooth function related to
the slowly varying sequence $\{ \alpha_n \}$. A straightforward computation shows that $Q_j$ is antisymmetric; based on this and  the equation for $\psi$ we obtain
\[
\begin{split}
  \Re \int_0^T \langle Q_j \psi, f \rangle dt & = \Re \int_0^T \langle Q_j \psi, 
(i \partial_t - \tilde{H}_{\lambda(t)}) \psi \rangle dt \\
  & =  \Im \int_0^T  \langle Q_j \psi, \partial_t \psi \rangle dt - 
\Re \int_0^T \langle Q_j \psi, \tilde{H}_{\lambda(t)} \psi \rangle dt \\
  & =  \frac12 \Im \int_0^T \partial_t \langle Q_j \psi, \psi \rangle dt -
 \Re \int_0^T \langle Q_j \psi, \tilde{H}_{\lambda(t)} \psi \rangle dt,
\end{split}
\]
which, by rearranging terms, becomes
\begin{equation} \label{idt}
- \Re \int_0^T \langle Q_j \psi, f \rangle dt  + \frac12 \Im \langle Q_j \psi, 
\psi \rangle |_0^T = \Re \int_0^T \langle Q_j \psi, \tilde{H}_{\lambda(t)} \psi \rangle dt.
\end{equation}
The right hand side can be expanded as follows
\[
\begin{split}
\Re \int_0^T \langle Q_j \psi, \tilde{H}_{\lambda(t)} \psi \rangle dt & = \Re \int_0^T \langle Q_j \psi, - \Delta \psi + \tilde{V}_{\lambda(t)} \psi  \rangle dt \\
& = \Re \int_0^T \langle Q_j \partial_r \psi, \partial_r \psi \rangle dt + \Re \int_0^T \langle [\partial_r,Q_j] \psi, \partial_r \psi \rangle dt
+ \Re \int_0^T \langle Q_j \psi, \tilde{V}_{\lambda(t)} \psi \rangle dt \\
& = \Re \int_0^T \langle [\partial_r,Q_j] \psi, \partial_r \psi \rangle dt
+ \frac12 \int_0^T \langle [\tilde{V}_{\lambda(t)},Q_j] \psi, \psi \rangle dt,
\end{split}
\]
where we have used twice the antisymmetry of $Q_j$. We now compute the commutators and start with the easier one,
\[
\begin{split}
\frac12[\tilde{V}_{\lambda(t)},Q_j] = - r \chi(2^j r) \partial_r \tilde{V}_{\lambda(t)} & =  \chi(2^j r) \frac{2}{r^2}
\left(1+8\frac{1+3 \lambda^4(t) r^4}{(1+\lambda(t)^4r^4)^2} \right) \\
& = \frac{2 \chi(2^{j} r)(9+26 r^{4}\lambda(t)^{4}+r^{8}\lambda(t)^{8}}{r^{2}(1+r^{4}\lambda(t)^{4})^{2}}> 0.
\end{split}
\]
The other commutator is
\[
[\partial_r, Q_j] = 2(2^j r \chi'(2^j r)+ \chi(2^j r)) \partial_r  + (3 \cdot 2^j \chi'(2^j r) + 2^{2j} r \chi''(2^j r)). 
\]
From these we obtain
\begin{equation}\label{commcalc1}
\begin{split}
& \Re \int_0^T \langle Q_j \psi, \tilde{H}_{\lambda(t)} \psi \rangle dt  = \frac12 \int_0^T \langle [\tilde{V}_{\lambda(t)},Q_j] \psi, \psi \rangle dt \\
+ &  \int_0^T \int_{0}^{\infty}  \left( 2[\chi(2^{j}r) + 2^{j} r \chi'(2^{j}r)] |\partial_{r} \psi|^{2} + [3   \cdot 2^{j}\chi'(2^{j}r) + r 4^{j} \chi''(2^{j}r)] 
\frac{\partial_{r}|\psi|^{2}}2 \right) rdr dt. 
\end{split}
\end{equation}
We impose the following condition on $\chi$
\begin{equation}\label{chicond}|\chi'(r)| + |r \chi''(r)| \leq \delta(\chi(r)+r \chi'(r))\end{equation}
for some sufficiently small $\delta >0$. This leads to
\[
\Re \int_0^T \langle Q_j \psi, \tilde{H}_{\lambda(t)} \psi \rangle dt  \geq \frac12 \int_0^T \langle [\tilde{V}_{\lambda(t)},Q_j] \psi, \psi \rangle dt + 
\frac12 \int_0^T \int_{0}^{\infty} (|\partial_{r}\psi|^{2} 
- 3 \cdot 2^{2j} \delta |\psi|^2)a_{j}(r) rdr dt,
\]
where $a_{j}(r) = \chi(2^{j}r) + 2^{j} r \chi'(2^{j}r)$. From \eqref{idt} and the previous estimates, we obtain
\begin{equation} \label{maincom}
\begin{split}
 & \int_0^T \langle [\tilde{V}_{\lambda(t)},Q_j] \psi, \psi \rangle dt +  \int_0^T \int_{0}^{\infty} (|\partial_{r}\psi|^{2} 
- 3 \cdot 2^{2j} \delta |\psi|^2)a_{j}(r) rdr dt \\
\leq & -2 \Re \int_0^T \langle Q_j \psi, f \rangle dt  +  \Im \langle Q_j \psi, \psi \rangle |_0^T. 
\end{split}
\end{equation}
The main idea of what follows next is to show that there is an appropriate choice of $\chi$, depending on the sequence $\{\alpha_n\}$ such that we can derive \eqref{ldp} from \eqref{maincom}. We first increase $\alpha_n$
the so that it remains slowly varying and, in addition, satisfies

\begin{equation} \label{ren}
\alpha_n=1, \ \ \ \mbox{for} \  n \leq n_0-j; \qquad \sum_{n \geq n_0-j} \alpha_n \approx 1.
\end{equation}
Here $n_0$ is a positive number to be chosen later. 

We claim that given a slowly varying sequence $\alpha_n$, satisfying \eqref{ren}, and $\delta > 0$ we can find
$\chi$ satisfying \eqref{chicond}, so that
\begin{equation} \label{chicond2}
a_j(r) \gtrsim \frac{\alpha_{n}}{1+2^{n+j}}, \qquad r \approx 2^n
\end{equation}
 and the following three fixed time bounds hold:
\begin{equation} \label{lek}
\| Q_j \psi \|_{L^2_r} \lesssim  \delta^{-1} \|\psi\|_{L^2_r}, \quad 
\| Q_j \psi\|_{X_{k,\alpha}} \lesssim \delta^{-1} \|\psi\|_{X_{k,\alpha}},
\end{equation}
\begin{equation}
2^{2j} 
\int_0^\infty a_j(r) |\psi|^2 rdr \lesssim \int_{0}^\infty a_j(r) |\partial_r \psi|^2 rdr + \frac12 \int_0^\infty [\tilde{V},Q] |\psi|^2 rdr.
\label{coercive}\end{equation}
These estimates are claimed for functions $\psi$ which are localized
at frequency $2^j$ in the frame $\tilde H_{\lambda(t)}$ and the constants involved do not depend on the value of $\lambda(t) \in (0,\infty)$ or the sequence $\{\alpha_n\}$ or $\delta$.

Going back to the estimate \eqref{maincom} and making use of \eqref{chicond2}, \eqref{lek} and \eqref{coercive}, 
we obtain
\begin{equation} \label{Xjle}
\|\psi\|_{X_{j,\alpha}}^{2} + 2^{-2j} \|\partial_{r} \psi\|_{X_{j,\alpha}}^{2} \lesssim \|\psi\|_{L^{\infty}_{t}L^{2}_{r}}^{2} +  \| f \|_{X_{j,\alpha}'+L^1_t L^2_r}^{2},
\end{equation}
when all terms are restricted to the time interval $[0,T]$, but
with the a constant independent of $T$. This implies \eqref{ldp}, which in turn was shown to imply our main claim \eqref{bllej}.  

It would be helpful to explain the role of the small constant $\delta$ in the above. $\delta$ needs to be small enough so that when taking into account \eqref{coercive}, the following holds true:
\[
\begin{split}
& \int_0^T \langle [\tilde{V}_{\lambda(t)},Q_j] \psi, \psi \rangle dt +  \int_0^T \int_{0}^{\infty} (|\partial_{r}\psi|^{2} 
- 3 \cdot 2^{2j} \delta |\psi|^2)a_{j}(r) rdr dt \\
\gtrsim & \ 2^{2j} \int_0^\infty a_j(r) |\psi|^2 rdr + \int_{0}^\infty a_j(r) |\partial_r \psi|^2 rdr. 
\end{split}
\]
Once this is achieved, $\delta$ simply becomes just another constant whose quantification is not necessary; in particular we could have ignored its quantification in \eqref{lek} or, as we already have done, its quantification in \eqref{Xjle}.

We now proceed with the construction of $\chi$ satisfying \eqref{chicond}, \eqref{chicond2}, \eqref{lek} and \eqref{coercive}.  Based on this, we construct a slowly varying function $\alpha$ such that
\[
\alpha(s) \approx \alpha_n \ \ \mbox{if} \ \ s \approx 2^n
\]
and with symbol regularity
\[
|\partial^l \alpha(s)| \lesssim_l s^{-l} \alpha(s), \qquad l \in \N.
\]
Due to the first condition in \eqref{ren} we can take $\alpha$ such that
$\alpha(s)=1$ for $s \leq 2^{n_0-j}$. We then construct the function $\chi$ by
\[
s \chi(s) = \int_0^s \alpha(2^{-j} \tau) h(\tau) d\tau,
\]
where $h$ is a smooth adapted variant of $ r^{-1}$,
namely $h(s)=1$ for $s \leq 2^{n_0}$ and $h(s) \approx 2^{n_0} s^{-1}$ for $s \geq 2^{n_0+1}$; in particular $h$ has symbol type estimates $|h^{(k)}(s)|\lesssim_k  2^{n_0} s^{-k-1}, x \geq 2^{n_0}$. With the observation that $\chi(s)=1$ for $s \leq 2^{n_0}$, one easily verifies the pointwise bounds
\begin{equation} \label{pchi}
\chi(s) \approx (1+2^{-n_0} s)^{-1},
\qquad |\chi^{(l)}(s)| \lesssim 2^{-ln_0} (1+2^{-n_0} s)^{-l-1}, \quad 1 \leq l \leq 4.
\end{equation}
Furthermore, we have 
\[
(s \chi(s))' = \alpha(2^{-j} s) h(s) \gtrsim (1+2^{-n_0} s)^{-1.1},
\qquad 
|(s \chi'(s))'| \lesssim 2^{-n_0} (1+2^{-n_0} s)^{-2}.
\]
It is a straightforward exercise to verify that $\chi$ satisfies
\eqref{chicond2}. A direct computation shows that
\[
\begin{split}
|\chi'(r)| + |r \chi''(r)| &\lesssim \frac{2^{-n_{0}}}{(1+2^{-n_{0}}r)^{2}} \lesssim  2^{-n_{0}}(r\chi(r))',
\end{split}
\]
where the constants involved are independent of $n_0$; thus by choosing $2^{-n_{0}} \sim \delta$ we obtain that $\chi$ satisfies \eqref{chicond}.

Next we seek to establish the first estimate in \eqref{lek} which requires an estimate on $\|Q_{j}\psi \|_{L^{2}_{r}}$. We have
\begin{equation} \label{Qjexp}
Q_{j}\psi = 2 r \chi(2^{j}r) \partial_{r} \psi + 2\chi(2^{j}r) \psi + r \partial_{r}(\chi(2^{j}r)) \psi.
\end{equation}
From this we obtain
\begin{equation} \label{leaux1}
\|Q_{j}\psi \|_{L^{2}_{r}} \lesssim \|r\chi(2^{j}r) \partial_{r}\psi\|_{L^{2}_r} +  \|\psi\|_{L^{2}_r}.
\end{equation}
Due to the frequency localization of $\psi$, we can write
\[
\psi(t) = \widetilde{P_{j}^{\lambda(t)}} \psi(t) = \int_{0}^{\infty} \widetilde{K_{j}}^{\lambda(t)}(r,s) \psi(t,s) sds,
\]
from which it follows that
\[
\partial_{r} \psi(t,r) = \int_{0}^{\infty} \partial_{r}\widetilde{K_{j}}^{\lambda(t)}(r,s) \psi(s) s ds. 
\]
Moreover, from Lemma \ref{p-lp} we obtain the following estimate:
\begin{equation} \label{drpjtermkernel}
|r \chi(2^{j}r) \partial_{r}\widetilde{K}_{j}^{\lambda(t)}(r,s)| \lesssim_N \frac{2^{2j} 2^{n_{0}}}{(1+2^{j}(r+s))(1+2^{j}|r-s|)^{N}}.
\end{equation}
This allows us to estimate as follows:
\begin{equation}
\begin{split}
 \|r \chi(2^{j}r) \partial_{r} \psi\|_{L^{2}_{r}} 
& \lesssim \sup_{s}  \|r \chi(2^{j}r)\partial_{r}\widetilde{K}_j^{\lambda(t)}(r,s)\|_{L^1_r}^\frac12 \cdot
\sup_{r}  \|r \chi(2^{j}r)\partial_{r}\widetilde{K}_j^{\lambda(t)}(r,s)\|_{L^1_r}^\frac12 \cdot \|\psi\|_{L^2} \\
& \lesssim 2^{n_0} \|\psi\|_{L^2}. 
\end{split}
\end{equation}
Inserting this estimate into \eqref{leaux1} leads to:
\[
\|Q_{j}\psi(t)\|_{L^{2}_{r}} \lesssim \delta^{-1} \|\psi(t)\|_{L^{2}_r},
\]
which concludes the first part of \eqref{lek}. We now turn to the second part of \eqref{lek}, namely establishing that
\[
\|Q_{j} \psi(t) \|_{X_{j,\alpha}}  \lesssim \delta^{-1} \| \psi(t)\|_{X_{j,\alpha}}.
\]
We start from the formula \eqref{Qjexp} for $Q_j\psi$. Using the rapid decay away from the diagonal in \eqref{drpjtermkernel}, and the fact that the weights $\alpha_k$ are slowly varying, we obtain
\[
\|r \chi(2^{j}r) \partial_{r} \psi\|_{X_{j,\alpha}}^{2} 
\lesssim  2^{2 n_{0}} \|\psi\|_{X_{j,\alpha}}^{2},
\]
from which it follows that 
\[ \|Q_{j}\psi\|_{X_{j,\alpha}}^{2} 
\lesssim  2^{2n_{0}} \|P_{j}\psi\|_{X_{j,\alpha}}^{2}.
\]
Recalling that $\delta \approx 2^{-n_0}$, this concludes our argument for \eqref{lek}. 

We are then left with verifying \eqref{coercive}; we start with
\[
\psi(t) = \widetilde{P}_{j}^{\lambda(t)} \psi(t) 
= \widetilde{P}_{j}^{\lambda(t)} \widetilde{H}_{\lambda(t)}^{-1} L_{\lambda(t)}L^{*}_{\lambda(t)} \psi 
= \int_{0}^{\infty}  K_{j}^{1,\lambda(t)}(r,s) \cdot L^{*}_{\lambda(t)} \psi (s) s ds. 
\]
From this we obtain:
\[
\sqrt{a_{j}(r)} 2^{j} \psi(t,r) = \int_{0}^{\infty}  K_7^{\lambda(t)}(s,r) L^{*}_{\lambda(t)} \psi(s) \sqrt{a_{j}(s)}s ds,
\]
where
\[
K_7^{\lambda(t)}(r,s)= \frac{\sqrt{a_{j}(r)}}{\sqrt{a_{j}(s)}} 2^{j} K_{j}^{1,\lambda(t)}(r,s).
\]
Using the symbol-type estimates on $\alpha$ and $h$, and noting that
\[
a_{j}(x) = \chi(2^{j}x) + 2^{j}x \chi'(2^{j}(x)) = \alpha(x) h(2^{j}x),
\]
we obtain $\frac{|a_{j}'(s)|}{a_{j}(s)} \lesssim s^{-1}$
which implies $|a_j(r)| \approx |a_j(s)|$ in the regime $r \sim s$. 
Based on this and the estimates on $K_j^{1,\lambda(t)}$ from \eqref{ker3}, we obtain
\[
|K_7^{\lambda(t)}(s,r)| \lesssim_N  \frac{2^{2j}}{(1+2^{j}(r+s)) (1+2^{j}|r-s|)^{N}}, \quad r \sim s. 
\]
If $r > 2s$ or $r \leq \frac{s}{2}$, then,
\[
\frac{a_{j}(r)}{a_{j}(s)} = \frac{\chi(2^{j}r)+2^{j}r \chi'(2^{j}r)}{a_{j}(s)} 
\lesssim  \frac{1}{(1+2^{-n_{0}+j}r)a_{j}(s)} 
\lesssim \frac{1}{(y\chi)'(2^{j}s)} 
\lesssim  (1+2^{j}s)^{1.1}.
\]
Using this we can estimate
\[
|K_{7}^{\lambda(t)}(r,s)| \lesssim_N \frac{\sqrt{a_{j}(r)}}{\sqrt{a_{j}(s)}} \frac{2^{2j}}{(1+2^{j}(r+s))(1+2^{j}|r-s|)^{N}} \lesssim \frac{2^{2j}}{(1+2^{j}(r+s))(1+2^{j}\text{max}\{r,s\})^{N-2}}.
\]
From this it follows that
\[\int_{0}^{\infty} a_{j}(r) 2^{2j} |\psi|^{2} r dr \lesssim  \int_{0}^{\infty} |L^{*}_{\lambda(t)} \psi|^{2}(r) a_{j}(r) r dr. 
\]
From the definition of $L^{*}_{\lambda(t)}$, see \eqref{Llambdadef}, it follows that
\[
|L^{*}_{\lambda(t)} u|^{2} \lesssim \frac{1}{r^{2}} |u(r)|^{2} +  |\partial_{r}u(r)|^{2}. 
\]
Since we have the straightforward inequality 
$|a_{j}(r)| \lesssim (1+2^{-n_{0}+j}r)^{-1} 
\lesssim  \chi(2^{j}r)$, we can estimate
\[
\begin{split}
a_{j}(r) |L^{*}_{\lambda(t)} \psi(r)| &\lesssim \frac{ |\psi(t,r)|^{2}}{r^{2}} a_{j}(r) +  |\partial_{r} \psi(t,r)|^{2} a_{j}(r)\\
&\lesssim (\frac{\chi(2^{j}r) |\psi|^{2}}{r^{2}}+|\partial_{r}\psi|^{2})a_{j}(r). 
\end{split}
\]
Therefore, we  have 
\[
\begin{split} 
\int_{0}^{\infty} a_{j}(r) 2^{2j}|\psi(r)|^{2} r dr & \lesssim    \int_{0}^{\infty}  |\psi|^{2} \frac{4 \chi(2^{j}r)}{r^{2}} \frac{(9+26 r^{4}\lambda(t)^{4}+r^{8}\lambda(t)^{8})}{(1+r^{4}\lambda(t)^{4})^{2}} r dr \\
&+  \int_{0}^{\infty} |\partial_{2} \psi|^{2} a_{j}(r) r dr,
\end{split}
\]
which is precisely \eqref{coercive}. This finishes the proof of \eqref{bllej}. To conclude the proof of the Lemma, we involve a simple energy argument: multiplying 
\eqref{wlin-eq} by $\bar{\psi}$, integrating with respect to $r$ and taking the imaginary part, gives
\[
\frac12 \frac{d}{dt} \|\psi(t) \|^2_{L^2_r} = \Im \la f, \psi \ra. 
\]
Integrating this on an arbitrary interval $[0,t_0]$ gives
\[
\|\psi(t_0) \|^2_{L^2_r} \leq \|\psi(0) \|^2_{L^2_r} + 2 \int_0^{t_0} |\Im \la f, \psi \ra| dt \leq \|\psi(0) \|^2_{L^2_r} + 2M^{-1} \|\psi\|_{LE_k}^2 
+2M \|f\|_{LE_k^*}^2.  
\]
Thus we obtain
\[
\|\psi\|^2_{L^\infty_t L^2_r} \leq \|\psi(0) \|^2_{L^2_r} + 2M^{-1} \|\psi\|_{LE_k}^2 
+2M \|f\|_{LE_k^*}^2.
\]
Using \eqref{bllej} for the term $\|\psi\|_{LE_k}$ and choosing $M$ large enough gives \eqref{bllej+} and this finishes the proof of our Lemma.

\end{proof}

\subsection{Full local energy decay}
Here we assemble the dyadic local energy bounds of the previous section into a full local energy bound. This does not require smallness of $\lambda'/\lambda^2$ 
in $L^2$. Instead, we will be able to track  the dependence  of the implicit constants on the above $L^2$ norm. This justifies introducing the notation (already mentioned in section \ref{LtHf})
\begin{equation} \label{Mfdef}
\Mf (T) := \|\frac{\lambda'}{\lambda^{2}}\|_{L^{2}_t[0,T]}.
\end{equation}

For some of the dyadic estimates we will also use the more refined 
quantities
\begin{equation} \label{Mjdef}
\Mf_j (T) = \|\chi_{\lambda = 2^j} \frac{\lambda'}{\lambda^{2}}\|_{L^{2}[0,T]},
\end{equation}
which measures the same $L^2$ norm but with a weight concentrated 
around the dyadic region $\{t\in [0,T]; \lambda(t) \approx 2^j\}$.  The two are related by
\[
\Mf^2(T) \approx \sum_j \Mf_j^2(T),
\]
where the constants used in $\approx$ are independent of $T$.

Our main well-posedness result concerning the linear $\tilde H$ equation
is as follows:
\begin{p1}\label{p:mainS} 
Assume that $\lambda'/\lambda^2 \in L^2[0,T]$. Then the equation 
\eqref{wlin-eq} is well-posed in $L^2$, and
the solution $\psi$ satisfies the following bound: 
\begin{equation} \label{lefinal+}
\begin{split} 
\|\psi\|_{\ell^2 (LE \cap L^\infty_t L^2_r) [0,T]} &\lesssim 
(1+\Mf(T)) \|\psi(0)\|_{L^{2}} + (1+\Mf(T))^2 \|f\|_{\ell^2 (LE^*+ L^1_t L^2_r) [0,T]}.
\end{split}
\end{equation}
\end{p1}
For clarity we also write this bound in an
expanded form,
\begin{equation}\label{lefinal+re}
\begin{split} 
 & \!\! \sum_{j} \|P_{j}\psi\|_{LE_{j} \cap L^\infty L^2}^{2}
+ 2^{-2j} \|\partial_r P_{j}\psi\|_{LE_{j} \cap L^\infty L^2}^{2} \\
\lesssim & 
\ (1+\Mf(T))^2 \|\psi(0)\|_{L^{2}}^{2} + 
(1+\Mf(T))^4 \sum_{j} \|P_{j}f\|_{L^{1}_{t}L^{2}_{r}+LE_j^*}^{2}.
\end{split}
\end{equation}
Here we have added also the bounds for $\partial_r P_j \psi$,
which follow almost freely from the argument due to the frequency localization.

\begin{proof} 
Since the equation \eqref{wlin-eq} coincides with its adjoint
equation and the time is reversible, a standard duality argument shows that $L^2$ well-posedness follows from the bound \eqref{lefinal+}. Hence we turn our attention to the proof of this bound.

 For each $j$, the functions $P_{j}^{\lambda}\psi$ solve the Cauchy problem
 \begin{equation} \label{Pjpsieq}
 \begin{cases} 
 i \partial_{t}P_{j}^\lambda\psi - \widetilde{H}_{\lambda(t)} P_{j}^\lambda \psi = P_{j} ^\lambda f + g_{j}\\
 P_{j}^\lambda \psi(0) = P_{j}^\lambda \psi_{0}
 \end{cases}
 \end{equation}
where $g_j$ arises from the commutator of $P_j$ with $\partial_t$,
\[
g_j = i [\partial_t, P_j^\lambda] \psi.
\]
 Hence our first task is to obtain good estimates for $g_j$. 
The above commutator is described in Lemma~\ref{dtP} in terms of the 
transference operator $\tcK$, which is in turn estimated in Lemma~\ref{l:mj-com}.  Here we will use these building blocks to 
prove the following 

\begin{l1}\label{l:gj}
The above commutators $g_j$ satisfy the bounds
\begin{equation}\label{e-to-le}
\| g_j\|_{LE_j^*} \lesssim \Mf_j(T)  \sum_k 
\chi_{k=j} \| \psi_k \|_{L^\infty_t L^2_r} 
\end{equation} 
with $\psi_k := P_k^\lambda \psi$, respectively 
\begin{equation}\label{le-to-e}
\| g_j\|_{L^1_t L^2_r} \lesssim \Mf_j(T)  \sum_k 
\chi_{k=j} \| \psi_k \|_{LE_k}.
\end{equation} 
In the above, all the space-time norms are restricted to the time interval $[0,T]$. 
\end{l1}

\begin{proof} Just as in the statement of the Lemma, in what follows below 
all the space-time norms are restricted to the time interval $[0,T]$; in order to keep the formualas compact we skip this from notation.

From Lemma~\ref{dtP} we have
\[
g_j = i \lambda' \FtH_\lambda^{-1}[ \tcK_\lambda,m_j] \FtH_\lambda \psi.
\]
We split this into
\[
g_j = \sum_k g_{jk}, \qquad g_{jk} = i \lambda' \FtH_\lambda^{-1}[ \tcK_\lambda,m_j]  \tm_k \FtH_\lambda \psi_k, 
\]
and then it remains to show that
\begin{equation}\label{gj-pieces}
\| g_{jk}\|_{LE_j^*} \lesssim \Mf_j(T)  
\chi_{k=j} \| \psi_k \|_{L^\infty_t L^2_r} ,
\qquad
\| g_{jk}\|_{L^1_t L^2_r} \lesssim \Mf_j(T)  \| \psi_k \|_{LE_k}. 
\end{equation}
These are dual bounds so it remains to prove the first one.
For that we need to estimate $g_{jk}$ in $L^2_t L^2_r(A_m)$, first for $m < j$
and then for $m \geq j$. In the first case we use \eqref{com-est} to obtain the fixed time bound
\[
\|g_{jk}(t)\|_{L^2_r(A_{<j})} \lesssim 2^{j} \chi_{\lambda = 2^j} \chi_{j=k} 
\frac{|\lambda'|}{\lambda^2} \| \psi_k(t)\|_{L^2},
\]
which after time integration yields
\[
\|g_{jk}\|_{L^2_t L^2_r(A_{<j})} \lesssim 2^{j} \Mf_j \| \psi_k\|_{L^\infty_t L^2_r}.
\]
In the second case we use \eqref{rcom-est} to obtain the fixed time bound
\[
\|g_{jk}(t)\|_{L^2_r(A_{m})} \lesssim 2^{-m} \chi_{\lambda = 2^j} \chi_{j=k}
\frac{|\lambda'|}{\lambda^2} \| \psi(t)\|_{L^2}, \qquad m \geq j,
\]
which after time integration yields
\[
\|g_{jk}\|_{L^2_t L^2_r(A_{m})} \lesssim 2^{-m} \Mf_j(T) \| \psi\|_{L^\infty L^2}.
\]
Combining the two cases we arrive at 
\[
\| g_{jk} \|_{LE_j^*} \lesssim 
\Mf_j(T)  \chi_{j=k} \| \psi_k\|_{L^\infty L^2}
\sum_{m \geq -j} 2^{-\frac{m+j}2}\lesssim \Mf_j(T) \chi_{j=k} \| \psi_k\|_{L^\infty_t L^2_r}, 
\]
as desired.

\end{proof}

Now we return to the proof of the proposition, and apply Lemma~\ref{l:LEj} to $P_j^{\lambda(t)} \psi$, using Lemma~\ref{l:gj} for $g_j$.
This yields
\[
\begin{aligned}
	\| P_j^{\lambda(t)} \psi\|_{LE_j \cap L^\infty_t L^2_r} \lesssim & \ \|P_j^{\lambda(0)} \psi(0)\|_{L^2}
	+ \| P_j^{\lambda(t)} f + g_j\|_{L^1L^2 + LE_j^*}
	\\
	\lesssim & \ \|P_j^{\lambda(0)} \psi(0)\|_{L^2}
	+ \| P_j^{\lambda(t)} f \|_{L^1L^2 + LE_j^*} + \|g_j\|_{LE_j^*}
	\\
	\lesssim & \ \|P_j^{\lambda(0)} \psi(0)\|_{L^2}
	+ \| P_j^{\lambda(t)} f \|_{L^1L^2 + LE_j^*} + \Mf_j(T) \| \psi \|_{L^\infty L^2}. 
\end{aligned}
\]
We square this bound and sum over $j \in \mathbb{Z}$ to arrive at  
\begin{equation}
	\label{lefinal-}\begin{split} \sum_{j} \|P_{j}^{\lambda(t)}\psi\|_{LE_{j} \cap L^\infty L^2}^{2} &\lesssim 
		\sum_{j} \|P_{j}^{\lambda(t)}f\|_{L^{1}_{t}L^{2}_{r}+LE_j^*}^{2}
		+ (1+ \Mf^2(T)) \|\psi\|_{L^\infty L^{2}}^{2}.
\end{split}\end{equation}
It remains to supplement this with an energy estimate for $\psi$.
We have 
\[
\frac{d}{dt} \| \psi(t)\|_{L^2}^2 = -2\Re \langle \psi, if \rangle
=  -2 \sum_{j,k \in \mathbb{Z}} \Re \langle P_k^{\lambda(t)} \psi, i P_j^{\lambda(t)} f \rangle.
\]
The terms in the last sum vanish unless $|k-j|\leq 3$, in which case
the $LE_k$ and $LE_j$ norms are equivalent by \eqref{lekpert}.
We integrate in time the last relation and use Cauchy-Schwarz to obtain
\[
\begin{aligned}
	\| \psi\|_{L^\infty L^2}^2 \lesssim & \ \|\psi(0)\|_{L^2}^2
	+ \sum_{|j-k| \leq 3} \| P_k^{\lambda(t)} \psi \|_{L^\infty L^2\cap LE_k} \| P_j^{\lambda(t)} f \|_{L^1 L^2 + LE_j^*}
	\\
	\lesssim & \ \|\psi(0)\|_{L^2}^2
	+ \left(\sum_{k} \| P_k^{\lambda(t)} \psi \|_{L^\infty L^2\cap LE_k}^2\right)^\frac12
	\left(\sum_j \| P_j^{\lambda(t)} f \|_{L^1 L^2 + LE_j^*}^2\right)^\frac12
\end{aligned}
\]
Finally, we insert this in \eqref{lefinal-} and apply Cauchy-Schwarz
one more time to arrive at \eqref{lefinal+}.

\end{proof}

\subsection{ Adding the Strichartz norms to the mix}\label{ss-str}
Here we start with the estimate \eqref{lefinal+}, written in the 
expanded form \eqref{lefinal+re}, and show that we can add in the 
weighted $L^4$, respectively $L^{\frac43}$ norms to arrive at 
the main bound \eqref{mainS} in Theorem~\ref{mainTl2}.

We do this in several steps, beginning with a Lemma which captures 
the essential one dimensional Strichartz estimate:
\begin{l1}\label{l4tlinftyrest}
Let $\psi$ be a solution to \eqref{wlin-eq}. Then the 
following estimate holds:
\begin{equation}\label{1d-Str}
 \sup_{j \geq -k} 2^{\frac{j}{2}} \|\psi\|_{L^{4}_{t}L^{\infty}_{r}(A_{j})}
\lesssim \|\psi \|_{L^\infty_t L^2_r \cap LE_k} + 2^{-k} \|\partial_r \psi \|_{L^\infty_t L^2_r \cap LE_k} + \|f\|_{N_k}.
 \end{equation}
\end{l1}
This lemma is designed for functions at frequency $2^k$, but does not actually assumes any frequency localization. We disregard for now the region $A_{<k}$, where matters will be simpler but will involve the frequency localization.
\begin{proof}
We begin by localizing the problem to each region $A_j$, using a suitable
cutoff function $\chi_j$. For $\psi_j = \chi_j \psi$ we can write a one 
dimensional Schr\"odinger equation
\[
(i \partial_t + \partial_r^2) \psi_j = f_j:= g_j + \chi_j f
\]
where the additional source term $g_j$ is given by
\[
g_j = 2 \chi'_j \partial_r \psi + \chi''_j \psi - {\chi_j} \frac{1}{r}\partial_r \psi - {\chi_j} \tilde V_\lambda \psi,
\]
and can be readily estimated by
\[
\| g_j\|_{LE_k^*} \lesssim  \|\psi \|_{LE_k} + 2^{-k} \|\partial_r \psi \|_{LE_k}. 
\]
Then \eqref{1d-Str} reduces to proving
\begin{equation}\label{1d-Str-loc}
  2^{\frac{j}{2}} \|\psi_j\|_{L^{4}_{t}L^{\infty}_{r}}
\lesssim \|\psi_j\|_{LE_k} + \|f_j\|_{N_k}
 \end{equation}
for $\psi_j$ localized in $A_j$.

\medskip

We  show that this last bound is a direct consequence of the 
one dimensional Strichartz estimate for $\psi$, which implies that 
\begin{equation}\label{1D-str}
\| \psi_j\|_{L^4_t L^\infty_r \cap L^\infty_r L^2_r(dr)} \lesssim 2^{\frac{k-j}2} \|\psi_j\|_{L^2_t L^2_r(dr)}
+ \|f_j\|_{2^\frac{k-j}2 L^2_t L^2_r(dr) + L^1_t L^2_r(dr) + L^{\frac43}_t L^1_r(dr)}.
\end{equation}
This is not exactly the classical Strichartz estimate, but it is derived from it in a straightforward manner. Precisely, the norms involved on the right-hand side terms only involve $L^p_t$ spaces with $p \leq 2$, so they are square summable with respect to time interval decompositions. This allows one to apply a second localization, now with respect to time,  on the time scale $2^{-k+j}$, with an error that may be included in the $L^2$ component of $f_j$. This reduces \eqref{1D-str} to the case when $\psi_j$ 
is localized in a time interval of size $2^{-k+j}$. But on this time scale the $L^2$ norm of $\psi$ bounds the averaged energy, while the $L^2L^2$ component of $f_j$ may be included in $L^1 L^2$ by H\"older's inequality.

The last inequality uses the one dimensional measure $dr$. Converting 
to the $rdr$ measure with $r \approx 2^j$  we arrive at
\begin{equation}
2^{\frac{j}2} \| \psi_j\|_{L^4_t L^\infty_r} \lesssim 2^{\frac{k-j}2} \|\psi_j\|_{L^2_{t,r}}
+ \|f_j\|_{2^\frac{k-j}2 L^2_{t,r} + L^1 L^2_r + 2^{-\frac{j}2} L^{\frac43}_t L^1_r}.
\end{equation}
To see that this implies \eqref{1d-Str-loc}, it suffices to verify that 
restricted to $A_j$ we have the uniform inclusion
\[
N_k \subset 2^\frac{k-j}2 L^2_{t,r} + L^1_t L^2_r + 2^{-\frac{j}2} L^{\frac43}_t L^1_r.
\]
It suffices to consider the $L^{\frac43}$ component of $N_k$, for which this is a straightforward interpolation computation.

\end{proof}

The corresponding bound inside $A_{<-k}$ is instead a Bernstein type inequality:

\begin{l1}
Let $\psi$ be a function which is localized at frequency $2^k$.
Then we have
\begin{equation}\label{1D-Bern}
2^{-\frac{k}{2}} \|\psi\|_{L^{4}_{t}L^{\infty}_{r}(A_{<-k})}
\lesssim \|\psi \|_{L^\infty_t L^2_r \cap LE_k}.
\end{equation}
\end{l1}

\begin{proof}
To capture the frequency localization we rewrite 
the above inequality in the form
\[
2^{-\frac{k}{2}} \|P_k^\lambda \psi\|_{L^{4}_{t}L^{\infty}_{r}(A_{<-k})}
\lesssim \|\psi \|_{L^\infty_t L^2_r \cap LE_k}. 
\]
Then the kernel bounds for $P_k^\lambda$ from \eqref{ker1} yield the following two 
bounds: 
\[
\|P_k^\lambda \psi\|_{L^{\infty}_{t,r}(A_{<-k})} \lesssim 2^k \|\psi \|_{L^\infty_t L^2_r}, \quad
\|P_k^\lambda \psi\|_{L^2_tL^{\infty}_{r}(A_{<-k})} \lesssim \|\psi \|_{LE_k}.
\]
Interpolating these to bound leads to the desired $L^4_t L^\infty_r$ bound.
\end{proof}

The next step is to combine the above two Lemmas in order to prove
a frequency localized $L^4_k$ bound:

\begin{l1}\label{l:SNk}
Let $\psi$ be a solution to \eqref{wlin-eq} which is localized at frequency $2^k$. Then we have
\begin{equation}\label{dyadic-Str}
\| \psi \|_{S_k} \lesssim \| \psi\|_{L^\infty_t L^2_r \cap LE_k} + \| f\|_{N_k}.
\end{equation}
\end{l1}

\begin{proof}
We need to estimate $\psi$ in $L^4_k$, which requires $L^4_{tr}$ bounds in each of the sets $A_j$. We have two cases:

i) If $j < -k$ then we use \eqref{1D-Bern} and H\"older's inequality in $r$.

ii) If $j > k$, then we first use the bounds in \eqref{ker2} for the kernel of $\partial_r P_k^\lambda$ to bound
$\partial_r \psi$,
\[
 2^{-k} \|\partial_r \psi \|_{L^\infty_t L^2_r \cap LE_k} \lesssim 
 \|\psi \|_{L^\infty_t L^2_r \cap LE_k}. 
\]
It then follows from \eqref{1d-Str} that
\[
 \sup_{j \geq -k} 2^{\frac{j}{2}} \|\psi\|_{L^{4}_{t}L^{\infty}_{r}(A_{j})}
\lesssim \|\psi \|_{L^\infty_t L^2_r \cap LE_k} + \|f\|_{N_k}.
 \]
This allows us to bound the $L^4_t L^\infty_r$ norm in  the left hand 
of \eqref{1d-Str}. Then we interpolate this estimate with the $L^\infty_t L^2_r$ energy bound to obtain an estimate in $L^6_{t,r}$, which is then interpolated 
with the $L^2_{t,r}$ bound from the local energy to obtain the desired conclusion; tracking the correct powers in this interpolation is left as an exercise.
\end{proof}

Finally, we begin to assemble the dyadic bounds so that we complete the proof of our main result in this section, Theorem~\ref{mainTl2}. Applying \eqref{dyadic-Str} to $P_k \psi$ 
and then Lemma~\ref{l:gj} we obtain
\[
\begin{aligned}
\| P_k^\lambda \psi \|_{S_k} \lesssim & \ \| P_k^\lambda \psi\|_{L^\infty_t L^2_r \cap LE_k} + \| P_k^\lambda f\|_{N_k} + \|g_k\|_{LE_k^*}
\\
\lesssim & \ \| P_k^\lambda \psi\|_{L^\infty_t L^2_r \cap LE_k} + \| P_k^\lambda f\|_{N_k} + \Mf_k(T) \|\psi\|_{L^\infty_t L^2_r}.
\end{aligned}
\]
Square summing these bounds yields
\begin{equation}\label{apriori-str}
\|\psi\|_{\ell^2 S} \lesssim (1+\Mf(T)) \|\psi \|_{\ell^2 (LE \cap L^\infty_t L^2_r)} + \| f\|_{\ell^2 N}.
\end{equation}

To continue we consider first the case when $f \in \ell^2 (LE^*+L^1L^2)$. 
Then we can apply Proposition~\ref{p:mainS} to conclude that 
\begin{equation}\label{direct}
\|\psi\|_{\ell^2 S} \lesssim (1+\Mf(T))^2 \|\psi_0 \|_{L^2} + (1+\Mf(T))^3\| f\|_{\ell^2 (LE^*+L^1L^2)}.
\end{equation}
From here, we claim that a duality argument 
yields
\begin{equation}\label{dual}
\|\psi\|_{\ell^2 (LE\cap L^\infty L^2)} \lesssim (1+\Mf(T))^2 \|\psi_0 \|_{L^2} + (1+\Mf(T))^3\| f\|_{\ell^2 N}.
\end{equation}
To see this, we pair $\psi$ and $f$ with a second solution $\tilde \psi$, $\tilde f$. The duality 
relation between the two in a time interval $[0,T]$ yields
\[
\langle\psi, \tpsi\rangle |_0^T
= \int_0^T \langle \psi, i\tilde f \rangle 
- \langle f, i\tilde \psi \rangle \, dt.
\]
This implies 
\[
\langle\psi(T), \tpsi(T)\rangle 
- \int_0^T \langle \psi, i\tilde f \rangle \, dt
\lesssim \| \psi(0)\|_{L^2} \|\tpsi(0)\|_{L^2}
+ \|f\|_{\ell^2 N} \| \tilde \psi \|_{\ell^2 S}, 
\]
which combined with \eqref{direct} applied for the backward problem yields
\[
\langle\psi(T), \tpsi(T)\rangle 
- \int_0^T \langle \psi, i\tilde f \rangle \, dt
\lesssim (\| \psi(0)\|_{L^2} + \|f\|_{\ell^2 N})
(\|\tilde \psi(T)\|_{L^2}
+  \| \tilde f \|_{\ell^2 (LE^*+L^1L^2)} ).
\]
This in turn gives \eqref{dual}.

Finally, combining \eqref{dual} with \eqref{apriori-str} we obtain the conclusion of Theorem~\ref{mainTl2}.

\subsection{The \texorpdfstring{$\ell^{1}$}{} Strichartz and local energy bounds}  
Here we turn our attention to the $l^1$ bounds in Theorem~\ref{mainTl1}, which we now prove. Precisely, we will show that

\begin{p1}
The following bound holds for solutions to the 
equation \eqref{wlin-eq}:
\begin{equation} \label{lejfinal++p}
\begin{split}  
\|\psi\|_{\ell^1 S} &\lesssim (1+\Mf(T))^2 \|\psi(0)\|_{\LX} + 
(1+\Mf(T))^4  \|f\|_{\ell^1 N}.
\end{split}
\end{equation}
In the above the space-time norms are restricted to the time interval $[0,T]$. 
\end{p1}

\begin{proof}
As in the case of the $\ell^2$ bound, our starting point is Lemma~\ref{l:LEj}, which applied to $P_j \psi$ yields  the dyadic bound
\begin{equation}
\label{lejj++}\begin{split}  \|P_{j}\psi\|_{S_j} &\lesssim  \|\psi_j(0)\|_{L^{2}} + 
 \|P_{j}f\|_{N_j}+ \|g_j\|_{LE^*_j}.
\end{split}\end{equation}
For $g_j$ we use Lemma~\ref{l:gj},
which implies that
\begin{equation} \label{gjLE3}
\|g_j\|_{LE^*_j} \lesssim \Mf_j(T) \sum_k \| \psi_k\|_{L^\infty_t L^2_r} \chi_{k=j}.
\end{equation}
Hence  we arrive at 
\begin{equation}
\label{lejjj++}\begin{split}  \|P_{j}\psi\|_{S_{j}} &\lesssim  \|\psi_j(0)\|_{L^{2}} +  \|P_{j}f\|_{N_j}
+  \Mf_j(T) \sum_k \chi_{k=j} \| P_k \psi\|_{L^\infty_t L^2_r}. 
\end{split}\end{equation}
Here we sum up with respect to $j$. For $\Mf_j$ we know that
\[
\sum_j \Mf_j^2(T) \lesssim \Mf^2(T).
\]
Then, using Cauchy-Schwarz and a convolution bound for the last term on the right, we obtain
\[
 \|\psi\|_{\ell^1 S} \lesssim  \|\psi(0)\|_{\LX} + 
 \|f\|_{\ell^1 N} + \Mf(T)  \| \psi\|_{\ell^2 LE}  
\]
Finally, using \eqref{lefinal+} for the last term we obtain \eqref{lejfinal++p}.

\end{proof}

\section{Bounds for the nonlinearity \texorpdfstring{$N(\psi)$}{}}
\label{s:nolinear}

An important intermediate step in the proof of our results is to have good estimates on the 
nonlinear term $N(\psi)$ appearing in the nonlinear equation \eqref{psieq2} for $\psi$.
This is our task in this section. In a nutshell, the estimates in this section will allow us to treat $N(\psi)$ in a perturbative fashion in two places: in the nonlinear analysis of the Schr\"odinger evolution \eqref{psieq2} and in quantifying its effect in the dynamics of the modulation parameters $\lambda$ and $\alpha$.  To serve these multiple goals, we include  several estimates on $N(\psi)$ in this section as follows:
\begin{itemize}
\item Lemma \ref{l:Wlest1} provides a simple estimate on $N(\psi)$ which suffices to treat the effect of $N(\psi)$ in \eqref{psieq2} as perturbative when closing the theory in $l^2S$. 

\item Lemma \ref{gNref} refines the earlier result, showing that even when  $\psi \in l^2S$, the nonlinearity  $N(\psi)$ belongs to $l^1 N$; this indicates that the nonlinear effect on the dynamic variable $\psi$ is perturbative in the stronger structure $l^1 S$, and allows us to close the 
$l^1 S$ bounds for $\psi$ even for large $l^1 L^2=\LX$ 
data.

\item The other results in Lemma \ref{l:Wlest23} and Lemma \ref{PkNref} are needed in order to quantify the effect of $N(\psi)$
in the dynamics of the modulation parameters $\lambda$ and $\alpha$.
\end{itemize}

For context, we begin by recalling the nonlinear equation for $\psi$ in \eqref{psieq2},
\begin{equation}
(i \partial_t - \tilde H_\lambda) \psi = N(\psi),
\qquad \psi(0) = \psi_0, 
\label{wnlin-eq1}\end{equation}
where 
\[
 N(\psi) =  W_\lambda \psi, \qquad W_\lambda = A_0 - 2 \frac{\delta^\lambda A_2}{r^2}-
\frac{1}r \Im({\psi}_2 \bar{\psi})
\]
with $A_2$ and $\psi_2$   uniquely determined by $\psi$, see
Proposition \ref{psito2}, and $\delta^\lambda A_2=A_2 - 2h_3^\lambda$. $A_0$ is
given by \eqref{aoef} which we recall for convenience
\[
A_0(r) =  -\frac12 |\psi|^2  + \frac{1}r \Im (\psi_2 \bar \psi)
+ [r\partial_r]^{-1}( |\psi|^2 - \frac{2}{r} \Im (\psi_2 \bar \psi)).
\]
Thus we can rewrite
\[
W_\lambda= -\frac12 |\psi|^2  - 2 \frac{\delta^\lambda A_2}{r^2}
+ [r\partial_r]^{-1}( |\psi|^2 - \frac{2}{r} \Im (\psi_2 \bar \psi)).
\]

Our goal here is to estimate the nonlinear term $N(\psi)=W_\lambda \psi$ in the above equation, which we will do in two different ways. 

We make the following convention: in the statements below whenever a space-time is involved, the time interval is restricted to either $[0,T]$ for some $T>0$ or to $[0,+\infty)$; the estimates are uniform with respect to any such choice. We will not indicate this restriction at the level of notation, in order to keep notation compact.  

The first result provides an estimate for  $W_\lambda \psi$ in the admissible dual Strichartz space $L^\frac43$, and will help us to interpret this term perturbatively when studying the evolution of $\psi$ in $l^2S$.

\begin{l1}\label{l:Wlest1}
The nonlinearity $N(\psi)$ satisfies
\begin{equation} \label{Wlest1}
    \| W_\lambda \psi \|_{L^\frac43_{t,r}} \les \|\psi\|_{L^4_{t,r}}
    (\| \frac{\psi}r \|_{L^2_{t,r}} + \|\psi\|_{L_{t,r}^4}^2).
\end{equation}
\end{l1}

\begin{proof} To keep the arguments below compact we use the notation $L^p:=L^p_{t,r}$.
The non-integral terms are estimated as follows:
\[
\|  -\frac12 |\psi|^2 \psi \|_{L^\frac43} \les \| \psi \|^3_{L^4}, \quad \| \frac{\delta^\lambda A_2}{r^2} \psi \|_{L^\frac43} \les \| \frac{\delta^\lambda A_2}{r^2} \|_{L^2} \| \psi \|_{L^4} \les  
\| \frac{\psi}r \|_{L^2} \|\psi\|_{L^4},
\]
where in the last inequality we used \eqref{psiLE22}. The integral terms are estimated using \eqref{rdrm} and \eqref{psiL44} as follows:
\[
\|  \psi \cdot  [r\partial_r]^{-1}( |\psi|^2) \|_{L^\frac43} \les \| \psi \|_{L^4} \|  [r\partial_r]^{-1}( |\psi|^2) \|_{L^2} 
\les \| \psi \|_{L^4} \| |\psi|^2 \|_{L^2} \les \| \psi \|^3_{L^4},
\]
\[
\|  \psi \cdot  [r\partial_r]^{-1}( \frac{\delta \psi_2}r \bar \psi) \|_{L^\frac43} 
\les \| \psi \|_{L^4} \|  [r\partial_r]^{-1}( \frac{\delta \psi_2}r \bar \psi) \|_{L^2} 
\les \| \psi \|_{L^4} \|  \frac{\delta \psi_2}r \bar \psi \|_{L^2} 
\les \| \psi \|_{L^4}^2  \|  \frac{\delta \psi_2}r  \|_{L^4} \les \| \psi \|_{L^4}^3,
\]
\[
\|  \psi \cdot  [r\partial_r]^{-1}( \frac{h_1^\lambda}r \bar \psi) \|_{L^\frac43} 
\les \| \psi \|_{L^4} \|  [r\partial_r]^{-1}( h_1^\lambda \frac{\bar \psi}r ) \|_{L^2} 
\les \| \psi \|_{L^4} \|  h_1^\lambda \frac{\bar \psi}r  \|_{L^2} 
\les \| \psi \|_{L^4}  \|  \frac{ \psi}r  \|_{L^2}.
\]
Adding all of the above estimates gives us \eqref{Wlest1}.

\end{proof}

The next estimate for $W_\lambda \psi$  will be useful
when considering its indirect contribution
to the modulation equations, and show that it yields
$L^1_t$ source terms there.

\begin{l1}\label{l:Wlest23}
The following bounds hold true:
\begin{equation} \label{Wlest3}
    \| W_\lambda \|_{L^2_{t,r}} \les \| \frac{\psi}r \|_{L^2_{t,r}} + \|\psi\|_{L^4_{t,r}}^2;
\end{equation}
\begin{equation} \label{Wlest2}
    \| \frac{W_\lambda \psi}r \|_{L^1_{t,r}} \les 
    \| \frac{\psi}r \|_{L^2_{t,r}} (\| \frac{\psi}r \|_{L^2_{t,r}} + \|\psi\|_{L^4_{t,r}}^2).
\end{equation}
\end{l1}

\begin{proof}
The estimate \eqref{Wlest2} is a direct consequence 
of \eqref{Wlest3}, so we focus on the latter.
The proof is fairly simple and it follows the same approach as the previous one, just that it relies more on the energy decay estimate. Indeed, the non-integral components are estimated as follows, using \eqref{psiLE22}:
\[
\|  |\psi|^2 \|_{L^2} \les \| \psi \|^2_{L^4}, \quad \| \frac{\delta^\lambda A_2}{r^2} \|_{L^2} \les \| \frac{\psi}r \|_{L^2},
\]
while the integral terms are estimated, using \eqref{rdrm} and \eqref{psiL44}, as follows
\[
\|  [r\partial_r]^{-1}( |\psi|^2) \|_{L^2} \les \| \psi \|^2_{L^4},
\quad \|  [r\partial_r]^{-1}( \frac{\delta \psi_2}r \bar \psi) \|_{L^2} 
 \les \|  \frac{\psi}r \|_{L^2} \| \psi \|^2_{L^4},
\]
\[
 \|  [r\partial_r]^{-1}( h_1^\lambda \frac{\bar \psi}r ) \|_{L^2} 
\les   \|  \frac{ \psi}r  \|_{L^2}.
\]
Obviously these estimates suffice to conclude with \eqref{Wlest3} and \eqref{Wlest2} follows from it. 

\end{proof}

The following result two results provide frequency localized refinement estimates for $N(\psi)$. The first one highlights an improvement for the nonlinearity $N(\psi)$ in that it gains the $l^1$ summability even in the context of the problem with data with $l^2$ structure; the precise statement is as follows.

\begin{l1} \label{gNref}
The following holds true: 
\begin{equation} \label{Nl1N}
\| N(\psi) \|_{l^1 N} \les \| \psi \|_{l^2S}^2. 
\end{equation}
\end{l1}

\begin{proof} It is obvious that \eqref{Nl1N} follows from the following frequency localized version:
\begin{equation} \label{PjN}
\|P_k^\lambda N(\psi) \|_{L^\frac43_j} \les \sum_{k_1,k_3} \tilde \chi_{k_1=k}  \tilde \chi_{k_3=k}
\| P_{k_1}^\lambda \psi \|_{S_{k_1}}
\| P_{k_3}^\lambda \psi \|_{S_{k_3}},
\end{equation}
where $\tilde \chi_{i=j} = 2^{-\frac{|i-j|}{10}}$. 
The rest of the proof is concerned with \eqref{PjN}. Recall that 
\[
N(\psi)= W_\lambda \cdot \psi, \quad W_\lambda= -\frac12 |\psi|^2  - 2 \frac{\delta^\lambda A_2}{r^2}+ [r\partial_r]^{-1}( |\psi|^2 - \frac{2}{r} \Im (\psi_2 \bar \psi)).
\]
We start with the cubic term. It follows directly from the definition of the spaces and the characterization of the kernel of $P_k^\lambda$ in \eqref{ker1} that
\[
\|P_k^\lambda (P_{k_1}^\lambda \psi P_{k_2}^\lambda \overline{\psi} P_{k_3}^\lambda \psi) \|_{L^\frac43_k} \les \tilde \chi_{k_1=k} \tilde \chi_{k_2=k} \tilde \chi_{k_3=k}
\| P_{k_1}^\lambda \psi \|_{L^4_{k_1}} \| P_{k_2}^\lambda \psi \|_{L^4_{k_2}}
\| P_{k_3}^\lambda \psi \|_{L^4_{k_3}}.
\]
 Next we use \eqref{rdrmw} with the choice $w(r)= (2^{k_1} r)^{-\frac19}, r \leq 2^{-k_1}$ and $(2^{k_1} r)^{\frac19}, r \geq 2^{-k_1}$ to obtain
\[
\| w [r\partial_r]^{-1}( P_{k_1}^\lambda \psi \cdot \overline \psi)\|_{L^2} \leq \| P_{k_1}^\lambda \psi\|_{L^4_{k_1}} \| \psi\|_{L^4}. 
\]
From this, the definition of the spaces and the characterization of the kernel of $P_k^\lambda$ in \eqref{ker1} it follows that
\[
\| P_k^\lambda ([r\partial_r]^{-1}( P_{k_1}^\lambda \psi \cdot \overline \psi) \cdot P_{k_3}^\lambda \psi )\|_{L^\frac43_k} \les \tilde \chi_{k_1=k}  \tilde \chi_{k_3=k}
\| P_{k_1}^\lambda \psi \|_{L^4_{k_1}} \| \psi \|_{L^4}
\| P_{k_3}^\lambda \psi \|_{L^4_{k_3}}. 
\]
From the definition of the space $LE_k$ and \eqref{point-le}, we obtain
\[
\|\frac{P_j^\lambda \psi}r\|_{L^2_t L^2(A_m)} \les 2^{-\frac{k+m}2} \| P_j^\lambda \psi \|_{LE_j}.
\]
From this it follows that with $w_1(r)= (2^{j} r)^{-\frac13}, r \leq 2^{-j}$ and $(2^{j} r)^{\frac19}, r \geq 2^{-j}$ we have
\[
\| w_1 \frac{P_j^\lambda \psi}r \|_{L^2_{t,r}} \les \| P_j^\lambda \psi \|_{LE_j}.
\]
Invoking \eqref{rdrmw} (and the trivial fact that $\| \psi_2 \|_{L^\infty} \leq 1$) gives us
\begin{equation} \label{LEjw}
\| w_1 [r\partial_r]^{-1}( \overline \psi_2 \frac{P_j^\lambda \psi}r) \|_{L^2_{t,r}} \les \| P_j^\lambda \psi \|_{LE_j}.
\end{equation}
Using this, the definition of the spaces and the characterization of the kernel of $P_k^\lambda$ in \eqref{ker1}, we obtain the following estimate
\[
\|P_k^\lambda( [r\partial_r]^{-1}( \overline \psi_2 \frac{P_{k_1}^\lambda \psi}r) P_{k_3}^\lambda \psi) \|_{L^\frac43_k} \les \tilde \chi_{k_1=k}  \tilde \chi_{k_3=k}
\| P_{k_1}^\lambda \psi \|_{LE_{k_1}} 
\| P_{k_3}^\lambda \psi \|_{L^4_{k_3}}. 
\]
We note that in all terms estimated above we took advantage of the fact that they contain (at least) two factors of $\psi$ which are easily amenable to frequency localization since it was important that we could use localization for at least two terms. This makes the  last term in the nonlinearity $\frac{\delta^\lambda A_2}{r^2} \cdot \psi$ more challenging; to be more precise, the challenging aspect is to bring some sort of frequency localization information on the term $\delta^\lambda A_2$. The transition from $\psi$ to $\delta^{\lambda,\alpha} \psi_2$ and $\delta^\lambda A_2$ is governed by the system \eqref{sysAp2} (which was studied in detail in the proof of Lemma \ref{psito2}); for convenience we recall the system here
\begin{equation}
\left\{ \begin{array}{l}
L_\lambda \delta^{\lambda,\alpha} \psi_2  = 2i h_3^\lambda \psi + \delta^\lambda A_2 \psi - \frac{1}r \delta^\lambda A_2 (2ie^{2 i \alpha}h_1^\lambda + \delta^{\lambda,\alpha} \psi_2), \cr \cr
L_\lambda \delta^\lambda A_2  = - 2 h_1^\lambda \Re{(e^{2i\alpha} \psi)} + \Im(\psi \overline{\delta^{\lambda,\alpha} \psi}_2) - \frac{(\delta^\lambda A_2)^2}r.
\end{array} \right.
\end{equation}
The idea is to localize $\psi$ in frequency and solve for the corresponding $\delta^{\lambda,\alpha} \psi_2$ and $\delta^\lambda A_2$; however the system is nonlinear, and this requires some care. A simple alternative is to force 
linearity in this system by letting $\delta_k^{\lambda,\alpha} \psi_2$ and $\delta^\lambda_k A_2$ be the solutions of the following system
\begin{equation}
\left\{ \begin{array}{l}
L_\lambda \delta_k^{\lambda,\alpha} \psi_2  = 2i h_3^\lambda P_k^\lambda \psi + \delta^\lambda A_2 P_k^\lambda \psi - \frac{1}r \delta^\lambda_k A_2 (2ie^{2 i \alpha}h_1^\lambda + \delta^{\lambda,\alpha} \psi_2), \cr \cr
L_\lambda \delta_k^\lambda A_2  = - 2 h_1^\lambda \Re{(e^{2i\alpha} P_k^\lambda \psi)} + \Im( P_k^\lambda \psi \overline{\delta^{\lambda,\alpha} \psi}_2) - \frac{\delta^\lambda A_2 \cdot \delta_k^\lambda A_2 }r,
\end{array} \right.
\end{equation}
with the same set of initial data as in Lemma \ref{psito2}
\[
\delta_k^{\lambda,\alpha} \psi_2(r_0)=0, \qquad \delta_k^{\lambda,\alpha} A_2(r_0)=0.
\]
Here we recall that $r_0 \approx \lambda^{-1}$. It is clear that this provides a linear decomposition of $\delta^{\lambda,\alpha} \psi_2$ and $\delta^\lambda A_2$ as follows
\[
\delta^{\lambda,\alpha} \psi_2 = \sum_{k \in \Z} \delta_k^{\lambda,\alpha} \psi_2, \qquad 
\delta^\lambda A_2 = \sum_{k \in \Z} \delta_k^{\lambda,\alpha} A_2. 
\]
Given that we already have the apriori knowledge $\| \delta^{\lambda,\alpha} \psi_2 \|_{L^\infty} + \| \delta^\lambda A_2  \|_{L^\infty} \les \| \psi\|_{L^2} \ll 1$ from Lemma \ref{psito2}, it suffices to analyse only the equation for $\delta^\lambda_k A_2$. This  can be rewritten in the integral form:
\[
 \delta^\lambda_k A_2 (r)=  h^\lambda_1(r) \int_{r_0}^r \frac1{h_1^\lambda}\left( - 2 h^\lambda_1 \Re{ (e^{2i\alpha} P_k^\lambda \psi)} + \Im(P_k^\lambda \psi \overline{\delta^{\lambda,\alpha} \psi}_2) - \frac{\delta^\lambda A_2 \cdot \delta_k^\lambda A_2}s \right) ds. 
\]
We record the following simple variant of \eqref{Linvest}:
\begin{equation} \label{Linvest2}
\| \frac{w}r h_1^\lambda \int_{r_0}^r \frac1{h_1^\lambda} f ds \|_{L^p} \les \|wf\|_{L^p}, \forall 1 \leq p \leq \infty,
\end{equation}
under the following assumptions on $w$: 

- $w: \R_{+} \rightarrow \R_{+}$ is slowly varying, that is $w(r_1) \approx w(r_2)$ for $2^{-1} r_2 \leq r_1 \leq 2r_2$;

- $r^{-1} w(r)$ is decreasing and $r^\frac12 w(r)$ is increasing.  

The second condition implies the slowly varying one, but we have stated them separately in order to also emphasize the first. 
The only role that $r_0$ plays in the above inequality is that is sits at the "height" of $h_1^\lambda$, that is $h_1^\lambda(r_0) \approx 1$ and $h_1^\lambda $ decays away from a neighborhood of size $\lambda^{-1}$ of $r_0$. 

Next we use \eqref{Linvest2} with the choice $w_{k_1}(r)= (2^{k_1} r)^{-\frac19}, r \leq 2^{-k_1}$ and $(2^{k_1} r)^{\frac19}, r \geq 2^{-k_1}$ to obtain
\[
\| \frac{w_{k_1}}r h^\lambda_1(r) \int_{r_0}^r \frac1{h_1^\lambda}\left( - 2 h^\lambda_1 \Re{ (e^{2i\alpha} P_{k_1}^\lambda \psi)} + \Im(P_{k_1}^\lambda \psi \overline{\delta^{\lambda,\alpha} \psi}_2) \right) ds \|_{L^4} \lesssim
\| w_{k_1} P_{k_1}^\lambda \psi\|_{L^4}
\les \| P_{k_1}^\lambda \psi\|_{L^4_{k_1}}. 
\]
This suggests that, at least at the linear level, the following holds true
\begin{equation} \label{deltaA2k}
\| w_{k_1} \frac{\delta^\lambda_{k_1} A_2}r \|_{L^4} \les \| P_{k_1}^\lambda \psi\|_{L^4_{k_1}}.
\end{equation}
The quadratic contributions (coming from the term $\frac{\delta^\lambda A_2 \cdot \delta_{k_1}^\lambda A_2}s$) are then estimated using the apriori assumption that $\|\delta^\lambda A_2\|_{L^\infty} \les \|\psi\|_{L^2} \ll 1$; the actual justification for \eqref{deltaA2k} can be done using a continuity/bootstrap argument on intervals $[r_0,R]$ and $[r,r_0]$ for any $r \leq r_0 \leq R$. 

Using \eqref{deltaA2k} and \eqref{LEjw}, as well as the definition of $L^\frac43_k$ and the characterizayion of the lernel of $P_k^\lambda$ in \eqref{ker1}, allows us to estimate
\[
\| P_k^\lambda (\frac{\delta^\lambda_{k_1} A_2}{r} \cdot \frac{P_{k_3} \psi}r) \|_{L^\frac43_k} \les 
\tilde \chi_{k_1=k}  \tilde \chi_{k_3=k}
\| P_{k_1}^\lambda \psi \|_{L^4_{k_1}} 
\| P_{k_3}^\lambda \psi \|_{LE_{k_3}}. 
\]
This finishes the proof of \eqref{PjN} and in turn, the proof of the Lemma.
\end{proof}

For technical reasons we also need the following result (which will be used in the proof of Lemma~\ref{l:PG}, precisely in estimating the term $e_k^2$ in Step 2).

\begin{l1} \label{PkNref}
The following holds true:
\begin{equation}\label{Npsi-21}
\| \frac{P_k N(\psi)}r\|_{L^2_t L^1_r} \lesssim   2^k \sum_{k'} 2^{-\frac{|k-k'|}{10}} \| P_{k'}^\lambda \psi \|_{S_{k'}}.
\end{equation}

\end{l1}

\begin{proof}
It is obvious that it suffices to establish the following estimate
\begin{equation} \label{Npsil2l1}
\| \frac{P_k^\lambda (W_\lambda \cdot P_{k'}^\lambda \psi)}{r}\|_{L^2_t L^1_r} \les 2^k 2^{-\frac{|k-k'|}{10}} \| P_{k'}^\lambda \psi \|_{S_{k'}}
\|\psi\|_{S}(1+ \|\psi\|_{S}).
\end{equation}
In order to prove \eqref{Npsil2l1}, we need the following technical result: 
\begin{equation} \label{Pkfl2l1}
\| \frac{P_k^\lambda f}r \|_{L^2_t L^1_r} \les 2^k \sum_{m \in \Z} 2^{-\frac{|m+k|}9} \|f \|_{L^2_t L^1_r(A_m)}. 
\end{equation}
The argument is fairly straightforward; we start with 
\[
\begin{split}
\| \frac{P_k^\lambda f}r \|_{L^2_t L^1_r (A_n)} & \approx 2^{-n} \| P_k^\lambda f  \|_{L^2_t L^1_r (A_n)} \\
& \les 2^{-n} \sum_{m \in \Z} \| \int \tilde K^\lambda_k(\cdot,s) \chi_{A_m}(s) f(s) sds  \|_{L^2_t L^1_r (A_n)} \\
& \les_N 2^{-n} \sum_{m \in \Z} \omega_{k,\lambda}(2^n)   \| \int \frac{2^{2k} \omega_{k,\lambda}(s) \chi_{A_m}(s) f(s)}{(1+2^k(2^m+2^n))(1+2^k|r-s|)^{-N}}  sds  \|_{L^2_t L^1_r (A_n)},
\end{split}
\]
where in passing to the last line we have used the bound \eqref{ker1} on the kernel $\tilde K^\lambda_k$; in what follows below we heavily rely on the bounds on $\omega_{k,\lambda}$ from \eqref{omdef}. 

If $n+k \leq 0$, then we can further bound the above quantity by
\[
\begin{split}
& 2^{n+2k} \sum_{m \leq -k} \omega_{k,\lambda}(2^n) \omega_{k,\lambda}(2^m)  \| f\|_{L^2_t L^1_r (A_m)} + 
2^{n+2k} \sum_{m \geq -k} \omega_{k,\lambda}(2^n)  2^{-N(k+m)} \| f\|_{L^2_t L^1_r (A_m)} \\
\les & 2^{n+2k} \sum_{m \in \Z} 2^{-|k+m|} \| f\|_{L^2_t L^1_r (A_m)}.
\end{split}
\]
If $n+k \geq 0$, then we can further bound the above quantity by
\[
\begin{split}
& 2^{n+2k} \sum_{m \leq -k}  \omega_{k,\lambda}(2^m) 2^{-N(k+n)} \| f\|_{L^2_t L^1_r (A_m)} +  2^{k} \sum_{m \geq -k} (2^{k+m}+2^{k+n})^{-1} \| f\|_{L^2_t L^1_r (A_m)} \\
\les & 2^{k} \sum_{m \in \Z} 2^{-\frac{|k+m|+|k+n|}2} \| f\|_{L^2_t L^1_r (A_m)}.
\end{split}
\]
Wrapping things up, we have succeeded to prove 
\[
\| \frac{P_k^\lambda f}r \|_{L^2_t L^1_r (A_n)} \les 2^{k} \sum_{m \in \Z} 2^{-\frac{|k+m|+|k+n|}2} \| f\|_{L^2_t L^1_r (A_m)},
\]
from which \eqref{Pkfl2l1} follows. 

Now we proceed with the proof of \eqref{Npsil2l1}. For convenience recall the formula for $W_\lambda$,
\[
 W_\lambda= -\frac12 |\psi|^2  - 2 \frac{\delta^\lambda A_2}{r^2}+ [r\partial_r]^{-1}( |\psi|^2 - \frac{2}{r} \Im (\psi_2 \bar \psi)).
\]
For the cubic component we note the following estimate
\[
\sup_{m \in \Z} 2^{\frac{|m+k'|}8} \| |\psi|^2 P_{k'}^\lambda \psi \|_{L^2_t L^1_r(A_m)} \les \| P_{k'}^\lambda \psi \|_{L^4_{k'}} \|\psi\|_{L^4} \| \psi \|_{L^\infty_t L^2_r},
\]
followed by
\[
\sum_{m \in \Z} 2^{\frac{|m+k'|}9} \| |\psi|^2 P_{k'}^\lambda \psi \|_{L^2_t L^1_r(A_m)} \les \| P_{k'}^\lambda \psi \|_{L^4_{k'}} \|\psi\|_{L^4} \| \psi \|_{L^\infty_t L^2_r}. 
\]
At this point we invoke \eqref{Pkfl2l1} to conclude that 
\[
\| \frac{P_k^\lambda(|\psi|^2 P_{k'}^\lambda \psi)}r \|_{L^2_t L^1_r} \les 2^k 2^{-\frac{|k-k'|}{10}} \| P_{k'}^\lambda \psi \|_{S_{k'}} \|\psi\|_{S}^2,
\]
which implies the claim made in \eqref{Npsil2l1} for the corresponding cubic component in $N(\psi)$. In order to extend the arguments above to the other components of $N(\psi)$, it is important to note that for the quadratic term $|\psi|^2$ all that we have used was  an estimate in $L^4_t L^\frac43_r$ (by invoking
the $L^4_{t,r}$ information on one of the $\psi$'s and $L^\infty_t L^2_r$ on the other $\psi$). 

Now we turn to the other components in $W_\lambda$. We first deal with the component which contains the operator $[r\partial_r]^{-1}$; from \eqref{rdrm} we know that the operator $[r\partial_r]^{-1}$ is bounded on $L^\frac43_r$, and, as a consequence, it is bounded on $L^4_t L^\frac43_r$ - as pointed earlier this is all that was needed for the term $|\psi|^2$ in the analysis of the cubic term above. From these considerations it follows that the term $P_{k'}^\lambda \psi \cdot [r\partial_r]^{-1} (|\psi|^2)$ can be treated in a similar fashion to the cubic one above. 

In analyzing the term $[r\partial_r]^{-1}(\frac{2}{r} \Im (\psi_2 \bar \psi))$ we write $\psi_2 = \delta^{\lambda,\alpha} \psi_2  + 2 e^{2i\alpha} h_1^\lambda$.  
The term $P_{k'}^\lambda \psi \cdot  [r\partial_r]^{-1} (\frac{\delta^{\lambda,\alpha} \psi_2}r \overline{\psi})$ is treated as above (similar to the cubic component) by using the $L^\infty_t L^2_r$ bound on $\frac{\delta^{\lambda,\alpha} \psi_2}r$ from \eqref{compsolh1-r} and the mapping properties of the operator $[r\partial_r]^{-1}$; the term $P_{k'}^\lambda \psi \cdot  [r\partial_r]^{-1} (\frac{h_1^\lambda}r \overline{\psi})$ is entirely similar since $\| \frac{h_1^\lambda}r \|_{L^2_r}= \| \frac{h_1}r \|_{L^2_r} \les 1$ and this bound is independent of time (thus uniform in time).  

From the above we conclude with the following estimate
\[
\| \frac{P_k^\lambda( P_{k'}^\lambda \psi \cdot [r\partial_r]^{-1}( |\psi|^2 - \frac{2}{r} \Im (\psi_2 \bar \psi)))}r \|_{L^2_t L^1_r} \les 2^k 2^{-\frac{|k-k'|}{10}} \| P_{k'}^\lambda \psi \|_{S_{k'}} \|\psi\|_{S}(1+ \|\psi\|_{S}),
\]
Finally we deal with the $\frac{\delta^\lambda A_2}{r^2}$ component in $W_\lambda$.
From \eqref{high-k} and \eqref{low-k}, it follows that
\[
\sum_{m \in \Z} 2^{\frac{|k'+m|}4}\| \frac{P_{k'}^\lambda \psi}r \|_{L^2_{t,r}(A_m)} \les \| P_{k'}^\lambda \psi\|_{LE_{k'}}. 
\]
From this and \eqref{compsolh1-r} we obtain that
\[
\sum_{m \in \Z} 2^{\frac{|k'+m|}4}\| \frac{\delta^\lambda A_2}r \cdot \frac{P_{k'}^\lambda \psi}r \|_{L^2_t L^1_r (A_m)} \les \| \frac{\delta^\lambda A_2}r \|_{L^\infty_t L^2_r} \| P_{k'}^\lambda \psi\|_{LE_{k'}} 
\les \| \psi \|_{L^\infty_t L^2_r} \| P_{k'}^\lambda \psi\|_{LE_{k'}}.  
\]
Combining this with \eqref{Pkfl2l1} gives the following:
\[
\| \frac{P_k^\lambda( P_{k'}^\lambda \psi \cdot \frac{\delta^\lambda A_2}{r^2})}r \|_{L^2_t L^1_r} \les 2^k 2^{-\frac{|k-k'|}{10}} \| P_{k'}^\lambda \psi \|_{S_{k'}} \|\psi\|_{S}(1+ \|\psi\|_{S}).
\]
This finishes the proof of \eqref{Npsil2l1}, and in turn of our Lemma. 
\end{proof}

\section{The modulation parameter dynamics - Part 1}
\label{s:mod}

The main goal of this section is to prove Theorems \ref{tmain-L} and \ref{tmain-L1}. Those theorems are mostly stated at the level of the map $u$, although the field $\psi$ appears in their conclusions, see \eqref{psiest-m} and \eqref{psiest-m1}. Below we restate those Theorems fully in terms of $\psi$, see Theorems \ref{PD-1} and \ref{PD-2}, and prove  the latter ones; Theorem \ref{tmain-L} and Theorem \ref{PD-1} are identical, while the equivalence between Theorem \ref{tmain-L1} and Theorem \ref{PD-2} follows directly from Proposition~\ref{uQpsi}.

As described in Section \ref{s:modulation-param}, $\dot H^1$ states for the Schr\"odinger Map flow are described via the differentiated $L^2$ field $\psi$ and the modulation parameters $\alpha$ and $\lambda$, which are chosen via the orthogonality condition \eqref{lell} (in conjunction with \eqref{lellf} in order to ensure uniqueness).

Dynamically, we expect the Schr\"odinger maps evolution
to be governed by a coupled system consisting of a Schr\"odinger type evolution for $\psi$ coupled with 
appropriate modulation equations for $\alpha$ and $\lambda$. The equation for $\psi$ has already been derived in  \eqref{psieq2}, and contains two types of nonlinear effects due to (i) the nonlinear terms on the 
right hand side of \eqref{psieq2}, and (ii) the time dependence of $\lambda$ in the left hand side of \eqref{psieq2}. The former are at least quadratic, and will play a perturbative role in our analysis. 
The effect of $\lambda$ is stronger, 
and was investigated in the last two sections; as seen there, the crucial control is that of $\|\frac{\lambda'}{\lambda^2}\|_{L^2}$.

In this section we begin investigating the other half of the story, namely the evolution of the modulation parameters $\alpha$ and $\lambda$. Our goal is twofold:

\begin{enumerate}
    \item We derive the modulation equation for $\lambda$ 
    and $\alpha$ as an ode system, separating the 
    contributions of $\psi$ into a leading order linear part and a perturbative quadratic part. In this analysis the essential control is that of the 
    local energy decay of $\psi$, primarily  $\|\frac{\psi}r\|_{L^2_{t,r}}$.

   \item We combine the bounds for the $\psi$ equation
   with the bounds for the modulation equation in order to simultaneously close the local energy bounds for $\psi$
   and the modulation equation bound for $\|\frac{\lambda'}{\lambda^2}\|_{L^2}$.
\end{enumerate}

The second part above is carried out on the full time of existence of the solutions, regardless of whether this time is finite or infinite. The result is stated as follows:

\begin{t1} \label{PD-1} Assume that we have a $2$-equivariant initial data $u_0 \in \dot H^1$ in the homotopy class of $Q^2$, and with energy 
below $8\pi+\delta^2$, and with $\delta$ sufficiently small.
Let $I_{max}=(T_{min},T_{max})$ be the maximal time of existence of the solution $u$ to the Schr\"odinger map flow \eqref{SM} with initial data $u_0$. 
Then the associated Coulomb gauge field $\psi$ and the associated parameter $\lambda$ (by the rule \eqref{lell}) satisfy
\begin{equation} \label{l1est}
    \|\frac{\lambda'}{\lambda^2} \|_{L^2_t(I_{max})}
    +  \|\frac{\alpha'}{\lambda} \|_{L^2_t(I_{max})}
    \les \| \psi(0) \|_{L^2},
\end{equation}
 and
\begin{equation} \label{psiest}
    \| \psi \|_{l^2 S(I_{max})} \les \| \psi(0) \|_{L^2}. 
\end{equation}
\end{t1}

Later on it will also be useful to have an $l^1$ version of the above theorem, 
more precisely of the estimate \eqref{psiest}:

\begin{t1}\label{PD-2}
Let $u$ be the solution $u$ to the Schr\"odinger map flow \eqref{SM} as in Theorem~\ref{PD-1}. Assume in addition that $\psi(0) \in \LX$. Then we have
\begin{equation} \label{psiest-l1}
    \| \psi \|_{l^1 S(I_{max})} \les \| \psi(0) \|_{\LX}. 
\end{equation}
\end{t1}
 A density argument using the local well-posedness result in Theorem~\ref{th2} shows that it suffices to prove   these results for $\dot H^1 \cap \dot H^2$ solutions. The forward in time and backward in time problems are similar, so for simplicity we drop the superscripts and we will work only on $[0,T_{max})$ in what follows.

This theorem provides two important pieces of information. Concerning the field $\psi$, it asserts that its dispersive properties (measured by local energy decay and Strichartz norms) persist on the maximal time of existence; in other words potential blow-up in finite or infinite time does not interfere with the dispersion. 

Concerning the parameter $\lambda$, the estimate \eqref{l1est}  prevents the scenario that $\lambda \rightarrow 0$ in finite time (this can also be seen as a consequence of the local result in Theorem~\ref{th2}).
However, it  does not prevent the scenarios 
$\lambda \rightarrow 0$ in infinite time, or
$\lambda \rightarrow \infty$ in finite or infinite time. In the next section we will refine our analysis of the ODE system so that we can investigate these other potential scenarios.

\subsection{The modulation equation}
In this section we derive the modulation equation for $\alpha, \lambda$ from the orthogonality condition \eqref{ldef}. From Corollary \eqref{cdif} we already know that $\lambda$ and $\alpha$ are differentiable in time on the maximal interval of existence.  Hence we can differentiate \eqref{ldef} with respect to time. The computations below are formally justified for more regular solutions as in Theorem~
\ref{th2}; such solutions satisfy in particular $\psi \in L^\infty H^1$. However, the final outcome, namely the modulation equations, are valid 
for all finite energy solutions via a density argument.

We differentiate \eqref{ldef} with respect to time and use \eqref{derhf} to obtain the following
\[
\begin{split}
0  =    & \ 
 \partial_t \la \psi_2-2ie^{2 i \alpha(t)} h_1(\lambda r), \varkappa(\lambda r) \ra 
\\  
= & \ 
 \la \partial_t \psi_2, \varkappa(\lambda r) \ra + 4 \alpha' \la e^{ 2i \alpha(t)} h_1(\lambda r), \varkappa(\lambda r)\ra \\
& + 4 i \lambda' \la e^{2 i \alpha(t)} r \frac{h_1(\lambda r) h_3(\lambda r)}{\lambda r}, \varkappa(\lambda r) \ra
+ \lambda' \la \psi_2-2ie^{ 2 i \alpha(t)} h_1(\lambda r), r \varkappa'(\lambda r) \ra \\
= & \  \la \partial_t \psi_2, \varkappa(\lambda r) \ra + 4  e^{2 i \alpha(t)} \left( \frac{\alpha'}{\lambda^2} \la h_1, \chi \ra 
+ i \frac{\lambda'}{\lambda^3} \la h_1 h_3, \varkappa \ra \right) \\
& + \frac{\lambda'}{\lambda^3} \la \psi_2(\lambda^{-1} r)-2ie^{2 i \alpha(t)} h_1(r), r \varkappa'(r) \ra.
\end{split}
\]
A closer look at the last expression in the sequence of equalities reveals the following:

\begin{itemize}
\item  the second term  $4  e^{2 i \alpha(t)} \left( \frac{\alpha'}{\lambda^2} \la h_1, \varkappa \ra + i \frac{\lambda'}{\lambda^3} \la h_1 h_3, \varkappa \ra \right)$ contains the primary dynamical information about the two real parameters $\alpha$ and $\lambda$ due to the nondegeneracy conditions \eqref{chi-nondeg}
\begin{equation}
    c_1 := \la h_1, \varkappa \ra \neq 0, \qquad  c_2:=\la h_1 h_3, \varkappa \ra \neq 0;
\end{equation}
in fact it is precisely the above computation (and later consideration) that justfies imposing these specific non-degenaracy conditions;  

\item the last term $\frac{\lambda'}{\lambda^3} \la \psi_2(\lambda^{-1} r)-2ie^{2 i \alpha(t)} h_1(r), r \varkappa'(r) \ra$ also contains some dynamical information,
but in quadratic form, and will be shown to be negligible/perturbative;

\item  the first term $\la \partial_t \psi_2, \varkappa(\lambda r) \ra$ will provide the leading linear source term in the 
modulation equations. Understanding this require additional analysis, which we present below. 
\end{itemize}

Using the covariant rules of calculus detailed in Section~\ref{seccoulomb}, we have
\begin{equation}\label{dt-psi2}
\partial_t \psi_2 = D_0 \psi_2 - i A_0 \psi_2= D_2 \psi_0 - i A_0 \psi_2 = i A_2 \psi_0 - i A_0 \psi_2.
\end{equation}
Our goal is now to reexpress the right hand side in terms of $\psi$ and the quantities $\delta^{\lambda} A_2$ and $\delta^{\lambda,\alpha} \psi_2$,
which measure the deviation from the soliton manifold. For $\psi$, we will seek to write its radial derivative in terms of $L_\lambda^*\psi$; this is because the operator $L_\lambda^*$ acts as a derivative in the $\tilde H_\lambda$ calculus, and thus has a favourable behavior in the low frequency regime.

Using the computations in Section~\ref{seccoulomb}, we have
\[
\begin{split}
- i \psi_0 & = D_1 \psi_1 + \frac{1}{r} \psi_1 + \frac{1}{r^2} D_2 \psi_2 \\
& = \partial_1 (\psi + i \frac{\psi_2}r) + \frac{1}{r}(\psi + i \frac{\psi_2}r)  + i \frac{1}{r^2} A_2 \psi_2 \\
& = \partial_1 \psi + \frac{1}{r} \psi + i \frac{\partial_r \psi_2}{r} - i \frac{\psi_2}{r^2} + i \frac{\psi_2}{r^2} + i \frac{1}{r^2} A_2 \psi_2 \\
& = \partial_r \psi + \frac{1}{r} \psi + \frac{i}r (i A_2 \psi) \\
& = \partial_r \psi + \frac{1-A_2}{r} \psi \\
& =- L_\lambda^* \psi + \frac{2h_3^\lambda-A_2}{r} \psi.
\end{split}
\]
Therefore
\[
i A_2 \psi_0 = 2 h_3^\lambda L^*_\lambda \psi + (A_2 - 2h^\lambda_3) L_\lambda^* \psi -  A_2 \frac{2h_3^\lambda-A_2}{r} \psi = 2 h_3^\lambda L^*_\lambda \psi + \delta^\lambda A_2 L_\lambda^* \psi +  A_2 \frac{\delta^\lambda A_2}{r} \psi.
\]
We now consider the second term on the right in \eqref{dt-psi2}.
Here we do not directly use the derivation in \eqref{a0bis} for $A_0$, but instead we seek to highlight again the component $L_\lambda^* \psi$ in the expression that gives $A_0$. For this purpose we rearrange the formula for $A_0$ as follows (this being a consequence of \eqref{curb} and the gauge condition $A_1=0$):
\[
\begin{split}
-\partial_r A_0 & = \Im(\psi_0 \overline{\psi}_1) = - \Re (i \psi_0 \overline{\psi}_1) =  \Re \left((-L_\lambda^* \psi - \frac{\delta^\lambda A_2}{r} \psi)(\overline{\psi} - i \frac{\overline{\psi_2}}{r})\right) \\
& = \Re (i L_\lambda^* \psi \frac{\overline{\psi_2}}r ) - \Re (L_\lambda^* \psi \cdot \overline{\psi}) 
- \Re (\frac{\delta^\lambda A_2}{r} \psi \overline{\psi}_1) \\
& = \Re( \frac{2 e^{- 2 i \alpha(t)}h_1^{\lambda}}r L_\lambda^* \psi ) + \Re (i L_\lambda^* \psi \frac{\overline{\delta^{\lambda,\alpha} \psi_2}}r ) 
- \Re (L_\lambda^* \psi \cdot \overline{\psi}) - \Re (\frac{\delta^\lambda A_2}{r} \psi \overline{\psi}_1);
\end{split}
\]
the third equality is justified by the computation above for $-i\psi_0$. 

From this we obtain the following representation for $A_0$:
\begin{equation}\label{A0-rep}
A_0(r)=\int_r^\infty \left( \Re( \frac{2 e^{- 2 i \alpha(t)}h_1^{\lambda}}r L_\lambda^* \psi ) + \Re (i L_\lambda^* \psi \frac{\overline{\delta^{\lambda,\alpha} \psi_2}}r ) 
- \Re (L_\lambda^* \psi \cdot \overline{\psi}) - \Re (\frac{\delta^\lambda A_2}{r} \psi \overline{\psi}_1) \right) dr,
\end{equation}
where the first term is linear in $\psi$ and the remaining terms are at least quadratic.

Bringing all the above computations together into the previous ODE for $\lambda$ and $\alpha$ (derived at the beginning of this section) gives us a first form of the modulation equation,
\begin{equation} \label{laeq}
 4  e^{2 i \alpha(t)} \left( \frac{\alpha'}{\lambda^2} c_1 + i \frac{\lambda'}{\lambda^3} c_2 \right)  = Lin   + Nlin,
\end{equation}
where $Lin$ collects all the linear terms:
\[
Lin=  - \la  2 e^{2 i \alpha}h_1^{\lambda} \int_r^\infty \Re( \frac{2 e^{- 2i \alpha}h_1^{\lambda}}r L_\lambda^* \psi ) dr, \varkappa^\lambda \ra  - \la 2h_3^\lambda L_\lambda^* \psi,  \varkappa^\lambda \ra,
\]
and $Nlin$ collects all nonlinear (quadratic or higher order) terms
\[
\begin{split}
Nlin & = - \frac{\lambda'}{\lambda^3} \la \psi_2(\lambda^{-1} r)-2ie^{2 i \alpha(t)} h_1(r), r \varkappa'(r) \ra -
\la \delta^\lambda A_2 \cdot L^*_\lambda \psi +  A_2 \frac{\delta^\lambda A_2}{r} \psi, \varkappa^\lambda  \ra \\
& - 2 e^{2 i \alpha(t)} \la \int_r^\infty \Re (i L^*_\lambda \psi \frac{\overline{\delta^{\lambda,\alpha}\psi_2}}{r}) - \Re (L^*_\lambda \psi \overline{\psi}) - \Re (\frac{\delta^\lambda A_2}{r} \psi \overline{\psi}_1) dr,  (h_1\varkappa)^\lambda \ra  \\
& + \la i A_0 \delta^{\lambda,\alpha} \psi_2, \varkappa^\lambda \ra.
\end{split}
\]
We have obtained the ODE system \eqref{laeq} that governs the evolution of the two parameters $\alpha$ and $\lambda$; this is a first version from which we derive an alternative version which is more amenable to estimates. We will prove that the contributions in $Nlin$ to the dynamics of $\lambda$ and $\alpha$ are perturbative. However the first term in $Nlin$ is special in that it contains the expression 
$\frac{\lambda'}{\lambda^3}$; this is why this term is treated slightly differently, and justifies the following separation: 
\[
Nlin = - \frac{\lambda'}{\lambda^3} \la \psi_2(\lambda^{-1} r)-2ie^{2 i \alpha(t)} h_1(r), r \varkappa'(r) \ra + \widetilde{Nlin}:= 
- \frac{\lambda'}{\lambda^3}  e + \widetilde{Nlin}.
\]
It will soon become clear that in \eqref{laeq}, the variable $\lambda$ is not scaled properly. For this reason we rewrite it as follows
\begin{equation} \label{laeq22}
 4  e^{2 i \alpha(t)} \left( \alpha'  c_1 + i \frac{\lambda'}{\lambda} c_2  \right)=  l(t) + q(t). 
\end{equation}
where 
\[
l(t)= \lambda^2  Lin, \quad q(t)=\lambda^2 Nlin.
\]
This is the main formulation of the ODE describing the parameter dynamics. Similar to $Nlin$, we split the $q$ terms as follows
\[
q= - \frac{\lambda'}{\lambda}  e + \tilde q,
\]
where $\tilde q(t)=\lambda^2 \widetilde{Nlin}$. 

We note that the form \eqref{laeq22} 
of the modulation equations will be used in Section \ref{s:mod2} in order to provide the most accurate dynamical information on $\alpha$ and $\ln \lambda$. However for the main result of this section, namely Theorem \ref{PD-1}, we will use a slight modification of it.

\subsection{Estimates for the modulation equations}
In this subsection we provide the main bounds for the source terms in the modulation equation,
which will be used in the proof of  Theorem \ref{PD-1}. 
For this purpose, we rewrite \eqref{laeq22} as follows:
\begin{equation} \label{laeq2}
 4  e^{2 i \alpha(t)} \left( c_1 \frac{\alpha'}{\lambda}   + i c_2 \frac{\lambda'}{\lambda^2}  + e^{-2 i \alpha(t)} \frac{\lambda'}{\lambda^2} e(t) \right)= \lambda^{-1} l(t) + \lambda^{-1} \tilde q(t). 
\end{equation}
We bound the entries $e,l,\tilde q$ in this equation as follows:
\begin{l1}\label{l:hqest}
  The functions $e,l,\tilde q$ in \eqref{laeq2} 
satisfy the following fixed time estimates:
\begin{equation} \label{hqest}
\begin{split}
    & |e(t)| \les \| \psi(t) \|_{L^2}, \quad
\lambda |e(t)| 
    \les \| \frac{\psi(t)}{r}\|_{L^2_r}, \quad
    \lambda^{-1} |l(t)| \les \| \frac{\psi(t)}{r}\|_{L^2_r}, \\
    & |\tilde q(t)| \les \| \frac{\psi(t)}{r}\|_{L^2_r}^2, \quad
    \lambda^{-1} |\tilde q(t)| \les \|\psi(t)\|_{L^2_r} \| \frac{\psi(t)}{r}\|_{L^2_r}. 
\end{split}
\end{equation}
\end{l1}

\begin{proof}[Proof of Lemma~\ref{l:hqest}]
We successively consider the bounds for $e$, $l$ and $\tilde q$.

\medskip

\emph{ a) The bounds for $e$.}
 From \eqref{compsolh1-r} it follows that  $\| \psi_2-2ie^{2 i \alpha(t)} h_1^\lambda \|_{L^\infty} \les \|\psi\|_{L^2}$,
 which directly implies that $|e(t)| \les \| \psi(t) \|_{L^2} \ll 1$. On the other hand we can bound  $e$ using the local energy norm of $\psi$, 
\[
\begin{split} 
 \lambda e(t) & =  \lambda \la \psi_2(\lambda^{-1} r)-2ie^{2 i \alpha(t)} h_1(r), r \varkappa'(r) \ra \\
&  =  \lambda^4 \la \delta^{\lambda,\alpha} \psi_2, r \varkappa'(\lambda r) \ra 
  =  \lambda^3 \la \delta^{\lambda,\alpha} \psi_2, (r \varkappa')^\lambda  \ra \\
&  =  \lambda \la \frac{\delta^{\lambda,\alpha} \psi_2}{r^2}, (r^3 \varkappa')^\lambda  \ra.
\end{split}
\]
Then using \eqref{psiLE2} we obtain
\[
\lambda |e(t)|  \les  \| \frac{\delta^{\lambda,\alpha} \psi_2}{r^2} \|_{L^2} \| \lambda  (r^3 \varkappa')^\lambda  \|_{L^2}  \les \| \frac{\psi}r \|_{L^2_{t,r}}, 
\]
as desired.

\medskip

\emph{ b) The bounds for $l$.}
We write $l(t)=l^1(t)+l^2(t)$ where
\begin{equation} \label{l1l2}
 l^1= - 2 \lambda^2 \la h_3^\lambda L_\lambda^* \psi , \varkappa^\lambda \ra, \quad l^2=- 4\lambda^2 \la  e^{2 i \alpha(t)}h_1^{\lambda(t)} \int_r^\infty \Re( \frac{ e^{- 2i \alpha(t)}h_1^{\lambda(t)}}s L_\lambda^* \psi ) ds, \varkappa^\lambda \ra. 
\end{equation}
A direct computation gives
\[
- \frac12 \lambda^{-1} l^1(t)= \lambda \la L_\lambda^* \psi,  h_3^\lambda \varkappa^\lambda \ra
= \lambda \la \psi,  L_\lambda (h_3^\lambda \varkappa^\lambda)  \ra =\lambda^2 \la \psi,  (L (h_3 \varkappa))^\lambda  \ra
=  \la \psi(\lambda^{-1} \cdot),  L (h_3 \varkappa)  \ra.
\]
This allows us to estimate
\[
| \lambda^{-1} l^1(t)| \les \| \psi(\lambda^{-1} \cdot) \|_{L^2(r \approx 1)} \approx \lambda^{-1} \| \frac{\psi(\lambda^{-1} \cdot)}{\lambda^{-1} \cdot} \|_{L^2(r \approx 1)} \approx \|\frac{\psi}r\|_{L^2(r \approx \lambda^{-1})} \les \|\frac{\psi}r\|_{L^2}. 
\]
Similarly, we compute:
\[
\begin{split}
-\frac{e^{-2 i \alpha}}4 \lambda^{-1} l^2 =  \lambda \la  \int_r^\infty \Re( \frac{ e^{- 2i \alpha}h_1^{\lambda}}s L_\lambda^* \psi ) ds, h_1^{\lambda} \varkappa^\lambda \ra =  \Re \left( e^{- 2 i \alpha(t)} \lambda \la  \int_r^\infty  \frac{h_1^{\lambda(t)}}s L_\lambda^* \psi  ds,  h_1^{\lambda} \varkappa^\lambda  \ra \right). 
\end{split}
\]
Defining 
\begin{equation} \label{g12def}
\gf_1(r) = \int_0^r h_1(s) \varkappa(s) sds, \qquad \gf_2=\frac{\gf_1 h_1}{r^2},
\end{equation}
this is further written as
\[
-\frac{e^{-2 i \alpha}}4 \lambda^{-1} l_2 =  \Re( e^{-2i \alpha(t)} \lambda \la   L_\lambda^* \psi ,  \gf_2^\lambda  \ra ), 
\]
thus resembling the previous expression modulo the phase plus a change in the test function $\gf_2$. Since $\gf_1(r)=0$ for $r \leq \frac12$ and $\gf_1(r)=c_1$ for  $r \geq 2$, it follows that $\gf_2(r)=0$ for $r \leq \frac12$ and $\gf_2(r)=\frac{c_1 h_1}{r^2}$ for $r \geq 2$; therefore the $\gf_2$ factor is similar to the $h_3 \varkappa$ factor appearing in the expression for $l^1$), the only difference being that it decays like $r^{-4}$ for large $r$, as opposed to being $0$. Thus we can estimate
\[
|\lambda \la   L_\lambda^* \psi ,  \gf_2^\lambda  \ra| =
|\lambda^2 \la \psi, (L\gf_2)^\lambda \ra|= | \la \psi(\lambda^{-1} \cdot), L\gf_2 \ra| \les \lambda^{-1} \| \frac{\psi(\lambda^{-1} r)}{\lambda^{-1} r} \|_{L^2} \| r L\gf_2\|_{L^2} \les \|\frac{\psi}r\|_{L^2}. 
\]
Based on the two estimates above on $\lambda^{-1} l^1$ and $\lambda^{-1} l^2$, we obtain the estimate on $\lambda^{-1} l$ as claimed in \eqref{hqest}.

\medskip

\emph{ c) The bounds for $\tilde q$.}
 There are many terms to be estimated here and we deal with them in the order in which they appear in the above expression of $Nlin$; recall that $\tilde q$ skips the first term in the expression of $Nlin$ since it collects only the terms which appear in $\widetilde{Nlin}$. 

\medskip

\emph{c1) The expression $\la \delta^\lambda A_2 \cdot L^*_\lambda \psi +  A_2 \frac{\delta^\lambda A_2}{r} \psi, \varkappa^\lambda  \ra$.}
Using \eqref{sysAp2}, the first term in this expression is 
\[
\begin{split}
\la \delta^\lambda A_2 \cdot L_\lambda^* \psi, \chi^\lambda  \ra & =   \la  L_\lambda^* \psi, \delta^\lambda A_2 \cdot \varkappa^\lambda \ra 
= \la  \psi,  L_\lambda (\delta^\lambda A_2 \cdot \varkappa^\lambda) \ra \\ 
& = \la  \psi,  L_\lambda \delta^\lambda A_2 \cdot \varkappa^\lambda \ra + \la  \psi,  \delta^\lambda A_2 \cdot \partial_r \varkappa^\lambda \ra \\
& = \la \psi , \left( - 2 h_1^\lambda \Re{(e^{2i\alpha} \psi)} + \Im(\psi \overline{\delta^{\lambda,\alpha} \psi}_2) - \frac{(\delta^\lambda A_2)^2}r \right) \varkappa^\lambda \ra
+ \la  \psi,  \delta^\lambda A_2 \cdot \frac1r  (r\varkappa')^\lambda \ra. \\
\end{split}
\]
Using \eqref{psiLE2} and \eqref{compsolh1-r} (of which we use the consequence $\| \delta^{\lambda,\alpha} \psi_2 \|_{L^\infty} + \|\delta^\lambda A_2\|_{L^\infty} \les \| \psi \|_{L^2} \ll 1$), this leads us to the following estimates:
\[
\begin{split}
|\la \delta^\lambda A_2 \cdot L_\lambda^* \psi, \varkappa^\lambda  \ra| & \les \la \frac{|\psi|^2}{r^2}, (1+ |\delta^{\lambda,\alpha} \psi_2|) r^2 \varkappa^\lambda \ra + 
\la \frac{|\psi|}r \cdot \frac{|\delta^\lambda A_2|}{r^2} \cdot |\delta^\lambda A_2|, r^2 \varkappa^\lambda \ra \\
& +  \la \frac{|\psi|}r \cdot \frac{|\delta^\lambda A_2|}{r^2}, r^2 |(r\varkappa')^\lambda| \ra  \les \lambda^{-2} (1+\| \psi\|_{L^2}) \|\frac{\psi}r\|^2_{L^2},
\end{split}
\]
and 
\[
|\la \delta^\lambda A_2 \cdot L_\lambda^* \psi, \varkappa^\lambda  \ra|  \les \lambda^{-1} (1+\| \psi\|_{L^2}) \| \psi\|_{L^2} \cdot \|\frac{\psi}r\|_{L^2}. 
\]
These are the estimates we were seeking for the first term; it is clear that the second term  $\la A_2 \frac{\delta^\lambda A_2}{r} \psi, \varkappa^\lambda  \ra$ is estimated in a similar manner (recall that $|A_2| \leq 2$).

\medskip

\emph{ c2) The expression $\la \int_r^\infty \Re (i L^*_\lambda \psi \frac{\overline{\delta^{\lambda,\alpha}\psi_2}}r ) dr, (h_1\chi)^\lambda \ra$.} 
Note that here we can harmlessly ignore the additional phase factor $2 e^{2i\alpha}$. We compute:
\[
\begin{split}
 \int_r^\infty L^*_\lambda \overline{\psi} \frac{\delta^{\lambda,\alpha} \psi_2}{s^2}  sds  & = -\overline{\psi(r)} \frac{\delta^{\lambda,\alpha} \psi_2}{r} + \int_r^\infty \overline{\psi} \cdot L_\lambda (\frac{\delta^{\lambda,\alpha} \psi_2}{r^2}) r dr \\
& = -\overline{\psi(r)} \frac{\delta^{\lambda,\alpha} \psi_2}{r} + \int_r^\infty \overline{\psi} \cdot \frac{L_\lambda  \delta^{\lambda,\alpha} \psi_2}{s^2} s ds - 2 \int_r^\infty \overline{\psi} \cdot \frac{ \delta^{\lambda,\alpha} \psi_2}{s^3} s ds \\
& =  -\overline{\psi(r)} \frac{\delta^{\lambda,\alpha} \psi_2}{r} - 2 \int_r^\infty \overline{\psi} \cdot \frac{ \delta^{\lambda,\alpha} \psi_2}{s^3} s ds \\
&  + \int_r^\infty \overline{\psi} \cdot \frac{2i h_3^\lambda \psi + \delta^\lambda A_2 \psi - \frac{1}s \delta^\lambda A_2 (2ie^{2 i \alpha}h_1^\lambda + \delta^{\lambda,\alpha} \psi_2)}{s^2} s ds,  \\
\end{split}
\]
where in passing to the last line we have used the system \eqref{sysAp2}. 

We estimate the term $\la \psi(r) \frac{\delta^{\lambda,\alpha} \psi_2}{r} , \varkappa^\lambda \ra$ just as we estimated the $\la  \psi,  \delta^\lambda A_2 \cdot \frac1r  (r\varkappa')^\lambda \ra$ term above. A direct estimate gives us:
\[
\| \frac1{r^2} \int_r^\infty \frac{\psi}s \cdot \frac{ \delta^{\lambda,\alpha} \psi_2}{s^2} s ds \|_{L^1(r \approx \lambda^{-1})} \les \| \frac{\psi}r\|_{L^2}  \cdot
\| \frac{\delta^{\lambda,\alpha} \psi_2}{r^2}\|_{L^2},
\]
and then using \eqref{psiLE2} we obtain:
\[
|\la \int_r^\infty \psi \cdot \frac{ \delta^{\lambda,\alpha} \psi_2}{s^3} s ds, \varkappa^\lambda \ra | \les \lambda^{-2} \| \frac{\psi}r\|_{L^2}^2.  
\]
Alternatively, we use \eqref{rdrmw} (with $p=1$ and $w=r^{-1}$) and  estimate
\[
\| \frac1{r} \int_r^\infty \frac{\psi}s \cdot \frac{ \delta^{\lambda,\alpha} \psi_2}{s} \frac1s s ds \|_{L^1} \les \| \frac{\psi}r\|_{L^2}  \cdot
\| \frac{\delta^{\lambda,\alpha} \psi_2}{r}\|_{L^2},
\]
from which, by using \eqref{compsolh1-r}, we obtain
\[
|\la \int_r^\infty \psi \cdot \frac{ \delta^{\lambda,\alpha} \psi_2}{s^3} s ds, \varkappa^\lambda \ra | \les \lambda^{-1} \|\psi\|_{L^2} \| \frac{\psi}r\|_{L^2}.  
\]
These are the expected estimates for the second term; finally the third term
\[
 \la \int_r^\infty \psi \cdot \frac{2i h_3^\lambda \psi + \delta^\lambda A_2 \psi - \frac{1}s \delta^\lambda A_2 (2ie^{2 i \alpha}h_1 + \delta^{\lambda,\alpha} \psi_2)}{s^2} s ds  , \varkappa^\lambda \ra
\]
is estimated in a similar fashion; the details are left as an exercise for the interested reader.

\medskip

\emph{c3) The expression $\la \int_r^\infty \Re (L^*_\lambda \psi \overline{\psi}), (h_1\varkappa)^\lambda \ra$.} Here we write: 
\[
\begin{split}
\la \int_r^\infty \Re (L^*_\lambda \psi \overline{\psi}) ds, (h_1\varkappa)^\lambda  \ra & = \la \int_r^\infty \Re ( (-\partial_s +\frac{2 h_3^\lambda-1}{s})\psi \cdot  \overline{\psi}) ds, (h_1\varkappa)^\lambda \ra \\
& = \la \int_r^\infty - \frac12 \partial_s  |\psi|^2 +  \frac{2 h_3^\lambda-1}{s} |\psi|^2 ds, (h_1\varkappa)^\lambda \ra \\
& =  \frac12 \la |\psi|^2,  (h_1\chi)^\lambda  \ra + \la \int_r^\infty   \frac{2 h_3^\lambda-1}{s} |\psi|^2 ds, (h_1\varkappa)^\lambda \ra. 
\end{split}
\]
The first term above is estimated as follows:
\[
\la |\psi|^2,  (h_1\varkappa)^\lambda   \ra \approx \int_{r \approx \lambda^{-1}} |\psi|^2 rdr \approx \lambda^{-2} \int_{r \approx \lambda^{-1}} 
\frac{|\psi|^2}{r^2} rdr,
\]
and this gives the correct contribution for the first $\tilde q$ bound.  Alternatively, the same term is estimated as follows:
\[
\la |\psi|^2,  (h_1\varkappa)^\lambda   \ra \approx \lambda^{-1} \int_{r \approx \lambda^{-1}} |\psi| \frac{|\psi|}{r} rdr \les \lambda^{-1}
\|\psi(t)\|_{L^2_r} \|\frac{\psi(t)}{r}\|_{L^2_r}^2,
\]
and this gives the correct contribution to the second $\tilde q$ bound. 

The second term above is estimated as follows:
\[
\begin{split}
|\la \int_r^\infty   \frac{2 h_3^\lambda-1}{s} |\psi|^2 ds, (h_1\varkappa)^\lambda  \ra| & = |\la \int_r^\infty   (2 h_3^\lambda-1) \frac{|\psi|^2}{s^2} sds, (h_1\varkappa)^\lambda  \ra| \\
& \les \| \frac{|\psi|^2}{r^2} \|_{L^1} \int (h_1\varkappa)(\lambda r) rdr \\
& \les \lambda^{-2} \| \frac{\psi}{r} \|^2_{L^2_r},
\end{split}
\]
and this gives the correct contribution to the first $\tilde q$ bound. We can also estimate this term as follows
\[
\begin{split}
 & \les |\la \int_r^\infty   \frac{|\psi|}{s} |\psi| ds, (h_1\varkappa)^\lambda  \ra| \les \| \frac{|\psi|}{r} \cdot |\psi| \|_{L^1_r} \int \frac{(h_1\varkappa)(\lambda r)}{r} rdr  \les \lambda^{-1} \| \frac{\psi(t)}{r} \|_{L^2_r} \|\psi(t)\|_{L^2},
\end{split}
\]
to obtain the correct contribution to the second $\tilde q$ bound.

\medskip

\emph{c4)  The expression $\la \int_r^\infty \Re (\frac{\delta^\lambda A_2}{r} \psi \overline{\psi}_1) ds, \varkappa^\lambda \ra$.} We write
\[
\frac{\delta^\lambda A_2}{r} \psi \overline{\psi}_1=  \frac{\delta^\lambda A_2}r \cdot \psi \cdot (\overline{\psi+i \frac{\psi_2}r}) 
=  \frac{\delta^\lambda A_2}r \cdot \psi \cdot (\overline{\psi}-i \frac{\overline{\psi_2}}r),
\]
and note that the previous estimates provide a template on how to obtain similar estimates for this term as well. The details are left as an exercise.

\medskip

\emph{ c5) The expression $ \la i A_0 \delta^{\lambda,\alpha} \psi_2, \varkappa^\lambda \ra$.} Based on the expression \eqref{A0-rep} for $A_0$, this term has two components: 

i) the contribution of the nonlinear component of $A_0$, which  combined with the trivial inequality $\| \delta^{\lambda,\alpha} \psi_2 \|_{L^\infty} \les \| \psi \|_{L^2}$
is estimated exactly as in the previous steps (c2)-(c4);

ii) the combination of the linear component of $A_0$ with the $\delta^{\lambda,\alpha} \psi_2$ factor, which in this context is also viewed as   a nonlinear contribution:
\[
\la \int_r^\infty \Re( \frac{2 e^{- 2i \alpha}h_1^{\lambda}}s L_\lambda^* \psi ) ds \cdot \delta^{\lambda,\alpha} \psi_2, \varkappa^\lambda \ra= \la \Re ( e^{- 2i \alpha} \int_r^\infty  \frac{2 h_1^{\lambda}}s L_\lambda^* \psi  ds) \cdot \delta^{\lambda,\alpha} \psi_2, \varkappa^\lambda \ra.
\]
We first integrate by parts
\[
\int_r^\infty  \frac{h_1^{\lambda}}{s} L_\lambda^* \psi  ds = - \frac{h_1^\lambda}r \psi +  \int_r^\infty \psi  L_\lambda(\frac{h_1^{\lambda}}{s^2}) s ds
= - \frac{h_1^\lambda}r \psi -2  \int_r^\infty \psi  \frac{h_1^{\lambda}}{s^3} s ds. 
\]
which is a-priori justified for regular $\psi$,
and then we estimate
\[
| \la  \frac{h_1^\lambda}r \psi \cdot \delta^{\lambda,\alpha} \psi_2, \varkappa^\lambda \ra |= 
| \la  \frac{h_1^\lambda}r \psi \cdot \frac{\delta^{\lambda,\alpha} \psi_2}{r^2}, r^2\varkappa^\lambda \ra | \les \lambda^{-2} \| \frac{\psi}r\|_{L^2} \| \frac{\delta^{\lambda,\alpha} \psi_2}{r^2} \|_{L^2} \les \lambda^{-2} \| \frac{\psi}r\|_{L^2}^2,
\]
as well as 
\[
| \la  \frac{h_1^\lambda}r \psi \cdot \delta^{\lambda,\alpha} \psi_2, \varkappa^\lambda \ra |= 
| \la  \frac{h_1^\lambda}r \psi \cdot \frac{\delta^{\lambda,\alpha} \psi_2}{r}, r \varkappa^\lambda \ra | \les \lambda^{-1} \| \frac{\psi}r\|_{L^2} \| \frac{\delta^{\lambda,\alpha} \psi_2}{r} \|_{L^2} \les \lambda^{-1} \| \frac{\psi}r\|_{L^2} \|\psi\|_{L^2}. 
\]
The second component above is estimated similarly in combination with the straightforward inequality:
\[
\| \int_r^\infty \psi  \frac{h_1^{\lambda}}{s^3} s ds \|_{L^2} \les \| \frac{\psi}r\|_{L^2}. 
\]
This finishes the proof of all our claims in \eqref{hqest}. 

\end{proof}

\subsection{Proof of Theorem \ref{PD-1}} 
The idea here is to carefully combine the bounds 
in Lemma~\ref{l:hqest} for the terms in the modulation equation with the local energy bounds for $\psi$.

Without restricting the generality of the argument, it suffices to establish the result on the forward in time maximal time of existence $I_{max}^+=[0,T_{max})$. 

\medskip

\emph{a) Modulation equation bounds.}
While the bounds in Lemma~\ref{l:hqest} are fixed time bounds, here we will primarily use them in an integrated fashion. Assume that we have space-time control of the local energy decay, i.e. that is we control the quantity  $\| \frac{\psi(t)}{r}\|_{L^2_t L^2_{r}(I \times \R)}$ for some time interval $I=[0,T]$ where $0 < T < T_{max}$. All the estimates involving space-time norms below are restricted to the time interval $I$.

The first bound on $e$ in \eqref{hqest} combined  
with the smallness of $\|\psi\|_{L^2}$ shows that $e$ 
is uniformly small,
\[
\| e\|_{L^\infty} \lesssim \|\psi\|_{L^\infty L^2} \ll 1.
\]
This shows that we can interpret $e$ perturbatively 
in \eqref{laeq2} and estimate using \eqref{hqest}  
\begin{equation} \label{lambdal2}
\|\frac{\alpha'}{\lambda}\|_{L^2_t(I)}+
\|\frac{\lambda'}{\lambda^2}\|_{L^2_t(I)} \les (1+\|\psi_0\|_{L^2_r}) \|\frac{\psi(t)}r\|_{L^2_t L^2_{r}(I)}.
\end{equation}

Incidentally, we remark, for later use, that this allows us to capture the perturbative nature of the term $ \frac{\lambda'}{\lambda} e(t)$ in \eqref{laeq2}, and thus justify why  it can be included in the perturbative term $q$. Indeed, the above with the second $e$ bound in \eqref{hqest} 
 yields
\begin{equation}\label{lpe-l2}
\| \frac{\lambda'}{\lambda} e(t) \|_{L^1_t} \les  \|\frac{\lambda'}{\lambda^2}\|_{L^2_t} \|\lambda e(t)\|_{L^2_t} 
\les  (1+\|\psi_0\|_{L^2_r}) \|\frac{\psi}r\|^2_{L^2},
\end{equation}
while using the first bound for $e$ in \eqref{hqest} gives
\begin{equation}\label{lpe-l1}
\| \frac{\lambda'}{\lambda^2} e(t) \|_{L^2_t} \les \|\frac{\lambda'}{\lambda^2}\|_{L^2_t} \| e(t)\|_{L^\infty_t} 
\les  (1+\|\psi_0\|_{L^2_r}) \|\psi\|_{L^\infty_t L^2_r} \|\frac{\psi}r\|_{L^2}. 
\end{equation}
This allows us to add the term $\frac{\lambda'}{\lambda} e(t)$ to the quadratic component $q(t)$ and extend the estimates on $\tilde q$ in \eqref{hqest} to $q$, namely 
\begin{equation} \label{hqest1}
     \| q \|_{L^1_t} \les \| \frac{\psi(t)}{r}\|_{L^2_t L^2_r}^2, \quad
    \| \lambda^{-1} q \|_{L^2_t} \les \|\psi \|_{L^\infty_t L^2_r} \| \frac{\psi}{r}\|_{L^2_{t,r}}. 
\end{equation}

This justifies the use of \eqref{laeq22} as the main formulation of the modulation equations, which will be used in all the analysis below. 

Applied in the modulation equation \eqref{laeq22}, the second bound above, combined with the second bound
for $l$ in \eqref{hqest},  yields
\begin{equation} \label{hqest2}
    \|\frac{\lambda'}{\lambda^2} \|_{L^2_t}
    +  \|\frac{\alpha'}{\lambda} \|_{L^2_t}
    \les \| \frac{\psi}{r}\|_{L^2_{t,r}}.
\end{equation}

\medskip

\emph{b) Bounds for $\psi$.}
Here we start with the local energy bound in Corollary~\ref{c:le}, which yields 
\[
\| \frac{\psi}r \|_{L^2_{t,r}} + \| \psi \|_{L^4_{t,r}} \les (1+\| \frac{\lambda'}{\lambda^2} \|_{L^2_t}) (\|\psi_0\|_{L^2} + \| W_\lambda \psi\|_{L^\frac43}).
\]
On the other hand from the bound \eqref{Wlest1}
for the nonlinearity we have 
\[
    \| W_\lambda \psi \|_{L^\frac43} \les \|\psi\|_{L^4}  (\| \frac{\psi}r \|_{L^2} + \|\psi\|_{L^4}^2).
\]
Combining the two we arrive at
\begin{equation}
 \label{hqest1a}
 \| \frac{\psi}r \|_{L^2_{t,r}} + \| \psi \|_{L^4_{t,r}} \les (1+\| \frac{\lambda'}{\lambda^2} \|_{L^2_t}) (\|\psi_0\|_{L^2}
+  (\| \frac{\psi}r \|_{L^2} + \|\psi\|_{L^4}^2)^2)
\end{equation}

\medskip

\emph{ c) A continuity argument.} 
So far, we have established the bounds  \eqref{hqest2}  and \eqref{hqest1a} in any interval $I$ of existence for the solution where all the norms involved are finite.
By Theorem~\ref{th2} we know the problem is well-posed at least on some small time interval $I=[0, \sigma \lambda_0^{-2}]$; further, we note that the theory developed in \cite{GuKaTs-2} gives control on the two norms $\| \frac{\psi}r \|_{L^2_{t,r}}$ and $\| \psi \|_{ L^4_{t,r}}$.

 By Theorem~\ref{th2} we know that the problem is well-posed at least on some small time interval $I=[0, \sigma \lambda_0^{-2}]$; further, we 
claim that on such a time interval the quantities 
$\| \psi \|_{ L^4_{t,r}(I \times [0,+\infty))}$ 
and $\| \frac{\psi}r \|_{L^2_{t,r}(I \times [0,+\infty))}$ are finite. 

Indeed the theory developed in \cite{GuKaTs-2} gives control on  $\| \psi \|_{ L^4_{t,r}}$ on the interval $I$. While there is no formal result stated in \cite{GuKaTs-2}, the arguments in that paper do give control on $\| \frac{\psi}r \|_{L^2_{t,r}}$ on $I$ as well. Alternatively, one could also argue based on a higher regularity argument.  At the level of maps in $\dot H^1 \cap \dot H^2$,   it follows from \eqref{Cextra} that $\psi \in L^2 \cap \dot H^1_e$. From this we obtain $\psi \in L^4_t$ and $\frac{\psi}r \in L^2_r$, which gives $\psi \in L^4_{t,r}$ and $\frac{\psi}r \in L^2_{t,r}$ locally in time and within a compact interval inside $[0,T_{max})$ where $\lambda(t)$ belongs to a compact interval inside $(0,+\infty)$ and we have uniform bounds for $\psi$ in $\dHde$.

Therefore the norms in  \eqref{hqest2}  and \eqref{hqest1a} are finite on any compact subinterval of the maximal existence interval $[0,T_{max})$.

Next we  combine the two bounds within a continuity argument in order
to arrive at the conclusion of the theorem.
For this we denote
\[
\Mf(T)=\| \frac{\lambda'}{\lambda^2} \|_{L^2_t[0,T]}
\qquad 
\MF(T) =\| \frac{\psi}r \|_{L^2_{t,r}([0,T] \times \R)} + \| \psi \|_{L^4_{t,r}([0,T] \times \R)}.
\]
From \eqref{hqest1}  and \eqref{hqest1a} we obtain
\[
\Mf(T) \leq C  \MF(T), \qquad \MF(T) \leq C (1+\Mf(T))(\delta+C\MF^2(T) + C\MF^3(T)). 
\]
From these relations we want to conclude that $\Mf(T),\MF(T) \lesssim \delta$ on the maximal interval of existence for the solutions $[0,T_{max})$.

Assuming that $\delta$ is small enough depending on the universal constant $C$, it is clear now that we can find a large enough constant $D$ such that if we know apriori that $\Mf(T),\MF(T) \leq 2D\delta$, then, based on the inequalities above, we can improve the result to $\Mf(T),\MF(T) \leq D\delta$. This allows us to run a standard continuity argument to obtain the desired result. 

For this we observe that 
both $\Mf$ and $\MF$ are continuous functions of $T$ for $T \in [0,T_{max})$,
and that as $T \to 0$ both $\Mf$ and $\MF$ have limit $0$. 
Denoting by $T_1$ the maximal time where $\Mf,\MF \leq 2D\delta$, our analysis above implies that 
$\Mf(T_1),\MF(T_1) \leq D\delta$. But by the continuity of $\Mf$, $\MF$, this contradicts the maximality of $T_1$ unless $T_1 = T_{max}$.

The above argument provides the result of the Theorem for the weaker structure when norms in $l^2 S(I^+_{max} \times \R)$ are replaced with 
$\| \frac{\psi}r \|_{L^2_{t,r} ( I^+_{max} \times \R)} + \| \psi \|_{ L^4_{t,r}( I^+_{max} \times \R)}$. However this weaker result is immediately upgraded to the full structure $l^2 S(I^+_{max})$ based on the result in Theorem~\ref{mainTl2} and the dual Strichartz bound \eqref{Wlest1} for the nonlinearity.

\subsection{Proof of Theorem \ref{PD-2}} 
Here we already have the bound \eqref{l1est} from Theorem~\ref{PD-1}. Then we can apply the linear bounds in Theorem~\ref{mainTl1} to obtain
\[
\|\psi\|_{\ell^1 S} \lesssim 
 \|\psi(0)\|_{\LX} + 
  \|N(\psi)\|_{\ell^1 N}.
\]
The nonlinear contribution is estimated using 
Lemma~\ref{gNref} to obtain
\[
\|\psi\|_{\ell^1 S} \lesssim 
 \|\psi(0)\|_{\LX} + 
  \|\psi\|_{\ell^2 S}.
\]
Finally we conclude using Theorem~\ref{PD-1}.

\section{The modulation parameter dynamics - Part 2} \label{s:mod2}
From the previous section we recall the ODE \eqref{laeq22} governing the dynamics of $\alpha$ and $\lambda$:
\begin{equation} \label{MODE}
 4  e^{2 i \alpha(t)} \left( \alpha'  c_1 + i \frac{\lambda'}{\lambda} c_2  \right)= l(t) + q(t),
\end{equation}
where, as expanded in \eqref{l1l2}, $l(t)=l^1(t)+l^2(t)$ with
\begin{equation} \label{ldec12}
 l^1= - 2 \lambda^2 \la h_3^\lambda L_\lambda^* \psi , \varkappa^\lambda \ra, \quad l^2=- 4\lambda^2 \la  e^{2 i \alpha(t)}h_1^{\lambda(t)} \int_r^\infty \Re( \frac{ e^{- 2i \alpha(t)}h_1^{\lambda(t)}}s L_\lambda^* \psi ) ds, \varkappa^\lambda \ra. 
\end{equation}
Here $q$ largely plays a perturbative role. Indeed, we have established (see \eqref{hqest1}) the following bound $\|q\|_{L^1_t} \les \|\frac{\psi}r\|_{L^2_{t,r}}^2$, where the right hand side is estimated directly in terms of the initial data on the maximal time of existence, see Theorem~\ref{PD-1}.

It is clear that the the component responsible for a potential finite time or infinite time blow-up of $\lambda$, or for the more complex dynamics of the set of parameters $(\lambda(t),\alpha(t))$ (as opposed to stability which would correspond to a small perturbation of the original state $(\lambda(0),\alpha(0))$), is the linear component $l(t)$; indeed if this was not present in the above ODE, then we would simply control $\alpha(t)$ and $\log \lambda(t)$ globally using the above bound for $q$ hence ruling out blow-up (both in finite and infinite time) and providing global bounds for both $\alpha$ and $\log \lambda$. The information we have so far on $l$, from Lemma~\ref{l:hqest} in the previous section, is not accurate enough  in order to 
clarify the potential scenarios for the behaviour of solutions to the ODE \eqref{MODE}, such as uniform bounds versus growth or even potential blow up in finite time. 

In this section we perform a refined analysis of the linear term $l(t)$, which will later be used to rule out the finite time blow-up scenario. This analysis has two main parts: an algebraic one that seeks to better understand the structure of the expressions $l^1$ and $l^2$ in \eqref{ldec12}, and an analytic one that provides effective estimates for the terms revealed by the algebraic part. The analysis carried in the algebraic part is inspired by the work of Gustafson, Nakanishi and Tsai \cite{GuNaTs}, although in our context it is carried with respect to different variables. The analytic part of our analysis is entirely new, and this is what allows us to obtain the second main result of this paper, the global well-posedness result of Theorem \ref{tmain-G}, as well as its refined $l^1$ Besov counterpart in  Theorem \ref{tmain-G1}.  
The main result of this section is the following:

\begin{t1} 
\label{t:mod-syst}
The system \eqref{MODE} can be rewritten as follows on the maximal 
time interval of existence $I = [0,T_{max})$:
\begin{equation} \label{mainaleq}
\begin{split}
\alpha' & = -\frac12 \Re \left( i e^{-2 i \alpha(t)} \partial_t ( \lc^1 + 2 \lc^2) \right) + \frac14 \Re Q(t),
\\ (\ln \lambda)' & = -\Re ( e^{-2 i \alpha(t)} \partial_t  \lc^1) + \frac12 \Im   Q(t). 
\end{split}
\end{equation}
Here $\lc^i: I \rightarrow \C$,  $i=1,2$ are continuous functions which satisfy the additional bounds 
\begin{equation} \label{lcorH12}
\| \lc^i \|_{\dot H^\frac12[I]} \lesssim \| \psi \|_{l^2S[I]}, \qquad \| \lc^i \|_{\dot B^\frac12_{2,1}[I]} \lesssim \| \psi \|_{l^1S[I]},
\end{equation}
while the function $Q$ obeys the following estimate:
\begin{equation} \label{q2est}
\|Q\|_{L^1_t[I]} \les \| \frac{\psi}r \|_{L^2[I]} (\| \frac{\psi}r \|_{L^2[I]} + \|\psi\|_{L^4[I]}^2).
\end{equation}
Further, assume $[0,T] \subset I$ is a compact interval and let $\lambda_T^{max}=max_{t \in [0,T]} \lambda(t)$. Then we have the following estimate:
\begin{equation} \label{lcorB12}
\| \lc^i \|_{\dot B^{\frac12}_{2,1}([0,T])} \lesssim \left(\ln ( 2+ ( \lambda_T^{max})^2  T)\right)^\frac12 \| \psi \|_{l^2S[I]}, \quad i=1,2.
\end{equation}

\end{t1}

\begin{r1}\label{r:norms}
Here the norm for the $\dot B^{\frac12}_{2,1}[I]$ quotient space of functions modulo constants is defined using extensions to $\R$,
\[
\| f \|_{\dot B^{\frac12}_{2,1}[I]} 
= \inf\{ \|f_{ext}\|_{\dot B^{\frac12}_{2,1}[\R]}
; \ f_{ext} = f \ \text{in } I\}.
\]
\end{r1}

We briefly comment on the system \eqref{mainaleq}, whose study is the main goal of the next section.
The $Q$ term in both ODE's is easily seen to play a perturbative role, so at leading order the dynamics is driven by $\Lambda^i$. Assuming these have been properly estimated,  what we have is a nonlinear ODE for the $\alpha$ variable and then an apparent linear ODE for $\lambda$;  the equation for $\lambda$ is not entirely linear since the term $\lc^1$ depends on $\lambda$, but this dependence can be incorporated into the estimates on $\lc^1$, thus rendering an essentially linear ODE. In some sense, this indicates that it is the equation for $\alpha$ should be our main focus, though ultimately 
is the growth in the $\lambda$ equation which we will want to control.

This section is organized as follows. In subsection \ref{lanalysis} we carry the formal algebraic computations for $l^1$ and $l^2$, the two components that appear in the linear term $l$, and this leads us to the new form of the modulation system in \eqref{mainaleq}.  The $L^1$ bound for $Q$ is proved in subsection \ref{s:Q}.
Finally in subsection~\ref{H12est} we turn our attention to the  expressions $\lc^1$ and $\lc^2$ appearing in \eqref{mainaleq}. These are formal objects which apriori are defined only if assume some additional spatial decay on the field $\psi$,  which apriori is  more than what our framework allows (which is $\psi(t) \in L^2_r$). The first step is to provide a rigorous justification of the fact that $\lc^1$ and $\lc^2$ are well-defined objects in our setting. Finally, we establish 
\eqref{lcorH12}; this is based on transferring local energy decay bounds from $\psi$ to $\dot H^\frac12$ bounds on $\lc^i, i=1,2$.

\subsection{An indepth analysis of the linear term \texorpdfstring{$l(t)$}{} and the derivation of \texorpdfstring{(\ref{mainaleq})}{}
} \label{lanalysis}

In this section we proceed with a further analysis of the linear term $l(t)$ appearing in \eqref{MODE}. This section draws inspiration and technical details from the work \cite{GuNaTs} as follows: the analysis of the operator $L_\lambda$, their inverse and adjoints has been already developed in Section 3 of \cite{GuNaTs} and the identification of the terms $\lc$ is equivalent to the identification of the ``normal form" in Section 9 in  \cite{GuNaTs}. However, for pure comparison reasons, we note that all of our computations are performed at the level of the gauge elements, while those appearing in \cite{GuNaTs} are performed at the level of the actual map.

In later computations in this section we will need to work with a right inverse of the operator $L_\lambda$. It is convenient to use a regularized version defined as follows
\[
R_{\lambda,\varphi} f = h_1^\lambda  \int_0^\infty \int_{r'}^r \frac1{h_1^\lambda(r'')} f(r'') dr'' \lambda^2  \varphi^\lambda(r') h_1^\lambda(r') r' dr',
\]
where $\varphi \in C_0^\infty((0,\infty);\R)$; in fact we assume that $\varphi$ is
supported in the interval $(\frac12,2)$.

Assuming also that we have the normalization $\la \varphi, h_1\ra=1$, a direct computation shows that the following relations hold true:
\[
L_\lambda R_{\lambda,\varphi} f=f, \qquad R_{\lambda,\varphi} L_\lambda f = f - h_1^\lambda \la f, \lambda^2 \varphi^\lambda \ra. 
\]
We also note that
\[
R_{\lambda,\varphi} = \lambda^{-1} D_\lambda R_{1,\varphi} D_{\lambda}^{-1}, \quad L_\lambda = \lambda D_\lambda L_1 D_{\lambda}^{-1}, 
\]
where $(D_\lambda f)(x)=f^\lambda(x)=f(\lambda x)$; the reason we use also the notation $D^\lambda f$ in addition to the standard $f^\lambda$ is that it is more streamlined in longer formulas just like the one above. 

We also need the adjoint of $R_{\lambda,\varphi}$. To compute it  we set $\lambda=1$ and use the following inversion formula:
\[
R_{1,\varphi} f(r)= \int_0^\infty \int_0^\infty  \frac{h_1(r)}{h_1(r'')} k(r,r',r'')   \varphi(r') h_1(r') r' f(r'') dr'' dr'= \int \int K(r,r',r'') f(r'') dr'' dr',
\]
where 
\[
k(r,r',r'')= 1_{\{r' \leq r'' \leq r\}} -1_{\{r \leq r'' \leq r'\}}. 
\]
 Then we have 
\[
\begin{split}
R^*_{1,\varphi} g(r'') & = \int_0^\infty \int_0^\infty K(r,r',r'') g(r) \frac{r}{r''} dr dr' \\
& = \frac{1}{r'' h_1(r'')} \int_0^\infty \int_0^\infty  h_1(r) k(r,r',r'') \varphi(r') h_1(r') r' r g(r) 
dr dr' \\
& = \frac{1}{r'' h_1(r'')} \left( \int_0^{r''} \int_{r''}^\infty  h_1(r) \varphi(r') h_1(r') r' r g(r) dr dr' -  \int_{r''}^\infty \int^{r''}_0 ... dr dr'  \right). \\
\end{split}
\]
This implies that 
\begin{equation} \label{Radj12}
R^{*}_{1,\varphi} g(r'')= -\frac{\la \varphi, h_1 \ra}{r'' h_1(r'')} \int_0^{r''}   h_1(r)  r g(r) dr, \qquad r'' \leq \frac12
\end{equation}
respectively
\begin{equation} \label{Radj2}
R^{*}_{1,\varphi} g(r'')=  \frac{\la \varphi, h_1 \ra}{r'' h_1(r'')}   \int_{r''}^\infty  h_1(r)  r g(r)
dr, \qquad r'' >2.
\end{equation}
In the remaining case  $r'' \in [\frac12,2]$ we may have contributions from both terms. It is easy to see that the general adjoint operator is given by
\[
 R_{\lambda,\varphi}^{*} = \lambda^{-1} D_\lambda R^{*}_{1,\varphi} D_{\lambda}^{-1}. 
\]
Finally, from the formula:
\[
R_{\lambda,\varphi} L^\lambda f = f  - h_1^\lambda \la f, \lambda^2 \varphi^\lambda \ra
\]
it follows that
\[
L_{\lambda}^* R_{\lambda,\varphi}^{*}  g = g -   \lambda^2 \varphi \la h_1^\lambda, g \ra. 
\]
Note that this last formula implies that in order to have a clean recovery formula in the sense $ L_{\lambda}^* R_{\lambda,\varphi}^{*}  g=g$, we need to free $g$ of the $h_1^\lambda$ mode, that is
\[
L_{\lambda}^* R_{\lambda,\varphi}^{*}  g=g \ \ \mbox{iff} \ \ \la g, h_1^\lambda\ra=0. 
\]

\vspace{.2in}

With the above formalism in place, we are now ready to proceed with our refined analysis of the linear term $l(t)=l^1(t)+l^2(t)$, see \ref{ldec12}. 
 We start with the analysis for $l^1$. Denoting 
 \[
 c=\frac{ \la h_3^\lambda \varkappa^\lambda, h_1^\lambda \ra}{ \la h_1^\lambda, h_1^\lambda \ra}=\frac{ \la h_3 \varkappa, h_1 \ra}{ \la h_1, h_1 \ra},
 \]
 we compute
\[
\begin{split}
\frac{l^1}{-2\lambda^2}  =  \la L_\lambda^* \psi,  h_3^\lambda \varkappa^\lambda  \ra & =  \la L_\lambda^* \psi,  h_3^\lambda \varkappa^\lambda -ch_1^\lambda  \ra \\
&=  \la L_\lambda^* \psi ,  L_\lambda^* R_{\lambda,\varphi}^{*}  (h_3 \varkappa -ch_1)^\lambda    \ra \\
& =  \la L_\lambda L_\lambda^* \psi ,  R_{\lambda,\varphi}^{*}  (h_3 \varkappa -ch_1)^\lambda   \ra  =  \la \tilde H_\lambda \psi ,  R_{\lambda,\varphi}^{*} (h_3 \varkappa -ch_1)^\lambda   \ra \\
& =  \la i \partial_t \psi - N(\psi),  R_{\lambda,\varphi}^{*} (h_3 \varkappa -ch_1)^\lambda    \ra \\
& =  i \la  \partial_t \psi ,  R_{\lambda,\varphi}^{*} (h_3 \varkappa -ch_1)^\lambda    \ra - \la N(\psi) ,  R_{\lambda,\varphi}^{*} (h_3 \varkappa -ch_1)^\lambda   \ra, 
\end{split} 
\]
where in the last two expressions we have used the Schr\"odinger map equation \eqref{psieq2}. The linear component $ i \la  \partial_t \psi,  R_{\lambda,\varphi}^*  (h_3 \varkappa  -ch_1)^\lambda  \ra$ above needs further work. We compute
\begin{equation} \label{R*h3chi}
\begin{split}
  R_{\lambda,\varphi}^{*}(h_3 \varkappa - ch_1)^\lambda = & \ \lambda^{-1} D_\lambda R_{1,\varphi}^{*} D_{\lambda}^{-1} (h_3 \varkappa - ch_1)^\lambda \\
 = & \ \lambda^{-1} D_\lambda R_{1,\varphi}^{*} (h_3 \varkappa - ch_1) = \lambda^{-1} D_\lambda G_1,
\end{split}
\end{equation}
where 
\[
G_1= R_{1,\varphi}^{*} (h_3 \varkappa - ch_1) = R^{1,*}_\varphi (h_3 \varkappa ) - c R_{1,\varphi}^{*} (h_1):=G_{1,weak}-cG.
\]

From \eqref{Radj12} and \eqref{Radj2} it follows $G_{1,weak}= R_{1,\varphi}^{*} (h_3 \varkappa)$ is compactly supported in the region $\frac12 \leq r'' \leq 2$, and it is also smooth,
\[
| \partial_r^\alpha R_{1,\varphi}^{*} (h_3 \varkappa)(r'')| \les_\alpha 1, \quad \alpha \geq 0;
\]
this is the reason we will refer to $G_{1,weak}$ as the weak contribution.

On the other hand $G$ is not compactly supported. In the regime $\frac12 \leq r'' \leq 2$ it satisfies $|\partial_r^\alpha G(r)| \les 1, \alpha \geq 0$. 
For $r \leq \frac12 $, from \eqref{Radj12} it follows that
\[
G(r)= -\frac{1}{r h_1(r)} \int_0^{r}   h_1^2(s)  s ds.
\]
A direct computation gives that for $r \leq \frac12$ we have
\[
G(r)=-\frac{(1+r^4)}{2r^3} (-\frac{r^2}{1+r^4} + \tan^{-1} r^2).
\]
From the above we obtain the pointwise bounds 
\[
 |\partial_r^\alpha G(r)| \les_\alpha r^{3-\alpha}, 
 \quad \alpha \geq 0, \quad r \leq \frac12 . 
\]
From \eqref{Radj12} it follows that for $r \geq 2$ we have 
\[
G(r)= \frac{1}{r h_1(r)} \int_r^{\infty}   h_1^2(s)  s ds.
\]
A direct computation gives $G(r)$ for $r>2$ as follows
\begin{equation}\label{glgr}
G(r) = \frac{\left(r^4+1\right) \left(\tan ^{-1}\left(\frac{1}{r^2}\right)+\frac{r^2}{r^4+1}\right) }{2 r^3}, \quad r >2;
\end{equation}
in particular we have $\frac{G(\frac{1}{r})}{r^{2}}=-G(r) + \frac{1}{r}, \quad r \geq 2$ and one could deduce either of the formulas in the regimes $r \leq \frac12$ and $r > 2$ from the other one.  From the above we obtain the pointwise bounds 
\[
| \partial_r^\alpha G(r)| \les_\alpha r^{-1-\alpha} , \quad \alpha \geq 0.
\]
A closer look at $rG$ shows that the asymptotic behavior of $rG$ at infinity is of the form 
\begin{equation} \label{rGass}
rG=1+ \frac13r^{-4} + O(r^{-8}).
\end{equation}
The following object $\tilde g=r(rG)'$ will occur later and providing estimates for it is helpful. From the above it follows that
\[
\tilde g = -\frac{4}{3 r^4}+O(r^{-8}),
\]
for $r$ large. Given the trivial bounds $|\tilde g(r)| \les 1$ for $r \les 1$, it follows that $\|\tilde g\|_{L^2} \les 1$.
Due to the good properties of the extra term $G_{1,weak}$, the above remains valid when $G$ is replaced by $G_1$ and $g$ by $\tilde g_1=r(rG_1)'$. Thus we record
\begin{equation} \label{tgiL2}
\|\tilde g_1\|_{L^2} \les 1.
\end{equation}
With the characterization of $G_1$ in place, we seek to better understand the term
\[
 \la i \partial_t \psi, \lambda^{-1} G^\lambda_1 \ra = i \lambda^{-2} \int_{0}^\infty \partial_t \psi(t,r) (r G_1)^\lambda(r) dr.
\]
A formal computation yields 
\[
\begin{split}
\lambda^2 \la i \partial_t \psi, \lambda^{-1} G_1^\lambda \ra
 & = i \int_{0}^\infty \partial_t \psi(t,r) (r G_1)^\lambda(r) dr \\
 & =  i \partial_t \int_{0}^\infty  \psi(t,r) (r G_1)^\lambda(r) dr -  i \int_{0}^\infty  \psi(t,r) \partial_t (r G_1)^\lambda(r) dr \\
 & =  i \partial_t \int_{0}^\infty  \psi(t,r) (r G_1)^\lambda(r) dr -  i \frac{\lambda'}{\lambda} \int_{0}^\infty  \psi(t,r) (r (r G_1)')^\lambda(r) dr. \\
\end{split}
\]
But the issue we face is that the integral
\[
\int_{0}^\infty  \psi(t,r) (r G_1)^{\lambda(t)}(r) dr
\]
is not finite for a generic $\psi \in L^2_r$ given the slow decay of $G_1$ described earlier; indeed if $\lambda(t)=1$ then $rG_1 = 1 + O(r^{-3})$ as $r \rightarrow \infty$, and the expression
\[
\int_1^\infty \psi(r) dr= \int_1^\infty \psi(r) \cdot \frac1r rdr
\]
may fail to be integrable since $\chi_{r \geq 1} r^{-1} \notin L^2$. Any additional decay on $\psi$ at infinity (in addition to $\psi \in L^2$) would fix this problem, but this would restrict the allowable initial data. 

To rectify this, the key observation is that the divergent part of the above integral comes
from the very low frequencies of $\psi$, which essentially do not change in time on compact time intervals. Based on this heuristic argument, we choose instead the first correction term to be
defined by
\begin{equation} \label{hcor1}
\lc^1(t)= \int_{0}^\infty  \psi(t,r) (r G_1)^{\lambda(t)}(r) - \psi(0,r) (r G_1)^{\lambda(0)}(r) dr.
\end{equation}
Since the term we just added is time independent, the earlier computation is not affected. Justifying that this modified integral is finite will be the subject of the following subsection.

For now we take for granted the fact that $\lc^1$ is well-defined and continue with our computations. Recalling what has been done above, we have established that
\[
l^1=-2 \lambda^2 \la  L_\lambda^* \psi, h_3^\lambda \varkappa^\lambda \ra= -2 i \partial_t \lc^1 + 2 q_1,
\]
where 
\[
q_1=\lambda \la N(\psi) ,  G_1^\lambda   \ra +  i \frac{\lambda'}{\lambda} \int_{0}^\infty  \psi(r) (r (r G_1)')^\lambda(r) dr. 
\]
We will prove at the end of this section that $q_1$ has good estimates in $L^1_t$, thus it can be absorbed into the better term $Q$ in \eqref{mainaleq}.

We turn our attention to the second linear contribution:
\[
 l^2 = - 4 \lambda^2 e^{2 i \alpha} \Re \left( e^{- 2 i \alpha(t)} \la  \int_r^\infty  \frac{h_1^{\lambda}}r L_\lambda^* \psi  dr,  h_1^{\lambda} \varkappa^\lambda  \ra \right). 
\]
In the previous section we have introduced (see \eqref{g12def} and the computations right after)
\[
\gf_1(r) = \int_0^r h_1(s) \varkappa(s) sds, \qquad \gf_2=\frac{\gf_1 h_1}{r^2},
\]
which allowed us to write
\[
-\frac{e^{-2 i \alpha}}4 l_2 =  \Re( e^{-2i \alpha(t)} \lambda^2 \la   L_\lambda^* \psi ,  \gf_2^\lambda  \ra ).
\]
 These functions have the following properties:

 i) $\gf_1(r)=\gf_2(r)=0$ for $r \leq \frac12$;

 ii) $\gf_1(r)=c_1$ for $r \geq 2$, $\gf_2(r) = O(r^{-4})$ as $r \rightarrow \infty$.

We then proceed with the analysis of the term $\la   L_\lambda^* \psi ,  \gf_2^\lambda  \ra$ just as we did earlier, writing
\[
\lambda^2 \la L_\lambda^* \psi,  \gf_2^\lambda \ra= i \partial_t \lc^2 + q_2,
\]
where 
\begin{equation} \label{hcor2}
\lc^2 (t)= \int_{0}^\infty  \psi(t,r) (r G_2)^{\lambda(t)}(r) - \psi(0,r) (r G_2)^{\lambda(0)}(r) dr, \quad G_2= R_{1,\varphi}^{*} (\gf_2 - c_3h_1),
\end{equation}
 and
\[
q_2=- \lambda \la N(\psi) ,  G_2^\lambda   \ra -  i \frac{\lambda'}{\lambda} \int_{0}^\infty  \psi(r) (r (r G_2)')^\lambda(r) dr. 
\]
Here we let 
\[
G_2=  R_{1,\varphi}^{*}  (\gf_2-c_3h_1)= R_{1,\varphi}^{*} \gf_2-c_3 G:= G_{2,weak} - c_3 G, \qquad c_3= \frac{ \la \gf_2, h_1 \ra}{ \la h_1, h_1 \ra},
\]
and note that $G_{2,weak}$ is the weak term, while $-c_3 G$ is the strong one, just as above. The only slight difference is that while $G_{1,weak}$ was compactly supported, $G_{2,weak}$ is not; instead $G_{2,weak}(r)=0$ for $r \leq \frac12$, $G_{2,weak} = O(r^{-4})$ as $r \rightarrow \infty$ and it is smooth. For all practical purposes $G_{2,weak}$ plays a similar role to $G_{1,weak}$.

This gives the following representation for $l_2$:
\[
l^2= -4 e^{2 i \alpha} \Re \left( e^{- 2 i \alpha(t)} i \partial_t \lc^2 \right)
 -4 e^{2 i \alpha} \Re \left( e^{- 2 i \alpha(t)} q_2 \right).
\]
At this time our ODE system \eqref{MODE} takes the following form:
\[
\begin{split}
4  e^{2 i \alpha(t)} ( c_1 \alpha'  + i c_2 \frac{\lambda'}{\lambda} )
 = -2 i \partial_t \lc^1 - 4 e^{2 i \alpha} \Re \left( e^{- 2 i \alpha} i \partial_t \lc^2  \right) + Q,
\end{split}
\]
where 
\[
Q=2q_1 - 4 e^{2 i \alpha} \Re \left( e^{- 2 i \alpha} q_2 \right) + q.
\]
 Taking the real and the imaginary part above we obtain the following system for $\alpha$ and $\lambda$:
\[
\begin{split}
c_1 \alpha' & = -\frac12 \Re \left( i e^{-2 i \alpha(t)} \partial_t ( \lc^1 + 2 \lc^2) \right) + \frac14 \Re (e^{-2 i \alpha(t)} Q(t)),
\\ 2c_2 (\ln \lambda)' & = -\Re ( e^{-2 i \alpha(t)} \partial_t  \lc^1) - \frac12 \Re  (i e^{-2 i \alpha(t)} Q(t)). 
\end{split}
\]
By simply relabeling $e^{-2 i \alpha(t)} Q(t)$ to be $Q(t)$, which has no effect on the $Q$ bounds, we obtain \eqref{mainaleq}.

\subsection{The \texorpdfstring{$L^1$}{} bounds for \texorpdfstring{$Q$}{}}\label{s:Q}
Our goal here is to prove the bound \eqref{q2est}, which asserts that $Q$ is small in $L^1_t$. 
We wll focus on the estimates for $q_1$; the ones for $q_2$ are similar while $q$ has already been estimated in \eqref{hqest1}.

Using \eqref{R*h3chi}, the characterization of $G_1$ and \eqref{Wlest2}, we obtain the following estimate:
\[
\begin{split}
\|\lambda^2 \la N(\psi) ,  R_\lambda^* \left( (h_3 \varkappa-ch_1) ^\lambda \right)  \ra \|_{L^1} & \les \|\frac{N(\psi)}r\|_{L^1} \|\lambda^2 r R_{\lambda,\varphi}^* \left( (h_3 \varkappa -ch_1) ^\lambda \right) \|_{L^\infty} \\
& \les \| \frac{\psi}r \|_{L^2} (\| \frac{\psi}r \|_{L^2} + \|\psi\|_{L^4}^2) \| (r G_1)^\lambda \|_{L^\infty} \\
& \les \| \frac{\psi}r \|_{L^2} (\| \frac{\psi}r \|_{L^2} + \|\psi\|_{L^4}^2). 
\end{split}
\]
This provides the desired estimate for the first component in $q_1$. We continue with the estimate for the second component in $q_1$.
Here we recall the estimate for $\tilde g_1$ from \eqref{tgiL2}; based on this we obtain
\[
| \int_{0}^\infty  \frac{\psi}r \tilde g_1^\lambda(r) r dr | \les \| \frac{\psi}r \|_{L^2_{t,r}} \| \tilde g_1^\lambda \|_{L^2_r} 
\les \| \frac{\psi}r \|_{L^2_{r}} \lambda^{-1} \|\tilde g_1\|_{L^2} 
\les \lambda^{-1} \| \frac{\psi}r \|_{L^2_{r}}. 
\]
From this we conclude with
\[
\| \frac{\lambda'}{\lambda} \int_{0}^\infty  \psi(r) (r (r G_1)')^\lambda(r) dr \|_{L^1_t} \les \|\frac{\lambda'}{\lambda^2}\|_{L^2_t} \| \lambda \int_{0}^\infty  \frac{\psi}r \tilde g_1^\lambda(r) r dr \|_{L^2_t} \les (1+\|\psi_0\|_{L^2_r}) \|\frac{\psi}r\|^2_{L^2}. 
\]
This provides the desired bound for the second component in $q_1$.

\subsection{The analysis of the correction terms \texorpdfstring{$\lc^i$}{}} \label{H12est}
The goal of this subsection is to justify that $\lc^i(t), i=1,2$ are well defined and continuous functions of $t$, and satisfy the bounds
\eqref{lcorH12}, \eqref{lcorB12}. We recall their definition from the previous section:
\[
\lc^i(t)= \int_{0}^\infty  \psi(t,r) (r G_i)^{\lambda(t)}(r) - \psi(0,r) (r G_i)^{\lambda(0)}(r) dr.
\]
The first observation is that these integrals are interpreted in a singular sense, as
\begin{equation}\label{lci-sharp}
\lc^i(t)= \lim_{R \rightarrow \infty} \int_{0}^R  \psi(t,r) (r G_i)^{\lambda(t)}(r) - \psi(0,r) (r G_i)^{\lambda(0)}(r) dr.
\end{equation}
For each $R$ the above integrals are well-defined but the existence of the limit is not obvious for $\psi \in C(L^2)$. However, for such 
$\psi$ we can harmlessly replace the sharp cutoff with a regularized 
cutoff, 
\begin{equation}\label{lci-mild}
\lc^i(t)= \lim_{R \rightarrow \infty} \int_{0}^\infty  \chi_{\leq 0}(r/R) \left[\psi(t,r) (r G_i)^{\lambda(t)}(r) - \psi(0,r) (r G_i)^{\lambda(0)}(r)\right] dr.
\end{equation}
Here we recall (from Section \ref{defnot}) that $\chi_{\leq 0}$ is smooth, supported in $(0,2)$ and is identically equal to $1$ on $(0,1)$ (we can use here any other function with these properties). 
 Indeed, for $\psi \in C(L^2)$ the difference of the integrals in \eqref{lci-sharp} and \eqref{lci-mild} is easily seen to converge to zero.

To understand this, we consider a dyadic decomposition in frequency for $\psi$, and the corresponding decomposition for $\lc^i$. Precisely, for every $k \in \Z$ we define
\begin{equation}\label{lcik}
\lc_k^i (t)= \int_0^\infty  P_k^{\lambda(t)} \psi(t) \cdot (rG_i)^{\lambda(t)} dr.
\end{equation}
also interpreted in a singular sense, first as  as $\lim_{R \rightarrow \infty} \int_0^R$ and then as 
\begin{equation}\label{lcikR}
\lc_k^i (t)= \lim_{R \to \infty} \int_0^\infty  P_k^{\lambda(t)} \psi(t) \cdot \chi_{\leq 0}(r/R)(rG_i)^{\lambda(t)} dr.
\end{equation}
With these notations we will show that
\begin{equation} \label{LsumLk}
\lc^i(t)= \sum_{k \in \Z} \lc_k^i (t) - \lc_k^i (0)
\end{equation}
as a uniformly convergent series. 
We then estimate $\lc^i$ by separately estimating the functions $\lc_k^i$,
which will be thought of as the generalized Littlewood-Paley pieces of $\lc^i$.

For now we make the observation that the limit in \eqref{lcikR}
does exist, uniformly in $t$.
To see that, we use the uniform $L^2$ bound for $\psi$ and the Fourier 
representation of $P_k^{\lambda(t)} \psi(t)$. Then it suffices to show
that the limit exists if we replace $P_k^{\lambda(t)} \psi(t)$ with a generalized eigenfunction $\psi_\xi^\lambda$, uniformly for $\xi$ in a fixed dyadic range and for $\lambda$ in a compact set. But this is a consequence of the asymptotic behavior of $\psi_\xi$ in Theorem~\ref{httransthm}.  

For fixed time $t$ we can insert an additional dyadic projection $\tilde P_k^{\lambda(t)}$ on $\psi$ in \eqref{lcikR}. This is  $L^2$ selfadjoint  so we can write
\[
 \int_0^\infty   P_k^{\lambda(t)} \psi(t) \cdot \chi_{\leq 0}(r/R)(rG_i)^{\lambda(t)} dr = \int_{0}^\infty P_k^{\lambda(t)} \psi(t) \tilde P_k^\lambda (\chi_{\leq 0}(r/R) \lambda G_i^\lambda) rdr.
\]
This is not a-priori justified without the $\chi$ cutoff as the function $G_i$ (and thus $\lambda G_i^\lambda$) does not belong to $L^2$ since it decays like $r^{-1}$ at $\infty$. However, it is justified with the cutoff 
$\chi_{\leq 0}(r/R)$ inserted, and  the limit as $R \to \infty$ is also well defined.  This will allow us to define 
$P_k^\lambda (\lambda G_i^\lambda)$ as an $L^2$ function, and thus have the representation
\begin{equation}\label{lcik2}
\lc_k^i (t)=  \int_0^\infty  P_k^{\lambda(t)} \psi(t) \cdot \tilde P_k^\lambda ( \lambda G_i)^{\lambda(t)} r dr.
\end{equation}
We will use this representation in order to obtain bounds for the functions
$\lc_k^i$. Our bounds will also hold uniformly with respect to the value $R$ used in the $\chi_{\leq 0}(r/R)$ cutoff, which is needed in order to establish uniform convergence in \eqref{LsumLk}. 
W will establish the bounds for $\lc_i^k$ in several steps:

\medskip

\bf Step 1: Uniform and $L^2$ bounds for  $\lc^i_k$.  \rm 
The key part of this step is to establish $L^2$ bounds for  the Littlewood-Paley pieces of $ ( \lambda G_i)^{\lambda(t)}$.
These have the form

\begin{l1}\label{l:PG}
The following bounds hold:
\begin{equation} \label{PkGest}
\| \tilde P_k^\lambda (\lambda G_i^\lambda) \|_{L^2} \les \la 2^k \lambda^{-1}\ra^{-4}, 
\end{equation}
\begin{equation} \label{PkGimproc}
\|  r  \tilde P_{k}^\lambda ( \lambda G_i^\lambda) \|_{L^2} 
\les 2^{-k} \la 2^k \lambda^{-1} \ra^{-4}. 
\end{equation}
\end{l1}
These will also hold uniformly with respect to the value $R$ used in the cutoff $\chi_{\leq 0}(r/R)$,
and all that is used is the symbol type behavior at infinity of $G_i$ with at least  $1/r$ decay. 

For technical reasons we need to provide similar estimate for the function $\tilde G_i=(rG_i)'$. From the properties of $G_i$ detailed in \eqref{lanalysis} and in particular \eqref{rGass}, it follows that 
$\tilde G_i$ is similar to $G_i$ except that it has better decay at infinity, that is $\tilde G_i = O(r^{-3})$. 

\begin{l1}\label{l:PG'}
The following holds true:
\begin{equation} \label{PkG'est}
\| P_k^\lambda (\lambda^2 \tilde G_i^\lambda) \|_{L^2} \les 2^k \la 2^k \lambda^{-1}\ra^{-4}. 
\end{equation}
\end{l1}

We remark that in order to get the estimate in \eqref{PkG'est} it suffices to have $O(r^{-2})$ decay at infinity, which is less than what $\tilde G_i$ has. Essentially the extra decay factor of $r^{-1}$ (over what $G_i$ has) gives us the extra gain of a factor of $2^k$ over \eqref{PkGest} and this is an improvement in the low frequency regime. 

The bounds in Lemma~\ref{l:PG} can be used in multiple ways. On one hand we can combine them with the uniform $L^2$ bounds on $\psi$, which leads us to 
a uniform bound for $\lc^i_k(t)$,
\begin{equation}\label{lcik-point}
|\lc^i_k(t)| \les \la 2^k \lambda(t)^{-1} \ra^{-4} \| P_k^{\lambda(t)} \psi(t)\|_{L^2}
\end{equation}
On the other hand we can combine them with local energy bounds for $\psi$, which yields an $L^2$ bound for $\lc^i_k(t)$,
\begin{equation}\label{lcik-l2}
\|\lc^i_k\|_{L^2(I)} \les 2^{-k}  \la 2^k (\lmax_I)^{-1}\ra^{-4}  \| P_k^{\lambda(t)} \psi\|_{LE_k}, \qquad \lmax_I = \sup_{t \in I} \lambda(t).
\end{equation}

We now return to the proof of the Lemmas:

\begin{proof}[Proof of Lemma~\ref{l:PG}]
 Using the symbol type behavior of $G_i$ at infinity, the Fourier
 transform $\mathcal{F}_{\tilde H_\lambda} (\lambda G_i^\lambda)$ can be easily defined as an improper integral by 
\[
(\mathcal{F}_{\tilde H_\lambda} \lambda G_i^\lambda) (\xi)= \lambda \int_0^\infty \psi^\lambda_\xi(r) G_i^\lambda(r) r dr 
= \lim_{R \to \infty} \lambda \int_0^\infty \chi_{\leq 0}(r/R) \psi^\lambda_\xi(r) G_i^\lambda(r) r dr.
\]
Omitting the $R$ truncation from here on, a direct rescaling 
yields
\[
\begin{split}
(\mathcal{F}_{\tilde H_\lambda} \lambda G_i^\lambda) (\xi) & = \lambda \int_0^\infty \psi^\lambda_\xi(r) G_i^\lambda(r) r dr \\
& =  \lambda \int_0^\infty \lambda^\frac12 \psi_{\lambda^{-1} \xi}(\lambda r)  G_i(\lambda r) r dr \\
& =   \int_0^\infty  \lambda^{-\frac12} \psi_{\lambda^{-1} \xi}(r) G_i(r) r dr .
\end{split}
\]
Thus we have 
\[
(\mathcal{F}_{\tilde H_\lambda} \lambda G_i^\lambda) (\xi)= \lambda^{-\frac12} f_i(\lambda^{-1} \xi), \qquad f_i(\eta)=\int_0^\infty   \psi_{\eta}(r) G_i(r) rdr,
\]
where the last integral defining $f_i$ is still defined as an improper integral. It remains to estimate the function $f_i$. This is done by pointwise estimates in the non-oscillatory regime $r \les \eta^{-1}$ and integration by parts in the oscillatory regime 
$r \ges \eta^{-1}$; using the properties of $G_i$ we obtain the following:
\[
\begin{split}
& |f_i(\eta)| \les \eta q(\eta)  \int_{r \les \eta^{-1}} rdr + O(\eta^{-\frac12}) 
\les \eta^{-\frac12}, \quad \eta \les 1; \\
& |f_i(\eta)|  \les \eta q(\eta) \int_{r \les \eta^{-1}} r^3 \cdot r^4 dr + O(\eta^{-\frac{9}2}) \les \eta^{-\frac{9}2}, \quad \eta \ges 1.
\end{split}
\]
A similar argument also yields bounds for the derivative of $f$,
namely 
\[
|\partial_\eta f_i(\eta)| \lesssim  \  \eta^{-\frac32} \la \eta \ra^{-4}.
\]
For the actual $\mathcal{F}_{\tilde H_\lambda} \lambda G_{i}^\lambda$,  the $f$ bounds translate into
\[
 |(\mathcal{F}_{\tilde H_\lambda} \lambda G_{i}^\lambda) (\xi)| \les \xi^{-\frac12} \la \lambda^{-1} \xi \ra^{-4}.
\]
respectively  
\[
 |\partial_\xi (\mathcal{F}_{\tilde H_\lambda} \lambda G_{i}^\lambda) (\xi)| \les  \xi^{-\frac32} \la \lambda^{-1} \xi \ra^{-4}.
\]

The first set of bounds directly imply \eqref{PkGest}. On the other hand for \eqref{PkGimproc} we use Lemma \ref{multrfourier} and the above computations to obtain:
\[
\| \lambda r  P_{k}^\lambda G^\lambda \|_{L^2_r} \les 
\| \frac{ \varphi_k(\xi) f_i(\lambda^{-1} \xi)}{\xi}\|_{L^2}
+ \| \partial_\xi (  \varphi_k(\xi) f_i(\lambda^{-1} \xi))\|_{L^2}
\]
It remains to show that each term on the right-hand side can be bounded, up to constants, by $2^{-k} \la 2^k \lambda^{-1} \ra^{-4}$. 
The estimate for the first term follows from \eqref{PkGest}; the same applies to the second term when $\partial_\xi$ is applied to the dyadic cutoff, while for the term containing $\partial_\xi f_i$ we 
use the bounds above for $\partial_\xi (\mathcal{F}_{\tilde H_\lambda} \lambda G_i^\lambda) (\xi)$.
\end{proof}

\begin{proof}[Proof of Lemma~\ref{l:PG'}]
  Just as in the previous lemma we compute 
\[
\begin{split}
(\mathcal{F}_{\tilde H_\lambda} \lambda^2 \tilde G_1^\lambda) (\xi)  =  \int  \lambda^{\frac12} \psi_{\lambda^{-1} \xi}(r) \tilde G_1(r) r dr.
\end{split}
\]
Thus with $\tilde f(\eta)=  \int  \lambda^{\frac12} \psi_{\eta}(r) \tilde G_1(r) rdr$,
we have $(\mathcal{F}_{\tilde H_\lambda} \lambda^2 \tilde G_1^\lambda) (\xi)= \tilde f(\lambda^{-1} \xi)$. 
A quick inspection of $\tilde G_1$ reveals that it has properties which are similar to those 
of $G_1$ plus the improved decay $\tilde G_1 = O(r^{-2})$, see \eqref{rGass}. By similar arguments as the ones used for estimating $f$ in Step 1, we obtain
\[
|\tilde f(\eta)| \les (\lambda \eta)^{\frac12},  \eta \les 1 \quad \mbox{and} \quad |\tilde f(\eta)| \les (\lambda \eta)^{\frac12} \eta^{-4}, \eta \ges 1.
\]
Recalling that $\lambda \eta=\xi$, these bounds translate into
\[
\begin{split}
& |(\mathcal{F}_{\tilde H_\lambda} \lambda^2 \tilde G_1^\lambda) (\xi)| \les \xi^{\frac12}, 
\quad \xi \les \lambda, \\
& |(\mathcal{F}_{\tilde H_\lambda} \lambda^2 \tilde G_1^\lambda) (\xi)| \les \xi^{\frac12} (\lambda^{-1} \xi)^{-4}, 
\quad \xi \ges \lambda.
\end{split}
\]
As a consequence of this we obtain the bound \eqref{PkG'est}.  
\end{proof}

\bigskip

\bf Step 2:  $L^2$ bounds for  $\partial_t \lc^i_k$.  \rm 
Our objective here is to prove that
\begin{equation}\label{dt-lci}
\|\partial_t \lc^i_k\|_{L^2(I)} \les 2^{k} \sum_{j} 2^{-\frac{|k-j|}{10}} \| P_j^{\lambda(t)} \psi\|_{S_j}.
\end{equation}

To prove this we begin with a simple computation,
\[
\begin{aligned}
   \partial_t \lc^i_k = & \  \la \partial_t P_k^\lambda \psi, (\lambda G_i)^\lambda\ra +  \la  P_k^\lambda \psi, \partial_t(\lambda G_i)^\lambda\ra
\\   
=  & \ i \la P_k^\lambda \psi, \tH_{\lambda} (\lambda G_i)^\lambda\ra
+ i \la P_k^\lambda  N(\psi) , (\lambda G_i)^\lambda\ra
+ \la  g_k, (\lambda G_i)^\lambda\ra
+  \la  P_k^\lambda \psi, \partial_t(\lambda G_i)^\lambda\ra
\\
:= & \ e^1_k + e^2_k + e^3_k + e^4_k,
\end{aligned}
\]
where we recall from \eqref{Pjpsieq} that $g_k = i [\partial_t, P_k^\lambda(t)] \psi$;  we now estimate separately each of the $e_k^l$
terms.
\medskip

\emph{ The bound for $e^1_k$.} This is identical to the proof of \eqref{lcik-l2}, as the operator $\tH_\lambda$ simply adds a $2^{2k}$ factor at frequency $2^{k}$; the details are left as an exercise. 

\medskip

\emph{ The bound for $e^2_k$.} This has the form
\[
e^2_k = \la  P_{k}^\lambda N(\psi), \lambda G_i^\lambda \ra
\]
For this we use the $L^2 L^1$ bound in \eqref{Npsi-21}: 
\[
\| e^2_k \|_{L^2_t} \les \| \la  \frac{P_{k}^\lambda N(\psi)}r, \lambda r G_i^\lambda \ra \|_{L^2_t} 
\les \| \frac{P_{k}^\lambda N(\psi)}r \|_{L^2_t L^1_r} \| \lambda r G_i^\lambda   \|_{L^\infty_{t,r}} \les 
2^k \sum_j 2^{-\frac{|k-j|}{10}} \| \psi_j\|_{S_j}. 
\]

\medskip

\emph{ The bound for $e^3_k$.} 
Here we look at $\la g_{k}, \lambda G_1^\lambda \ra$ and establish the estimate \eqref{dt-lci} on $[0,T]$.  From Lemma~\ref{dtP} we have
\[
g_k = i \lambda' \FtH_\lambda^{-1}[ \tcK_\lambda,m_k] \FtH_\lambda \psi. 
\]
We split this into
\[
g_k = \sum_{l,j} g_{lkj}, \qquad g_{lkj} = i \lambda' \FtH_\lambda^{-1} m_l [ \tcK_\lambda,m_k]  \tm_j \FtH_\lambda \psi_j,
\]
where the summands vanish unless either $l=k+O(1)$
or $j=k+O(1)$,
and are estimated at fixed time by Lemma~\ref{l:mj-com} as 
\begin{equation} \label{glkj}
\| g_{lkj}\|_{L^2_{r}} \lesssim 2^k \frac{|\lambda'|}{\lambda^2} \chi_{2^j = \lambda}  \chi_{k=l} \chi_{k=j} \| \psi_j \|_{L^\infty L^2}.
\end{equation}
Next we write
\[
\la g_{lkj}, \lambda G^\lambda \ra =\la  g_{lkj}, \lambda  \tP_l^\lambda G^\lambda \ra,
\]
and, using \eqref{PkGest}, we estimate again at fixed time
\[
\begin{split}
|\la g_{lkj}, \lambda G^\lambda \ra| & \les \| g_{lkj}\|_{L^2_{r}} \la 2^l \lambda^{-1}\ra^{-4} \\
& \les  2^k \frac{|\lambda'|}{\lambda^2} \chi_{2^k = \lambda}  \chi_{k=l} \chi_{k=j} 
\la 2^l \lambda^{-1}\ra^{-4}
\| \psi_j \|_{L^\infty L^2}.
\end{split}
\]
After $j$ and $l$ summation we arrive at 
\[
\begin{split}
|\la g_{k}, \lambda G^\lambda \ra| & \les 
 2^k \frac{|\lambda'|}{\lambda^2}  \sum_j \chi_{2^j = \lambda} \chi_{k=j} 
\| \psi_j \|_{L^\infty L^2}.
\end{split}
\]
Finally we take the $L^2_t$ norm  
to get a final bound
\[
\begin{split}
\|\la g_{k}, \lambda G_1^\lambda \ra\|_{L^2_t} \les 2^k \sum_{j} M_j \chi_{k=j} \| P_j^\lambda \psi\|_{L^\infty_t L^2_r} .
\end{split}
\]
\medskip
\emph{ The bound for $e^4_k$.}
Here we first compute  
\[
\partial_t [\lambda  G_1^\lambda] =  \lambda'(  G_1 + r G'_1)^\lambda = \frac{\lambda'}{\lambda^2}
\lambda^2 \tilde G_1^\lambda, \qquad \tilde G_1= (rG_1)'.
\]
 Using \eqref{PkG'est} we obtain
\[
\begin{split}
\| \frac{\lambda'}{\lambda^2} \la P_{k}^\lambda \psi , \lambda^2 \tilde G_1^\lambda  \ra| \|_{L^2_t} 
& \les \| \frac{\lambda'}{\lambda^2}  \|_{L^2_t} \sup_t \| P_{k}^\lambda \psi  \|_{L^2_r} \| \tilde P_k^\lambda (\lambda^2 \tilde G_1)^\lambda\|_{L^2}  \\
& \les  \| \frac{\lambda'}{\lambda^2}  \|_{L^2_t} 2^k   \la 2^k (\lmax_I)^{-1}\ra^{-4}\| P_k^\lambda \psi \|_{L^2}. 
\end{split}
\]
as needed.

\medskip

\textbf{ Step 3: The dyadic summation of the $\lc^i_k$ bounds.}
Here we complete the proof of Theorem~\ref{t:mod-syst}. Our starting point consists of 
the bounds \eqref{lcik-l2} and \eqref{dt-lci}, which hold in any time interval $[0,T]$ of existence for the solution. Here we may assume without loss of generality that $T \gtrsim \lambda(0)^{-2}$,
which is the existence time for the local well-posedness result.
To use these bounds in order to estimate $\Lambda^i$ in \eqref{LsumLk}, we separate the dyadic frequency indices $k$  as follows:
\begin{enumerate}
\item low frequencies, $2^{k}  \lesssim T^{-\frac12}$. 
\item medium frequencies, $T^{-\frac12} \lesssim 2^{k} \lesssim   \lmax_T$.
\item high frequencies, $\lmax_T \lesssim 2^{k}$.
\end{enumerate}

\medskip

\emph{ Low frequencies:}
Here we directly use H\"older's inequality and \eqref{dt-lci} to obtain a uniform bound
\begin{equation} \label{corlf}
| \lc_k^i(t)- \lc_k^i(0)| \les (2^{2k} t)^\frac12 \sum_{j} \chi_{j=k} \| P_j \psi\|_{S_j}, 
\end{equation}
as well as an integrated form of this,
\begin{equation} \label{corlf-L2}
\| \lc_k^i(t)- \lc_k^i(0)\|_{L^2[0,T]} \les 2^{k} T\sum_{j} \chi_{j=k}  \| P_j \psi\|_{S_j} . 
\end{equation}

\medskip

\emph{ Medium and high frequencies:}
In this case our time interval is long enough so we can interpolate between 
the $L^2$ and $\dot H^1$ bound for $\lc_k^i$, which yields
\begin{equation} \label{corlfa}
\| \lc_k^i\|_{L^\infty} \les  \la 2^k (\lmax_T)^{-1}\ra^{-2} \sum_{j} \chi_{j=k} \| P_j \psi\|_{S_j}. 
\end{equation}
where the distinction is that we have extra decay for high frequencies.

\bigskip

Given the above bounds with decay both at low and at high frequencies, it is clear that we have uniform convergence in \eqref{LsumLk}, with a uniform overall bound obtained by applying H\"older's inequality for the intermediate frequencies, namely 
\begin{equation}
\|\lc^i\|_{L^\infty} \lesssim (\ln (T (\lmax_T)^2))^\frac12 \| u\|_{l^2 S}.
\end{equation}
If instead we assume $l^1$ summability then this becomes a uniform bound, 
\begin{equation}
\|\lc^i\|_{L^\infty} \lesssim \| u\|_{l^1 S}.
\end{equation}

\medskip

We now consider the $\dot H^{\frac12}$ bound for $\lc^i$. One difficulty we face is that we are in a bounded interval $[0,T]$. But this can be easily bypassed by a standard extension and truncation
argument. Here the constants do not matter, but we have to consider them carefully when we truncate.
The best strategy is to separate two cases, exactly as above:

\begin{enumerate}
    \item Low frequency: here we take the function $\lc_k^i - \lc_k^i(0)$ in $[0,T]$, 
which we simply extend by reflection to $[0,2T]$ and then by $0$ outside this interval. This leaves the bounds \eqref{corlf-L2} and \eqref{dt-lci} unchanged.

\item Medium and high frequencies: here take the function $\lc_k^i$ in $[0,T]$,
which we extend by reflection to $[-T,2T]$, and then truncate it outside $[0,T]$.
This leaves the bounds \eqref{lcik-l2} and \eqref{dt-lci} unchanged.
\end{enumerate}

After this, we estimate the dyadic $L^2$ norms by interpolating between the $L^2$ and the $\dot H^1$
norms. For the low frequency part the two norms balance exactly at time frequency $T^{-1}$, and we 
obtain (using $\PP_m$ for Littlewood-Paley projections in time)
\[
\| \PP_m (\lc_k^i(t) - \lc_k^i(0))_{ext} \|_{\dot H^{\frac12}} \lesssim (2^k T^\frac12) \chi_{m = - \log_2 T} \sum_{j} \chi_{j=k} \| P_j \psi\|_{S_j}. 
\]
For the medium frequency part the two norms balance exactly at time frequency $2^{2k}$, and we 
obtain 
\[
\| \PP_m (\lc_k^i)_{ext} \|_{\dot H^{\frac12}} \lesssim  \chi_{m = 2k} \sum_{j} \chi_{j=k} \| P_j \psi\|_{S_j}. 
\]
Finally for the high frequency part the two norms balance exactly at time frequency $2^{m_k} = 2^{2k}(2^k (\lmax_T)^{-1})^4$, and we 
obtain 
\[
\| \PP_m (\lc_k^i)_{ext} \|_{\dot H^{\frac12}} \lesssim  \chi_{m = m_k} 
(2^k (\lmax_T)^{-1})^{-2}
\sum_{j} \chi_{j=k} \| P_j \psi\|_{S_j}. 
\]
Since all constants in front of the sums are $\lesssim 1$ and we have in all cases off-diagonal decay, the summation with respect to $m$ is straightforward, and we obtain immediately the bounds in \eqref{lcorH12}. For \eqref{lcorB12} we observe that we have additional decay above both for the low and for the high frequencies, so it suffices to bound the contributions of intermediate $k$. 
But there we simply convert the $l^2$ norm to $l^1$ using H\"older's inequality, noting that 
the number of intermediate dyadic regions is about $\ln (T (\lmax_T)^2)$.

\section{An abstract ode result}
\label{s:ode}

In this section we consider the solvability question for 
nonlinear ode systems such as our modulation equations
\eqref{mainaleq}. As written there, the two components are partially uncoupled, in that it suffices to solve first the $\alpha$ equation, and then the $\lambda$ equation is a direct integration. To keep the notations simple in this section we consider a more general vector valued model of the form 
\begin{equation} \label{mainz}
\alpha'=  \Re( n(\alpha)  f'(t)) + g(t), \quad \alpha(0)=\alpha_0,
\end{equation}
where $n$ is assumed to be a globally
$C^3$ function, $\|n\|_{C^3} < +\infty$, while $f$ and $g$ are taken in the spaces
\[
f \in \dot B^{\frac12}_{2,1}[\R], \qquad g \in L^1(\R).
\]
For this problem we look for solutions in the space
\begin{equation}
Z=(\dot B^{\frac12}_{2,1}[\R] + \dot W^{1,1}[\R]) \cap L^\infty \subset C^0,
\end{equation}
where the $L^\infty$ component of the norm 
simply has the goal of controlling constants.
We successively consider the small data case and then the large data case. For the latter we will also use the same space in an interval $I$, either bounded or unbounded,
\begin{equation}
Z[I]=(\dot B^{\frac12}_{2,1}[I] + \dot W^{1,1}[I]) \cap L^\infty \subset C^0.
\end{equation}

\subsection{Iterating small Besov data}
Our goal here is to show the following
\begin{p1} \label{mainzL}
There exists $\epsilon >0$ such that if 
\begin{equation}
\|f\|_{\dot B^{\frac12}_{2,1}[\R]} \leq \epsilon,  \qquad 
\| g \|_{L^1[\R]} \leq \epsilon,
\end{equation}
then \eqref{mainz} has a unique global solution $\alpha$
with 
\begin{equation}
\|\alpha-\alpha_0\|_{Z}\les \epsilon. 
\end{equation}
\end{p1}
We remark that as a direct consequence, the same result holds also in any interval $I$.
\begin{proof} 
We will use a fixed point argument in the space 
$Z$. We also denote by $Z'$ the space of derivatives of functions in $Z$, namely
\[
Z' = \dot B^{-\frac12,1}[\R] +  L^{1}[\R].
\]
For these spaces we have the following properties:
\begin{l1}
a) $Z$ is an algebra.

b) We have the bilinear bound $ Z \cdot Z' \subset Z'$. 

c) For a $C^3$ function $n$ with $n(0) = 0$ we have the Moser inequality
\[
\| n(f)\|_{Z} \lesssim  \| f \|_{Z} +\|f\|_Z^3.
\]
\end{l1}

\begin{proof}
a) By the Leibniz rule, part (a) may be seen as a direct consequence of part (b). 

\medskip

b) Consider a product $v w'$ with $v, w \in Z$. For $v$ and $w$ we need to consider the two components of the $Z$ norm.  So
 we let $v=v_1+v_2, w=w_1+w_2$ with 
 $v_1,w_1 \in \dot B^{\frac12}_{2,1}(\R)$ and $v_2,w_2 \in \dot W^{1,1}(\R)$, 
 so that 
 \[
 \|v_1\|_{\dot B^{\frac12}_{2,1}} + \|v_2\|_{\dot W^{1,1}} \les \|v\|_{Z}, \qquad \|w_1\|_{\dot B^{\frac12}_{2,1}} + \|w_2\|_{\dot W^{1,1}} \les \|w\|_{Z}.
 \]
The contribution of $w_2$ to $vw'$ is easily estimated  in $L^1$, so it remains to consider the product $vw'_1$. We expand it as follows
\[
v w_1' = \sum_k P_{<k-10} v P_k w_1' + 
\sum_k P_{k} v P_{<k-10} w_1' + \sum_{|k_1-k_2| < 10} P_{k_1} v_1 P_{k_2} w_1'
+ \sum_{|k_1-k_2| < 10} P_{k_1} v_2 P_{k_2} w_1',
\]
and we consider each sum separately.
The summands in the first sum are localized  at frequency $2^k$, so it is convenient to bound the sum in the Besov norm,
\[
\left\| \sum_k P_{<k-10} v P_k w_1'  \right\|_{\dot B^{-\frac12,1}}
\lesssim \sum_k 2^{-\frac{k}2} \| P_{<k-10} v\|_{L^\infty} \| P_k w_1'\|_{L^2}
\lesssim \|v\|_{L^\infty} \|w_1\|_{\dot B^{\frac12}_{2,1}}.
\]
The second sum is similar, 
\[
\left\|  \sum_k P_{k} v P_{<k-10} w_1' \right\|_{\dot B^{-\frac12,1}}
\lesssim \sum_k 2^{-\frac{k}2} \| P_{k} v\|_{L^\infty} \| P_{<k-10} w_1'\|_{L^2}
\lesssim \|v\|_{L^\infty} \|w_1\|_{\dot B^{\frac12}_{2,1}}.
\]
But the remaining two sums are instead estimated in $L^1$:
\[
\| \sum_{|k_1-k_2| < 10} P_{k_1} v_1 P_{k_2} w_1' \|_{L^1}
\lesssim \sum_{|k_1-k_2| < 10} 2^{k_2} \| P_{k_1} v_1 \|_{L^2} \|P_{k_2} w_1 \|_{L^2}
\lesssim  \|v_1\|_{\dot B^{\frac12}_{2,1}} \|w_1\|_{\dot B^{\frac12}_{2,1}}.
\]
respectively
\[
\|\sum_{|k_1-k_2| < 10} P_{k_1} v_2 P_{k_2} w_1'  \|_{L^1} 
\lesssim \sum_{|k_1-k_2| < 10} \|  P_{k_1} v_2\|_{L^1} \| P_{k_2} w_1'  \|_{L^\infty} \lesssim  \|v_2'\|_{L^1} \|w_1\|_{\dot B^{\frac12}_{2,1}}.
\]

\medskip 

c) Let $f = f_1+f_2$ where 
$f_1 \in \dot B^{\frac12}_{2,1}$ and $f_2 \in \dot W^{1,1}$.
The easier case is when $f_1=0$, where the Moser bound is a simple application of chain rule. To deal with $f_1$, we consider a continuous Littlewood-Paley expansion
\[
f_1 = \int_{-\infty}^\infty P_h f_1 dh.
\]
For instance, for $h \in \R$, we can define $P_{\leq h}$ to be the zero order multiplier whose symbol is $\chi_{\leq 0}(\frac{\xi}{2^h})$; recall from Section \ref{defnot} that $\chi_{\leq 0}$ is smooth, supported in $(0,2)$ and is identically equal to $1$ on $(0,1)$ . Then we let $P_h f_1=\frac{d}{dh} P_{\leq h} f_1$ and record that
\[
f_{1,\leq k} := P_{\leq k} f_1 = \int_{-\infty}^k P_h f_1 dh;
\]
here $k \in \R$ (not to be confused with the standard choice in $\Z$). Accordingly we can define the operators $P_{< h}, P_{\geq h}, P_{> h}$. 
We record the following basic inequalities:
\begin{equation} \label{basicPh}
\|2^{-h} P_{\leq h} f' \|_{L^p} \leq \|f\|_{L^p} , \quad 
\|2^h P_{\geq h} f \|_{L^p} + \| P_{\geq h} f' \|_{L^p} \leq \|f'\|_{L^p}.
\end{equation}
which hold for any $1 \leq p \leq \infty$. The proof of these estimates is simply based on estimating kernels in $L^1$ and it relies on the smoothness of the frequency cut-off used (that is $\chi_{\leq 0}$); the details are left as an exercise.  

Based on the calculus above, we obtain the following expansion 
\[
\begin{aligned}
n(f) = & \ n(f_2) + \int_{-\infty}^\infty \frac{d}{dh} n (f_{1,<h} + f_2) dh
\\
= \ & n(f_2) + \int_{-\infty}^\infty  n' (f_{1,<h} + f_2) P_h f_1 dh
\\
= \ & n(f_2) + \int_{-\infty}^\infty  (n'(f_{1,<h} + f_2) -n' (f_{1,<h} + f_{2,<h})) P_h f_1 
+ n' (f_{1,<h} + f_{2,<h}) P_h f_1 
dh 
\\
:= \ & n(f_2) + \int_{-\infty}^\infty (g_1^h + g_2^h) dh := n(f_2) + g_1 + g_2.
\end{aligned}
\]
Here we estimate $g_1^h$ and $g_2^h$ separately. We place $g_1^h$ in $\dot W^{1,1}$, for which we compute
\[
\begin{aligned}
(g_1^h)' = & \ (n''(f_{1,<h} + f_2)(f'_{1,<h} + f'_2) -n'' (f_{1,<h} + f_{2,<h}) 
(f'_{1,<h} + f'_{2,<h}))
P_h f_1 \\  & \ + 
(n'(f_{1,<h} + f_2) -n' (f_{1,<h} + f_{2,<h})) P_h f_1'. 
\end{aligned}
\]
Estimating this in $L^1$ we obtain
\[
\begin{aligned}
\| (g_1^h)'\|_{L^1} \lesssim & \  \|f_2 -f_{2,<h}\|_{L^1} \|f'_{1,<h} + f'_{2,<h}\|_{L^\infty} \| P_h f_1\|_{L^\infty}
+ \|f'_2 -f'_{2,<h}\|_{L^1} \| P_h f_1\|_{L^\infty} \\ 
& +\| f_2-f_{2,< h} \|_{L^1} \| P_h f'_1\|_{L^\infty} \\ 
\lesssim & \| 2^h f_{2,\geq h}\|_{L^1} \|2^{-h} (f'_{1,<h} + f'_{2,<h})\|_{L^\infty} \| P_h f_1\|_{L^\infty}
+ \|f'_{2,\geq h}\|_{L^1} \| P_h f_1\|_{L^\infty} \\
& + \| 2^h f_{2,\geq h}\|_{L^1}  \| 2^{-h} P_h f'_1\|_{L^\infty}
\\
\lesssim & \ \| f_2' \|_{L^1} (1+\|f_1\|_{L^\infty}+ \|f_2\|_{L^\infty}) \| P_h f_1\|_{L^\infty}  
\\
\lesssim & \ \| f_2' \|_{L^1} (1+\|f\|_{Z}) \| P_h f_1\|_{{\dot B^{\frac12}_{2,1}}},
\end{aligned}
\]
where we have used the basic estimates in \eqref{basicPh} and  the $P_h f_1$ factor was bounded in terms of the Besov norm.

After  $h$ integration this yields
\[
\| g'_1\|_{L^1} \lesssim \| f_2'\|_{L^1} \|f_1\|_{\dot B^{\frac12}_{2,1}}(1+\|f\|_{Z}).
\]

On the other hand $g_2^h$ may be placed in the Besov space, by estimating it in 
$L^2$ and in $\dot H^1$, as follows:
\[
\| g_2^h\|_{L^2} \lesssim \|  P_h f_1\|_{L^2} \lesssim 2^{-\frac{h}{2}}
  \|  P_h f_1\|_{\dot B^{\frac12}_{2,1}}
\]
respectively
\[
\| (g_2^h)'\|_{L^2} \lesssim \|  P_h f'_1\|_{L^2} + 
\|f'_{1,<h} + f'_{2,<h}\|_{L^\infty} 
\|  P_h f_1\|_{L^2} 
\lesssim 2^{\frac{h}{2}}(1+\|f\|_{Z})
  \|  P_h f_1\|_{\dot B^{\frac12}_{2,1}}
\]
Combining the two we obtain 
\[
\| g_2^h\|_{\dot B^{\frac12}_{2,1}} \lesssim  (1+\|f\|_{Z})
  \| P_h f_1\|_{\dot B^{\frac12}_{2,1}}.
\]
This concludes the argument for the Lemma. 
\end{proof}

This lemma allows us to run a fixed point argument as follows. Given $\alpha \in Z$, we define
\[
N(\alpha) = \Re(n(\alpha f') + g)
\]
Then \eqref{mainz} is equivalent to the integral form
\[
\alpha(t) = \int_0^t N(\alpha)(s) ds + \alpha_0. 
\]
We solve this using the contraction principle in a ball in $Z$ of size $C \epsilon$, where $C$ is a universal large constant. To achieve this it suffices to show that the operator $u \to N(u)$ 
maps $B_Z(C\epsilon)$ to $B_{Z'}(\epsilon)$ with a small Lipschitz constant.
But this is a direct consequence of the lemma above.

\end{proof}

\subsection{Iterating large Besov data}

Here we consider the same ode \eqref{mainz}, but we allow $f$ and $g$ to be large in $\dot B^{\frac12}_{2,1}$, respectively $L^1$. 

\begin{p1} \label{mainzLL}
For  $M \gtrsim 1$, assume that  
\begin{equation}
\|f\|_{\dot B^{\frac12}_{2,1}(\R)} \leq M,  \qquad 
\| g \|_{L^1(\R)} \leq M.
\end{equation}
then \eqref{mainz} has a unique global solution $\alpha$
with 
\begin{equation}\label{large-f}
\|\alpha-\alpha_0\|_{Z}\les M^2. 
\end{equation}
\end{p1}

We remark that as a direct consequence, the same result holds also in any interval $I$.

\begin{proof}
The main ingredient of the proof is the following divisibility lemma for the 
$\dot B^{\frac12}_{2,1} (\R)$ norm:

\begin{l1}\label{l:division}
a) Given an arbitrary partition of $\R$ into intervals $I_j=[a_j,b_j]$,
we have the divisibility bounds
\begin{equation}\label{B-divisible}
\sum_{j} \|f\|^2_{\dot B^{\frac12}_{2,1} [I_j]} \les \|f\|^2_{\dot B^{\frac12}_{2,1} [\R]}, 
\end{equation}
respectively 
\begin{equation}\label{W-divisible}
\sum_{j} \|f\|^2_{\dot W^{1,1} [I_j]} \les \|f\|^2_{\dot W^{1,1} [\R]}. 
\end{equation}

b) For a converse bound, we have 
\begin{equation}\label{B-triangle}
\|f\|_{\dot B^{\frac12}_{2,1} (\R)+ \dot W^{1,1}(\R)} \lesssim \sum_{j} \|f\|_{\dot B^{\frac12}_{2,1}(I_j)+\dot W^{1,1}[I_j]}. 
\end{equation}
\end{l1}

\begin{proof}
a) The bound \eqref{W-divisible} is obvious, so we focus on \eqref{B-divisible}. Since on the right of \eqref{B-divisible} we have an $\ell^1$-Besov norm, without any restriction in generality 
we can assume that $f$ is localized at a single dyadic frequency. By scaling,
we can assume that this frequency is $1$. Normalizing we can also assume 
that $\|f\|_{L^2} = 1$ therefore 
\[
\| \partial_x f\|_{L^2} + \| \partial_x^2 f\|_{L^2} \lesssim 1.
\]
We measure the $L^2$ norm of $f$ in unit size intervals $[k,k+1]$ and 
bound it using a frequency envelope $(c_k)$, so that the following properties hold:

\begin{enumerate}[label=\roman*)]
\item $\| f \|_{L^2[k,k+1]} +\| f' \|_{L^2[k,k+1]} 
\leq c_k, \forall k$;

\item $(c_k)$ is  slowly varying in the following sense: $c_i \les \la i-j \ra^2 c_j, \forall i,j$;

\item $\sum_k c_k^2 \les 1.$
\end{enumerate}
In particular by Sobolev embeddings we will also have 
\[
\| f \|_{L^\infty[k,k+1]} +
\| f' \|_{L^\infty[k,k+1]} 
+\| f' \|_{L^2[k,k+1]} 
\les c_k.
\]

For an interval $I$ we denote 
\[
c_I^2 = \sum_{I \cap [k,k+1] \neq\emptyset} c_k^2.
\]

To measure the Besov norm of the function $f$ in an interval $I$ we use 
the interpolation inequality
\[
\|f\|_{\dot B^{\frac12}_{2,1} [I]}^2
\lesssim \| f-f_I\|_{L^2[I]} \|f'\|_{L^2[I]},
\]
where $f_I=|I|^{-1} \int_I f$ is the average of $f$ on $I$. This is easily seen in $\R$, and then it can be transferred to an interval $I$ 
using a suitable extension which preserves the size of norms on the right (e.g. one even reflection followed by a constant extension).

Then from Poincare's inequality we obtain
\[
\| f\|_{\dot B^{\frac12}_{2,1} [I]}^2
\lesssim |I| \|f'\|_{L^2[I]}^2.
\]
To prove \eqref{B-divisible} we separate in intervals $I_j$ into two classes:

\medskip

i) long intervals, $|I_j| \geq 1$.
Here we use the above interpolation inequality to write
\[
\|f\|_{B^{\frac12}_{2,1} [I_j]}^2
\lesssim c_{I_j}^2.
\]
Then after summation we get 
\[
\sum_{I_j \ long} \|f\|_{B^{\frac12}_{2,1} [I_j]}^2 
\lesssim \sum_{I_j \ long} c_{I_j}^2 \lesssim 1,
\]
since each $c_k$ is counted at most twice.

ii) Short intervals, $|I_j| \leq 1$.
Then we use the second bound above
to write
\[
\|f\|_{B^{\frac12}_{2,1} [I_j]}^2
\lesssim |I_j| c_{I_j}^2.
\]
Here we may have many subintervals $I_j$
intersecting a given unit interval $[k,k+1]$, but the sum of their lengths
is bounded by $3$. Therefore after summation we obtain
\[
\sum_{I_j \ short} \|f\|_{B^{\frac12}_{2,1} [I_j]}^2 
\lesssim \sum_{I_j \ short} |I_j| c_{I_j}^2 \lesssim \sum_k c_k^2 \lesssim 1.
\]
This concludes the proof  of part (a) of the lemma.
\bigskip

b) The result would be trivial for the $\dot W^{1,1}$ norm. In order to deal with the larger  space $(\dot B^{\frac12}_{2,1} + \dot W^{1,1}) [I_j]$, it suffices to prove the following Lemma for a single interval:

\begin{l1}\label{l:better-split}
Let $I \subset \R$ be any interval and $f \in (\dot B^{\frac12}_{2,1} + \dot W^{1,1})[I]$. Then there exists a decomposition $f = f_1+f_2$ in $I$ where 
\[
\| f_1\|_{\dot B^{\frac12}_{2,1} (\R)} + \|f_2\|_{ \dot W^{1,1} [I]}
\lesssim \|f\|_{(\dot B^{\frac12}_{2,1} + \dot W^{1,1})[I]}
\]
with $f_1$ supported in $I$.    
\end{l1}
\begin{proof}
Rescaling we can assume that $I = [0,1]$. It suffices to consider the case 
when we have $f \in \dot B^{\frac12}_{2,1} (\R)$. This means that there exists an extension,
still denoted by $f$, so that 
\[
\|f\|_{\dot B^{\frac12}_{2,1} (\R)} \lesssim \|f\|_{\dot B^{\frac12}_{2,1}(I)}.
\]
We consider a dyadic decomposition  $f = \sum_{j \in \Z} P_j f$, and construct 
$f_1$ by truncating the dyadic pieces appropriately. Precisely, we set
\[
f_{1} = \sum_{j \geq 0} f_{1,j} : = \sum_{j > 0 } \chi_j P_j f, 
\]
where the cutoff functions $\chi_j$ have support in $I=[0,1]$ and equal $1$ 
in $[2^{-j}, 1-2^{-j}] \subset I$, and satisfy $|\chi_j'| \lesssim 2^{-j}$.
Then we can estimate $f_1$ as follows:
\[
\|f_1\|_{\dot B^{\frac12}_{2,1} (\R)} \lesssim \sum_{j > 0} 
2^{\frac{j}2} \| f_{1,j} \|_{L^2} + 2^{-\frac{j}2} \| f'_{1,j} \|_{L^2}
\lesssim \sum_{j > 0} 
2^{\frac{j}2} \| P_j f\|_{L^2} + 2^{-\frac{j}2} \| P_j f' \|_{L^2}
\lesssim \|f\|_{\dot B^{\frac12}_{2,1} (\R)}.
\]
On the other hand $f_2$ is given by 
\[
f_2 = f_{2,lo}+ f_{2,hi}:=  \sum_{j \leq 0} P_j f + \sum_{j >0} (1-\chi_j) P_j f,
\]
and we bound its derivative in $L^1$ by estimating separately the two sums as follows:
\[
\| f'_{2,lo}\|_{L^1[I]} \lesssim \sum_{j \leq 0} \| P_j f'\|_{L^\infty}
\lesssim \sum_{j \leq 0} 2^{j} \| P_j f\|_{L^\infty} \lesssim \|f\|_{\dot B^{\frac12}_{2,1} (\R)},
\]
respectively, using H\"older's inequality, 
\[
\| f'_{2,hi}\|_{L^1[I]} \lesssim \sum_{j > 0} 2^{-j}  \| ((1-\chi_j) P_j f)'\|_{L^\infty}\lesssim \sum_{j > 0} 2^{-j} ( \| P_j f'\|_{L^\infty} + 2^{j} \|P_j f\|_{L^\infty}) \lesssim \|f\|_{\dot B^{\frac12}_{2,1} (\R)}
\]
which is exactly as needed.
\end{proof}
Once we have Lemma~\ref{l:better-split}, the conclusion of part (b) of Lemma~\ref{l:division}
follows directly from the triangle inequality.
\end{proof}

Now we get back to the proof of the proposition. Lemma~\ref{l:division} allows us to split $\R$ into $\les \frac{M^2}{\epsilon^2}$ intervals $I_j$ with the property that $\|f\|_{\dot B^{\frac12}_{2,1}[I_j]+ \dot W^{1,1}[I_j]}\leq \epsilon$ (to do so we use that the Besov norm $\|f\|_{\dot B^{\frac12}_{2,1}[I_j]+ \dot W^{1,1}[I]}$ depends continuously on the interval $I$; since we can identify two nearby intervals via scaling and translation, this is a consequence of the fact that the 
scaling and translation groups are continuous in $\dot B^{\frac12}_{2,1}[\R]+ \dot W^{1,1}[\R]$).
Then we successively apply  the small data result in Proposition~\ref{mainzL} on each of these subintervals. We arrive 
at a solution $\alpha$ which is continuous in time and satisfies 
\begin{equation}
\| \alpha \|_{(\dot B^{\frac12}_{2,1} + \dot W^{1,1}) [I_j]}   \lesssim \epsilon.
\end{equation}
To obtain a global bound from here, we use the triangle type inequality
\eqref{B-triangle}, and the desired bound \eqref{large-f} follows.
\end{proof}

\section{The final bootstrap}
\label{s:final}

\subsection{Proof of Theorem \ref{tmain-G}}

Here we bring all the elements together to prove our second main result in this paper, that is lack of finite time blow-up and some control on the growth rate. The idea is the following: we start from the system \eqref{mainaleq}, which we recall here:
\begin{equation} \label{mainaleq2}
\begin{split}
\alpha' & = -\frac12 \Re \left( i e^{-2 i \alpha(t)} \partial_t ( \lc^1 + 2 \lc^2) \right) + \frac14 \Re Q(t),
\\ (\ln \lambda)' & = -\Re ( e^{-2 i \alpha(t)} \partial_t  \lc^1) + \frac12 \Im   Q(t). 
\end{split}
\end{equation}
The main inputs in this system are $\lc^i$ which have estimates in $\dot B^{\frac12}_{2,1}$ as given by \eqref{lcorB12}; the inputs $Q$ are small in $L^1$ by \eqref{qest}, and easily manageable. Thus we can employ the ode theory developed in the previous section in order to obtain information on $\alpha$ and $\lambda$, particularly about their rate of growth in time. There is a major conflict in this analysis, which we need to resolve: on one hand, the bounds in \eqref{lcorB12} depend on how large $\lambda(t)$ becomes on the time interval $I$ (through the factor $\lambda_T^{max}$ there), while, on the other hand, these bounds are then being used to estimate precisely the growth of $\lambda(t)$. Below we show how to resolve this 
conflict, and prove the growth bound \eqref{lamdagb} for $\lambda$. This in turn shows that $\lambda$ cannot reach infinity nor zero in finite time, and thus that the solutions are global in time.

We now turn to the details of the strategy described above. Assume that we work on an time interval $I=[0,T] \subset [0,T_{max})$ where our solution exists; in particular 
$\lambda^{max}_T=\max \{\lambda(t): t \in I \}$ is well defined. We recall that \eqref{lcorB12} provides a bound on the Besov bound of $\lc^i$ as follows
\[
\| \lc^i \|_{\dot B^{\frac12}_{2,1}([0,T])} \les  \ln ( 2+ ( \lambda_T^{max})^2  T)^\frac12 \| \psi_0 \|_{L^2}, \quad i=1,2.
\]
while \eqref{qest} shows that 
\[
\|Q\|_{L^1} \lesssim \| \psi_0 \|_{L^2}^2
\]
 Then applying the large data ode result in Proposition~\ref{mainzLL} to the system \eqref{mainaleq2} we obtain the bound
\[
\| \alpha(t) -\alpha(0)\|_{Z} 
+ \| \ln \lambda(t) - \ln \lambda (0)\|_{Z} \les\ln ( 2+ ( \lambda_T^{max})^2  T) \| \psi_0 \|_{L^2}^2.
\]

For simplicity we normalize $\lambda(0) = 1$. 
Then the above bound implies in particular that
\[
\ln \lambda_T^{max} \leq C\ln ( 2+ ( \lambda_T^{max})^2  T) \| \psi_0 \|^2_{L^2}
\leq  C \| \psi_0 \|^2_{L^2} \ln (2+T) +  C \| \psi_0 \|^2_{L^2} \ln \lambda_T^{max}. 
\]
By taking $\|\psi_0\|_{L^2} \leq \delta \ll 1$ small enough, we can absorb the last term and conclude that 
\[
\ln \lambda_T^{max} \leq \| \psi_0 \|^2_{L^2} \ln (2+T),
\]
which implies the bound
\[
\lambda_T^{max} \les e^{C\delta^2 \ln (2+T) } = (2+T)^{C \delta^2}. 
\]

The bound from below on $\lambda(t)$ is much simpler and follows directly from the local well-posedness result in Theorem~\ref{th2}. Precisely,
if we had $|\lambda(T)^2 T| \ll 1$ then Theorem~\ref{th2} would imply that $\lambda(T) \approx \lambda(0)$. This cannot be if $T \gg 1$,
so we conclude that 
\[
\lambda(T) \gtrsim \la T \ra^{-\frac12}. 
\]

\subsection{Proof of Theorem~\ref{tmain-G1}}
We recall that here we work with data which,
at the level of $\psi$, has smallness in $L^2$
but also the $l^1$ Besov bound 
\[
\| \psi(0)\|_{\LX} \lesssim M.
\]
By Theorem~\ref{tmain-L1}, this in turn yields
a global bound for the solution, namely
\[
\|\psi\|_{l^1 S} \lesssim M.
\]
This last bound allows us to apply the second part of \eqref{lcorH12} in Theorem~\ref{t:mod-syst}, which yields an $l^1$ Besov bound 
for $\Lambda^1$, $\Lambda^2$ in the modulation system \eqref{mainaleq2},
\[
\| \Lambda^{1,2}\|_{B^{\frac12}_{2,1}} \lesssim M.
\]
On the other hand for $Q$ we still have the favourable small $L^1$ bound as in the previous proof. 

Now we can directly use Proposition~\ref{mainzLL} for  the modulation system \eqref{mainaleq2}, which yields the global $\alpha$ 
bound
\[
\|\alpha\|_{
\dot B^{\frac12}_{2,1} + \dot W^{1,1}}
+ \|\ln \lambda\|_{
\dot B^{\frac12}_{2,1} + \dot W^{1,1}}
\les M^2. 
\]
Finally this implies the desired bound \eqref{lamdagb1}, completing the proof of the 
the theorem.

\subsection{Proof of Theorem~\ref{tstable}} Here we establish the stability result, simply by concatenating the required building blocks.

On one hand for $\psi$ we can use Theorem~\ref{tmain-L}, which as a corollary 
of \eqref{psiest-m1} yields 
\[
\|\psi \|_{L^\infty \LX} \lesssim \gamma,
\]
after which we return to $u$ via Proposition~\ref{uQpsi} to get 
\[
\|u - Q_{\alpha(t),\lambda(t)}\|_{L^\infty \bX} \lesssim \gamma.
\]

On the other hand for the modulation parameters we can use the small Besov data result in Proposition~\ref{mainzL}, which yields
\[
\|\alpha\|_{
\dot B^{\frac12}_{2,1} + \dot W^{1,1}}
+ \|\ln \lambda\|_{
\dot B^{\frac12}_{2,1} + \dot W^{1,1}}
\les \gamma. 
\]
and the difference bound
\[
 | \ln \lambda(t) - \ln \lambda(0)| + |\alpha(t) - \alpha(0)| \les \gamma.
\]

Combining the two we obtain  the bound \eqref{tsolution} and conclude the proof of the theorem.

\bibliographystyle{amsplain} \bibliography{SM-refs}

\end{document}